\documentclass[11pt,english,a4paper]{amsart}
\title{Hypercyclicity of Toeplitz operators with smooth symbols}

\date{}

\author[Fricain]{Emmanuel Fricain}
 \address{Univ. Lille, CNRS, UMR 8524 - Laboratoire Paul Painlevé, F-59000 Lille, France}
 \email{emmanuel.fricain@univ-lille.fr}

\author[Grivaux]{Sophie Grivaux}
\address{Univ. Lille, CNRS, UMR 8524 - Laboratoire Paul Painlevé, F-59000 Lille, France}
\email{sophie.grivaux@univ-lille.fr}

\author[Ostermann]{Ma\"eva Ostermann}
\address{Univ. Lille, CNRS, UMR 8524 - Laboratoire Paul Painlevé, F-59000 Lille, France}
\email{maeva.ostermann@univ-lille.fr}

\thanks{The authors were supported in part by the Labex CEMPI (ANR-11-LABX-0007-01). This work was supported in part by the project COMOP of the French National Research Agency (grant ANR-24-CE40-0892-01). The authors acknowledge the support of the CDP C$^2$EMPI, as well as of the French State under the France-2030 program, the University of Lille, the Initiative of Excellence of the University of Lille, and the European Metropolis of Lille for their funding and support of the R-CDP-24-004-C2EMPI project. The third author also acknowledges the support of the CNRS}


\usepackage[notransparent]{svg} 

\usepackage{amssymb,amsmath,amsthm,amsrefs,newlfont,enumerate,array,graphicx,tikz,tkz-tab,xcolor,bbold,yhmath}

\usepackage[colorlinks, citecolor=blue, linkcolor=red, pdfstartview=FitB]{hyperref}

\usepackage[noabbrev]{cleveref} 
\usepackage[shortlabels]{enumitem}
\makeatletter
\newcommand{\customlabel}[2]{\def\@currentlabel{#2}\label{#1}}
\makeatother

\usepackage{float}
\usepackage{subcaption}

\numberwithin{equation}{section}

\newtheorem{theorem}{Theorem}[section]
\newtheorem*{theorem*}{Theorem}
\newtheorem{lemma}[theorem]{Lemma}
\newtheorem{corollary}[theorem]{Corollary}
\newtheorem{fact}[theorem]{Fact}
\newtheorem{proposition}[theorem]{Proposition}

\newtheorem{question}[theorem]{Question}

\theoremstyle{definition}
\newtheorem{definition}[theorem]{Definition}
\newtheorem{remark}[theorem]{Remark}
\newtheorem{example}[theorem]{Example}
\newtheorem{notation}[theorem]{Notation}
 
\DeclareMathOperator{\dist}{dist}
\DeclareMathOperator{\w}{wind}
\DeclareMathOperator{\spa}{span}
\DeclareMathOperator{\Ran}{Ran}
\newcommand*{\ordot}[1]{\ifx\hfuzz#1\hfuzz \,\cdot\,\else#1\fi}
\newcommand*{\dual}[3][]{#1\langle{\ordot{#2}}\kern1pt #1|\kern1pt {\ordot{#3}}#1\rangle_{p,q}}
\newcommand*{\dualg}[3][]{#1\left[{\ordot{#2}}\kern1pt #1|\kern1pt {\ordot{#3}}#1\right]_{p,q}}

\def\D{\ensuremath{\mathbb D}}

\def\T{\ensuremath{\mathbb T}}
\def\R{\ensuremath{\mathbb R}}
\def\Z{\ensuremath{\mathbb Z}}

\def\C{\ensuremath{\mathbb C}}

\let\Re\relax\DeclareMathOperator{\Re}{Re}

\begin{document}

\maketitle

\begin{abstract}
This paper is devoted to the study of the dynamics of Toeplitz operators $T_F$ with smooth symbols $F$ on the Hardy spaces of the unit disk $H^p$, $p>1$. Building on a model theory for Toeplitz operators on $H^2$ developed by Yakubovich in the 90's, we carry out an in-depth study of hypercyclicity properties of such operators. Under some rather general smoothness assumptions on the symbol, we provide some necessary/sufficient/necessary and sufficient conditions for $T_F$ to be hypercyclic on $H^p$. In particular, we extend previous results on the subject by Baranov-Lishanskii and Abakumov-Baranov-Charpentier-Lishanskii. We also study some other dynamical properties for this class of operators.
\end{abstract}

\tableofcontents

\section{Introduction and main results}\label{Section-intro}
Our aim in this paper is to investigate Toeplitz operators on the Hardy space $H^p$ of the open unit disk $\mathbb{D}$, where $1<p<+\infty$, from the point of view of linear dynamics.

Given a function $F\in L^\infty(\mathbb T)$, where $\mathbb T$ denotes the unit circle in $\mathbb C$, the Toeplitz operator $T_F$ on $H^p$ is defined by
\[T_Fu=P_+(Fu)\quad\text{for every}~u\in H^p.\]
Here $P_+$ denotes the canonical (Riesz) projection from $L^p(\mathbb T)$ onto $H^p$, and it is well-known that it is bounded when $1<p<+\infty$. Then the Toeplitz operator $T_F$ is bounded on $H^p$, which we denote by $T_F\in\mathcal B(H^p)$. The class of Toeplitz operators plays a prominent role in operator theory, partly because of their numerous applications to various other domains such as complex analysis, theory of orthogonal polynomials, probability theory, information and control theory, mathematical physics, etc. We refer the reader to one of the references \cites{BottcherSilbermann1990,Nikolski-2020,Hartmann-2016,Basor-Gohberg-1994,Brown-Halmos,Nikolski-2002-1,Nikolski-2002-2,Douglas-1972} for an in depth study of Toeplitz operators from various points of view. See also Appendix \ref{Section:Rappels} for some useful reminders on Toeplitz operators and the Riesz projection $P_+$.
\par\smallskip
Our focus here will be on the study of the hypercyclicity of Toeplitz operators on $H^p$ spaces, as well as  of related properties such as chaos, frequent hypercyclicity, ergodicity...
Definitions as well as related concepts will be presented in the forthcoming sections and in  \Cref{Section:Rappels}. We refer to books \cite{BayartMatheron2009} and \cite{GrosseErdmannPeris2011} for comprehensive accounts on linear dynamics, together with a thorough presentation of linear dynamical systems, both from the topological and from the measurable point of view. It may be useful to mention that the main difficulty when studying the dynamics of Toeplitz operators is that, in general, explicit formulas for the powers of the Toeplitz operator are not available, except when the symbol is analytic or anti-analytic. 

Given a bounded operator $T$ acting on a (real or complex) separable Banach space $X$, $T$ is said to be \textit{hypercyclic} on $X$ if it admits a vector $x\in X$ with a dense orbit $\{T^nx\,;\,n\ge0\}$. 
Such a vector $x$ is called a hypercyclic vector for $T$. Hypercyclicity is clearly a reinforcement of the classical notion of cyclicity, where one requires the existence of a vector $x\in X$ such that the linear space of its orbit is dense in $X$. 
A hypercyclic operator $T\in \mathcal{B}(X)$ satisfies the following spectral properties: every connected component of its spectrum $\sigma(T)$ intersects the unit circle, and the point spectrum of its adjoint is empty. In particular, it easily follows that a Toeplitz operator with an analytic symbol is never hypercyclic.  
On the other hand, many operators can be shown to be hypercyclic, thanks to arguments based on the study of their eigenvectors. 
Indeed, a well-known criterion due to Godefroy and Shapiro \cite{GodefroyShapiro1991} states the following: 
\par\smallskip
\emph{Suppose that the two subspaces of $X$
\begin{equation}\label{eq:critere-godefroy-shapiro-subspaces}
H_-(T)~=~\overline{\spa}\,[\ker(T-\lambda)\,;\,|\lambda|<1]
\quad\text{and}\quad H_+(T)~=~\overline{\spa}\,[\ker(T-\lambda)\,;\,|\lambda|>1]
\end{equation}
are equal to $X$.
Then $T$ is hypercyclic on $X$.}
\par\smallskip

This criterion will be used repeatedly in the paper to show that for large classes of symbols $F\in L^\infty(\mathbb T)$, the associated Toeplitz operator $T_F$ on $H^p$ is hypercyclic.
\par\smallskip
The study of Toeplitz operators from the point of view of linear dynamics began in the seminal work of Godefroy and Shapiro \cite{GodefroyShapiro1991}, where they gave a necessary and sufficient condition for the adjoint of a multiplication operator on $H^2$ to be hypercyclic, thus characterizing hypercyclic Toeplitz operators with anti-analytic symbols: if $F\in H^\infty$ and $F$ is not constant, then $T_{\overline{F}}\in \mathcal{B}(H^2)$ is hypercyclic if and only if $F(\mathbb D)\cap\mathbb T\neq\varnothing$.
The next step was done by Shkarin, who considered in \cite{Shkarin2012} the case where 
\[F(e^{i\theta})=ae^{-i\theta}+b+c e^{i\theta} \quad\text{for every }e^{i\theta}\in \mathbb T,\] and proved that $T_F$ is hypercyclic on $H^2$ if and only if $|a|>|c|$ and the bounded component of $\mathbb C\setminus F(\mathbb T)$ intersects both $\mathbb D$ and $\mathbb C\setminus \overline{\mathbb D}$. 
Baranov and Lishanskii investigated next in \cite{BaranovLishanskii2016} the more general case where $F$ has the form
\[F(e^{i\theta})=P(e^{-i\theta})+\varphi(e^{i\theta})\quad \text{for a.e. }e^{i\theta}\in\mathbb T,\]
where $P$ is an analytic polynomial and $\varphi$ belongs to $H^\infty$.
They provided some necessary conditions (of a spectral kind) for $T_F$ to be hypercyclic on $H^2$, as well as some sufficient conditions, based on the study of the eigenvectors of $T_F$ and on the Godefroy-Shapiro Criterion. 
The same line of approach was taken in the subsequent work \cite{AbakumovBaranovCharpentierLishanskii2021} of Abakumov, Baranov, Charpentier and Lishanskii, where the authors extended results of \cite{BaranovLishanskii2016} to the case of more general symbols $F$ of the form 
\[F(e^{i\theta})=R(e^{-i\theta})+\varphi(e^{i\theta})\quad \text{for a.e. }e^{i\theta}\in\mathbb T,\]
where $R$ is a rational function without poles in $\overline{\mathbb D}$ and $\varphi$ belongs to $H^\infty$. 
The novel feature of the approach taken in \cite{AbakumovBaranovCharpentierLishanskii2021} is the use of deep results of Solomyak \cite{Solomyak1987} providing necessary and sufficient conditions for finite sets of functions in $H^2$ to be cyclic for some analytic Toeplitz operators on $H^2$. 
Taken in combination with the description of some eigenvectors of $T_F$, this yields substantial extensions of certain results from \cite{BaranovLishanskii2016}.
\par\smallskip

Our work is a further contribution to the study of hypercyclicity of Toeplitz operators. 
Our approach here is somewhat different from the ones taken in the previous works mentioned above, although we will use too properties of eigenvectors and the Godefroy-Shapiro Criterion in order to obtain sufficient conditions for $T_F$ to be hypercyclic. 
So as to exhibit suitable families of spanning eigenvectors for $T_F$, we rely on deep constructions by Yakubovich of model operators for Toeplitz operators on $H^2$ with smooth symbol. It should be mentioned that the model theory for Toeplitz operators with smooth symbols has a rich history. See for instance \cites{Clark-1980-1,Clark-1980-2,Clark-1981,Clark-1987,Clark-1990,Clark-Morrel-1978,Duren-1961,Peller-1986,Wang-1984}. In a series of papers culminating in the works \cite{Yakubovich1991} and \cite{Yakubovich1996} (see also \cites{Yakubovich1989,Yakubovich1993}), Yakubovich showed that for a very large class of positively wound smooth symbols $F$ on $\mathbb T$, the operator $T_F$ is similar to a direct sum of multiplication operators by the independent variable $z$ on certain closed subspaces of Smirnov spaces $E^2(\Omega)$, where the sets $\Omega$ are suitable unions of connected components of $\sigma(T_F)\setminus F(\mathbb T)$.
In the case where the symbol $F$ is negatively wound, in an analogous way $T_F$ is similar to the direct sum of the adjoints of these multiplication operators on these subspaces of the spaces $E^2(\Omega)$. 
In particular, under these conditions, $T_F$ admits an $H^\infty$ functional calculus on the interior of the spectrum $\sigma(T_F)$ of $T_F$. Moreover, using this model, Yakubovich proved in \cite{Yakubovich1991} and \cite{Yakubovich1996} that whenever $F$ is sufficiently smooth and negatively wound, $T_F$ is cyclic on $H^2$.
The proof of this last property relies on the fact that the eigenvectors of $T_F$ associated to eigenvalues $\lambda\in\sigma(T_F)\setminus F(\mathbb T)$ span a dense subspace of $H^2$.
\par\smallskip
It is not difficult to see (see \Cref{Prop:OrientationForHC} below) that, for smooth symbols at least, a necessary condition for the hypercyclicity of $T_F$ is that $F$ be negatively wound. 
We will suppose that it is the case in the rest of this introduction.

Thus, in order to decide whether the Godefroy-Shapiro Criterion can be applied to $T_F$ (proving then hypercyclicity), we need to find conditions on a subset $A$ of $\sigma(T_F)\setminus F(\mathbb T)$ implying that the vector space
\[\spa\,[\ker(T_F-\lambda)\,;\,\lambda\in A]\]
is dense in $H^p$.
In other words, given any $x\in H^q$ (where $q$ is the conjugate exponent of $p$), we need to be able to decide whether $x=0$ as soon as it vanishes on the eigenspaces $\ker(T_F-\lambda)$ for all $\lambda\in A$. 
An easy observation, based on the analyticity of the eigenvector fields and the Uniqueness Principle for analytic functions (see \Cref{Prop:EigenvectorDenseYaku} for details), shows that if $x$ vanishes on $\ker(T_F-\lambda)$ for all $\lambda\in A\cap\Omega$, then $x$ vanishes on $\ker(T_F-\lambda)$ for all $\lambda\in \Omega$, where $\Omega$ is any connected component of $\sigma(T_F)\setminus F(\mathbb T)$ such that $A\cap \Omega$ has an accumulation point in $\Omega$. Recall that we do know (as mentioned above, it is a consequence of results from  \cite{Yakubovich1991} and \cite{Yakubovich1996}) that if for every connected component $\Omega$ of $\sigma(T_F)\setminus F(\mathbb T)$ and for every $\lambda\in \Omega$, $x$ vanishes on $\ker(T_F-\lambda)$, then $x=0$. So things boil down to the following kind of question, which we state here informally:
\begin{question}\label{Question1}
Let $\Omega$ be a connected component of $\sigma(T_F)\setminus F(\mathbb T)$ and let $x\in H^q$ be such that 
\[x\text{ vanishes on }\ker(T_F-\lambda)~\text{for every}~\lambda\in\Omega.\]
Let $\Omega'$ be another connected component of $\sigma(T_F)\setminus  F(\mathbb T)$. Under which conditions on $\Omega$ and $\Omega'$ is it true that \[x \text{ vanishes on } \ker(T_F-\lambda)~\text{for every}~\lambda\in\Omega'?\]
\end{question}
We will investigate this question in some depth, finding rather general conditions on pairs $(\Omega,\Omega')$ of connected components of $\sigma(T_F)\setminus F(\mathbb T)$ implying an affirmative answer to \Cref{Question1}. 
This will allow us to derive necessary and sufficient geometrical conditions implying the hypercyclicity of $T_F$ under fairly general conditions on the (negatively wound) symbol $F$. 
Moreover, this approach allows us to treat the case where $T_F$ acts on $H^p,~p>1$ rather than  the Hilbertian case $p=2$ only.
\par\smallskip
Our conditions are stated essentially in terms of smoothness of the symbol $F$, which is at least required to belong to a class $C^{1+\varepsilon}(\T)$, for some suitable $\varepsilon >0$ (see \Cref{Section-Riesz} for the definition of $C^{1+\varepsilon}(\T)$) and of geometrical properties of the set of connected components of $\sigma(T_F)\setminus F(\mathbb T)$. Rather than state here formally some of our results, which would require some preparation and notation, we prefer to present in an informal way some cases where we are able to characterize hypercyclicity of $T_F$, along with a picture of a situation where we are in the case considered and $T_F$ is hypercyclic.
\par\medskip
Let $p>1$ and $q$ its conjugate exponent, i.e. $\frac{1}{p}+\frac{1}{q}=1$. Consider the following conditions on the symbol $F$:
\par\smallskip
\begin{enumerate}[(H1)]
    \item\label{H1}$F$ belongs to the class $C^{1+\varepsilon}(\mathbb T)$ for some $\varepsilon>\max(1/p,1/q)$, and its derivative $F'$ does not vanish on $\mathbb T$;
    \par\smallskip
   \item\label{H2} the curve $F(\mathbb T)$ self-intersects a finite number of times, i.e. there exist real numbers $\theta_0<\theta_1<\dots<\theta_m=\theta_0+2\pi$ such that, letting $\alpha_j$ be the open arc $\alpha_j=\{e^{i\theta};\,\theta_j<\theta<\theta_{j+1}\}$, we have
     \begin{enumerate}[(a)]
       \item $F$ is injective on  each arc $\alpha_j,~0\le j \le m-1$;
        \item for every $i\neq j,~0\le i,j\le m-1$, the sets $F(\alpha_i)$ and $F(\alpha_j)$ are disjoint;
    \end{enumerate}
    \par\smallskip
    \item\label{H3}for every $\lambda\in\mathbb C\setminus F(\mathbb T)$, $\w_F(\lambda)\le0$, where $\w_F(\lambda)$ denotes the winding number of the curve $F(\mathbb T)$ around $\lambda$;
    \par\smallskip
    \item\label{H4}for every $0\le j\le m-1$, $F_{|\overline{\alpha}_j}$ admits an analytic extension to a neighborhood of the closed subarc $\overline{\alpha}_j=\{e^{i\theta};\,\theta_j\le\theta\le\theta_{j+1}\}$, where  $\theta_0,\dots,\theta_m$ are given by \ref{H2}. 
 \end{enumerate}
             \par\medskip
Under the assumption \ref{H2}, let $\mathcal O$ be the set of the self-intersection points of $F(\mathbb T)$ (i.e. the set of points that have at least two preimages under $F$). Note that there exists a collection of subarcs $\alpha_0,\dots,\alpha_{m-1}$ given by \ref{H2} such that $\mathcal O$ is exactly the set of the extremities of all the arcs $F(\alpha_j)$. Then the function $F$ admits an inverse $F^{-1}$ defined on $F(\mathbb T)\setminus \mathcal O$, and we denote by $\zeta$ the function $1/F^{-1}$ on $F(\mathbb T)\setminus \mathcal O$. Note also that under the assumption \ref{H4}, given two adjacent arcs $\alpha_j$ and $\alpha_k$, the corresponding analytic extensions of $F_{|\overline{\alpha}_j}$ and $F_{|\overline{\alpha}_k}$ may differ. 
 \par\smallskip

The smoothness and injectivity conditions \ref{H1} and \ref{H2} are those from \cite{Yakubovich1991} in the case $p=2$.
A more general condition than \ref{H2}, allowing the sets $F(\alpha_i)$ and $F(\alpha_j)$ to coincide for some indices $i\not = j$, is provided
in \cite{Yakubovich1996}. Some of our results also hold under this weaker assumption, but for simplicity's sake, we begin by restricting ourselves here to the case where \ref{H2} holds, i.e.  to the case where the curve $F(\mathbb T)$ has only a finite number of self-intersection points. The weaker assumptions of \cite{Yakubovich1996} will be considered in \Cref{section:JFA}.
\par\medskip
Let $F$ satisfies the first three assumptions \ref{H1}, \ref{H2} and \ref{H3}. 
Recall that in this case
\[\sigma(T_F)~=~\{\lambda\in\mathbb C\setminus F(\mathbb T)\,;\,\w_F(\lambda)<0\}\cup F(\mathbb T).\]
See \Cref{Lemma:spectral-properties-toeplitz-operators}. We denote by $\mathcal C$ the set of all connected components of $\sigma(T_F)\setminus F(\mathbb T)$.
For every $\Omega\in\mathcal C$, let $\w_F(\Omega)$ be the common value of $\w_F(\lambda),~\lambda\in \Omega.$
We will prove (see \Cref{Th:CNforHC-Intersection}) that a necessary condition for $T_F$, acting on $H^p$, to be hypercyclic, is that $\mathbb T$ intersects every connected component of the interior
of $\sigma(T_F)$.
Here are now some examples of necessary/sufficient/necessary and sufficient conditions that we obtain for the hypercyclicity of $T_F$ on $H^p$:
\par\smallskip
\begin{enumerate}
\item[\textbf{Case 1:}] If $\Omega\cap\mathbb T\neq\varnothing$ for \emph{every} $\Omega\in\mathcal C$, then $T_F$ is hypercyclic on $H^p$ (\Cref{Th:SimpleConseqHc}). See \Cref{Fig1}.
\begin{figure}[ht]
\includegraphics[page=1,scale=.9]{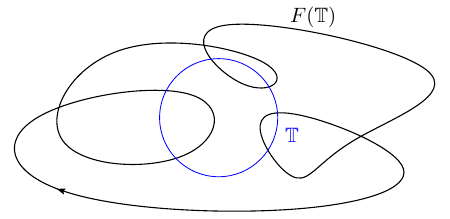}
            \caption{}
            \label{Fig1}
    \end{figure}
\item[\textbf{Case 2:}] If $\mathcal C$ has only one element, i.e. $\sigma(T_F)\setminus F(\mathbb{T})$ is connected,  then $T_F$ is hypercyclic on $H^p$ if and only if $\overset{\circ}{\sigma(T_F)}\cap \mathbb T\neq\varnothing$ (\Cref{Th:ConnexeHC}). See \Cref{Fig2}.
\begin{figure}[ht]
\includegraphics[page=2,scale=.9]{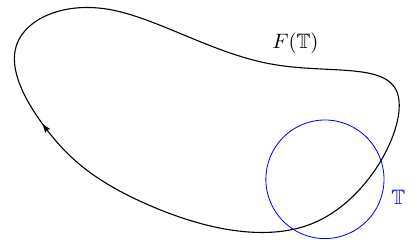}
        \caption{}
        \label{Fig2}
    \end{figure}
\item[\textbf{Case 3:}] If $\overline{\Omega}\cap\overline{\Omega'}$ is finite for all distinct elements $\Omega,\Omega'$ of $\mathcal C$, then $T_F$ is hypercyclic on $H^p$ if and only if $\Omega\cap\mathbb T\neq\varnothing$ for every $\Omega\in\mathcal C$ (\Cref{Th:CNS-IndMax1}). See \Cref{Fig3}.
\begin{figure}[ht]
\includegraphics[page=3,scale=.9]{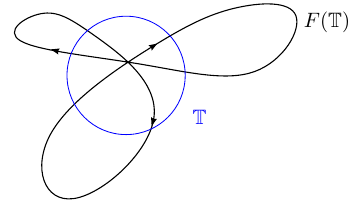}
    \caption{}
    \label{Fig3}
    \end{figure}
    \item[\textbf{Case 4:}] Say that $\Omega\in\mathcal C$ is a maximal component of $\sigma(T_F)\setminus F(\mathbb T)$ if the following holds: for every $\Omega'\in\mathcal C$ which is adjacent to $\Omega$, $|\w_F(\Omega')|<|\w_F(\Omega)|$.\\
    If $\mathbb T$ intersects all the maximal components of $\sigma(T_F)\setminus F(\mathbb T)$, then $T_F$ is hypercyclic on $H^p$ (\Cref{Th:HcInstersCompMax}). See \Cref{Fig4}.
\begin{figure}[ht]
\includegraphics[page=4,scale=.85]{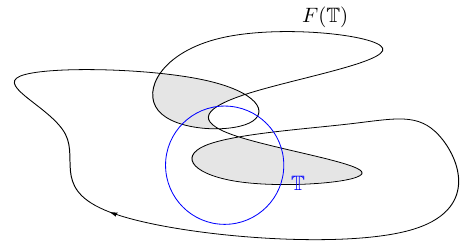}
            \caption{}
            \label{Fig4}
    \end{figure}
    \item[\textbf{Case 5:}] Suppose that $F$ satisfies the additional assumption \ref{H4} and that the following holds:
for every pair $(\Omega,\Omega')$ of adjacent connected components of $\sigma(T_F)\setminus F(\mathbb T)$ with $|\w_F(\Omega)|<|\w_F(\Omega')|$, one of the following properties holds:
\begin{enumerate}[(a)]
\item  there exists a self intersection point $\lambda_0$ of the curve $F(\mathbb T)$ and an open neighborhood $V$ of $\lambda_0$ such that: $V\cap\mathcal{O}=\{\lambda_0\}$, $V\cap\partial \Omega'=V\cap\partial \Omega\cap \partial\Omega'$,
$\lambda_0$ belongs to $\partial\Omega\cap\partial\Omega'\cap V$, and the restriction of the function $\zeta$ to $(\partial\Omega\cap\partial\Omega'\cap V)\setminus\{\lambda_0\}$ cannot be extended continuously at the point $\lambda_0$;
\item $\partial \Omega'$ is a Jordan curve and there exists a closed arc $\Gamma$ of $\mathbb T$ such that $\partial \Omega'=F(\Gamma)$ and $F$ has an analytic extension to a neighborhood of $\Gamma$.
\end{enumerate}
    Then $T_F$ is hypercyclic on $H^p$ if and only if $\mathbb T$ intersects every component of $\overset{\circ}{\sigma(T_F)}$ (\Cref{Th:CNShc2}). See \Cref{Fig5}.
\begin{figure}[ht]
\includegraphics[page=5]{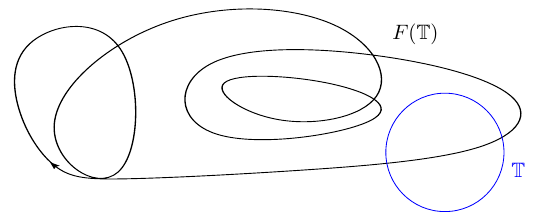}
            \caption{}
            \label{Fig5}
    \end{figure}
\end{enumerate}
Our results allow us to retrieve and generalize previous results from \cite{Shkarin2012}, \cite{BaranovLishanskii2016} and \cite{AbakumovBaranovCharpentierLishanskii2021} on hypercyclicity of Toeplitz operators on $H^2$, up to one important point: while our results apply to general smooth functions $F$ on $\mathbb T$, those of  \cite{AbakumovBaranovCharpentierLishanskii2021} apply to symbols of the form $F(e^{i\theta})=R(e^{-i\theta})+\varphi(e^{i\theta})$ where $R$ is a rational function without poles in $\overline{\mathbb D}$ and $\varphi$ belongs to the disk algebra $A(\mathbb D)$, the space of functions which are analytic on $\mathbb D$ and admit a continuous extension to $\overline{\mathbb D}$. 
While the restriction of $R$ to $\mathbb T$ is as smooth as we wish, that of $\varphi$ is only continuous, so that the restriction of $F$ to $\mathbb T$ may be only continuous. See \Cref{Section:ABCL} for details.
\par\smallskip
The paper is organized as follows:  in \Cref{Section:Yakubovich-court}, we first present Yakubovich's results on models for Toeplitz operators with negatively wound symbols on $H^2$, and state their extension to the $H^p$ setting that  will be required in the sequel of the paper. In \Cref{Section:CNforHc}, we present some necessary conditions for the hypercyclicity of Toeplitz operators on $H^p$. \Cref{Section:PassageFrontiere,Section:ApplicationCShc} are devoted to the proofs of our main results, providing sufficient as well as necessary and sufficient conditions for the hypercyclicity of these operators. The more general case where the symbol is not necessarily required to be injective on the circle minus a finite number of points (i.e. the setting of \cite{Yakubovich1996}) is discussed then in \Cref{section:JFA}.
In \Cref{Section:ABCL}, we compare our results and approach to those taken by Abakumov, Baranov, Charpentier and Lishanskii in \cite{AbakumovBaranovCharpentierLishanskii2021}. 
\Cref{Section:AutresResult+Questions} contains some generalizations and further results, as well as  some open questions arising from our study.
In particular, we investigate in \Cref{Section:AutresResult+Questions} other notions in linear dynamics such as supercyclicity, chaos, frequent hypercyclicity and ergodicity with respect to an invariant ergodic measure in the setting of Toeplitz operators 

The paper also contains two appendices: the first one, \Cref{Section:Rappels}, presents some reminders on topics and tools which are repeatedly used in the paper: basics on the Fredholm theory of Toeplitz operators, Carleson measure, Smirnov spaces, quasiconformal maps, Privalov's type theorems, and, lastly, linear dynamics. 
The second one, \Cref{Section:Yakubovich_demoHp}, contains a full proof of the extension to the $H^p$ case of the main results from \cite{Yakubovich1991} on Toeplitz operators.
Our approach follows closely that of \cite{Yakubovich1991}, and we do not make any claim for originality here. 
However, since these results are extremely beautiful, but not so easy to read in their original presentation, we thought it worthwhile to provide a detailed account of them, and of their extension to the $H^p$ setting.

\par\smallskip
\textbf{Acknowledgments.} We are grateful to Dmitry Yakubovich for stimulating discussions on the topic of this paper, and in particular for pointing out to us the simplification in the proof of \Cref{Prop:Clef_Ext-Int} given in \Cref{Remark-Yakubovich}. We also thank Dmitry Khavinson for providing us with useful references on Smirnov spaces. Finally, we warmly thank the referee for his/her tremendous work and careful reading, which helped us improve the first version of our manuscript and to correct a few inaccuracies.

\section{Model theory for Toeplitz operators on $H^p$ spaces with smooth symbols}\label{Section:Yakubovich-court}
In this section we present the notation, setting and main results from \cite{Yakubovich1991} and \cite{Yakubovich1996} (in the $H^2$-case) which are crucial to our approach to hypercyclicity properties of Toeplitz operators on $H^p,p>1$. We will restrict ourselves here to the bare essentials, see \Cref{Section:Yakubovich_demoHp} for further details, explanations and proofs.

For $1<p<+\infty$, let $H^p$ denote the Hardy space of all analytic functions $u$ on the open unit disk $\mathbb D$ such that 
\[\|u\|_{H^p}~:=~\sup_{0<r<1}M_p(u,r)<+\infty\quad\text{where}~M_p(u,r)=\left(\int_0^{2\pi}|u(re^{i\theta})|^p\frac{\mathrm d\theta}{2\pi}\right)^{1/p}.\]
A function $u\in H^p$ has non tangential boundary values  $u^*$ almost everywhere on $\mathbb T$. We will often still denote this boundary value as $u$. It is well-known that $\|u\|_{H^p}=\|u^*\|_{L^p(\mathbb T)}$. The dual of $H^p$ is canonically identified to $H^q$, where $q$ is the conjugate exponent of $p$ and the duality is given by the following formula:
\begin{equation}\label{eq:duality-defn}
\dual{x}{y}~=~\frac{1}{2\pi}\int_0^{2\pi}x(e^{i\theta})y(e^{-i\theta})\,\mathrm{d}\theta,
\end{equation}
where $x\in H^p$ and $y\in H^q$.
This duality bracket is linear on both sides; we keep this convention even when $p=2$, so that in particular adjoints of operators on $H^2$  must be understood as Banach space adjoints, and not Hilbert space adjoints.
\par\smallskip
We denote by $P_+$ the Riesz projection from $L^p(\mathbb T)$ into $H^p$ defined by 
\[P_+u(z)~=~\frac{1}{2i\pi}\int_{\mathbb T}\frac{u(\tau)}{\tau-z}\,\mathrm d\tau\quad\text{for }u\in L^p(\mathbb T) \,\text{ and } z\in\mathbb D.\]

\subsection{Some spectral theory}
Given $F\in L^\infty(\mathbb T)$, the Toeplitz operator $T_F$ on $H^p$ defined by $T_Fu=P_+(Fu),~u\in H^p$, is bounded on $H^p$.  With the duality between $H^p$ and $H^q$ given by \eqref{eq:duality-defn}, we have $T_F^*=T_f$ where the symbol $f$ is defined as $f(z)=F(1/z)$, $z\in\T$ (see \Cref{lem:adjoint-banachique}). Suppose that $F$ is continuous on $\mathbb T$. Then $T_F-\lambda$ is a Fredholm operator on $H^p$ if and only if $\lambda\notin F(\mathbb T)$; in this case, its Fredholm index is equal to $-\w_F(\lambda)$ and we can describe the spectrum of $T_F$ as 
\[\sigma(T_F)~=~\{\lambda\in\mathbb C\setminus F(\mathbb T)\,;\,\w_F(\lambda)\neq0\}\cup F(\mathbb T),\]
where $\w_F(\lambda)$ is the winding number of the curve $F(\mathbb T)$ with respect to the point $\lambda$. See \Cref{Lemma:spectral-properties-toeplitz-operators}. 
Note that if $F$ is at least $C^1$ on $\mathbb T$ (which is the case in the whole of this paper),  the winding number is given by
\[\w_F(\lambda)~=~\frac1{2i\pi}\int_{\mathbb T}\frac{F'(\tau)}{F(\tau)-\lambda}\,\mathrm d\tau~=~\frac1{2\pi}\int_0^{2\pi}\frac{F'(e^{i\theta})}{F(e^{i\theta})-\lambda}e^{i\theta}\,\mathrm d\theta.\]
By continuity of the winding number, the function $\w_F$ is constant on each connected component $\Omega$ of $\sigma(T_F)\setminus F(\mathbb T)$. We denote this common value by $\w_F(\Omega)$.

\subsection{Eigenvectors} Suppose now that $\w_F$ takes only non positive values on $\mathbb C\setminus F(\mathbb T)$, i.e. that assumption \ref{H3} from the introduction is satisfied. It then follows from the $H^p$ version of the Coburn Theorem (\Cref{Th:Coburn}) and \Cref{Lemma:spectral-properties-toeplitz-operators}  that for every $\lambda\in\mathbb C\setminus F(\mathbb T)$,
\[\ker(T_F^*-\lambda)~=~\{0\}\quad\text{and}\quad\dim\,\ker(T_F-\lambda)~=~-\w_F(\lambda).\]
\par\smallskip
Henceforward, we assume that $F$ satisfies assumptions \ref{H1}, \ref{H2} and \ref{H3}. In order to represent $T_F$ as the adjoint of a multiplication operator by $\lambda$ on a certain space of analytic functions, Yakubovich provided in \cite{Yakubovich1991} an explicit expression of a spanning family of elements of $\ker(T_F-\lambda),~\lambda\notin F(\mathbb T)$. Let $\lambda\in\mathbb C\setminus F(\mathbb T)$, and consider the function $\phi_\lambda$ defined on $\mathbb T$ by 
\[\phi_\lambda(\tau)~=~\tau^{-\w_F(\lambda)}(F(\tau)-\lambda)\quad\text{for }\tau\in \mathbb T.\]
Since $\phi_\lambda$ is of class $C^{1+\varepsilon}$ and does not vanish on $\mathbb T$, and since $\w_{\phi_\lambda}(0)=0$, one can define a logarithm $\log\phi_\lambda$ of $\phi_\lambda$ on $\mathbb T$ that is $C^{1+\varepsilon}$ on $\mathbb T$, and set
\begin{equation}\label{eq:vecteur-propre-Yaku33434ZD}
F_\lambda^+~=~\exp(P_+(\log \phi_\lambda)).
\end{equation}
More details on the construction on $F_\lambda^+$ are contained in {\Cref{Section:Yakubovich_demoHp}}.
The functions $F_\lambda^+$ and $1/F_\lambda^+$ both belong to $A(\mathbb D)$ and for every connected component $\Omega$ of $\mathbb C\setminus F(\mathbb T)$ and  every $z\in \overline{\mathbb D}$, the map $\lambda\mapsto F_\lambda^+(z)$ is analytic on $\Omega$ and continuous on $\overline{\Omega}$ (see \Cref{Lem:analytic-en-dehors-de-la-courbe} and \Cref{Lemme-propriete-varphi} for details).
\par\smallskip
Let $N~=~\max\{|\w_F(\lambda)|\,;\,\lambda\notin F(\mathbb T)\}$ and for each $j=0,\dots,N-1$, let
\begin{equation}\label{eq:def-Omegajplus}
\Omega_j^+~=~\{\lambda\notin F(\mathbb T)\,;\,|\w_F(\lambda)|>j\}
\end{equation}
be the set of points $\lambda$ in $\mathbb C\setminus F(\mathbb T)$ where the (negative) winding number of $F$ is strictly less than $-j$. We have
\[\Omega_{N-1}^+  \subseteq\Omega_{N-2}^+  \subseteq\dots  \subseteq \Omega_{1}^+  \subseteq\Omega_{0}^+,\]
and $\sigma(T_F)\setminus F(\mathbb T)$ (which is the set of all $\lambda$'s in $\sigma(T_F)$ such that $T_F-\lambda$ is Fredholm) coincides with $\Omega_0^+$. Note also that 
\[\partial\Omega_{N-1}^+  \subseteq\partial\Omega_{N-2}^+  \subseteq\dots  \subseteq \partial\Omega_{1}^+  \subseteq\partial\Omega_{0}^+=F(\mathbb T).\]
For every $j=0,\dots,N-1$ and every $\lambda\in\Omega_j^+$, set 
\begin{equation}\label{eve}
  h_{\lambda,j}(z)~=~z^j\frac{F_\lambda^+(0)}{F_\lambda^+(z)}\quad\text{for any }z\in\mathbb D.  
\end{equation}
These functions belong to $A(\mathbb D)$, hence to $H^p$, and it can be checked that $(T_F-\lambda)h_{\lambda,j}=0$ for every $\lambda\in\Omega_j^+$. So, we have
\[\ker(T_F-\lambda)~=~\spa\big[h_{\lambda,j}\,;\,0\le j<|\w_F(\lambda)|\big]\]
for every $\lambda\in\sigma(T_F)\setminus F(\mathbb T)$ (see \Cref{lem:description-eigenvectors}). Moreover, the map $\lambda\longmapsto h_{\lambda,j}$ is analytic from $\Omega_j^+$ into $A(\mathbb D)$ and, for a fixed $z\in\mathbb D$, the function $\lambda\longmapsto h_{\lambda,j}(z)$ belongs to the Smirnov space $E^q(\Omega_j^+)$, where $q$ is the conjugate exponent of $p$
(see \Cref{lemme-Uborne}). 
\par\smallskip
The reader may wonder why we normalize the eigenvectors in \Cref{eve} by the constant $F_\lambda^+(0)$. In fact, in the construction of the model, we would like that the isomorphism $U$ which implements the model satisfies that $(U1)_0\equiv 1$ (see \Cref{defn-model-U-section2-5}) and, since $(U1)_0(\lambda)=h_{\lambda,0}(0)$, this imposes this constant. With this property for $U$, it turns out that one can obtain an explicit formula for the inverse of $U$. See \cite{Yakubovich1991} for more explanations on this point. 

\subsection{Interior and exterior boundary values}\label{Subsection:IntExt}
\begin{figure}[ht]
\includegraphics[page=6,scale=1]{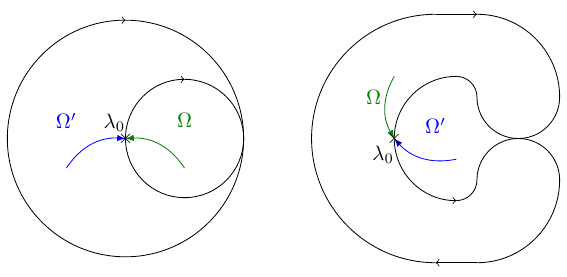}
        \caption{}\label{Fig2.6}
    \end{figure}

Assumption \ref{H2} implies that the curve $F(\mathbb T)$ has only a finite number of points of self-intersection which is the set $\cal O$. Whenever $\gamma$ is a subarc of $F(\mathbb T)$ containing no point of $\cal O$, $\gamma$ is included in the boundary of exactly two connected components $\Omega$ and $\Omega'$ of $\mathbb C\setminus F(\mathbb T)$, and we have
\[|\w_F(\Omega)-\w_F(\Omega')|=1.\]
If $|\w_F(\Omega)|>|\w_F(\Omega')|$, we say that $\Omega$ is the \textit{interior component}  and $\Omega'$ is the \textit{exterior component} (with respect to $\gamma$). Now, let $u$ be continuous on $\mathbb{C}\setminus F(\mathbb T)$. For $\lambda_0\in F(\mathbb T)\setminus\mathcal O$, let $I_{\lambda_0}$ be a small subarc of $F(\mathbb T)\setminus\mathcal O$ containing $\lambda_0$. We define (when they exist) the following non-tangential limits, called respectively the interior and exterior boundary values of $u$ at $\lambda_0$:
\[
u^{int}(\lambda_0)~=~\lim_{\substack{\lambda\to\lambda_0\\
\lambda\in \Omega}}u(\lambda)\quad\text{and}\quad u^{ext}(\lambda_0)~=~\lim_{\substack{\lambda\to\lambda_0\\
\lambda\in \Omega'}}u(\lambda),
\]
where $\Omega$ and $\Omega'$ are respectively the interior and exterior connected components with respect to $I_{\lambda_0}$. Functions in a Smirnov space of a domain $\Omega$ of $\C$ (with a rectifiable boundary $\Gamma=\partial\Omega$) admit non-tangential limits almost everywhere on $\Gamma$ (see \Cref{smirnov}). If $\Omega$ and $\Omega'$ are two adjacent domains along an arc $\gamma$, and if $u$ belongs to some Smirnov space $E^r(\Omega\cup\Omega')$, then the interior and exterior boundary values $u^{int}$ and $u^{ext}$ exist almost everywhere on $\gamma$ (see \Cref{Fig2.6}).

With these notations, suppose that $F$ satisfies \ref{H1}, \ref{H2} and \ref{H3}. Since $F$ is invertible from $\mathbb T\setminus F^{-1}(\cal O)$ onto $F(\mathbb T)\setminus\cal O$, one can define the map
$\zeta=1/F^{-1}$ on $F(\mathbb T)\setminus\cal O$. In particular, $\zeta$ is defined almost everywhere on $F(\mathbb T)$. 
Then (see \Cref{cor:condition-aux-bords-hlambdaj}) the eigenvectors $h_{\lambda,j}$ of $T_F$ given by (\ref{eve}) satisfy the following relation: for every $z\in\mathbb D$, we have
\begin{equation}
    h_{\lambda,j}^{int}(z)-\zeta(\lambda) h_{\lambda,j+1}^{int}(z)=h_{\lambda,j}^{ext}(z)~\text{ for almost all $\lambda\in \partial \Omega_{j+1}^+$}.
    \label{Eq:RelationEigenvectors}
\end{equation}

\subsection{Boundary conditions in some Smirnov spaces}
Let $r>1$. On the direct sum $\bigoplus_{0\le j\le N-1} E^r(\Omega_j^+)$,  we consider the following norm:
\[\left\|(u_j)_{0\le j\le N-1}\right\|~=~\left(\sum_{j=0}^{N-1}\|u_j\|_{E^r(\Omega_j^+)}^r\right)^{1/r}\quad\text{for every }(u_j)_{0\le j\le N-1}\in \bigoplus_{j=0}^{N-1}E^r(\Omega_j^+).\]
Note that for each connected component $\Omega$ of $\Omega_j^+$, the non-tangential limit of the function $u_j$ exists almost everywhere on the boundary of $\Omega$ and that for every $0\le j\le N-1$, the non tangential limits $u_j^{int},u_{j+1}^{int}$ and $u_j^{ext}$ exist almost everywhere on $\partial \Omega_{j+1}^+$. 
\par\smallskip
Let $E^r_F$ be the closed subspace of $\bigoplus E^r(\Omega_j^+)$ formed by the $N$-tuples  $(u_j)_{0\le j\le N-1}$ in $\bigoplus E^r(\Omega_j^+)$ satisfying, for all $0\le j< N-1$,
\begin{equation}\label{eq:boundary-conditions-ErF}
u_j^{int}-\zeta u_{j+1}^{int}~=~u_j^{ext}\quad\text{a.e. on}~\partial \Omega_{j+1}^+.
\end{equation}
Note that \Cref{Eq:RelationEigenvectors} implies that for all $z\in\mathbb D$, the $N$-tuple $(h_{\,\cdot,0}(z),\dots,h_{\,\cdot,N-1}(z))$ belongs to $E^r_F$.
Note also 
that this subspace is an invariant subspace for the multiplication operator $M_\lambda:\bigoplus E^r(\Omega_j^+)\to \bigoplus E^r(\Omega_j^+)$ defined by
\[(M_\lambda u)_j~=~\lambda u_j(\lambda)\quad\text{for every }u=(u_i)_{0\le N-1}\in \bigoplus_{j=0}^{N-1} E^r(\Omega_j^+).\]
More generally let  $v\in H^\infty(\overset{\circ}{\sigma(T_F)})$ be a bounded analytic function on the interior 
of the spectrum of $T_F$. Then the space $E^r_F$ is invariant by the multiplication operator $M_{v}$ on $\bigoplus
E^r(\Omega_j^+)$ defined by $M_v(u_j)_{0\le j\le N-1}=(vu_j)_{0\le j\le N-1}$, where $u_j\in E^r(\Omega_j^+)$ for every ${0\le j\le N-1}$.
 This means that $H^\infty(\overset{\circ}{\sigma(T_F)})$ is contained in the multiplier algebra of $E^r_F$.
\subsection{The main result of Yakubovich}\label{subsection-Yaku-Operator-U34ZZE111}
After all this preparation, we are now ready to define the operator $U$ that gives $M_\lambda$ as a model operator for $T_F^*$.
\par\smallskip
Let $p>1$, and let $q$ be the conjugate exponent of $p$ (i.e. $1/p+1/q=1$). Suppose that the symbol $F$ of the Toeplitz operator $T_F\in\cal B(H^p)$ satisfies \ref{H1}, \ref{H2} and \ref{H3}. Let $h_{\lambda,j}$, $0\le j\le N-1$, be given by (\ref{eve}). For every function $g\in H^q$, define $Ug=((Ug)_j)_{0\le j\le N-1}$ by setting
\begin{equation}\label{defn-model-U-section2-5}
(Ug)_j(\lambda)~=~\dual{h_{\lambda,j}}{g}~\textrm{ for every }\lambda\in\Omega_j^+.
\end{equation}

Note that, since $h_{\lambda,j}$ is an eigenvector of $T_F$ associated to the eigenvalue $\lambda$, we have for every $g\in H^q$ and every $0\le j\le N-1$ that
\[\dual{h_{\lambda,j}}{T_F^*g}~=~\dual{T_Fh_{\lambda,j}}{g}~=~\lambda\dual{h_{\lambda,j}}{g}\quad \text{ for every }\lambda\in\Omega_j^+.\]
In other words, \[U(T_F^*g)~=~M_\lambda Ug \quad \text{ for every } g\in H^q.\]
Since the functions $h_{\lambda,j}$, $0\le j\le N-1$, form a basis of the eigenspace $\ker(T_F-\lambda)$, it follows that

\begin{fact}\label{basic-fact}
 For any $g\in H^q$,
$g$ vanishes on $\ker(T_F-\lambda)$ if and only if  $(Ug)_j(\lambda)=0$ for every $0\le j<|\w_F(\lambda)|$.
\end{fact}

We are now ready to give the $H^p$ version of the main result of Yakubovich in \cite{Yakubovich1991}.

\begin{theorem}\label{T:Yakbovich_Hp}
The operator $U$ defined above is a linear isomorphism from $H^q$ onto $E_F^q$. Moreover, we have
\[T_F^*~=~U^{-1}M_\lambda U.\]
In other words, the following diagram commutes:
\begin{center}
\includegraphics[page=7,scale=1.1]{figures.pdf}
\end{center}
\end{theorem}

This model is built from eigenvectors of $T_F$, and the relation $UT_F^*=M_\lambda U$ means nothing else than this. See the introductions of the papers \cites{Yakubovich1991, Yakubovich1993} for insightful comments on this construction.
A detailed proof of \Cref{T:Yakbovich_Hp} will be given in \Cref{Section:Yakubovich_demoHp}. 
\par\smallskip
Note that \Cref{T:Yakbovich_Hp} is stated in \cite{Yakubovich1991} for $p=2$ and for any value $\varepsilon>0$ of the parameter $\varepsilon$ appearing in condition \ref{H1}. Here we prove it only for $\varepsilon >\max(1/p, 1/q)$, so that we have to require that $\epsilon >1/2$ in the case where $p=2$. However, D. Yakubovich commented to us in a private communication that there is an inaccuracy in his statement in \cite{Yakubovich1991} and that he is  also able to prove his result only for $\varepsilon >1/2$.
\par\smallskip
We conclude this section with some important consequences of \Cref{T:Yakbovich_Hp}. 
\subsection{The functional calculus}
Since $H^\infty(\overset{\circ}{\sigma(T_F)})$ is contained in the multiplier algebra of $E^q_F$, a first important consequence of \Cref{T:Yakbovich_Hp} is the existence of an $H^\infty$ functional calculus for $T_F$ on the interior of the spectrum of $T_F$.
\begin{corollary}[Yakubovich \cite{Yakubovich1991}]\label{Cor:HinftyCalculus}
Suppose that $F$ satisfies \ref{H1}, \ref{H2} and \ref{H3}. Then $T_F$ admits an $H^\infty$ functional calculus on $\overset{\circ}{\sigma(T_F)}$, i.e. there exists a constant $C>0$ such that
\[\|u(T_F)\|~\le~ C\sup\{|u(\lambda)|\,;\, \lambda\in \overset{\circ}{\sigma(T_F)}\}\]
for every function $u\in H^\infty(\overset{\circ}{\sigma(T_F)})$.
\end{corollary}

Indeed, let $u\in H^\infty(\overset{\circ}{\sigma(T_F)})=H^\infty(\overset{\circ}{\sigma(T_F^*)})$. By \Cref{T:Yakbovich_Hp}, we can set $u(T_F^*)=U^{-1}M_{u(\lambda)}U$. This defines a functional calculus for $T_F^*$ on 
$H^\infty(\overset{\circ}{\sigma(T_F)})$, and
$$\|u(T_F^*)\|~\le~ \|U\|\, \|U^{-1}\|\, \sup\{|u(\lambda)|\,;\, \lambda\in \overset{\circ}{\sigma(T_F)}\}.$$ See \Cref{model}.
 Set now $u(T_F)=u(T_F^*)^*$ for every $u\in H^\infty(\overset{\circ}{\sigma(T_F)})$; this defines an $H^\infty$ functional calculus for $T_F$ on the interior of ${\sigma(T_F)}$, and \Cref{Cor:HinftyCalculus} follows.

\subsection{Spanning eigenvectors}
The second important consequence of \Cref{T:Yakbovich_Hp} is that eigenvectors of $T_F$ span a dense subset of $H^p$. 
Before we justify this statement, we set a notation which will be used repeatedly throughout this work.

\begin{definition}\label{def-importantes}

For every subset $A$ of $\mathbb C$, set \begin{equation}\label{eq:definition-sep}
H_A(T_F)~:=~\overline{\spa}\,[\ker(T_F-\lambda)\,;\,\lambda\in A].
\end{equation}
\end{definition}
Let  $\mathcal C$ denote the set of all connected components of $\sigma(T_F)\setminus F(\mathbb T)$. 
Here is now a result which will be crucial in the sequel.
\begin{proposition}\label{Prop:EigenvectorDenseYaku}
Let $p>1$ and let $F$ satisfy \ref{H1},\ref{H2} and \ref{H3}.
\begin{enumerate}
    \item Let $\Omega\in\mathcal C$, and let $A$ be a subset of $\mathbb C$. If $A\cap \Omega$ has an accumulation point in $\Omega$, then
    \[\overline{\spa}\,[\ker(T_F-\lambda)\,;\,\lambda\in A\cap\Omega]~=~\overline{\spa}\,[\ker(T_F-\lambda)\,;\,\lambda\in \Omega].\]
    \item  We have
\begin{equation*}
      \overline{\spa}\,[H_\Omega(T_F)\,;\,\Omega\in\cal C]~=~H^p.
  \end{equation*}
\end{enumerate}
\end{proposition}
\begin{proof} Given $\Omega\in\mathcal{C}$, observe that the functions $\lambda\mapsto \dual{h_{\lambda,j}}{g}$, $0\le j<|\w_F(\Omega)|$, are analytic on $\Omega$ (they belong to $E^p(\Omega)$) and span a dense subspace of $H_\Omega(T_F)$.  Then the uniqueness principle for analytic functions yields (1).
In order to prove (2), suppose that $g\in H^p$ vanishes on all eigenspaces $\ker(T_F-\lambda),\lambda\in\Omega,\, \Omega\in\mathcal{C}$. Then we have, for every $0\le j\le N-1$, 
\[\dual{h_{\lambda,j}}{g}~=~0\quad\text{ for every}~\lambda\in\Omega_j^+.\]
This means exactly that $Ug=0$, and since $U$ is injective by \Cref{T:Yakbovich_Hp}, $g=0$. 
\end{proof}

\par\smallskip
This consequence of \Cref{T:Yakbovich_Hp} is used in \cite{Yakubovich1991} to derive the cyclicity of $T_F$ under the assumptions \ref{H1}, \ref{H2} and \ref{H3}. Indeed, since $\ker(T_F-\lambda)^*=\varnothing$ for every $\lambda\in\sigma(T_F)\setminus F(\T)$, it follows from \Cref{Prop:EigenvectorDenseYaku} and
\Cref{lem:cyclicite-sophie-nikolski} that $T_F$ is cyclic.
In order to investigate the hypercyclicity of $T_F$, an important part of our work will be to refine \Cref{Prop:EigenvectorDenseYaku}, to as to be able to conclude (under suitable conditions on the symbol $F$), that
\begin{align*}
H_-(T_F)&~=~\overline{\spa}\,[\ker(T_F-\lambda)\,;\,|\lambda|<1]~=~H^p;\\
H_+(T_F)&~=~\overline{\spa}\,[\ker(T_F-\lambda)\,;\,|\lambda|>1]~=~H^p.
\end{align*}
These two properties will allow us to apply the Godefroy-Shapiro Criterion, and thus to deduce that $T_F$ is hypercyclic on $H^p$. 
This will be the main goal of our \Cref{Section:PassageFrontiere}. Applications will be presented in \Cref{Section:ApplicationCShc}. Meanwhile, we present some necessary conditions (of a spectral nature) for the hypercyclicity of $T_F$.

\section{Necessary conditions for hypercyclicity of Toeplitz operators}\label{Section:CNforHc}
\subsection{A condition on the winding number}
A first observation is that for a sufficiently smooth symbol $F$, $F$ has to be negatively wound for $T_F$ to have a chance of being hypercyclic.
\begin{proposition}\label{Prop:OrientationForHC}
Let $F$ be continuous on $\mathbb T$ and let $p>1$. If there exists a point $\lambda_0\in\mathbb C\setminus F(\mathbb T)$ with $\w_F(\lambda_0)>0$, then $T_F$ is not hypercyclic on $H^p,p>1$.
\end{proposition}
\begin{proof}
Recall (see \Cref{Lemma:spectral-properties-toeplitz-operators}) that, since $\lambda_0\notin F(\mathbb T)$,  $T_F-\lambda_0$ is a Fredholm operator, and its Fredholm index satisfies
\[j(T_F-\lambda_0)~=~\dim(\ker(T_F-\lambda_0))-\dim(\ker(T_F^*-\lambda_0))~=~-\w_F(\lambda_0).\]
Thus \[\dim(\ker(T_F^*-\lambda_0))~=~\w_F(\lambda_0)+\dim(\ker(T_F-\lambda_0))~\ge~ \w_F(\lambda_0)>0.\]
Then $\lambda_0$ is an eigenvalue for $T_F^*$. But the adjoint of a hypercyclic operator must have empty point spectrum, and hence $T_F$ cannot be hypercyclic. 
\end{proof}
Assumption \ref{H3} is thus natural in our context.
\par\smallskip
\subsection{Intersection with connected components of the interior of the spectrum}
It is also well-known (see \Cref{Subsection:LinearDyn}) that whenever $T$ is a hypercyclic operator on a complex separable Banach space $X$, every connected component of the spectrum of $T$ must intersect the unit circle $\mathbb T$. It turns out that, in our current setting, a stronger property must hold.
\begin{theorem}\label{Th:CNforHC-Intersection}
Let $p>1$ and let $F$ satisfy assumption \ref{H1}, \ref{H2} and \ref{H3}. If $T_F$ is hypercyclic on $H^p$, then every connected component of the interior 
of the spectrum of $T_F$ must intersect $\mathbb T$.
\end{theorem}

\begin{figure}[ht]
\includegraphics[page=8,scale=1]{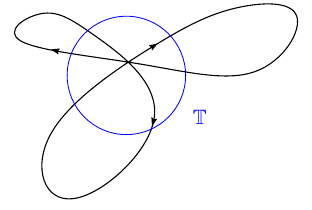}
            \caption{}
            \label{Fig8}
    \end{figure}

For instance, in the situation represented in \Cref{Fig8}, \Cref{Th:CNforHC-Intersection} states that for $T_F$ to be hypercyclic, the unit circle must intersect all the ``petals'' of the flower-like domain $\overset{\circ}{\sigma(T_F)}$ - which is much stronger that requiring that $\sigma(T_F)$, which is connected here, intersects $\mathbb T$. 

    \Cref{Th:CNforHC-Intersection} follows immediately from \Cref{Prop:CNgeneraHC} below, which has some independent interest, and the fact that under assumptions \ref{H1}, \ref{H2} and \ref{H3}, $T_F$ has an $H^\infty(\overset{\circ}{\sigma(T_F)})$-functional calculus (\Cref{Cor:HinftyCalculus}).
    
    \begin{proposition}\label{Prop:CNgeneraHC}
    Let $X$ be a complex separable Banach space, and let $T\in \mathcal{B}(X)$ be hypercyclic.
     Suppose that $T$ admits an $H^\infty$-functional calculus on a certain open subset $O$ of $\mathbb C$, and that all connected components of $O$ intersect $\sigma(T)$.
    Then every connected component of $O$ intersects $\mathbb T$.
    \end{proposition}
    
\begin{proof}
Let $O_1$ be a connected component of $O$, and let $O_2=O\setminus O_1$. Define a function $u_1$ on $O$ by setting $u_1(\lambda)=1$ if $\lambda\in O_1$ and $u_1(\lambda)=0$ if $\lambda\in O_2$. Set $u_2=1-u_1$. Then $u_1$ and $u_2$ are elements of $H^\infty(O)$. Set $P_i=u_i(T),~i=1,2$. The $P_i$'s are well-defined bounded projections on $X$, with $P_1P_2=P_2P_1=0$ and $P_1+P_2=I$. Denote by $M_i$ the range of $P_i,i=1,2$. Then $M_i$ is a closed $T$-invariant subspace of $X$ and $X=M_1\oplus M_2$. If $O$ has at least two connected components, $\sigma(T)$ intersects both $O_1$ and $O_2$, $\sigma(P_i)=\sigma(u_i(T))=u_i(\sigma(T))=\{0,1\}$, and then $P_i$ is a non-trivial projection. So $M_i\neq \{0\}$ and $M_i\neq X$. If $O$ is connected, then $O_1=O$ and $P_1=I$.
Denote by $T_i$ the operator induced by $T$ on $M_i$, $i=1,2$. 
Then
\begin{fact}\label{spectre}
We have
 $\sigma(T_1)=\overline{O_1}\cap \sigma(T)$.
\end{fact}
\begin{proof}[Proof of \Cref{spectre}]Since $\sigma(T)=\sigma(T_1)\cup\sigma(T_2)$, we clearly have  $\sigma(T_1)  \subseteq \sigma(T)$. If $\lambda\notin \overline{O_1}$, define a function $f$ by setting $f(z)=\frac{1}{z-\lambda}u_1(z)$, $z\in O$. Then $f\in H^{\infty}(O)$, and the operator $S:=f(T)$ leaves $M_1$ invariant and satisfies $S(T-\lambda)x_1=(T-\lambda)Sx_1=x_1$ for every $x_1\in M_1$. So $T_1-\lambda$ is invertible on $M_1$, and $\lambda\notin\sigma(T_1)$. Thus $\sigma(T_1)  \subseteq \overline{O_1}\cap \sigma(T)$.
Suppose now that $\lambda\in O_1\cap \sigma(T)$. Then $\lambda\notin\overline{O_2}$ (because $O_1$ is open and $O_1\cap O_2=\varnothing$), and by the argument above $\lambda\notin\sigma(T_2)$. Since $\lambda\in\sigma(T)$, necessarily $\lambda\in \sigma(T_1)$, so that $O_1\cap \sigma(T)  \subseteq\sigma(T_1)$. Hence $\overline{O_1}\cap \sigma(T)  \subseteq\sigma(T_1)$ and this shows that $\sigma(T_1)=\overline{O_1}\cap \sigma(T)$.
\end{proof}
Moreover, since the operator $T$ admits an $H^\infty(O)$-functional calculus, $T_1$ admits an $H^\infty(O_1)$-functional calculus. Let $u\in H^\infty(O_1)$. Then $u$ can be extended into a function $ \widetilde u\in H^\infty(O)$ by setting $ \widetilde u(z)=0$ for $z\in O_2$. For $x_1\in M_1$, define 
\begin{equation}\label{eq:functional-calculus-somme-directe}
u(T_1)x_1~=~ \widetilde{u}(T)x_1.    
\end{equation}
Then $u(T_1)$ is a bounded operator on $M_1$ and we have 
\[\|u(T_1)\|~\leq~ \| \widetilde u(T)\|~\le~ C\| \widetilde u\|_{\infty,O}~=~C \|u\|_{\infty,O_1}\]
for a positive constant $C$ independent of $u$. Then we easily see that \Cref{eq:functional-calculus-somme-directe} defines an $H^\infty(O_1)$-functional calculus for $T_1$. 
Since $T$ is a hypercyclic operator on $X$, the operator $T_1$ is also hypercyclic on $M_1$ (see \Cref{lemme-HC-direct-sum}). 

Let us finally show that it is impossible to have $O_1 \subseteq\overline{\mathbb D}$ nor $O_1 \subseteq\mathbb C\setminus\mathbb D$.
If $O_1 \subseteq\overline{\mathbb D}$, then
\[\|T_1^n\|~\le~ C\|z^n\|_{\infty,O_1}~\le~ C\quad\text{for every}~n\ge0,\]
so that $T_1$ is a power-bounded operator on $M_1$. This stands in contradiction with the hypercyclicity of $T_1$.
If now $O_1 \subseteq\mathbb C\setminus\mathbb D$, then since $\sigma(T_1)=\overline{O_1}\cap\sigma(T_F)$, $T_1$ is invertible, and so $T_1^{-1}$ is also hypercyclic (see \Cref{Subsection:LinearDyn}). Since $0\notin O_i$, by the functional calculus again, we have that 
\[\|T_i^{-n}\|~\le~ C\|z^{-n}\|_{\infty,O_i}~\le~ C\quad \text{for every}~n\ge0,\]
contradicting the hypercyclicity of $T_1^{-1}$.
So $O_1\cap (\mathbb C\setminus\overline{\mathbb D})\neq\varnothing$ and $O_1\cap\mathbb D\neq\varnothing$, i.e. $O_1\cap\mathbb T\neq\varnothing$.
\end{proof}

\section{Proving the hypercyclicity of $T_F$: from a connected component  to another}\label{Section:PassageFrontiere}
We remind the reader that $\mathcal C$ denotes the set of all connected components of $\sigma(T)\setminus F(\mathbb T)$; we also recall that the spaces $H_-(T)$ and $H_+(T)$ are defined in \Cref{eq:critere-godefroy-shapiro-subspaces}, and the spaces $H_\Omega(T)$ for $\Omega\in\mathcal C$ in \Cref{eq:definition-sep}. 
Let us begin with the easiest case where $T_F$ can be shown to be hypercyclic.

\begin{theorem}\label{Th:SimpleConseqHc}
Let $p>1$ and let $F$ satisfy \ref{H1}, \ref{H2} and \ref{H3}. Suppose that  
\begin{equation}\label{C1}
\Omega\cap \mathbb T~\neq~\varnothing \quad\text{for every}~ \Omega\in\cal C.
\end{equation}
Then $T_F$ is hypercyclic on $H^p$.
\end{theorem}
\Cref{Fig8423} illustrates an example of a situation where $\sigma(T_F)\setminus F(\mathbb T)$ has two connected components and \Cref{Th:SimpleConseqHc} applies: 
\begin{figure}[ht]
\includegraphics[page=9,scale=1.1]{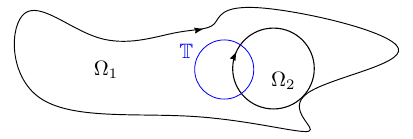}
    \caption{}\label{Fig8423}
\end{figure}
\begin{proof}
    \Cref{Th:SimpleConseqHc} is an immediate consequence of \Cref{Prop:EigenvectorDenseYaku}. Indeed, since $\Omega\cap\mathbb T\neq\varnothing$ for every $\Omega\in\cal C$, $\Omega\cap\mathbb D$ and $\Omega\cap (\mathbb C\setminus\overline{\mathbb D})$ both have accumulation points in $\Omega$ for every $\Omega\in\cal C$ and then, for every $\Omega\in\cal C$,
    \[\overline{\spa}\big[\ker(T_F-\lambda)\,;\,\lambda\in\Omega\cap\mathbb D\big]~=~\overline{\spa}\big[\ker(T_F-\lambda)\,;\,\lambda\in\Omega\cap(\mathbb C\setminus\overline{\mathbb D})\big]=H_\Omega(T_F).\]
    Hence $H_\Omega(T_F)\subseteq H_-(T_F)$ and $H_\Omega(T_F)\subseteq H_+(T_F)$ for every $\Omega\in\cal C$. So by \Cref{Prop:EigenvectorDenseYaku}, we deduce that $H_+(T_F)=H_-(T_F)=H^p$, and the conclusion follows from the Godefroy-Shapiro Criterion. 
\end{proof}
\par\smallskip
What happens now if, in \Cref{Fig8423}, $\mathbb T$ intersects only one of these two connected components $\Omega_1$ and $\Omega_2$? There are at least two kinds of situations which have to be considered.
\par\smallskip
\begin{enumerate}
    \item[\textbf{Case a:}] the case where $\mathbb T$ intersects the component $\Omega_2=\{\lambda\notin F(\mathbb T)\,;\,|\w_F(\lambda)|=2\}$, but not the component $\Omega_1=\{\lambda\notin F(\mathbb T)\,;\,|\w_F(\lambda)|=1\}$ as is \Cref{Fig937}. 
\begin{figure}[ht]
\includegraphics[page=10,scale=1.1]{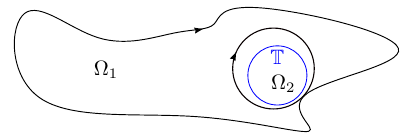}
            \caption{}\label{Fig937}
    \end{figure}
\item[\textbf{Case b:}] the case where $\mathbb T$ intersects $\Omega_1$ but not $\Omega_2$ as in \Cref{Fig10478}.
\begin{figure}[ht]
\includegraphics[page=11,scale=1.1]{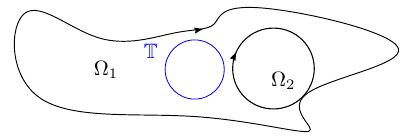}
            \caption{}\label{Fig10478}
    \end{figure}
\end{enumerate}
In the first case, since $\mathbb C\setminus\overline{\mathbb D}$ has accumulation points in both $\Omega_1$ and $\Omega_2$, according to \Cref{Prop:EigenvectorDenseYaku}, we have $H_+(T_F)=H^p$. However, we only know a priori, by \Cref{Prop:EigenvectorDenseYaku} again, that $H_-(T_F)$ contains the subspace $H_{\Omega_2}(T_F)$, and in order to be able to conclude that $H_-(T_F)=H^p$, we would need to have also that 
\[H_{\Omega_1}(T_F)  ~\subseteq~ H_-(T_F).\]
This would obviously be true provided that \[H_{\Omega_1}(T_F)  ~\subseteq~ H_{\Omega_2}(T_F)\]
and then \Cref{Prop:EigenvectorDenseYaku} would imply that $H_-(T_F)=H^p$.
\par\medskip
In the second case, the easy fact is that $H_+(T_F)=H^p$ and it is true that $H_{\Omega_1}(T_F) \subseteq H_-(T_F)$. So what we would need in order to obtain that $H_-(T_F)=H^p$ is that 
\[H_{\Omega_2}(T_F)  ~\subseteq ~H_{\Omega_1}(T_F).\]

\begin{figure}[ht]
\includegraphics[page=12,scale=1.1]{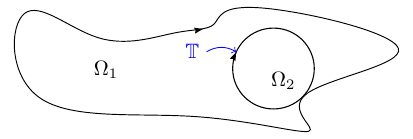}
            \caption{}\label{Fig9485}
            \end{figure}

Observe that one can also have a mix of these two situations, where the unit circle is contained in the boundary of both $\Omega_1$ and $\Omega_2$ as in \Cref{Fig9485}.

Going back to our general situation, let us say that two connected components $\Omega,\Omega'\in\cal C$ are \textit{adjacent} is there exists an arc $\gamma$ with non-empty interior contained in $\partial\Omega\cap\partial\Omega'$. In this case, the winding numbers $\w_F(\Omega)$ and $\w_F(\Omega')$ differ by $1$ exactly. The main problem we are facing is the following:

\begin{question}\label{Q:InclusionEspaceEngVP}
    Let $\Omega$ and $\Omega'$ be two adjacent components of $\sigma(T_F)\setminus F(\mathbb T)$. When is it true that 
  \[H_{\Omega'}(T_F)  ~\subseteq~ H_{\Omega}(T_F)\,?\]
\end{question}

 As illustrated by the example in the next subsection, the difficulty of this question varies, depending on whether $|\w_F(\Omega)|>|\w_F(\Omega')|$ (meaning that $\Omega$ is the interior component and $\Omega'$ is the exterior component with respect to $\gamma$) or $|\w_F(\Omega)|<|\w_F(\Omega')|$ (meaning that $\Omega'$ is the interior component and $\Omega$ is the exterior component). 
\par\smallskip
In our study of \Cref{Q:InclusionEspaceEngVP}, the boundary conditions in  \Cref{eq:boundary-conditions-ErF} defining the subspace $E^q_F=\mathrm{Im}(U)$ of $\bigoplus E^q(\Omega^+_j)$ (where $U$ is the operator defined in \Cref{Section:Yakubovich-court}) will play a crucial role. As a consequence, we will often use the following uniqueness property of functions belonging to Smirnov spaces (which is also recalled in \Cref{smirnov}):
\par\smallskip
\begin{fact}\label{fait-uniqueness-property}
 Let $\Omega$ be a bounded simply connected domain and let $r\ge1$. Let $f\in E^r(\Omega)$. Then $f$ admits a non-tangential limit $f^*$ almost everywhere on $\partial \Omega.$ If $f^*=0$ on a subset of $\partial\Omega$ of positive measure, then $f\equiv0$ on $\Omega$.
\end{fact}

\subsection{Example of two circles}\label{SubSection:Example2Circles} We consider here the following function:
\[F(e^{i\theta})~=~\begin{cases}
-1+2e^{-i3\theta/2}&\text{if}~0\le \theta<4\pi/3\\
e^{-3i\theta}&\text{if}~4\pi/3\le \theta<2\pi.
\end{cases}\]
 \par
\begin{figure}[ht]
\includegraphics[page=13,scale=1.1]{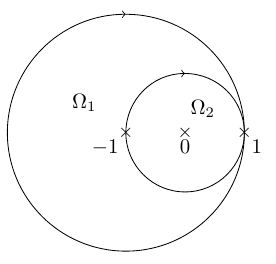}
\caption{}
\label{Fig13}
\end{figure}

\smallskip  
The image of $\mathbb T$ under $F$ is represented in \Cref{Fig13}.
 \par\smallskip           
Our aim here is to show, using the explicit expression of the function $\zeta=1/F^{-1}$ on $\partial\Omega_2\setminus\{1\}$ which is available in this case, that $T_F$ is hypercyclic on $H^p$. 
 \par\smallskip  
Observe that 
\[F(e^{-i\theta})~=~\begin{cases}
-1+2e^{i3\theta/2}&\text{if}~-4\pi/3<\theta\leq 0\\
e^{3i\theta}&\text{if}~-2\pi<\theta\leq -4\pi/3,
\end{cases}\]
which is equivalent to say that 
\[F(e^{-i\theta})~=~\begin{cases}
-1+2e^{i3\theta/2}&\text{if}~-4\pi/3<\theta\leq 0\\
e^{3i\theta}&\text{if}~0<\theta\leq 2\pi/3.
\end{cases}\]
 \par\smallskip  
Observe now that $\partial\Omega_2\setminus\{1\}=\{F(e^{-i\theta}):0<\theta<2\pi/3\}$. Hence, for $\lambda\in\partial\Omega_2\setminus\{1\}$, let $\theta(\lambda)\in (0,2\pi/3)$ be such that $\zeta(\lambda)=e^{i\theta(\lambda)}$. We have $\lambda=e^{3i\theta(\lambda)}$, and since $3\theta(\lambda)\in (0,2\pi)$, we get that $\theta(\lambda)=\frac{1}{3}\arg_{(0,2\pi)}(\lambda)$. So we deduce that 
\[\zeta(\lambda)~=~\exp\left[\frac {i}{3}\arg_{(0,2\pi)}(\lambda)\right]\quad\text{for every }\lambda\in\partial\Omega_2\setminus\{1\}.\]
 \par\smallskip  
Let $g\in H^q$ and write $Ug\in  E^q_F$ as $Ug=(u,v)\in E^q(\Omega_0^+)\oplus E^p(\Omega_1^+)$, where $\Omega_0^+=\Omega_1\cup\Omega_2$ and $\Omega_1^+=\Omega_2$. Then $u$ and $v$ satisfy $u^{int}-\zeta v^{int}=u^{ext}$ almost everywhere on $\partial \Omega_2$ (see \Cref{eq:boundary-conditions-ErF}). Note that, with respect to $\mathbb T=\partial\Omega_2$, $\Omega_1$ is the exterior component and $\Omega_2$ is the interior component.
\par\medskip
$\bullet$ Suppose that $g$ vanishes on $H_{-}(T_F)$. Since $\Omega_2=\mathbb D$, we have
    \[u(\lambda)~=~\dual{h_{\lambda,0}}{g}~=~0\quad \text{and}\quad v(\lambda)~=~\dual{h_{\lambda,1}}{g}~=~0  \quad\text{for every } \lambda\in\Omega_2.\]
So $u^{int}=v^{int}=0$ almost everywhere on $\partial\Omega_2$, and thus the boundary relation becomes 
    \[u^{ext}~=~0\quad\text{a.e. on}~\partial\Omega_2.\] 
    Since $u_{|\Omega_1}\in E^q(\Omega_1)$ and its non-tangential limit is zero on a subset of $\partial\Omega_1$ with positive measure, we deduce that $u=0$ on $\Omega_1$. This means that 
    $g$ vanishes on $H_{\Omega_1}(T_F)$.
    So $g$ vanishes on $H_{\Omega_1}(T_F)$ and $H_{\Omega_2}(T_F)$, and by \Cref{Prop:EigenvectorDenseYaku}, we deduce that $g=0$ and that $H_-(T_F)=H^p$.
\par\medskip
$\bullet$  Suppose now that $g$ vanishes on $H_{+}(T_F)$. Since $\Omega_1 \subseteq \mathbb C\setminus\overline{\mathbb D}$, $g$ vanishes on $H_{\Omega_1}(T_F)$, i.e. 
    \[u(\lambda)~=~\dual{h_{\lambda,0}}{g}=0\quad \text{for every } \lambda\in\Omega_1.\]
    Then $u^{ext}=0$ almost everywhere on $\partial\Omega_2$ and the boundary relation becomes 
    \[u^{int}-\zeta v^{int}~=~0\quad\text{a.e. on}~\partial\Omega_2.\]
The function $\zeta$ admits an analytic extension to $\Omega_2\setminus[0,1]$, given by
\begin{equation}
\zeta(\lambda)~=~\exp\left[\frac 13\log_{(0,2\pi)}(\lambda)\right],\quad\text{for every }\lambda\in\Omega_2\setminus[0,1],
\end{equation}
where $\log_{(0,2\pi)}$ is an analytic branch of the logarithm with imaginary part in $(0,2\pi)$.
Note that this extension is such that
\begin{equation}\label{discontinuite}
\lim_{\underset{y>0}{y\to0}}\zeta(x+iy)~\neq~ \lim_{\underset{y<0}{y\to0}}\zeta(x+iy)\quad \text{for every } x\in]0,1[.
\end{equation}
Since $\zeta$ is analytic and bounded on $\Omega_2\setminus[0,1]$, the function $\zeta v$ belongs to $E^q(\Omega_2\setminus[0,1])$, and since we have $u^{int}=\zeta v^{int}$ almost everywhere on $\partial\Omega_2$, which is a subset of $\partial(\Omega_2\setminus[0,1])$ of positive measure, we deduce that $u=\zeta v$ on $\Omega_2\setminus[0,1]$.
\par\smallskip
If we suppose now that $v$ is not identically zero on $\Omega_2$, there exists a point $x_0\in(0,1)$ such that $v$ does not vanish on a neighborhood of $x$. Thus 
\[\lim_{\underset{y>0}{y\to0}}\zeta(x_0+iy)~=~\lim_{\underset{y>0}{y\to0}}\frac{u(x_0+y)}{v(x_0+iy)}~=~\lim_{\underset{y<0}{y\to0}}\frac{u(x_0+y)}{v(x_0+iy)}~=~\lim_{\underset{y<0}{y\to0}}\zeta(x_0+iy),\]
and this contradicts \Cref{discontinuite}. So $v=0$ on $\Omega_2$ and thus $u=0$ on $\Omega_2$ as well. This means that 
$g$ vanishes on $H_{\Omega_2}(T_F)$,
 so   $g=0$ and thus $H_+(T_F)=H^p$.
\par\smallskip
Hence we have shown that $H_{-}(T_F)=H_{+}(T_F)=H^p$ and, by the Godefroy-Shapiro Criterion, the operator $T_F$ is hypercyclic on $H^p$. 
\par\smallskip
This example highlights the fact that when we need to go from a connected component to an adjacent one with a lower winding number, the situation is quite simple, whereas going to an adjacent component with a higher winding number becomes substantially more difficult. 

\subsection{From the interior to the exterior}
The easiest case is the first one, when $\Omega$ is the exterior component of $\gamma$ and $\Omega'$ is the interior one. In this case, \Cref{Q:InclusionEspaceEngVP} always admits an affirmative answer.

\begin{theorem}\label{Th:Simple_Int-Ext}
Let $p>1$. Suppose that $F$ satisfies \ref{H1}, \ref{H2} and \ref{H3} and let $\Omega,\Omega'\in\cal C$ be two adjacent components of $\sigma(T_F)\setminus F(\mathbb T)$ along an arc $\gamma$, such that $|\w_F(\Omega)|>|\w_F(\Omega')|$.
Then $H_{\Omega'}(T_F)  \subseteq H_\Omega(T_F)$.
\end{theorem}

\begin{proof}
Suppose that $g\in H^q$ vanishes on $H_\Omega(T_F)$, i.e. that $\dual{h_{\lambda,j}}{g}=0$ for every $\lambda\in\Omega$ and every $0\le j<|\w_F(\Omega)|$. Set $j_0=|\w_F(\Omega)|-1\ge0$, so that $\Omega \subseteq\Omega_{j_0}^+$. 
Denoting by $U$ the operator from \Cref{T:Yakbovich_Hp} and setting $Ug=(u_j)_{0\le j\le N-1}$ with $u_j\in E^q(\Omega_j^+)$ for $0\le j\le N-1$, our assumption on $g$ can be rewritten as $u_j=0$ on $\Omega$ for every $0\le j\le j_0$.
\par\smallskip
The assumption $|\w_F(\Omega)|>|\w_F(\Omega')|$ means that $\Omega$ is the interior component of $\gamma$ and $\Omega'$ is the exterior one. It follows that
\[u_j^{int}~=~0\quad\text{a.e. on}~\partial\Omega~\text{for every}~0\le j\le j_0,\]
so that, for every $0\le j<j_0$, $u_j^{int}=0$ and $u_{j+1}^{int}=0$ almost everywhere on $\partial\Omega$.
\par\smallskip
The functions $u_j$ satisfy the boundary relations
\[u_j^{int}-\zeta u_{j+1}^{int}~=~u_j^{ext}\quad\text{a.e. on}~\partial\Omega_{j+1}^+~\text{for every}~0\le j< N-1.\]
Since $|\w_F(\Omega)|>j_0\ge j+1$ for every $0\le j< j_0$, it follows that $\Omega \subseteq \Omega_{j+1}^+$ and 
\[u_j^{int}-\zeta u_{j+1}^{int}~=~u_j^{ext}\quad\text{a.e. on}~\partial\Omega~\text{for every }0\le j< j_0,\]
so that $u_j^{ext}=0$ almost everywhere on $\gamma  \subseteq \partial \Omega$. Now $\gamma \subseteq\partial \Omega'$ is a subset of positive measure of $\partial \Omega'$ and ${u_j}_{|\Omega'}\in E^q(\Omega')$ for every $0\le j< j_0$, so $u_j=0$ on $\Omega'$ for every $0\le j< j_0$. Since $|\w_F(\Omega')|=|\w_F(\Omega)|-1=j_0$,
\[u_j(\lambda)~=~\dual{h_{\lambda,j}}{g}~=~0\quad \text{for every}~\lambda\in\Omega'~\text{and every}~0\le j<|\w_F(\Omega')|,\]
i.e. $g$ vanishes on all the eigenspaces $\ker(T_F-\lambda)$ with $\lambda\in\Omega'$. So $g$ vanishes on $H_{\Omega'}(T_F)$, and \Cref{Th:Simple_Int-Ext} is proved.
\end{proof}

\subsection{From the exterior to the interior}\label{subsection:Ext-to-Int}
Thus one can always ``go up'' from one component of $\sigma(T_F)\setminus F(\mathbb T)$ to an adjacent one. We now need to ``go down'', i.e. to deal with the case where $|\w_F(\Omega)|<|\w_F(\Omega')|$. Here the situation is more intricate, and we are able to show that $H_{\Omega'}(T_F)  \subseteq H_\Omega(T_F)$ only under stronger assumptions, concerning both the smoothness of the function $F$ and a certain geometric property \ref{cond(G)} of the boundaries of the domains $\Omega$ and $\Omega'$. 
\par\medskip
Let $\Omega,\Omega'\in\mathcal C$ be two adjacent components of $\sigma(T_F)\setminus F(\mathbb T)$ along an arc $\gamma$ such that $|\w_F(\Omega)|<|\w_F(\Omega')|$. We say that the pair $(\Omega,\Omega')$ satisfies the property \ref{cond(G)} if the following holds:
\par\smallskip
\begin{enumerate}[(P)]
    \item \label{cond(G)} there exists a self-intersection point $\lambda_0$ of the curve $F(\mathbb T)$ and an open neighborhood $V$ of $\lambda_0$ such that $V\cap\mathcal{O}=\{\lambda_0\}$ and \par\smallskip  
\begin{enumerate}[(i)]
    \item $\Omega$ is the only adjacent component of $\Omega'$ in a neighborhood of the point $\lambda_0$, i.e. $V\cap\partial \Omega'=V\cap\partial \Omega\cap \partial\Omega'$;
    \par\smallskip
    \item$\lambda_0\in\partial\Omega\cap\partial\Omega'$, and the restriction of the function $\zeta$ to $(\partial\Omega\cap\partial\Omega'\cap V)\setminus\{\lambda_0\}$ cannot be extended continuously at the point $\lambda_0$;
    \par\smallskip
    \item there exists a simple closed arc $\Delta$ with $\Delta\cap F(\mathbb T)=\{\lambda_0\}$ such that $\Delta$ separates $\Omega'\cap V$ into two connected components and $\zeta_{|(\partial\Omega'\cap V)\setminus\{\lambda_0\}}$ admits a bounded analytic extension to $(V\cap \Omega')\setminus\Delta$.
\end{enumerate}
\end{enumerate}
\par
\begin{figure}[ht]
    \includegraphics[page=53,scale=1.15]{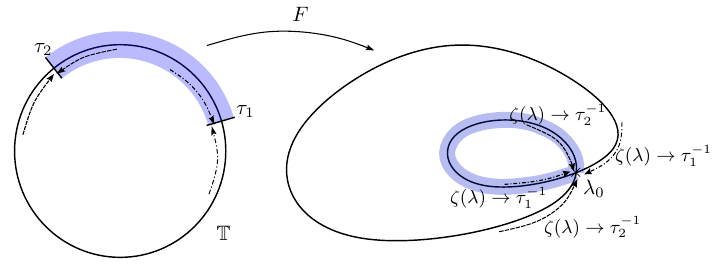}
    \caption{}\label{Fig136589}
\end{figure}

\medskip
The meaning of the geometric assumption \ref{cond(G)} may be somewhat difficult to grasp. Let us first explain how one can see whether an intersection point $\lambda_0$ on $\partial\Omega\cap\partial\Omega'$ satisfies or not properties (i) and (ii) of \ref{cond(G)}. These two properties involve the behavior of $\zeta$ in a neighborhood of $\lambda_0$. Recall that $\zeta$ is constructed from $F$ (we have indeed $\zeta=1/F^{-1}$ on $F(\mathbb T)\setminus\mathcal O$), so in each direction of the curve near $\lambda_0$, $\zeta$ admits a limit. For example, in \Cref{Fig136589}, in two different directions, the limits are $1/\tau_1$, and in the other two directions, the limits are $1/\tau_2$.
\par
\begin{figure}[ht]
\begin{subfigure}[b]{0.32\textwidth}
\includegraphics[page=25,scale=.87]{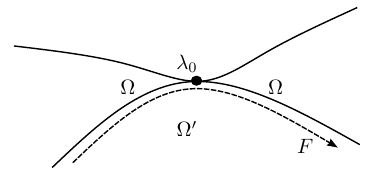}
\caption{}\label{Fig14a}
\end{subfigure}~
\begin{subfigure}[b]{0.32\textwidth}
\includegraphics[page=26,scale=.87]{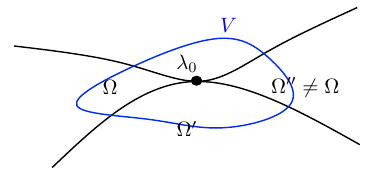}
\caption{}\label{Fig14b}
\end{subfigure}~
\begin{subfigure}[b]{0.32\textwidth}
\includegraphics[page=27,scale=.87]{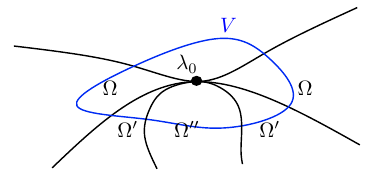}
\caption{}\label{Fig14c}
\end{subfigure}
    \caption{}
\end{figure}

\smallskip
For $\lambda_0$ to satisfy properties (i) and (ii) of \ref{cond(G)}, we require that on $\partial\Omega\cap\partial\Omega'$, it is possible to approach $\lambda_0$ from two different directions, and that the limits of $\zeta$ in these two directions are not the same. As a result, it will not be possible to extend $\zeta$ to $\partial\Omega\cap\partial\Omega'$ in such a way that it is continuous at $\lambda_0$. This means that we are not in the situation of Figure {\sc\ref{Fig14a}}.
We also require that there is a neighborhood $V$ of $\lambda_0$ such that $V\cap\partial\Omega'=V\cap\partial\Omega\cap\partial\Omega'$, i.e. $\Omega$ is the only adjacent component of $\Omega'$ in a neighborhood of $\lambda_0$. So, if for example the intersection at the point $\lambda_0$ is simple, then we cannot be in the situation described in Figure {\sc\ref{Fig14b}}.  
If the intersection at the point $\lambda_0$ is not simple, then there is no component $\Omega''$ inside $\Omega'$ which would also have $\lambda_0$ as a boundary point, i.e. we are not in the situation in Figure {\sc\ref{Fig14c}}.

\par
\begin{figure}[ht]
\includegraphics[page=28,scale=1.1]{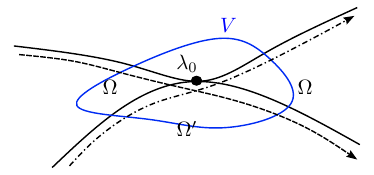}
\caption{}\label{Fig15909t}
\end{figure}

\medskip
An intersection point on $\partial\Omega\cap\partial\Omega'$ which satisfies the assumptions (i) and (ii) in \ref{cond(G)} is a point where we have the type of situation represented in \Cref{Fig15909t}.

Observe that when $F$ satisfies assumptions \ref{H1}, \ref{H2} and \ref{H4}, the property (iii) in \ref{cond(G)} is automatically satisfied. Indeed, we have:\par\smallskip
\begin{fact}\label{Fait:H4=>(P)iii}
Let  $F$ satisfy \ref{H1}, \ref{H2} and \ref{H4}. Let $\lambda_0\in\cal O$. Then for every sufficiently small connected open neighborhood $V$ of $\lambda_0$ satisfying $V\cap \mathcal O=\{\lambda_0\}$ and every connected component $\Sigma$ of $V\setminus F(\mathbb T)$, there exists a simple closed arc $\Delta$ with $\Delta\cap F(\mathbb T)=\{\lambda_0\}$ which separates $\Sigma$ into two connected components such that $\zeta_{|(\partial \Sigma\cap V)\setminus\{\lambda_0\}}$ has a bounded analytic extension to $\Sigma\setminus\Delta$.
\end{fact}
\par\smallskip
Note that examples such as the ones considered in \Cref{SubSection:Example2Circles} satisfy the hypothesis \ref{H1}, \ref{H2} and \ref{H4} and thus also (iii) of property \ref{cond(G)}. Note also that \Cref{Fait:H4=>(P)iii} holds under a weaker hypothesis than \ref{H4}: indeed, as we will see in the proof, we only use the fact that $F$ has an analytic extension to a neighborhood of the extremities of the $\alpha_j$'s. 

\begin{proof}Let $\theta_0<\theta_1<\dots<\theta_m=\theta_0+2\pi$ given by \ref{H2} and let $\alpha_j$ the open arcs defined in \ref{H2}. We have
\[\overline{\alpha}_j~=~\{e^{i\theta}\,;\,\theta_j\le \theta\le \theta_{j+1}\}.\]
\par\smallskip
Since $F$ satisfies \ref{H4}, for every $0\le j< m$, there exists  a neighborhood $U_j$ of the closed arc $\overline{\alpha}_j$ such that  $F_{|\overline{\alpha}_j}$ admits an analytic extension $F_j$ to $U_j$. Since $F'(e^{i\theta})\neq0$ for every $e^{i\theta}\in\mathbb T$, $F_j'(e^{i\theta})\neq 0$ for every point $e^{i\theta}\in \overline{\alpha}_j$, and hence $F_j$ is a local biholomorphism at every point $e^{i\theta}\in\overline{\alpha}_j$, including the extremities $e^{i\theta_j}$ and $e^{i\theta_{j+1}}$ of $\overline{\alpha}_j$. 
\par\smallskip
Let now $\lambda_0\in \mathcal O$, and let 
\[
K_{\lambda_0}=\{1\le k\le m\,;\,\lambda_0=F(e^{i\theta_k})\}.
\]
For each $k\in K_{\lambda_0}$, $e^{i\theta_k}$ is an extremity of the two arcs 
$\alpha_k$ and $\alpha_{k-1}$, with $\alpha_0=\alpha_m$. Let $U_{k,0},U_{k,1}$ and 
$V_{k,0},V_{k,1}$ be connected open neighborhoods of $e^{i\theta_k}$ and $\lambda_0$ 
respectively such that $F_k:U_{k,0}\overset{\sim}{\longrightarrow} V_{k,0}$ and 
$F_{k-1}:U_{k,1}\overset{\sim}{\longrightarrow} V_{k,1}$ are biholomorphisms, with 
$F_0=F_m$. We denote by $\xi_{k,0}:V_{k,0}\overset{\sim}{\longrightarrow} U_{k,0}$ 
and $\xi_{k,1}:V_{k,1}\overset{\sim}{\longrightarrow} U_{k,1}$ their respective 
inverses. We have $F_k(\xi_{k,0}(\lambda))=\lambda$ for every $\lambda\in V_{k,0}$ 
and $F_{k-1}(\xi_{k,1}(\lambda))=\lambda$ for every $\lambda\in V_{k,1}$. Suppose 
now that $\lambda\in F(\overline{\alpha}_k\cap U_{k,0})$, i.e. that $\lambda=F(e^{i\theta})$ for some $e^{i\theta}\in \overline \alpha_k\cap U_{k,0}$. Since the restriction of $F_k$ to 
$\overline \alpha_k$ coincides with $F$, we have $F_k(e^{i\theta})=\lambda$. But 
$F_k(\xi_{k,0}(\lambda))=\lambda$, $F_k$ is injective on $U_{k,0}$, and 
$\xi_{k,0}(\lambda)$ and $e^{i\theta}$ both belong to $U_{k,0}$, and thus 
$\xi_{k,0}(\lambda)=e^{i\theta}$. It follows that $\xi_{k,0}$ and $1/\zeta$ 
coincide on $F(\alpha_k\cap U_{k,0})$. In the same way, 
$\xi_{k,1}$ and $1/\zeta$ coincide on $F(\alpha_{k-1}\cap 
U_{k,1})$ (see \Cref{Fig:H4=>iii}).

\begin{figure}[ht]
    \includegraphics[page=62]{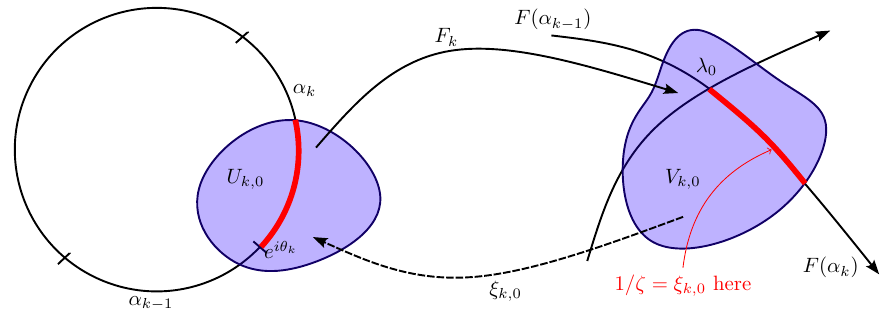}
    \caption{}
    \label{Fig:H4=>iii}
\end{figure}

Let now $V=\cap_{k\in K_{\lambda_0}}(V_{k,0}\cap V_{k,1})$, which is a connected open 
neighborhood of $\lambda_0$. Let $\Sigma$ be a connected component of $V\setminus 
F(\mathbb T)$. Making $V$ smaller if necessary, we can assume that  
$\partial\Sigma\cap\cal O=\{\lambda_0\}$ and that $\partial\Sigma\cap V$ consists of two arcs intersecting at the point $\lambda_0$, as in \Cref{Fig:H4=>iii}. 

Let $\Delta$ be any simple closed arc 
(for instance a segment) with $\Delta\cap F(\mathbb T)=\{\lambda_0\}$ which separates $\Sigma$ in two 
components $\Sigma_1,\Sigma_2$. Let now $(k_1,\varepsilon_1)$ and $(k_2,\varepsilon_2)$ with $k_1,k_2\in 
K_{\lambda_0}$ and $\varepsilon_1,\varepsilon_2\in\{0,1\}$ be such that 
$\partial\Sigma_1 \cap F(\mathbb T)\subset F(\overline \alpha_{k_1-\epsilon_1}\cap 
U_{k_1,\varepsilon_1})$ and $\partial\Sigma_2 \cap F(\mathbb T)\subset 
F(\overline{\alpha}_{k_2-\epsilon_2}\cap U_{k_2,\varepsilon_2})$. Then 
$\xi_{k_1,\varepsilon_1}$ and $1/\zeta$ coincide on $(\partial\Sigma_1\cap F(\mathbb 
T))\setminus\{\lambda_0\}$, while $\xi_{k_2,\varepsilon_2}$ and $1/\zeta$ coincide on $(\partial\Sigma_2\cap F(\mathbb T))\setminus\{\lambda_0\}$. It follows that the 
function $\widetilde\zeta$ defined on $\Sigma\setminus\Delta$ by setting $\widetilde\zeta=1/\widetilde\xi$, where
\[\widetilde{\xi}(\lambda)~=~
\begin{cases}
\xi_{k_1,\varepsilon_1}(\lambda)&\text{if}~\lambda\in\Sigma_1\\
\xi_{k_2,\varepsilon_2}(\lambda)&\text{if}~\lambda\in\Sigma_2\\
\end{cases}\]
is a bounded analytic extension to $\Sigma\setminus\Delta$ of the restriction of $\zeta$ to $\partial\Sigma\setminus\{\lambda_0\}$.
\end{proof}

In some circumstances, we are able to handle situations where \ref{cond(G)} does not hold. We will get back to this at the end of this section, and also provide some examples in \Cref{Section:ApplicationCShc}.

\begin{theorem}\label{Th:Complique_Ext-Int}Let $p>1$.
Suppose that $F$ satisfies \ref{H1}, \ref{H2} and \ref{H3} and let $\Omega,\Omega'\in\cal C$ be two adjacent components of $\sigma(T_F)\setminus F(\mathbb T)$ along an arc $\gamma$. Suppose additionally that the pair $(\Omega,\Omega')$ satisfies the property \ref{cond(G)}. Then $H_{\Omega'}(T_F)= H_\Omega(T_F)$.
\end{theorem}

\begin{proof}
    The beginning of the proof of \Cref{Th:Complique_Ext-Int} is the same as that of \Cref{Th:Simple_Int-Ext}. Let $g\in H^q$ which vanishes on $H_\Omega(T_F).$ Set $j_0=|\w_F(\Omega)|-1$ and define $Ug=(u_j)_{0\le j\le N-1}$ with $u_j\in E^q(\Omega_j^+)$, $0\le j\le N-1$. We know that
    $u_j=0~\text{on}~\Omega~\text{for every}~0\le j\le j_0.$
    Now the assumption $|\w_F(\Omega)|<|\w_F(\Omega')|$ means that $\Omega'$ is the interior component of $\gamma$ and $\Omega$ is the exterior component of $\gamma$. Let $0\leq j\leq j_0$. 
    Since  $|\w_F(\Omega')|=|\w_F(\Omega)|+1=j_0+2>j_0+1$, $\Omega' \subseteq\Omega_{j+1}^+$, and we have $\partial\Omega\cap\partial\Omega'  \subseteq \partial\Omega_{j+1}^+$ for every $0\leq j\leq j_0$. So, according to \Cref{eq:boundary-conditions-ErF}, we have 
    \[u_j^{int}-\zeta u_{j+1}^{int}~=~u_j^{ext}\quad\text{a.e. on}~\partial\Omega\cap\partial\Omega'~\text{for every}~0\le j\le j_0.\]
    Now $u_j^{ext}=0$ almost everywhere on $\partial\Omega\cap\partial\Omega'$ for every $0\le j\le j_0$, and hence  \[u_j^{int}-\zeta u_{j+1}^{int}~=~0\quad\text{a.e. on}~\partial\Omega'\cap\partial \Omega~\text{for every}~0\le j\le j_0.\]
    We would like to be able to deduce from this relation that $u_j=0$ on $\Omega'$ for every $j$ with $0\le j\le |\w_F(\Omega')|-1=j_0+1$. It is here that the assumptions \ref{cond(G)} come into play. We state separately as \Cref{Prop:Clef_Ext-Int} the result we will need in order to conclude the proof.

\begin{proposition}\label{Prop:Clef_Ext-Int}Let $p>1$.
Suppose that $F$ satisfies \ref{H1}, \ref{H2} and \ref{H3}, and let $\zeta=1/F^{-1}$ on $F(\T)\setminus\mathcal{O}$.
Fix $\Omega_0\in\cal C$. Suppose that there exists $\lambda_0\in\partial \Omega_0$ which is a point of self-intersection of the curve $F(\mathbb T)$ 
such that for every sufficiently small connected open neighborhood  $V$ of $\lambda_0$,
there exists a connected component $\Sigma$ of $V\cap\Omega_0$ such that\par\smallskip
\begin{enumerate}[(a)]
\item the restriction of $\zeta$ to $\partial\Sigma\setminus\{\lambda_0\}$ cannot be extended continuously at the point $\lambda_0$;
\par\smallskip
\item there exists a simple closed arc $\Delta$ with $\Delta\cap F(\mathbb T)=\{\lambda_0\}$ such that $\Delta$ separates $\Sigma$ into two connected components and $\zeta_{|(\partial\Sigma\cap V)\setminus\{\lambda_0\}}$ admits a bounded analytic extension to $\Sigma\setminus\Delta$.
\end{enumerate}
Let $r>1$, and let $u,v\in E^r(\Omega_0)$. If $u-\zeta v=0$ almost everywhere on $\partial \Omega_0\cap V$, then $u=v=0.$ 
\end{proposition}
\par\smallskip
Note that the condition on $\lambda_0$ in \Cref{Prop:Clef_Ext-Int} implies that for every connected neighborhood $V$ of $\lambda_0$, we can find a connected component $\Sigma$ of $V\cap\Omega_0$ such that $\partial\Sigma\cap F(\mathbb T)$ cannot be written as the image by $F$ of exactly one arc in $\mathbb T$. Also note the following: if the assumptions of \Cref{Th:Complique_Ext-Int} are satisfied, if we set $\Omega_0=\Omega'$ and if $V$ is any sufficiently small connected neighborhood of the point $\lambda_0$, then $V\cap\Omega_0$ is always a connected set, so that the assumptions of \Cref{Prop:Clef_Ext-Int} are satisfied with $\Sigma=V\cap\Omega_0$. 
\par\smallskip
Taking \Cref{Prop:Clef_Ext-Int} for granted, let us explain how we conclude the proof of \Cref{Th:Complique_Ext-Int}.
For every $0\le j\le j_0$, the functions $u_j$ and $u_{j+1}$ belong to $E^q(\Omega_j^+)$ and  $E^q(\Omega_{j+1}^+)$ respectively.
Since $|\w_F(\Omega')|>j_0+1\ge j+1$ for every $0\le j \le j_0$, we have $\Omega' \subseteq\Omega_{j+1}^+ \subseteq\Omega_j^+$, and thus $u_j$ and $u_{j+1}$ both belong to $E^q(\Omega')$.
We apply \Cref{Prop:Clef_Ext-Int} to $\Omega_0=\Omega'$, $\lambda_0$ and $V$ given by the assumption \ref{cond(G)}, $\Sigma=\Omega_0\cap V$, $u=u_j$ and $v=u_{j+1}$.
We have $u-\zeta v=0$ almost everywhere on $\partial\Omega\cap\partial\Omega'$, hence $u-\zeta v=0$ almost everywhere on $\partial\Omega\cap\partial\Omega'\cap V$. 
Now by assumption \ref{cond(G)}, $\partial\Omega\cap\partial\Omega'\cap V=\partial\Omega'\cap V$.
So $u-\zeta v=0$ almost everywhere on $\partial\Omega'\cap V$. \Cref{Prop:Clef_Ext-Int} applies, and yields that $u=v=0$ on $\Omega'$.
\par\smallskip
We have proved that $u_j$ and $u_{j+1}$ vanish on $\Omega'$ for every $0\le j\le j_0$, which is the conclusion we were looking for: $g$ vanishes on $\ker(T_F-\lambda)$ for every $\lambda\in\Omega'$, and so $H_{\Omega'}(T_F)  \subseteq H_{\Omega}(T_F)$. The other inclusion is given by \Cref{Th:Simple_Int-Ext}.
\end{proof}
It remains to prove \Cref{Prop:Clef_Ext-Int}.

\begin{proof}[Proof of \Cref{Prop:Clef_Ext-Int}]
In agreement with assumption \ref{H2}, let $\theta_0<\theta_1<\dots<\theta_{m}=\theta_0+2\pi$ satisfying \ref{H2}. Let \[ \alpha_j~=~\{e^{i\theta}\,;\,\theta_j< \theta<\theta_{j+1}\}.\]Without loss of generality, for all $0\le j\le m-1$, we assume that no point of $\mathcal O$ belongs to $F(\alpha_j)$. 
We define $\gamma_j$ to be the image by $F$ of the closure of $\alpha_j$, i.e. $\gamma_j=F(\overline \alpha_j)$, so that $$F(\mathbb T)~=\bigcup_{1\le j\le m} \gamma_j.$$ The self-intersection points of the curve $F(\mathbb T)$ are to be found among the points $F(e^{i\theta_j}),~j=0,\dots,m-1$. Whenever $\lambda\in F(\mathbb T)$ is a self intersection point of the curve, we set
\[J_\lambda~:=~\big\{j\in\{0,\dots,m-1\}\,;\,\lambda=F(e^{i\theta_j})\big\}.\] Then $J_\lambda$
contains at least two elements. If we enumerate $J_{\lambda}$ as $J_{\lambda}=\{j_1,\dots,j_s\}$ for some $s>1$, where $j_1<\dots<j_s$, then for each $1\le k\le s$, $\lambda$ belongs to the arcs $\gamma_{j_k}$ and $\gamma_{j_{k+1}}$, where $\gamma_{m}=\gamma_0$. See \Cref{Figdk02iu4}.
\begin{figure}[ht]
\includegraphics[page=16,scale=1]{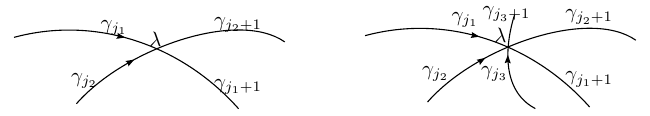}
\caption{}\label{Figdk02iu4}
\end{figure}
 
Let $\Omega_0,~\lambda_0,~V$ and $\Sigma$ be given satisfying the assumptions of  \Cref{Prop:Clef_Ext-Int}.
We know that the restriction of the function $\zeta$ to $\partial\Sigma$ cannot be extended continuously at the point $\lambda_0$. Reducing $V$ if necessary, it follows that there exist two integers $1\le k\neq l\le s$ and $\varepsilon_k,\varepsilon_l\in\{0,1\}$ such that 
\begin{equation}\label{discontinuite-bis}
\lim_{\underset{\lambda\in\gamma_{j_k+\varepsilon_k}}{\lambda\to\lambda_0}}\zeta(\lambda)~=~ e^{i\theta_{j_k}}\quad\text{and}\quad\lim_{\underset{\lambda\in\gamma_{j_l+\varepsilon_l}}{\lambda\to\lambda_0}}\zeta(\lambda)~=~ e^{i\theta_{j_l}},
\end{equation}
\(\partial\Sigma\cap V~=~(\gamma_{j_k+\varepsilon_k}\cup\gamma_{j_l+\varepsilon_l})\cap V\)
and 
\((\gamma_{j_k+\varepsilon_k}\cap V)\cap(\gamma_{j_l+\varepsilon_l}\cap V)~=~\{\lambda_0\}.\)
What it means is best seen on a picture as in \Cref{Fig17lfa}  with for example $k=1,l=3,\varepsilon_1=0$ and $\varepsilon_3=1$:
\begin{figure}[ht]
\includegraphics[page=17,scale=1]{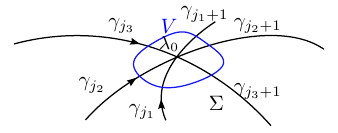}
\caption{}
\label{Fig17lfa}
\end{figure}

Since $u-\zeta v=0$ almost everywhere on $\partial\Omega_0\cap V$, $u-\zeta v=0$ almost everywhere on $\partial\Sigma\cap V$. Let $\Delta$ be a simple closed arc given by the hypothesis of \Cref{Prop:Clef_Ext-Int} and containing the point $\lambda_0$, and let $\Delta'=\overline{\Sigma}\cap\Delta$. Reducing again $V$ if necessary, we can assume that 
$F(\T)\cap\Sigma$ consists of two arcs intersecting at the point $\lambda_0$, that
$\Sigma\setminus \Delta'$ has exactly two components and that $\zeta$ has a bounded analytic extension to $\Sigma\setminus\Delta'$.
We denote by $D_k$ the component containing the set $V\cap\gamma_{j_k+\varepsilon_k}\cap\Sigma$, and by $D_l$ the component containing the set $V\cap\gamma_{j_l+\varepsilon_l}\cap\Sigma$. See \Cref{Fig4.18}.

\begin{figure}[ht]
 \includegraphics[page=59,scale=1]{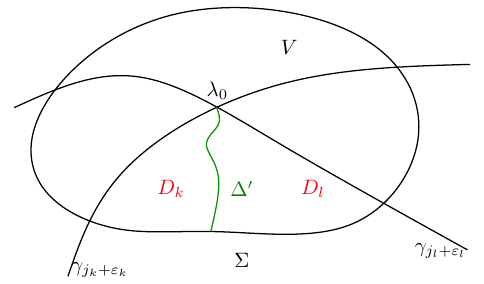}
\caption{}\label{Fig4.18}
     \end{figure}
     Since $u,v\in E^r(\Sigma)$ and $\zeta$ is bounded on $\Sigma\setminus\Delta'$, $u-\zeta v\in E^r(D_k)$. Our assumption that $u-\zeta v=0$ almost everywhere on $\partial\Sigma\cap V$ implies that $u-\zeta v=0$ almost everywhere on $\gamma_{j_k+\varepsilon_k}\cap \overline{D_k}$, which is a subset of $\partial D_k$ with positive measure. Hence $u-\zeta v=0$ on $D_k$ by the uniqueness property of Smirnov functions (see \Cref{fait-uniqueness-property}). In the same way, $u-\zeta v=0$ on $D_l$.
     
Suppose that $v$ is not identically $0$ on $\Sigma$. Then the meromorphic function $u/v$ coincide with $\zeta$ on $D_k\cup D_l=\Sigma\setminus \Delta'$. Since $\zeta$ is bounded on $\Sigma\setminus \Delta'$, it follows that $u/v$  is also bounded on $\Sigma\setminus \Delta'$. Thus $u/v$ has no pole on $\Delta'$, i.e. $u/v\in Hol(\Sigma)$, and so $\zeta$ admits a bounded analytic extension to $\Sigma$.  But by (\ref{discontinuite-bis}), this extension satisfies
\begin{equation}\label{discontinuite2}
 \lim_{\underset{\lambda\in\gamma_{j_k+\varepsilon_k}}{\lambda\to\lambda_0}}\zeta(\lambda)~=~ e^{i\theta_{j_k}}\quad\text{and}\quad\lim_{\underset{\lambda\in\gamma_{j_l+\varepsilon_l}}{\lambda\to\lambda_0}}\zeta(\lambda)~=~ e^{i\theta_{j_l}}.   
\end{equation}
Since $\theta_{j_k}\neq\theta_{j_l}$, this gives a contradiction with \Cref{Lindelof-Jordan}. So $v=0$ on $\Sigma$ and thus on $\Omega_0$ and so we also have $u=0$ on $\Omega_0$.
\end{proof}

\begin{remark}\label{Remark-Yakubovich} 
It was kindly pointed out to us by D. Yakubovich that 
if $F$ satisfies \ref{H4} (and in particular if $F$ has an analytic extension to a neighborhood of $\mathbb T$),  the proof of \Cref{Prop:Clef_Ext-Int} can be simplified, avoiding the use of Lindelöf's Theorem. Indeed, let $\gamma_1,\gamma_2$ be the connected components of $(\partial\Sigma\cap V)\setminus\{\lambda_0\}$ (we may assume that $\lambda_0$ is the only point of $\mathcal O$ which belongs to the two closed arcs $\overline{\gamma}_j$, $j=1,2$). 
Then there exists a disk $D$ centered in $\lambda_0$ and two analytic functions $\zeta_1$ and $\zeta_2$ on $D$ such that $\zeta_1$ coincides with $\zeta$ on $\gamma_1\cap D$ and $\zeta_2$ coincides with $\zeta$ on $\gamma_2\cap D$ (see the proof of \Cref{Fait:H4=>(P)iii}).
Let $W=\Sigma\cap D$. Then $\zeta_1$ and $\zeta_2$ belong to $H^\infty(W)$, and so $u-\zeta_1 v$ and $u-\zeta_2 v$ belong to $E^q(W)$. Moreover, $u-\zeta_1 v=0$ almost everywhere on $\gamma_1$ and $u-\zeta_2 v=0$  almost everywhere on $\gamma_2$. Since $\gamma_1$ and $\gamma_2$ are sub-arcs of $\partial W$ with positive measure, we have $u-\zeta_1 v=u-\zeta_2 v=0$ on $W$, so that $\zeta_1 v=\zeta_2 v$ on $W$. Since it is supposed in \Cref{Prop:Clef_Ext-Int} that the restriction of $\zeta$ to $\partial\Sigma$ 
cannot be extended continuously at the point $\lambda_0$, we have $\zeta_1(\lambda_0)\neq\zeta_2(\lambda_0)$. Hence there exists $\eta>0$ such that for every $z\in W\cap D(\lambda_0,\eta)$, $\zeta_1(z)\neq \zeta_2(z)$. It follows that $v=0$ on $W$, and hence that $u=0$ on $W$.
By the uniqueness principle, $u=v=0$ on $\Omega_0$, which concludes the proof of \Cref{Prop:Clef_Ext-Int} under the assumption \ref{H4}. 
\end{remark}

\begin{figure}[ht]
\includegraphics[page=19,scale=1]{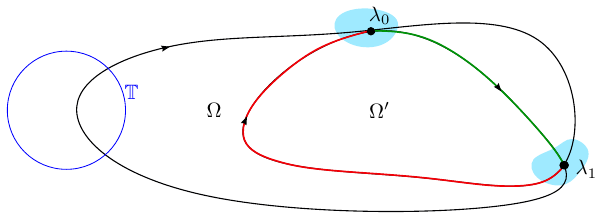}
\caption{}\label{Fig19acb}\end{figure}

We conclude this subsection with a modified version of \Cref{Prop:Clef_Ext-Int} which is tailored to deal  with the kind of configuration given in \Cref{Fig19acb}, where \Cref{Th:Complique_Ext-Int} cannot be applied.

In \Cref{Fig19acb}, for any open neighborhood $V$ of one of the two self intersection points $\lambda_0$ and $\lambda_1$ of the curve $F(\mathbb T),~V\cap\partial\Omega\cap\partial\Omega'\varsubsetneq V\cap \partial\Omega'$. When running through the proof of \Cref{Th:Complique_Ext-Int} for such an example, we arrive to the equation $u-\zeta v=0$ almost everywhere on $\partial\Omega\cap\partial\Omega'\cap V$. In order to be able to apply \Cref{Prop:Clef_Ext-Int}, we need to be able to ensure that $u-\zeta v=0$ almost everywhere on $\partial\Omega'\cap V$, and since $V\cap\partial\Omega\cap\partial\Omega'\neq V\cap \partial\Omega'$, it may not be necessarily the case.
\par\medskip
However, the following consequence of \Cref{Prop:Clef_Ext-Int} will allow us to conclude that $u-\zeta v=0$ almost everywhere on $\partial \Omega'\cap V$.
\par\smallskip
\begin{proposition}\label{Prop:Ext-Int_AutreCas}Let $p>1$ and $q$ be its conjugate component.
Suppose that $F$ satisfies \ref{H1}, \ref{H2}, and \ref{H3}. Let $\Omega_0\in\mathcal C$ be such that its boundary $\partial \Omega_0$ is a Jordan curve. Let  $\Gamma$ be a closed arc of $\mathbb T$ such that $\partial \Omega_0=F(\Gamma)$, and suppose that $F$ has an analytic extension to a neighborhood of $\Gamma$. Let $u$ and $v$ belong to $E^q(\Omega_0)$. If $u-\zeta v=0$ on a subset $Z$ of $\partial\Omega_0$ of positive measure, then $u=v=0$.
\end{proposition}

\begin{proof}
Let $\Gamma=\{e^{i\theta}\,;\,\alpha\le\theta\le\beta\},\, \alpha<\beta<\alpha+2\pi$, be such that $\partial\Omega_0=F(\Gamma)$ and let $\lambda_0=F(e^{i\alpha})$.
For every connected open neighborhood $V$ of $\lambda_0$, $\Sigma=V\cap\Omega_0$ is connected and the restriction of $\zeta$ to $\partial\Sigma$ cannot be extended continuously at the point $\lambda_0$.
As in the proof of \Cref{Prop:Clef_Ext-Int}, we obtain an open neighborhood $C$ of $\partial\Omega_0$ and a closed simple arc $\Delta$  with non empty interior having $\lambda_0$ as one of its extremities such that $\zeta$ has a bounded analytic extension to $C\setminus \Delta$.  Without loss of generality, we can  also suppose that $\partial (C\cap \Omega_0)\cap\Delta$ has exactly two points, which are the extremities of $\Delta$ (see \Cref{FIG23} below).
\begin{figure}[ht]
\includegraphics[page=20,scale=1.2]{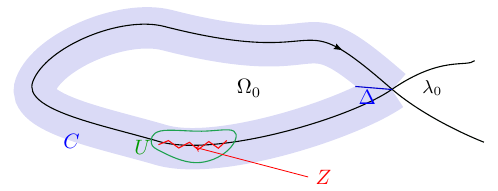}
            \caption{}\label{FIG23}\end{figure}

Replacing if necessary $Z$ by a smaller subset of $\partial \Omega_0$ with positive measure, we can assume that there exists a simply connected domain $U$ with $Z \subseteq U$ and $\overline{U} \subseteq (C\setminus\Delta)$.
Since the functions $u$ and $\zeta v$ both belong to $E^q(U\cap\Omega_0)$, and have boundary limits which coincide on $Z$, which is a subset of $\partial(U\cap\Omega_0)$ of positive measure, $u-\zeta v=0$ on $U\cap \Omega_0$. It follows that $u-\zeta v=0$ on $(C\cap\Omega_0)\setminus\Delta$ and hence $u-\zeta v=0$ almost everywhere on $\partial\Omega_0$. Thus the assumption of \Cref{Prop:Clef_Ext-Int} is satisfied. Hence $u=v=0$ and we are done.
\end{proof}
\par\smallskip
Note that if there is another self-intersection point $\lambda_1$ on $\partial\Omega_0$, then $\zeta$ will have an analytic extension to $V_{\lambda_1}\cap\Omega_0$, where $V_{\lambda_1}$ is some neighborhood of $\lambda_1$. Then, we can transfer on $\partial\Omega_0$ the fact that $u-\zeta v=0$ from one side of $\lambda_1$ to the other side. This observation will be used in \Cref{Subsection:MemeInd} and is the key point of the proof of \Cref{Th:MemeInd}.

As a consequence, we obtain the following variation of \Cref{Th:Complique_Ext-Int}:

\begin{theorem}\label{Th:Complique-Jordan_version}Let $p>1$.
Suppose that $F$ satisfy \ref{H1}, \ref{H2} and \ref{H3}, and let $\Omega,\Omega'\in\cal C$ be two adjacent components of $\sigma(T_F)\setminus F(\mathbb T)$ such that $|\w_F(\Omega)|<|\w_F(\Omega')|$. Suppose that $\partial\Omega'$ is a Jordan curve and there exists a closed arc $\Gamma$ of $\mathbb T$ such that $\partial \Omega'=F(\Gamma)$ and $F$ has an analytic extension to a neighborhood of $\Gamma$. Then $H_{\Omega'}(T_F)= H_\Omega(T_F)$.
\end{theorem}
\subsection{A case where the two components have the same winding number}\label{Subsection:MemeInd}
\begin{figure}[ht]
    \includegraphics[page=55,scale=1]{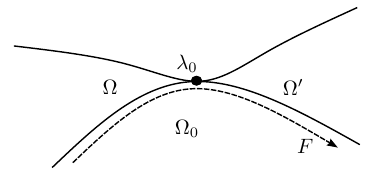}
    \caption{}\label{Fig21efj}
\end{figure}

Here is a last configuration that we are able to deal with. It is somewhat different from the preceding ones since the two connected components $\Omega$ and $\Omega'$ involved are not adjacent.
\begin{theorem}\label{Th:MemeInd}Let $p>1$.
Suppose that $F$ satisfies \ref{H1}, \ref{H2} and \ref{H3}, and let $\Omega,\Omega'\in\cal C$ be two connected components of $\sigma(T_F)\setminus F(\mathbb T)$ with $\w_F(\Omega)=\w_F(\Omega')$. Suppose that there exists a connected component $\Omega_0$ with $|\w_F(\Omega_0)|>|\w_F(\Omega)|$ and a point $\lambda_0\in\partial\Omega\cap\partial\Omega_0\cap \partial\Omega'$ such that
\begin{itemize}
    \item $\Omega$ and $\Omega'$ are both adjacent to $\Omega_0$ and $\partial\Omega\cap\partial\Omega_0$ and $\partial\Omega'\cap\partial\Omega_0$ both contain an arc with non empty interior having $\lambda_0$ as an extremity;
    \item the restriction of the function $\zeta$ to the boundary of $\Omega_0$ admits a bounded analytic extension to a neighborhood of the point $\lambda_0$.
\end{itemize}
 Then 
\[H_\Omega(T_F~)=~H_{\Omega'}(T_F).\]
\end{theorem}
\Cref{Fig21efj} illustrates a situation where \Cref{Th:MemeInd} can be applied. 
\begin{proof}
Let $V$ be a connected open neighborhood $V$ of $\lambda_0$ such that $\zeta$ admits a bounded analytic extension to $V$.
Also, since $\Omega$ and $\Omega'$ are both adjacent to $\Omega_0$ and $\partial\Omega\cap\partial\Omega_0$ and $\partial\Omega'\cap\partial\Omega_0$ both contain an arc with non empty interior having $\lambda_0$ as an extremity, it follows that both $V\cap \partial\Omega\cap\partial\Omega_0$ and $V\cap\partial\Omega'\cap\partial\Omega_0$ contain $\lambda_0$ and have positive measure.
Suppose that $g\in H^q$ vanishes on $H_\Omega(T_F)$. Writing $Ug=(u_j)_{0\le j\le N-1}$ with $u_j\in E^q(\Omega_j^+)$, we have $u_j=0$ on $\Omega$ for every $0\le j<|\w_F(\Omega)|.$ Since $|\w_F(\Omega_0)|>j+1$, $\Omega_0 \subseteq\Omega_{j+1}^+$, and thus
\[u_j^{int}-\zeta u_{j+1}^{int}~=~u_j^{ext}\quad\text{a.e. on }~ \partial\Omega_0\cap(\partial\Omega\cup\partial\Omega') \]
and $u_j^{ext}=0$ almost everywhere on $\partial\Omega\cap\partial\Omega_0$. So $u_j^{int}-\zeta u_{j+1}^{int}=0$ almost everywhere on $\partial\Omega\cap\partial\Omega_0$.
Now, the function $u_j- \zeta u_{j+1}$ belongs to $E^q(\Omega_0\cap V)$ and vanishes on $\partial\Omega_0\cap\partial\Omega\cap V$, which is a subset of positive measure of the boundary of $\Omega_0\cap V$. Hence $u_j-\zeta u_{j+1}$ is identically zero on $\Omega_0\cap V$, and in particular $u_j^{int}-\zeta u_{j+1}^{int}=0$ almost everywhere on $V\cap\partial\Omega_0\cap\partial\Omega'$. Since $u_j^{int}-\zeta u_{j+1}^{int}=u_j^{ext}$ almost everywhere on $\partial\Omega_0\cap\partial\Omega'$, this yields that $u_j^{ext}=0$ almost everywhere on $V\cap \partial\Omega_0\cap\partial\Omega'$, which is a subset of positive measure of $\partial\Omega'$, and so $u_j=0$ on $\Omega'.$
\par\smallskip
This being true for every $0\le j<|\w_F(\Omega)|=|\w_F(\Omega')|,$ it follows that $g$ vanishes on $H_{\Omega'}(T_F)$. Hence $H_{\Omega'}(T_F)  \subseteq H_\Omega(T_F)$. In the same way, $H_{\Omega}(T_F)  \subseteq H_{\Omega'}(T_F)$, and \Cref{Th:MemeInd} is proved.
\end{proof}

\subsection{Applying \Cref{Th:Simple_Int-Ext,Th:Complique_Ext-Int}}
We conclude this section by outlining the kind of argument  which will be used repeatedly in the proofs of our forthcoming results when applying our Theorems \ref{Th:Simple_Int-Ext} and \ref{Th:Complique_Ext-Int}.
Recall that $H_-(T_F)=\overline{\spa}\,[\ker(T_F-\lambda)\,;\,|\lambda|<1]$, that $H_+(T_F)=\overline{\spa}\,[\ker(T_F-\lambda)\,;\,|\lambda|>1]$, and that we have 
\begin{equation}\label{Eq:DensiteVP}
\overline{\spa}\,[H_\Omega(T_F)\,;\,\Omega\in\cal C]~=~H^p
\end{equation}
by \Cref{Prop:EigenvectorDenseYaku}.
Here is now the key argument which will be applied once we have managed to show, using either \Cref{Th:Simple_Int-Ext} or \Cref{Th:Complique_Ext-Int}, that $H_{\Omega'}(T_F) \subseteq H_\Omega(T_F)$ for two adjacent components $\Omega$ and $\Omega'$.

\begin{proposition}\label{PropositionP}
    Let $\Omega,\Omega'\in\cal C$ be two adjacent components. Suppose that
    \[H_{\Omega'}(T_F)~ \subseteq~ H_\Omega(T_F).\]
    \begin{enumerate}
        \item If
    \(\Omega\cap\mathbb D\neq\varnothing,\)
    then $H_{\Omega'}(T_F) \subseteq H_\Omega(T_F)\subseteq H_-(T_F)$.
    \par\smallskip
    \item If
    \(\Omega\cap(\C\setminus\overline{\mathbb D})\neq\varnothing,\)
    then $H_{\Omega'}(T_F) \subseteq H_\Omega(T_F)\subseteq H_+(T_F)$.
    \end{enumerate}
\end{proposition}

\begin{proof}
 If $\Omega\cap\mathbb D\neq \varnothing$, then $\Omega\cap\mathbb D$ has an accumulation point in $\Omega$. Thus, by \Cref{Prop:EigenvectorDenseYaku}, 
\begin{equation*}
\overline{\spa}\,[\ker(T_F-\lambda)\,;\,\lambda\in\Omega\cap \mathbb D]~=~H_\Omega(T_F), 
\end{equation*}
so that
    \(H_\Omega(T_F) \subseteq H_-(T_F).\)
This proves assertion (1). The proof of assertion (2) is exactly similar.
\end{proof}
Suitable assumptions will ensure that we can apply  \Cref{PropositionP} to pairs of components $(\Omega,\Omega')$, where $\Omega'$ varies of all of $\cal C$. By \Cref{Eq:DensiteVP}, this will yield that $H_+(T_F)=H_-(T_F)=H^p$, and will enable us to apply the Godefroy-Shapiro Criterion.

\section{Applications and examples}\label{Section:ApplicationCShc}
We already gave in \Cref{Section:CNforHc} a necessary condition for the hypercyclicity of $T_F$ under assumptions \ref{H1}, \ref{H2} and \ref{H3} on the symbol $F$: the unit circle $\mathbb T$ has to intersect every connected component of the interior of the spectrum of $T_F$ (\Cref{Th:CNforHC-Intersection}). Then, at the beginning of \Cref{Section:PassageFrontiere}, we stated a condition implying, in a rather straightforward manner, that $T_F$ is hypercyclic (\Cref{Th:SimpleConseqHc}): if
$\Omega\cap\mathbb T\neq\varnothing$ for \textit{every} $\Omega\in \cal C$, then $T_F$ is hypercyclic on $H^p$. In this section, we explore more deeply sufficient conditions for hypercyclicity, the general philosophy being the following: if we only know that $\Omega\cap\mathbb T\neq\varnothing$ for \textit{some} $\Omega\in \cal C$, can we deduce that $T_F$ is hypercyclic?

\subsection{The one-component case}
\begin{figure}[ht]
\includegraphics[page=21,scale=.9]{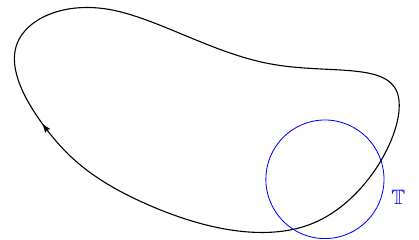}
\caption{}
\label{Fig22ljglkj}
\end{figure}

By \Cref{Th:CNforHC-Intersection,Th:SimpleConseqHc}, if $\sigma(T_F)\setminus F(\mathbb T)$ is connected, we obtain immediately a necessary and sufficient condition for the hypercyclicity of $T_F$.
\begin{theorem}\label{Th:ConnexeHC}Let $p>1$ and let $F$ satisfy \ref{H1}, \ref{H2} and \ref{H3}. Suppose that the set $\sigma(T_F)\setminus F(\mathbb T)$ is connected. Then, the following are equivalent:
\begin{enumerate}
    \item $T_F$ is hypercyclic on $H^p$;
    \item $\overset{\circ}{\sigma(T_F)}\cap\mathbb T\neq\varnothing$.
\end{enumerate}
\end{theorem}
\Cref{Fig22ljglkj} illustrates a situation where \Cref{Th:ConnexeHC} can be applied.
  
    \Cref{Th:ConnexeHC} allows us for instance to retrieve and extend to the $H^p$ setting the pioneering result of Shkarin, who characterized in \cite{Shkarin2012} the symbols $F$ of the form $F(e^{i\theta})=ae^{-i\theta}+b+ce^{i\theta},~a,b,c\in\mathbb C$, inducing a hypercyclic Toeplitz operator $T_F$ on $H^2$.
    \begin{theorem}\label{Th:ShkarinHp}
    Let $F(e^{i\theta})=ae^{-i\theta}+b+ce^{i\theta},~a,b,c\in\mathbb C$ with $|a|\neq|c|$, and let $p>1$. The following are equivalent:
    \begin{enumerate}
        \item $T_F$ is hypercyclic on $H^p$
        \item $|a|>|c|$ and $\cal E\cap \mathbb T\neq\varnothing$, where $\cal E$ is the interior of the (non-degenerate) elliptic curve $F(\mathbb T).$
    \end{enumerate}
    \end{theorem}
\begin{proof}
    Let us first show that (1)$\implies$(2). 
    Since $|a|\neq|c|$, $F$ satisfies \ref{H1} and \ref{H2}, and the curve $F(\mathbb T)$ is an ellipse with non-empty interior $\cal E$. 
    \begin{fact}\label{Fait:OrientationEllipse}
        We have $\w_F(\cal E)<0$ when $|a|>|c|$ and $\w_F(\cal E)>0$ when $|a|<|c|$.
    \end{fact}
    \begin{proof}Let $a=|a|e^{i\alpha}$ and $c=|c|e^{i\gamma}$, and consider the map $ \widetilde F$ defined by $ \widetilde F(e^{i\theta})=|a|e^{-i\theta}+|c|e^{i\theta}$. Then 
    \[F(e^{i\theta})~=~b+e^{i(\alpha+\gamma)/2} \widetilde F(e^{i(\theta+ (\gamma-\alpha)/2})\quad\text{ for every } e^{i\theta}\in\mathbb T.\]
    In particular $F(\mathbb T)$ is the image of $ \widetilde F(\mathbb T)$ by an affine transformation,  and thus these two curves have the same orientation. It is now easy to see geometrically how the orientation of the elliptic curve $ \widetilde F(\mathbb T)$ depends on $|a|$ and $|c|$: note that $ \widetilde F(1)=|a|+|c|$, $ \widetilde F(-1)=-(|a|+|c|)$ and $ \widetilde F(\mathbb T)$ is negatively oriented if and only if $|c|-|a|=\mathrm{Im}( \widetilde F(i))<0$.
Hence $F(\mathbb T)$ is negatively oriented (which is equivalent to the property that $\w_F(\cal E)<0$) if and only if $|a|>|c|$.
    \end{proof}
 By   \Cref{Fait:OrientationEllipse}, \Cref{Prop:OrientationForHC} and \Cref{Th:CNforHC-Intersection}, (2) is a necessary condition for the hypercyclicity of $T_F$ on $H^p$. The converse implication is a direct consequence of \Cref{Th:ConnexeHC}.
\end{proof}

\begin{remark}
When $|a|=|c|$, it was observed by Shkarin in  \cite{Shkarin2012} that $T_F$ is a normal operator on $H^2$, so that it cannot be hypercyclic on $H^2$. Thus in the case where $p=2$, the assumption that $|a|\neq|c|$ can be dispensed with in \Cref{Th:ShkarinHp}. This observation extends to the case where $p>2$. Indeed, if $T_F$ were hypercyclic on $H^p$, then it would be hypercyclic on $H^2$ as well (see \Cref{rem-p-r} for details). We do not know what happens in the case $1<p<2$. Note that when $|a|=|c|$, the curve $F(\mathbb T)$ is a segment and $\sigma(T_F)=F(\mathbb T)$. 
\end{remark}

\subsection{When there are no adjacent components}
\begin{figure}[ht]
\includegraphics[page=22,scale=1]{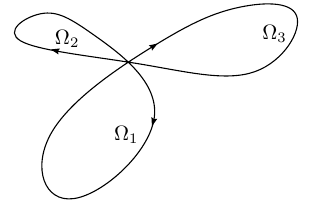}
                \caption{}
                \label{figure-petale}\end{figure}

In this subsection, we obtain a characterization of the hypercyclicity of $T_F$ when $F(\mathbb T)$ has the kind of form given in \Cref{figure-petale}.

    Here the interior of $\sigma(T_F)$ consists of the three disjoint connected components $\Omega_1,\Omega_2$ and $\Omega_3$ of $\sigma(T_F)\setminus F(\mathbb T)$. 
    Such symbols are called ``symbols with loops" by Clark, and studied for instance in \cites{Clark-1980-1, Clark-Toeplitz-2, Clark-quasisim}.
    In this kind of configuration, we have:
    \begin{theorem}\label{Th:CNS-IndMax1} Let $p>1$ and let $F$ satisfy \ref{H1}, \ref{H2} and \ref{H3}. Suppose that, for all pairs $(\Omega,\Omega')$ of distinct connected component of $\sigma(T_F)\setminus F(\mathbb T)$, the set $\overline{\Omega}\cap\overline{\Omega'}$ is finite. Then the following assertions are equivalent:
    \begin{enumerate}
        \item $T_F$ is hypercyclic on $H^p$;
        \item $\Omega\cap \mathbb T\neq\varnothing$ for every connected component $\Omega$ of $\sigma(T_F)\setminus F(\mathbb T)$.
    \end{enumerate}
    \end{theorem}
    
    \begin{proof}
    Our assumption implies that every $\Omega\in\cal C$ is a connected component of $\overset{\circ}{\sigma(T_F)}$, so the implication (1)$\implies$(2) is given by \Cref{Th:CNforHC-Intersection} and the implication (2)$\implies$(1) follows immediately from \Cref{Th:SimpleConseqHc}.
    \end{proof}
 
\subsection{A result involving maximal components}\label{Subsection:CompMax}
Let $F$ satisfy \ref{H1}, \ref{H2} and \ref{H3}. We say that $\Omega\in\cal C$ is a \textit{maximal component} of $\sigma(T_F)\setminus F(\mathbb T)$ if for every $\Omega'\in \cal C$ which is adjacent to $\Omega$, $|\w_F(\Omega')|<|\w_F(\Omega)|$.
In \Cref{Figure25} below, the maximal components are $\Omega_2$ and $\Omega_3$, and $\w_F(\Omega_2)=-2$ while $\w_F(\Omega_3)=-3$.
\begin{figure}[ht]
\includegraphics[page=23,scale=1]{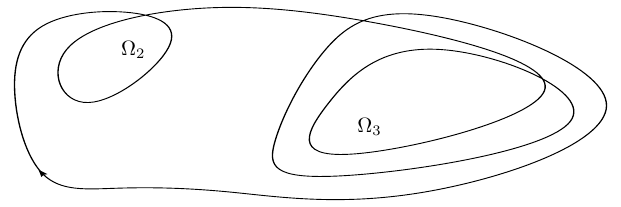}
            \caption{}\label{Figure25}\end{figure}
\par\smallskip
On the other hand, in Figure~\ref{figure-petale} above, the maximal components are $\Omega_1,\Omega_2,\Omega_3$ because there are no adjacent components to any of them. 
Using \Cref{Th:Simple_Int-Ext}, it is possible to weaken the assumptions of \Cref{Th:SimpleConseqHc} by requiring only that every maximal component of $\sigma(T_F)\setminus F(\mathbb T)$ intersects $\mathbb T$.

\begin{theorem}\label{Th:HcInstersCompMax}
Let $p>1$ and $F$ satisfy \ref{H1}, \ref{H2} and \ref{H3}. If $\Omega\cap\mathbb T\neq\varnothing$ for every maximal component $\Omega$ of $\sigma(T_F)\setminus F(\mathbb T)$, then $T_F$ is hypercyclic on $H^p$.
\end{theorem}

\begin{proof}
Let $\Omega\in\cal C$. We claim:
\begin{fact}\label{Fait:PassageVersCompMax}There exists a sequence $(\Omega_n)_{0\le n\le r}$ of elements of $\cal C$ such that:
\begin{itemize}
    \item $\Omega_0=\Omega$;
    \item for each $0\le n< r$, $\Omega_n$ and $\Omega_{n+1}$ are adjacent and $|\w_F(\Omega_n)|<|\w_F(\Omega_{n+1})|$;
    \item $\Omega_r$ is a maximal component of $\sigma(T_F)\setminus F(\mathbb T)$.
\end{itemize}

\end{fact}
\begin{proof}[Proof of \Cref{Fait:PassageVersCompMax}]
Starting from $\Omega_0=\Omega$, we choose, if possible, $\Omega_1$ adjacent to $\Omega_0$ such that $|\w_F(\Omega_1)|>|\w_F(\Omega_{0})|$. If there is no such $\Omega_1\in\mathcal{C}$, then $\Omega_0$ is already a maximal component. In the opposite situation we choose, if possible, $\Omega_2$ adjacent to $\Omega_1$ such that $|\w_F(\Omega_2)|>|\w_F(\Omega_{1})|$, etc. The process stops after finitely many steps, since $\sigma(T_F)\setminus F(\mathbb T)$ has a finite number of connected components. If the process stops after we have chosen $\Omega_r$, i.e. if there is no $\Omega_{r+1}\in\mathcal{C}$ adjacent to
$\Omega_r$ with $|\w_F(\Omega_{r+1})|>|\w_F(\Omega_{r})|$, this means that $\Omega_r$ is a maximal component.
\end{proof}
\par\smallskip
Starting from a given component $\Omega\in\cal C$, let $(\Omega_n)_{0\le n\le r}$ be given by \Cref{Fait:PassageVersCompMax}. Since $\Omega_n$ and $\Omega_{n+1}$ are adjacent and $|\w_F(\Omega_n)|<|\w_F(\Omega_{n+1})|$, \Cref{Th:Simple_Int-Ext} applies and
\[H_{\Omega_n}(T_F)~ \subseteq~ H_{\Omega_{n+1}}(T_F)\quad\text{for every}~0\le n<r.\]
In particular, $H_\Omega(T_F)\subseteq H_{\Omega_r}(T_F)$.
Since $\Omega_r$ is a maximal component of $\sigma(T_F)\setminus F(\mathbb T)$, by hypothesis, $\Omega_r\cap\mathbb T\neq\varnothing$, and thus by \Cref{PropositionP}, we have 
\[
H_\Omega(T_F) ~\subseteq~ H_{\Omega_{r}}(T_F) ~\subseteq~ H_+(T_F)\quad\text{and}\quad H_\Omega(T_F) ~\subseteq~ H_{\Omega_{r}}(T_F) ~\subseteq~ H_-(T_F).
\]
This being true for \textit{all} connected components $\Omega\in\cal C$, it follows from \Cref{Prop:EigenvectorDenseYaku} that $H_+(T_F)=H_-(T_F)=H^p$. By the Godefroy-Shapiro Criterion, $T_F$ is hypercyclic on $H^p$.
\end{proof}
\par\smallskip
When $\Omega_0$ is connected to $\Omega_r$ by a chain of adjacent components $(\Omega_n)_{0\le n\le r}$ satisfying $|\w_F(\Omega_n)|<|\w_F(\Omega_{n+1})|$ for every $0\le n<r$, we write
$$\Omega_0 ~\overset{I.W.}{\longrightarrow}~\Omega_{r}.$$

Note that when $\Omega_0 \overset{I.W.}{\longrightarrow}\Omega_{r}$, then \Cref{Th:Simple_Int-Ext} implies that 
\begin{equation}\label{IWcondition-inclusion-spectrale}
H_{\Omega_0}(T_F) ~\subseteq~ H_{\Omega_r}(T_F).
\end{equation}
\par\smallskip
In spirit, the proof of \Cref{Th:HcInstersCompMax} can be summarized as follows: knowing that a certain function $u\in E_F^q$ vanishes on a certain maximal component $\Omega_r\in\mathcal{C}$, we need to deduce that $u$ vanishes on another component $\Omega_0$ with $\Omega_0 \overset{I.W.}{\longrightarrow}\Omega_{r}$. 
Letting $\gamma_n$ be an arc contained in the common boundaries of $\partial \Omega_n$ and $\partial\Omega_{n+1},~0\le n<r$, where $(\Omega_n)_{0\le n\le r}$ is a chain of adjacent components connecting $\Omega_0$ and $\Omega_r$,  we thus see that the proof of \Cref{Th:Simple_Int-Ext} shows exactly the following: using the boundaries relations in \Cref{eq:boundary-conditions-ErF} of $E_F^q$, the fact that $u$ vanishes on $\Omega_r$ goes over from the interior component of $\gamma_{r-1}$ (i.e. $\Omega_r$) to its exterior component, $\Omega_{r-1}$, then from the interior component of $\gamma_{r-2}$ (i.e. $\Omega_{r-1}$) to its exterior component $\Omega_{r-2}$, etc. until we reach $\Omega_0$.
 \par\smallskip
\begin{example}
For instance, here is what the proof of \Cref{Th:Simple_Int-Ext} tells us  in the situation described in \Cref{Fig25-1}:
\begin{figure}[ht]
\includegraphics[page=24,scale=1]{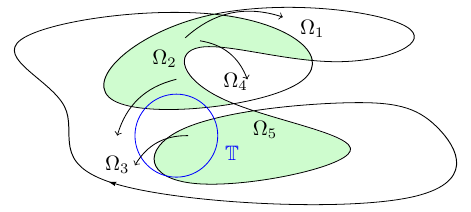}
\caption{}\label{Fig25-1}
\end{figure}

Note that $\Omega_0^+=\Omega_1\cup\Omega_2\cup\Omega_3\cup\Omega_4\cup\Omega_5$ and $\Omega_1^+=\Omega_2\cup\Omega_5$, and the maximal components are $\Omega_2$ and $\Omega_5$. Starting from $u_0\in E^q(\Omega_0^+)$ and $u_1\in E^q(\Omega_1^+)$, and knowing that $u_0=u_1=0$ on $\Omega_2\cup\Omega_5$, we wish to deduce that $u_0$ is identically zero, using the boundary relation $u_0^{int}-\zeta u_1^{int}=u_0^{ext}$ almost everywhere on $\partial\Omega_1^+=\partial\Omega_2\cup\partial\Omega_5$. The proof of \Cref{Th:Simple_Int-Ext} goes as follows:
\begin{enumerate}[i)]
    \item since $u_0=u_1=0$ on $\Omega_2$, $u_0=0$ on $\Omega_1$;
    \item since $u_0=u_1=0$ on $\Omega_2$, $u_0=0$ on $\Omega_3$ (alternatively, we could have used that since $u_0=u_1=0$ on $\Omega_5$, $u_0=0$ on $\Omega_3$);
    \item since $u_0=u_1=0$ on $\Omega_2$, $u_0=0$ on $\Omega_4$.
\end{enumerate}
Hence $u_0=0$ on $\Omega_1\cup\Omega_3\cup\Omega_4$, and we are done.
\end{example}
The assumptions of
\Cref{Th:HcInstersCompMax} can be weakened to yield the following result, which is similar in spirit to \cite[Th. 1.5]{AbakumovBaranovCharpentierLishanskii2021}; see also the forthcoming \Cref{th 7.4}. 

\begin{figure}[ht]
    \includegraphics[page=54,scale=1]{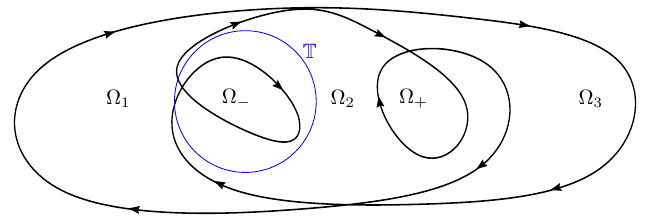}
    \caption{}\label{FIG26dke}
\end{figure}

\begin{theorem}\label{Th:HcInstersCompMax-faible}
Let $p>1$, and let $F$ satisfy assumptions \ref{H1}, \ref{H2}, \ref{H3} and \ref{H4}. 
Suppose that $\mathbb T$ intersects every maximal component $\Omega_0\in\mathcal{C}$  with $\w_F(\Omega_0)=-1$, and that for every non-maximal component $\Omega\in\mathcal{C}$, there exist two maximal components
$\Omega_{+}$ and $\Omega_{-}$ in $\mathcal{C}$ such that
\begin{itemize}
    \item $\Omega_-\cap\mathbb D\neq\varnothing \text{ and } \Omega \overset{I.W.}{\longrightarrow}\Omega_{-}$;
    \item $\Omega_+\cap(\mathbb C\setminus\overline{\mathbb D})\neq\varnothing  \text{ and }   \Omega \overset{I.W.}{\longrightarrow}\Omega_{+}$.
\end{itemize}
Then $T_F$ is hypercyclic on $H^p$.
\end{theorem}
Note that when a maximal component $\Omega_0\in\cal C$ satisfies $\w_F(\Omega_0)=-1$, it is a connected component of the interior of ${\sigma(T_F)}$, and so the condition $\Omega_0\cap\mathbb T$ is necessary for the hypercyclicity of $T_F$ by \Cref{Th:CNforHC-Intersection}.
Note also that, for a function $F$ satisfying assumptions \ref{H1}, \ref{H2}, \ref{H3} and \ref{H4}, the condition of \Cref{Th:HcInstersCompMax-faible} is weaker than the condition of \Cref{Th:HcInstersCompMax}. Indeed, in \Cref{Th:HcInstersCompMax}, we require that every maximal connected component meets $\mathbb T$, and thus intersects both $\mathbb D$ and $\mathbb C\setminus\overline{\mathbb D}$. It is not the case in the situation described in \Cref{FIG26dke}. 

On this figure, we see that the maximal connected components are $\Omega_-$ and $\Omega_+$ and that $|\w_F(\Omega_+)|=|\w_F(\Omega_-)|=3$. We have  $\Omega_- \subseteq \mathbb D$ and $\Omega_+ \subseteq\mathbb C\setminus\overline{\mathbb D}$. Thus the hypothesis of \Cref{Th:HcInstersCompMax-faible} are satisfied, because we can go up from each one of the non-maximal connected components $\Omega_1,\Omega_2,\Omega_3$ to $\Omega_+$ and $\Omega_-$ via adjacent components, that is $\Omega_k\overset{I.W.}{\longrightarrow}\Omega_{\pm}$, $k=1,2,3$.

\begin{proof}[Proof of \Cref{Th:HcInstersCompMax-faible}]
Let $g\in H^q$ vanish on $H_-(T_F)$, and let $\Omega$ be a non-maximal component of $\sigma(T_F)\setminus F(\mathbb T)$. According to our hypothesis, there exists a maximal component $\Omega_-$ such that $\Omega_-\cap\mathbb D\neq\varnothing$ and 
$\Omega \overset{I.W.}{\longrightarrow}\Omega_{-}$. Then, by \Cref{IWcondition-inclusion-spectrale}, $H_\Omega(T_F) \subseteq H_{\Omega_-}(T_F)$ and since $H_{\Omega_-}(T_F)\subseteq H_-(T_F)$ by \Cref{Prop:EigenvectorDenseYaku}, 
\[
H_{\Omega}(T_F) ~\subseteq~ H_-(T_F).
\]
Thus $g$ vanishes on $H_\Omega(T_F)$ for every non-maximal component $\Omega\in\mathcal C$.
Also, our hypothesis implies that $g$ vanishes on $H_\Omega(T_F)$ for all maximal components $\Omega\in\mathcal{C}$ with $\w_F(\Omega)=-1$.
\par\smallskip
Let $Ug=(u_j)_{0\le j\le N-1}$, where $U$ is the operator defined in \Cref{Section:Yakubovich-court} and $u_j\in E^q(\Omega_j^+)$ for every $0\le j\le N-1$. This means that for every non-maximal  $\Omega\in\mathcal{C}$ such that $j<|\w_F(\Omega)|$, $u_j=0$ on $\Omega$, and moreover $u_0=0$ on $\Omega$  for any maximal component $\Omega\in\mathcal{C}$ with $\w_F(\Omega)=-1$. It remains to prove that for every maximal component $\Omega\in\mathcal{C}$ with $|\w_F(\Omega)|\ge2$ and for every $0\le j<|\w_F(\Omega)|$, we have $u_j=0$ on $\Omega$.

Let thus $\Omega\in\mathcal{C}$ be a maximal component with $|\w_F(\Omega)|\ge2$. Then every adjacent component $\Omega'$ of $\Omega$ is non-maximal. So for every $0\le j<|\w_F(\Omega)|-1$ and for every adjacent component $\Omega'$ of $\Omega$, $u_j=0$ on $\Omega'$. 
The boundary conditions on the functions $u_j$ applied on $\partial\Omega$ imply that for every $0\le j< |\w_F(\Omega)|-1$, $u_j-\zeta u_{j+1}\in E^q(\Omega)$ enjoys the following property: for every adjacent component $\Omega'$ of $\Omega$, $u_j-\zeta u_{j+1}=0$ almost everywhere on $\partial\Omega'\cap\partial \Omega$. So $u_j-\zeta u_{j+1}=0$ almost everywhere on $\partial\Omega$. By \Cref{Fait:H4=>(P)iii} and \Cref{Prop:Clef_Ext-Int} applied to any self-intersection point $\lambda_0\in\partial\Omega$ where $\zeta$ is discontinuous (such a point does exist since $|\w_F(\Omega)|\ge2$), 
$u_j=u_{j+1}=0$ on $\Omega$.
We have thus proved that $g$ vanishes on $H_\Omega(T_F)$.
\par\smallskip
In conclusion, we have shown that $g$ vanishes on $H_\Omega(T_F)$ for all $\Omega\in\mathcal{C}$, and so $g=0$ by \Cref{Prop:EigenvectorDenseYaku}. Hence $H_-(T_F)=H^p$.
Similarly, we have $H_+(T_F)=H^p$, and thus, by the Godefroy-Shapiro Criterion, $T_F$ is hypercyclic on $H^p$.
\end{proof}
Observe the following: whenever $\Omega\in\cal C$ is a maximal component such that $|\w_F(\Omega)|\ge2$, there exists a self-intersection point $\lambda_0\in\partial\Omega$ such that assumption (a) in \Cref{Prop:Clef_Ext-Int} is satisfied. Assumption \ref{H4} in \Cref{Th:HcInstersCompMax-faible} is here to ensure, by \Cref{Fait:H4=>(P)iii}, that assumption (b) in \Cref{Prop:Clef_Ext-Int} is also satisfied.
\par\smallskip
In the next subsection, we rely on \Cref{Th:Complique_Ext-Int} to present a general situation where it is possible to extend the nullity from the exterior component of a boundary arc to the interior component (which is harder than going from the interior component to the exterior component).
\subsection{A case with a necessary and sufficient condition}
Here is the main theorem of this section, which provides a necessary and sufficient condition for the hypercyclicity of $T_F$ under  fairly general assumptions on $F$. Recall that the function $F$ admits an inverse $F^{-1}$, defined on $F(\T)\setminus\mathcal{O}$, where  $\mathcal{O}$ is the set of self-intersection points of the curve, and that $\zeta=1/F^{-1}$ on $F(\T)\setminus\mathcal{O}$. 

\begin{theorem}\label{Th:CNShc1}
Let $p>1$, and let $F$ satisfy \ref{H1}, \ref{H2} and \ref{H3}. Suppose that every pair $(\Omega,\Omega')$ of adjacent connected component of $\sigma(T_F)\setminus F(\mathbb T)$ with $|\w_F(\Omega)|<|\w_F(\Omega')|$ satisfies the property \ref{cond(G)}. 
\par\smallskip
Then the following assertions are equivalent:
\begin{enumerate}
    \item $T_F$ is hypercyclic on $H^p$;\par\smallskip
    \item $\Theta\cap\mathbb T\neq\varnothing$ for every connected component $\Theta$ of $\overset{\circ}{\sigma(T_F)}$.
\end{enumerate}
\end{theorem}
\par\smallskip
Note that, apart from the regularity hypothesis on the function $F$ given by the assumption (iii) in \ref{cond(G)}, that was introduced in \Cref{subsection:Ext-to-Int}, all the assumptions of this theorem are purely geometric.
The following example illustrates how \Cref{Th:CNShc1} can be applied in order to deduce the hypercyclicity (or non-hypercyclicity) of an operator whose symbol satisfies the assumption \ref{cond(G)}.

\begin{example}
Let $F$ satisfy \ref{H1}, \ref{H2}, \ref{H3}  and \ref{H4}, and suppose that the curve $F(\mathbb T)$ is represented as in \Cref{FIG30}.
\begin{figure}[ht]
\includegraphics[page=29,scale=1]{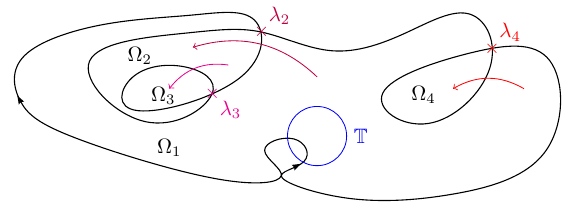}
\caption{}\label{FIG30}\end{figure}
                
    Then this function $F$ satisfies \ref{cond(G)}, and the intersection points that we use in \ref{cond(G)} are  $\lambda_2$ for the pair $(\Omega_1,\Omega_2)$, $\lambda_3$ for the pair $(\Omega_2,\Omega_3)$ and  $\lambda_4$ for the pair $(\Omega_1,\Omega_4)$.
  Since $\mathbb T$ intersects $\Omega_1$, $\mathbb T$ intersects $\overset{\circ}{\sigma(T_F)}$ which is connected in this case. So, by \Cref{Th:CNShc1}, $T_F$ is hypercyclic on $H^p$ for every $p>1$. 
\end{example}
Let us now get back to the proof of \Cref{Th:CNShc1}.
As can be expected, the role of assumption \ref{cond(G)} is to allow us to propagate a condition of the form ``$u=0$'' from $\Omega$, which is the exterior component of a suitable arc $\gamma \subseteq\partial\Omega\cap\partial\Omega'$, to its interior component $\Omega'$, using the point $\lambda_0$ and \Cref{Th:Complique_Ext-Int}. 

\begin{proof}[Proof of \Cref{Th:CNShc1}]
    Thanks to \Cref{Th:CNforHC-Intersection}, we only need to prove that (2)$\implies$(1). Let $\Theta$ be a connected component of $\overset{\circ}{\sigma(T_F)}$.
\par\smallskip
    Our strategy of proof will be to connect any two components $\Omega, \Omega'\in\cal C$ with $\Omega, \Omega' \subseteq \Theta$ via a finite sequence of adjacent  components of $\mathcal C$ contained  in $\Theta$. Given $\Omega,\Omega'\in\cal C$, we write $\Omega\longleftrightarrow\Omega'$ if there exists a finite sequence $(\Omega_n)_{0\le n\le r}$ of elements of $\cal C$  with $\Omega_n \subseteq \Theta$, $0\leq n\leq r$ such that $\Omega_0=\Omega$, $\Omega_r=\Omega'$ and $\Omega_n$ is adjacent to $\Omega_{n+1}$ for each $0\le n<r$. 
    \begin{lemma}\label{Lemma51}
    For any $\Omega,\Omega'\in \cal C$ with $\Omega,\Omega' \subseteq \Theta$, we have $\Omega\longleftrightarrow\Omega'$.
    \end{lemma}
    \begin{proof}[Proof of \Cref{Lemma51}]
    Choose $\lambda\in \Omega$ and $\lambda'\in\Omega'$. Since $\Theta\setminus \mathcal O$ is a connected open subset of $\overset{\circ}{\sigma(T_F)}$, it is path-wise connected.
    Let $\gamma:[0,1]\longmapsto \Theta\setminus\mathcal O$ be a continuous map such that $\gamma(0)=\lambda$ and $\gamma(1)=\lambda'$. Without loss of generality, we can suppose that the set 
    \[
    Z~:=~\{t\in [0,1]\,;\, \gamma(t)\in F(\mathbb T)\}
    \]
    consists of finitely many points $0<t_0<t_1<\dots<t_r<1$, and that $\gamma(t_n)\in F(\mathbb T)\setminus\mathcal O$ for every $0\le n\le r<1$. Set $t_{-1}=0$ and $t_{r+1}=1$. For every $0\le n\le r$ and for $t_n<t<t_{n+1}$, the point $\gamma(t)$ remains in a given connected component of $\sigma(T_F)\setminus F(\mathbb T)$, which we call $\Omega_n$. Also $\Omega_0=\Omega$, $\Omega_r=\Omega'$ and $\gamma(0)\in\Omega_0$, $\gamma(1)\in\Omega_r$. Since $\partial\Omega_n\cap\partial\Omega_{n+1}$ contains the point $t_n$, which is not a self-intersection point of the curve $F(\mathbb T)$, $\Omega_n$ and $\Omega_{n+1}$ are adjacent for every $0\leq n<r$. It follows that $\Omega\longleftrightarrow\Omega'$.
    \end{proof}
    
    The next step in the proof of \Cref{Th:CNShc1} is
    \begin{lemma}\label{Lemma52}
    If $\Omega,\Omega'\in\cal C$ with $\Omega,\Omega' \subseteq \Theta$, then $H_\Omega(T_F)=H_{\Omega'}(T_F)$.
    \end{lemma}
    
    \begin{proof}[Proof of \Cref{Lemma52}]
    By \Cref{Lemma51}, we have $\Omega\longleftrightarrow\Omega'$, so let $(\Omega_n)_{0\le n \le r}$ a finite sequence of adjacent components such that $\Omega_0=\Omega$ and $\Omega_r=\Omega'$. For any $0\le n<r$, the pair $(\Omega_n,\Omega_{n+1})$ satisfies the property \ref{cond(G)}, so by \Cref{Th:Complique_Ext-Int} we deduce that $H_{\Omega_n}(T_F)=H_{\Omega_{n+1}}(T_F)$ and thus
    \[H_{\Omega}(T_F)~=~H_{\Omega_0}(T_F)~=~\dots~=~H_{\Omega_r}(T_F)~=~H_{\Omega'}(T_F).\qedhere\]
    \end{proof}
    
    Since $\Theta\cap\T\neq\varnothing$, there exists $\Omega_0\in\cal C$, $\Omega_0 \subseteq \Theta$, such that $\Omega_0\cap(\mathbb C\setminus\overline{\mathbb D})\neq\varnothing$. Then, according to \Cref{Prop:EigenvectorDenseYaku}, $H_{\Omega_0}(T_F) \subseteq H_+(T_F)$.
      By \Cref{Lemma52}, for every $\Omega\in\mathcal C$, $\Omega \subseteq \Theta$, we have  $H_\Omega(T_F)=H_{\Omega_0}(T_F) \subseteq H_+(T_F)$ and thus 
    \[\overline{\spa}\,\big[H_\Omega(T_F),\Omega\in\cal C,\Omega \subseteq \Theta\big] ~\subseteq~ H_+(T_F).\]
    Since this reasoning applies to all connected components $\Theta$ of $\overset{\circ}{\sigma(T_F)}$, it follows from \Cref{Prop:EigenvectorDenseYaku} that $H_+(T_F)=H^p$. In the same way, connecting any $\Omega\in\mathcal{C}$ with $\Omega  \subseteq \Theta$ to a fixed component $\Omega_1  \subseteq \Theta$ such that $\Omega_1\cap\D\neq\varnothing$, we obtain that $H_-(T_F)= H^p$. By the Godefroy-Shapiro Criterion, $T_F$ is hypercyclic on $H^p$.
\end{proof}

There are some natural situations where the hypothesis of \Cref{Th:CNShc1} are not satisfied, and where we need other tools in order to be able to decide whether the operator $T_F$ is hypercyclic or not. \Cref{figure:inclusionellipse-cercle} represents an example of a curve for which assumption \ref{cond(G)} of \Cref{Th:Complique_Ext-Int} is not satisfied:
\begin{figure}[ht]
\includegraphics[page=30,scale=1]{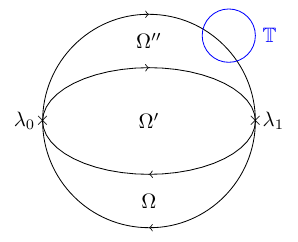}
\caption{}
\label{figure:inclusionellipse-cercle}
\end{figure}
            
Indeed, in the situation of \Cref{figure:inclusionellipse-cercle}, we cannot find a neighborhood $V$ of $\lambda_0$ such that $V\cap \partial\Omega'=V\cap\partial\Omega\cap\partial\Omega'$ (and similarly for $\lambda_1$), so that property (i) of \ref{cond(G)} does not hold. Note however that for any continuous parametrization $F$ of the curve (given by \Cref{figure:inclusionellipse-cercle}), the ellipse is the image of a single closed subarc of $\mathbb T$. So the following more general version of \Cref{Th:CNShc1} solves cases like this one.

\begin{theorem}\label{Th:CNShc2}
Let $p>1$, and let $F$ satisfy \ref{H1}, \ref{H2} and \ref{H3}. Suppose that, for every pair $(\Omega,\Omega')$ of adjacent connected components of $\sigma(T_F)\setminus F(\mathbb T)$ with $|\w_F(\Omega)|<|\w_F(\Omega')|$, one of the following properties holds:
\begin{enumerate}[(a)]
    \item \label{(aa)} the pair $(\Omega,\Omega')$ satisfies the condition \ref{cond(G)}.
    \item \label{(bb)} $\partial \Omega'$ is a Jordan curve and there exists a closed arc $\Gamma$ of $\mathbb T$ such that $\partial \Omega'=F(\Gamma)$ and $F$ has an analytic extension to a neighborhood of $\Gamma$.
\end{enumerate}
Then the following assertions are equivalent:
\begin{enumerate}
    \item $T_F$ is hypercyclic on $H^p$;
    \item $O\cap\mathbb T\neq\varnothing$ for every connected component $O$ of $\overset{\circ}{\sigma(T_F)}$.
\end{enumerate}    
\end{theorem}

\begin{proof}
The proof is exactly the same as that of \Cref{Th:CNShc1}, except for the proof of \Cref{Lemma52}. To show that $H_{\Omega_n}(T_F)=H_{\Omega_{n+1}}(T_F)$, we use either \Cref{Th:Complique_Ext-Int} if we are in case (a) or \Cref{Th:Complique-Jordan_version} if we are in case (b).
\end{proof}

\begin{example}
As an illustration, we describe on \Cref{FIG32} below the connections between adjacent components which allow us to apply \Cref{Th:CNShc2} in this case, and to show that $T_F$ is hypercyclic on $H^p$.
\begin{figure}[ht]
\includegraphics[page=31,scale=.9]{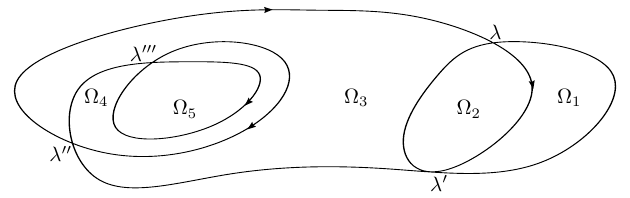}
 \caption{}\label{FIG32}\end{figure}
 
    The pairs $(\Omega,\Omega')$ of adjacent  components with $|\w_F(\Omega)|<|\w_F(\Omega')|$ are
    \begin{itemize}
        \item $(\Omega_1,\Omega_2)$, where \ref{(bb)} applies (but not \ref{(aa)});
        \item $(\Omega_3,\Omega_2)$, where \ref{(bb)} applies (but not \ref{(aa)});
        \item $(\Omega_3,\Omega_4)$, where \ref{(aa)} applies using the point $\lambda_0=\lambda''$ (but the point $\lambda'''$ does not satisfy assumption \ref{(aa)}, and \ref{(bb)} does not apply either);
        \item $(\Omega_4,\Omega_5)$, where either \ref{(bb)} or \ref{(aa)} applies using the point $\lambda_0=\lambda'''$.
    \end{itemize}
\end{example}
        
\begin{example}\label{ex 5.13}
    Still, there are examples where the philosophy of \Cref{Th:CNShc2} applies although its assumptions are not satisfied.  For instance, consider a function $F$ satisfying \ref{H1}, \ref{H2}, \ref{H3} and analytic on a neighborhood of $\mathbb T$ such that $F(\mathbb T)$ is represented as in \Cref{Fig30-98}:
\begin{figure}[ht]
\includegraphics[page=32,scale=.9]{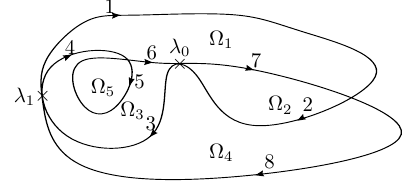}
\caption{}\label{Fig30-98}
\end{figure}

    Because of the components $\Omega_2$ and $\Omega_3$, a function $F$ such that $F(\mathbb T)$ has this representation does not satisfy the hypothesis of \Cref{Th:CNShc2}. Nonetheless, it is still true that $T_F$ is hypercyclic on $H^p$ if and only if $\mathbb T$ intersects the interior of ${\sigma(T_F)}$.
    Indeed, let $g\in H^q$ and suppose for instance that $g$ vanishes  on $H_{\Omega_1}(T_F)$. Since the restriction of $\zeta$ to $\partial\Omega_3$ has a continuous extension to the point $\lambda_1$, we can apply  \Cref{Th:MemeInd} to deduce that $g$ also vanishes on $H_{\Omega_4}(T_F)$. We can now apply \Cref{Prop:Clef_Ext-Int} to $\Omega_2$ and $\Omega_3$, using the point $\lambda_0$. Indeed, let $U:H^q\to E^q(\Omega_0^+)\oplus E^q(\Omega_1^+)\oplus E^q(\Omega_2^+)$ be the operator defined in \Cref{Section:Yakubovich-court} where $\Omega_0^+=\Omega_1\cup\Omega_2\cup\Omega_3\cup\Omega_4\cup\Omega_5$,  $\Omega_1^+=\Omega_2\cup\Omega_3\cup\Omega_5$, and $\Omega_2^+=\Omega_5$ and write $Ug=(u_0,u_1,u_2)$ with 
\[
u_j(\lambda)~=~\dual{h_{\lambda,j}}{g}\quad \text{ for every } \lambda\in\Omega_j^+, j=0,1,2.
\]
Since $g$ vanishes on $H_{\Omega_1}(T_F)\cup H_{\Omega_4}(T_F)$, we know that $u_0=0$ on $\Omega_1\cup\Omega_4$. Then, according to the boundary relation in \Cref{eq:boundary-conditions-ErF}, we get 
\begin{align*}
&u_0^{int}-\zeta u_1^{int}~=~0\quad\text{a.e. on }\partial\Omega_3\cap(\partial\Omega_1\cup\partial\Omega_4)  ;\\
&u_0^{int}-\zeta u_1^{int}~=~0\quad\text{a.e. on }\partial\Omega_2.
\end{align*}
Applying \Cref{Prop:Clef_Ext-Int} to the point $\lambda_0$ yields that $u_0=u_1=0$ on $\Omega_2\cup\Omega_3$, and thus $g$ vanishes on $H_{\Omega_2}(T_F)\cup H_{\Omega_3}(T_F)$. Finally, since $\partial \Omega_5$ is a Jordan curve and the image by $F$ of a single closed subarc of $\mathbb T$, we can apply \Cref{Th:Complique-Jordan_version} and get that $H_{\Omega_5}(T_F) \subseteq H_{\Omega_3}(T_F)$. Thus $g$ vanishes on $H_{\Omega_5}(T_F)$.  Hence $g$ vanishes on $H_{\Omega_i}(T_F)$ for every $1\le i \le 5$, and \Cref{Prop:EigenvectorDenseYaku} implies that $g=0$. 
\par\smallskip
If at the beginning, $g$ vanishes  on another space $H_{\Omega_j}(T_F)$, using similar arguments as well as \Cref{Th:Simple_Int-Ext}, we can also conclude that $g=0$. Therefore one can conclude from the Godefroy-Shapiro Criterion that $T_F$ is hypercyclic on $H^p$ as soon as $\mathbb T$ intersects the interior of ${\sigma(T_F)}$, which is connected in this case.
\end{example}
    \subsection{The case where $\max|\w_F(\Omega)|=2$}
    In the case where $\max\{|\w_F(\Omega)|\,;\,\Omega\in\cal C\}=2$, some of our results above take a simpler form.
    \begin{theorem}\label{Th:Case_wind2}
    Let $p>1$, and let $F$ satisfy \ref{H1}, \ref{H2} and \ref{H3}. Suppose that $\max\{|\w_F(\Omega)|\,;\,\Omega\in\cal C\}=2$, that for every $\Omega\in\cal C$ with $|\w_F(\Omega)|=2$, $\partial \Omega$ is a Jordan arc which is the image by $F$ of a single closed subarc $\Gamma$ of $\mathbb T$, and that $F$ has an analytic extension to a neighborhood of $\Gamma$. Then the following assertions are equivalent:
    \begin{enumerate}
    \item $T_F$ is hypercyclic on $H^p$;
    \item $O\cap\mathbb T\neq\varnothing$ for every connected component $O$ of $\overset{\circ}{\sigma(T_F)}$.
\end{enumerate} 
    \end{theorem}
    \begin{figure}[ht]
\includegraphics[page=33,scale=.85]{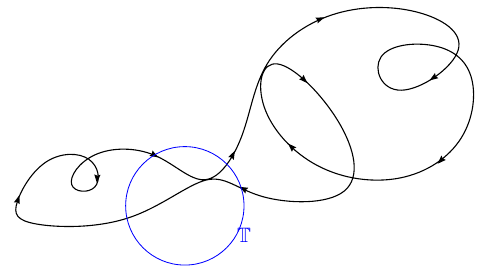}
\caption{}
\label{Fig31kglc}
\end{figure}
    \Cref{Fig31kglc} illustrates a situation where \Cref{Th:Case_wind2} can be applied.
 \Cref{Th:Case_wind2} also applies to the case considered in \Cref{SubSection:Example2Circles}, where $F(\mathbb T)$ consists of two tangent circles, the inner one being the unit circle (see \Cref{Fig13}). We will get back to this example in \Cref{Section:AutresResult+Questions}.
  Note that \Cref{Th:Case_wind2} is an immediate consequence of \Cref{Th:CNShc2} since any pair $(\Omega,\Omega')$ of elements of $\cal C$ with $|\w_F(\Omega)|<|\w_F(\Omega')|$ is such that $|\w_F(\Omega)|=1$ and $|\w_F(\Omega')|=2$, and then by hypothesis, $\partial\Omega'$ is a Jordan curve which is the image by $F$ of a single closed subarc of $\mathbb T$. Hence we are in case (b) of \Cref{Th:CNShc2}.
  \par
\begin{figure}[ht]
\includegraphics[page=34,scale=.9]{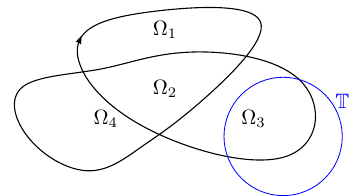}
 \caption{}\label{figure:dsfsdfs1323}
 \end{figure}

\smallskip
 \Cref{figure:dsfsdfs1323} illustrates a situation where $\max\{|\w_F(\Omega)|\,;\,\Omega\in\cal C\}=2$ but \Cref{Th:Case_wind2}  does not apply.
    
\begin{figure}[ht]
    \includegraphics[page=57,scale=1.1]{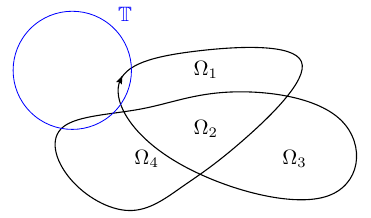}
    \caption{}\label{Fig33kfjdbk}
\end{figure}

If $\mathbb T$ does not intersect the component $\Omega_2$ of $\sigma(T_F)\setminus F(\T)$ with $|\w_F(\Omega_2)|=2$, but  intersects only one of the other components $\Omega$ with $|\w_F(\Omega)|=1$, we are unable to conclude that $T_F$ is hypercyclic.
However if for instance $\mathbb T\cap \Omega_1\neq\varnothing$ and $\mathbb T\cap\Omega_4\neq\varnothing$ as in \Cref{Fig33kfjdbk}, then $T_F$ is hypercyclic on $H^p$. 

Indeed, let $g$ vanish on $H_-(T_F)$. Then $g$ vanishes on $\Omega_1$ and $\Omega_4$. Let $Ug=(u_0,u_1)$ where $U$ is the operator defined in \Cref{Section:Yakubovich-court}. So $u_0=0$ on $\Omega_1\cup\Omega_4$. Thanks to the boundary conditions given by \Cref{eq:boundary-conditions-ErF}, we obtain that $u_0^{int}-\zeta u_1^{int}=0$ on $\partial \Omega_2\cap V$ (where $V$ is a small neighborhood of $\lambda_0$ where $\partial\Omega_1\cap\partial \Omega_4=\{\lambda_0\}$), and we can apply \Cref{Prop:Clef_Ext-Int} to conclude that $u_0=u_1=0$ on $\Omega_2$. So $T_F$ is hypercyclic on $H^p$.
    
\section{An extension to a more general setting}\label{section:JFA}

In this section, we will extend some results of the previous sections, replacing \ref{H2} by a more general assumption \ref{H2'} which allows the point $F(e^{i\theta})$ to travel several times along certain portions of the curve $F(\mathbb{T})$. This weaker assumption is introduced by Yakubovich in the paper \cite{Yakubovich1996}, and a model theory for Toeplitz operators on $H^2$ is developed here, under assumptions \ref{H1}, \ref{H2'} and \ref{H3}.  We will not provide the full details of this approach here, and we will simply state the results from \cite{Yakubovich1996} that will be needed in order to generalize some of the theorems from the previous sections to this more general setting.
In particular, we will consider only the case of Toeplitz operators acting on $H^2$ in this section (although it can reasonably be conjectured that the results can be extended to the $H^p$ setting). We refer to the paper \cite{Yakubovich1996} for all the details of the construction.

\subsection{Yakubovich's model theory for Toeplitz operators}\label{subsection:JFA-model}
Let $F:\mathbb{T}\longrightarrow\mathbb{C}$ be a symbol satisfying the assumptions \ref{H1}, \ref{H2'} and \ref{H3}, where 
the condition \ref{H2'} runs as follows:
\par\smallskip
\begin{enumerate}[(H2')]
    \item \label{H2'} Suppose that $\mathbb C\setminus F(\mathbb T)$ has a finite number of connected components, and there exist real numbers $\theta_0<\theta_1<\dots<\theta_m=\theta_0+2\pi$ such that, letting $\alpha_j$ be the open arc $\alpha_j=\{e^{i\theta};\,\theta_j<\theta<\theta_{j+1}\}$, we have\par\smallskip
     \begin{enumerate}[(a)]
       \item $F$ is injective on  each arc $\alpha_j,~0\le j \le m-1$;\par\smallskip
        \item for every $i\neq j,~0\le i,j\le m-1$, $F(\alpha_i)=F(\alpha_j)$ or  the sets $F(\alpha_i)$ and $F(\alpha_j)$ are disjoint.
    \end{enumerate}
\end{enumerate}\par
\begin{figure}[ht]
\includegraphics[page=35,scale=.8]{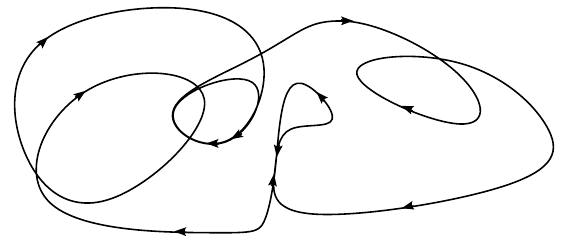}
\caption{}
\label{Fig34kfjn}
\end{figure}

\medskip
Under the assumption \ref{H2'}, let $\mathcal O$ be  the set the extremities of the subarcs $F(\alpha_j)$ where $\alpha_1,\dots,\alpha_k$ are given by  \ref{H2'}.
\Cref{Fig34kfjn} is a picture illustrating this assumption \ref{H2'}.
\par\smallskip
In order to develop a model theory for operators satisfying these weaker assumptions, it is required in \cite{Yakubovich1996} that ``the point $F(e^{i\theta})$ does not travel too many times in one way along a portion of the curve $F(\mathbb{T})$".
In order to explain precisely what is meant by this assumption, we recall some definitions and notations concerning the interior and exterior components with respect to a subarc of $F(\mathbb T)\setminus\mathcal O$.
\par\smallskip   
Let $\gamma$ be a non-oriented subarc of $F(\mathbb T)$ containing no point of $\cal O$. Suppose first that $\gamma$ is included in the boundary of exactly two connected components $\Omega$ and $\Omega'$ of $\mathbb C\setminus F(\mathbb T)$.
If $|\w_F(\Omega)|>|\w_F(\Omega')|$, we say that $\Omega$ is the \textit{interior component with respect to $\gamma$} and $\Omega'$ is the \textit{exterior component}. If $|\w_F(\Omega)|=|\w_F(\Omega')|$, interior and exterior connected components can be chosen interchangeably. Recall that the interior and exterior boundary values of a function $u$ at $\lambda_0\in\gamma$, denoted by $u^{int}(\lambda_0)$ and $u^{ext}(\lambda_0)$, are the non-tangential limits of $u$ at $\lambda_0$ through the interior and exterior component when these limits exist, i.e. 
\[
u^{int}(\lambda_0)~=~\lim_{\substack{\lambda\to\lambda_0\\
\lambda\in \Omega}}u(\lambda)\quad\text{and}\quad u^{ext}(\lambda_0)~=~\lim_{\substack{\lambda\to\lambda_0\\
\lambda\in \Omega'}}u(\lambda),
\]
where $\Omega$ and $\Omega'$ are respectively the interior and exterior connected components with respect to $\gamma$.
When $\gamma$ is included in the boundary of exactly one component $\Omega$, then we say that $\Omega$ is both the interior and the exterior component. Now we need to define the interior and exterior boundary values at a point $\lambda_0\in\gamma$ in this case: 
let $V$ be a neighborhood of $\lambda_0$ such that $V\setminus\gamma \subseteq\Omega$ has exactly two connected components $V_1$ and $V_2$. Then we note
\[
u^{int}(\lambda_0)~=~\lim_{\substack{\lambda\to\lambda_0\\
\lambda\in V_1}}u(\lambda)\quad\text{and}\quad u^{ext}(\lambda_0)~=~\lim_{\substack{\lambda\to\lambda_0\\
\lambda\in V_2}}u(\lambda)
\]
whenever these two non-tangential limits exist. The open sets $V_1$ and $V_2$ which are used to define $u^{int}(\lambda_0)$ and $u^{ext}(\lambda_0)$ can be chosen interchangeably.
\par\smallskip   
For every point $\lambda\in F(\mathbb T)\setminus\cal O$, let $\gamma  \subseteq F(\mathbb T)\setminus\cal O$ be a non-trivial curve such that $\lambda\in\gamma$. Let $\Omega$ and $\Omega'$ be respectively the interior and exterior components with respect to $\gamma$. Then set
\[w_i(\lambda)~:=~\w_F(\Omega)\quad \text{ and }\quad w_e(\lambda)~:=~\w_F(\Omega').\]
The functions $w_i$ and $w_e$ are in fact respectively the interior and exterior limits of the function $\w_F$ on the arc $\gamma$.
As it is now allowed to travel the arc $\gamma$ in both directions, 
let $n_i(\lambda)$ and $n_e(\lambda)$ denote the number of times that the point 
$F(e^{i\theta})$ travels along the arc $\gamma$ in each one of the two directions, in such a way that $n_i(\lambda)\ge n_e(\lambda)$. Recall now the following geometrical interpretation of the winding number of a curve at a point. Adding $1$ to the winding number can be interpreted as turning one more time around the point while keeping it on the left. In other words, saying that we travel $n_i(\lambda)$ times in one direction and $n_e(\lambda)$ times in the other direction on the curve $F(\T)$ in the vicinity of the point $\lambda$ is exactly equivalent to saying that we need to add $n_i(\lambda)$ and subtract $n_e(\lambda)$ to the winding number of one component (either interior of exterior) at the point $\lambda$ in order to obtain the winding number of the other one; in other words,
\begin{equation}\label{lien-ni-ne}
n_i(\lambda)-n_e(\lambda)~=~|w_i(\lambda)-w_e(\lambda)|~=~|w_i(\lambda)|-|w_e(\lambda)|
\end{equation}
or, equivalently, 
\begin{equation}\label{lien-ni-ne-bis}
|w_i(\lambda)|-n_i(\lambda)~=~|w_e(\lambda)|-n_e(\lambda).
\end{equation}
We refer for instance to \cite[p. 286]{Yger2001} for this geometrical interpretation of the winding number of a curve at a point, which will be also used in the proof of \Cref{lemme:relation-wind-int-ext}.
\par
\begin{figure}[ht]
\includegraphics[scale=.85,page=63]{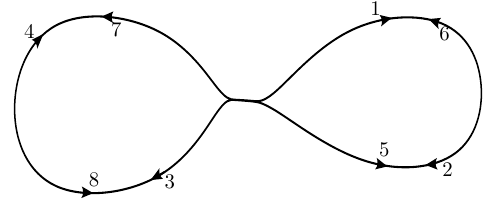}
    \caption{}
    \label{fig:spectrecourbe}
\end{figure}

\smallskip
In this new situation it is also possible to define an operator $U$ in much the same way as in \Cref{Section:Yakubovich-court}, but this time, the boundary condition induced by \Cref{remark-conditions-bords-vect-H2prime3434A4} depends of 
$n_i(\lambda)+1$ many functions $h_{\lambda,j}$ and  $n_e(\lambda)+1$ many functions  $h_{\lambda,j}$. 
It turns out that the following additional assumption is needed \cite[Lemma 4.1]{Yakubovich1996}:
\begin{equation}\label{Eq:9034ujoijw}
|w_i(\lambda)|-n_i(\lambda)~=~|w_e(\lambda)|-n_e(\lambda)~\ge~0\quad \text{ a.e. on}~F(\mathbb T),
\end{equation}
so that the $h_{\lambda,j}$'s involved in the boundary relation induced by \Cref{remark-conditions-bords-vect-H2prime3434A4} are eigenvectors of $T_F$. Note that this additional assumption asserts that $\sigma(T_F)\setminus F(\mathbb T)$ is non empty. Indeed, without this assumption $\sigma(T_F)\setminus F(\mathbb T)$ can be empty even if $F$ satisfies the assumptions \ref{H1}, \ref{H2'} and \ref{H3}. \Cref{fig:spectrecourbe} is a picture illustrating such an example.

Some authors have obtained models for certain classes of smooth symbols on $\mathbb T$ which do not satisfy \eqref{Eq:9034ujoijw}. In particular for unimodular symbols $F$, it was shown by Clark in \cite{Clark-1987} that whenever $F$ has the form $F= \omega \varphi/\bar\varphi$, where $\omega$ is an inner function and $\varphi$ and $1/\varphi$ belong to $H^\infty(\mathbb D)$ (in particular $\w_F(\lambda)\ge0$ for every $\lambda\notin F(\mathbb T)$), $T_F$ is similar to an isometry. Sufficiently smooth unimodular symbols with $\w_F(\lambda)\ge0$ for every $\lambda\notin F(\mathbb T)$ are of this form. For a certain subclass of rational functions, this isometry was described by Gamal' in \cite{Gamal2008}. For a more general class of continuous unimodular symbols with small sets of singularities which are smooth outside the set of singularities, she proved in \cite{Gamal2010} that $T_F$ is similar to a unilateral shift of finite multiplicity.
\par\smallskip
Under the additional assumption given by \eqref{Eq:9034ujoijw}, Yakubovich was able to define in \cite{Yakubovich1996} a model operator similar to the one obtained under the assumption \ref{H2} in \cite{Yakubovich1991} (see \Cref{T:Yakbovich_Hp}). We do not state this result in detail, but restrict ourselves to a rather informal statement of what is precisely needed for our purposes:

\begin{theorem}[Yakubovich \cite{Yakubovich1996}]\label{Th:YakubovichJFA}
Let $F$ satisfy the assumptions \ref{H1}, \ref{H2'} and \ref{H3}, and suppose moreover that 
\begin{equation}\label{JFA_cond_nb_passage}
|w_i(\lambda)|-n_i(\lambda)~=~|w_e(\lambda)|-n_e(\lambda)~\ge~0\quad\text{a.e. on}~F(\mathbb T).
\end{equation}
Then the operator 
$U$ of \Cref{T:Yakbovich_Hp} can be replaced by a bounded operator $ \widetilde U$, defined on $H^2$ and taking values in a vector-valued (Hilbertian) Smirnov space of functions on the connected components of $\sigma(T_F)\setminus F(\mathbb{T})$, satisfying similar properties to the ones of \Cref{T:Yakbovich_Hp}: 
\begin{itemize}
\item for every $\lambda\in \sigma(T_F)\setminus F(\mathbb{T})$, we have $ \widetilde UT_F^*g(\lambda)=\lambda  \widetilde Ug(\lambda)$. Also,  $ \widetilde{U}g(\lambda)=0$ if and only if $g$ vanishes on  $\ker(T_F-\lambda)$;
\item  the operator $ \widetilde U:H^2\to \Ran( \widetilde U)$ is invertible;
\item the image $\Ran( \widetilde U)$ of $ \widetilde{U}$ is the subspace of the vector-valued Smirnov space of functions on the connected components of $\sigma(T_F)\setminus F(\mathbb{T})$ satisfying the boundary conditions of \cite[Eq. (2.1)]{Yakubovich1996}.
\end{itemize}
\end{theorem}
\subsection{A necessary condition for hypercyclicity}
\Cref{Th:YakubovichJFA} implies in particular that, in the same way as what happened under the stronger assumption  \ref{H2}, the operator $T_F$ still has an $H^\infty$ functional calculus:

\begin{corollary}[Yakubovich \cite{Yakubovich1996}]\label{Cor:HinftyCalculusJFA}
Let $F$ satisfy the assumptions of \Cref{Th:YakubovichJFA}.
Then $T_F$ admits an $H^\infty$ functional calculus on the interior of ${\sigma(T_F)}$, and there exists a constant $C>0$ such that
\[\|u(T_F)\|~\le~ C\sup\left\{|u(\lambda)|\,;\, \lambda\in \overset{\circ}{\sigma(T_F)}\right\}\]
for every function $u\in H^\infty(\overset{\circ}{\sigma(T_F)})$.
\end{corollary}
Since \Cref{Cor:HinftyCalculusJFA} provides us with an $H^\infty$ functional calculus in this more general situation, \Cref{Th:CNforHC-Intersection} can be extended to this context. We have:

\begin{theorem}\label{Th:CNforHC-Intersection-JFA}
Let $F$ satisfy the assumptions of  \Cref{Th:YakubovichJFA}.
If $T_F$ is hypercyclic on $H^2$, then every connected component of the interior
of the spectrum of $T_F$ must intersect $\mathbb T$.
\end{theorem}

\subsection{Sufficient conditions for hypercyclicity} 
In this subsection, we will not really work with the operator $\widetilde{U}$ but rather with the operator $U$ constructed in \Cref{subsection-Yaku-Operator-U34ZZE111}, which satisfies similar properties to those proved when \ref{H2} holds under the more general hypothesis \ref{H2'}. The main difference between the two settings concerns the range of this operator, which is more difficult to describe and to use in the general case where only \ref{H2'} is satisfied. Let us briefly recall the construction. In \Cref{subsection-Yaku-Operator-U34ZZE111}, we considered the operator $U$ defined as 
\[
(Ug)_j(\lambda)~=~\dual{h_{\lambda,j}}{g}\quad \text{ for every } g\in H^q,\, \lambda\in\Omega_j^+,~\text{and}~ 0\le j<N,
\]
where the functions $h_{\lambda,j}$ defined by  $h_{\lambda,j}(z)=z^j F_{\lambda}^+(0)/F_\lambda^+(z)$, for $z\in\mathbb D$ and $0\le j<N$, form a basis of the eigenspace $\ker(T_F-\lambda)$. This construction remains valid when we replace \ref{H2} by \ref{H2'}. In this case, the operator $U$ satisfies similar properties to the ones that it enjoys under assumption \ref{H2}. More precisely, we have the following:
\begin{theorem}
Let $F$ satisfy the assumptions of \Cref{Th:YakubovichJFA}. Then $U$ is a bounded and injective operator from $H^2$ into $\bigoplus_{0\le j<N} E^2(\Omega_j^+)$.
\end{theorem}

\begin{proof}
The proof of the boundedness of $U$ under the hypothesis \ref{H2'} is carried out in \Cref{lemme-Uborne}. The injectivity of $U$ follows from the injectivity of $ \widetilde U$. Indeed, since the functions $h_{\lambda,j}$, $0\leq j<|\w(\lambda)|$, form a basis of $\ker(T_F-\lambda)$ for every $\lambda\in\sigma(T_F)\setminus F(\T)$, we know that $g\perp\ker(T_F-\lambda)$ if and only if $(Ug)_j(\lambda)=0$ for every $0\leq j<|\w_F(\lambda)|$ (see \Cref{basic-fact}). According to \Cref{Th:YakubovichJFA}, $g\perp\ker(T_F-\lambda)$ if and only if $ \widetilde Ug(\lambda)=0$. Hence if $Ug=0$, we also have $ \widetilde Ug=0$, and 
then $g=0$ by \Cref{Th:YakubovichJFA}. 
\end{proof}

In particular, we have the following analogue of \Cref{Prop:EigenvectorDenseYaku}:

\begin{proposition}\label{Prop:EigenvectorDenseYakuJFA}
Let $F$ satisfy the assumptions of \Cref{Th:YakubovichJFA}.
\begin{enumerate}
    \item Let $\Omega\in\mathcal C$, and let $A$ be a subset of $\mathbb C$. If $A\cap \Omega$ has an accumulation point in $\Omega$, then
    \[\overline{\spa}\,[\ker(T_F-\lambda)\,;\,\lambda\in A\cap\Omega]~=~\overline{\spa}\,[\ker(T_F-\lambda)\,;\,\lambda\in \Omega].\]
    \item  We have
\begin{equation*}
      \overline{\spa}\,[H_\Omega(T_F)\,;\,\Omega\in\cal C]~=~H^2.
  \end{equation*}
\end{enumerate}
\end{proposition}
\par\smallskip
In the rest of this section, we give a version, under the assumption \ref{H2'}, of some of the results of  \Cref{Section:CNforHc,Section:ApplicationCShc}
which are also true in this context. 

\begin{figure}[ht]
\includegraphics[page=36,scale=1]{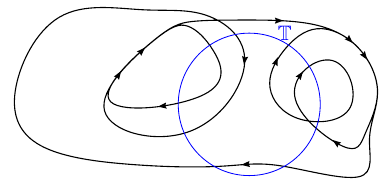}
\caption{}
\label{Fig35lkm}
\end{figure}

First, since \Cref{Prop:EigenvectorDenseYaku} generalizes here into \Cref{Prop:EigenvectorDenseYakuJFA},  we can extend \Cref{Th:SimpleConseqHc} (which was a direct consequence of \Cref{Prop:EigenvectorDenseYaku}) to the case where only \ref{H2'} is satisfied. This is the content of the next theorem:

\begin{theorem}\label{Th:SimpleConseqHc-JFA}
Let $F$ satisfy the assumptions of  \Cref{Th:YakubovichJFA}.
Suppose that
\begin{equation}
\Omega\cap \mathbb T~\neq~\varnothing \quad\text{for every}~ \Omega\in\cal C.
\end{equation}
Then $T_F$ is hypercyclic on $H^2$.
\end{theorem}
\Cref{Fig35lkm} illustrates a situation where \Cref{Th:SimpleConseqHc-JFA} can be applied.

  As already said, the difficulty when working  with assumption \ref{H2'} lies with the description of the range of $U$. In \Cref{remark-conditions-bords-vect-H2prime3434A4}, we show that when \ref{H2'} holds, the eigenvectors $h_{\lambda,j}$ satisfy a boundary equation of the following form: for every $z\in\mathbb D$,
\[
\sum_{p=0}^{n_{i}(\lambda)}a_p(\lambda)h_{\lambda,k+p}^{int}(z)~=~\sum_{\ell=0}^{n_{e}(\lambda)}b_\ell(\lambda)h_{\lambda,k+\ell}^{ext}(z)\quad \text{for a.e. } \lambda\in\partial\Omega_n^+
\]
for every $n\ge k+\max(n_i(\lambda),n_e(\lambda))$, where the quantities $a_p(\lambda)$ and $b_l(\lambda)$ are defined by \Cref{eq:polynome1-condition-aux-bords,eq2:polynome1-condition-aux-bords} respectively. Now, using the boundedness of the operator $U$ (which is still true under hypothesis \ref{H2'}), we see that if, given $g\in H^q$, we write $Ug$ as $Ug=(u_{j})_{0\leq j\leq N-1}$, then the functions $u_j$ satisfy the boundary relations  
\begin{equation}\label{eq:condition-bords-H2prime-73343TT}
\sum_{p=0}^{n_{i}(\lambda)}a_p(\lambda) u_{k+p}^{int}(\lambda)~=~\sum_{\ell=0}^{n_{e}(\lambda)}b_\ell(\lambda) u_{k+\ell}^{ext}(\lambda)\quad \text{for a.e. } \lambda\in\partial\Omega_n^+
\end{equation}
for every $n\ge k+\max(n_i(\lambda),n_e(\lambda))$.

A particular case of \ref{H2'} where the assumption \eqref{JFA_cond_nb_passage} of \Cref{Th:YakubovichJFA} is satisfied is when $n_e=0$ almost everywhere on $F(\mathbb{T})$, which means that the point $F(e^{i\theta})$ travels over every subarc of $F(\mathbb T)\setminus\cal O$ in one direction only. Under this condition, 
\Cref{eq:condition-bords-H2prime-73343TT} implies that the range of $U$ is contained in a subspace of $\bigoplus E^2(\Omega_j^+)$ consisting of $N$-tuples $(u_j)_{0\le j\le N-1}\in \bigoplus_{0\le j\le N-1} E^2(\Omega_j^+)$ satisfying a boundary relation of the type 
\[u_{j}^{ext}(\lambda)~=~a_0(\lambda)\,u_{j}^{int}(\lambda)+...+a_{n_i(\lambda)}(\lambda)\,u^{int}_{j+n_i(\lambda)}(\lambda)\quad\text{for a.e. }\lambda\in \partial\Omega_{j+n_i(\lambda)}^+.\]
In particular, if $u=Ug$ vanishes on the interior component of an arc $\gamma$ included in $F(\mathbb{T})\setminus\mathcal{O}$, then $u$ vanishes also on the exterior component of $\gamma$, and so \Cref{Th:Simple_Int-Ext} can be extended when $F$ satisfies \ref{H1}, \ref{H2'} and \ref{H3} and $n_e=0$ almost everywhere on $F(\mathbb T)$. Thus \Cref{Th:HcInstersCompMax} can be extended as follows:

\begin{theorem}\label{Th:HcInstersCompMax_JFA}
Let $F$ satisfy \ref{H1}, \ref{H2'} and \ref{H3}. Suppose moreover that $n_e=0~a.e.$ on $F(\mathbb{T})$. If $\Omega\cap\mathbb T\neq\varnothing$ for every maximal component $\Omega$ of $\sigma(T_F)\setminus F(\mathbb T)$, then $T_F$ is hypercyclic on $H^2$.
\end{theorem}

\begin{figure}[ht]
\includegraphics[page=37,scale=.95]{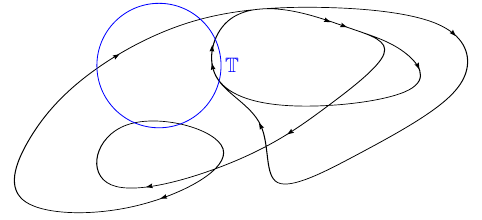}
\caption{}\label{FIG38}\end{figure}

\Cref{FIG38} gives an example of a function $F$ to which \Cref{Th:HcInstersCompMax_JFA} applies but not \Cref{Th:SimpleConseqHc-JFA}.

 Although going ``from the interior component to the exterior component of a boundary arc" is still possible under assumption \ref{H2'} combined with the hypothesis that $n_e=0$ almost everywhere on $F(\mathbb{T})$, we do not know of general conditions, possibly similar to those of Theorems \ref{Th:Complique_Ext-Int} or \ref{Th:Complique-Jordan_version},
implying that one can go ``from the exterior component to the interior component of a boundary arc". Still, we are able to handle the following interesting example.
\begin{example}\label{Exemple:2cercle-JFA}
One of the simplest situations where we would need to go from the exterior component of a boundary arc to the interior, in order to conclude that \(T_{F}\) is hypercyclic, is the one given by a curve as in \Cref{Fig37kll}.
\begin{figure}[ht]
\includegraphics[page=38,scale=.9]{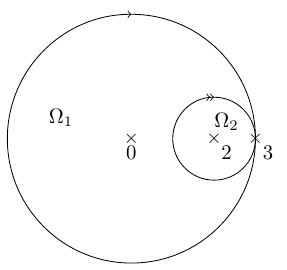}
            \caption{}
            \label{Fig37kll}
            \end{figure}
            
\noindent It is similar to the curve considered in \Cref{SubSection:Example2Circles}, but the inner circle is traveled twice. Set
\[F(e^{i\theta })~=~
\begin{cases}
3\,e^{-5i\theta /3}\quad &\text{if}\quad  0\le\theta <6\pi /5\\
2+e^{-5i\theta }\quad &\text{if}\quad  6\pi /5\le \theta <2\pi.
\end{cases}
\]
In this case, the boundary relations allow us to prove directly that \(T_{F}\) is hypercyclic on \(H^{2}\). To this aim, we proceed in a similar fashion to what we did in \Cref{SubSection:Example2Circles}: for every $\lambda\in\partial\Omega_2\setminus\{3\}$, there exist $\theta_1(\lambda)\in(0,2\pi/5)$ and $\theta_2(\lambda)\in(2\pi/5,4\pi/5)$ such that
\[\lambda~=~F(e^{-i\theta_1(\lambda)})~=~2+e^{5i\theta_1(\lambda)}\quad\text{and}\quad\lambda~=~F(e^{-i\theta_2(\lambda)})~=~2+e^{5i\theta_2(\lambda)}.\]
So let $\zeta_1(\lambda)=e^{i\theta_1(\lambda)}$ and $\zeta_2(\lambda)=e^{i\theta_2(\lambda)}$. We have
\begin{align*}
 \zeta _{1}(\lambda )&~=~\exp\left(\dfrac{i}{5}\arg_{(0 ,2\pi )}\,(\lambda -2) \right)\quad\text{for every }\lambda \in\partial \Omega _{2}\setminus\{3\};\\
 \zeta _{2}(\lambda )&~=~\exp\left(\dfrac{1}{5}\arg_{(2\pi ,4\pi )}\,(\lambda -2) \right)\quad\text{for every }\lambda \in\partial \Omega _{2}\setminus\{3\}.
\end{align*}
Let \(g\in H^{2}\). Then write \(U{g}=(u_{0},u_{1},u_{2})\), where \(u_{0}\in E^{2}(\Omega _{1}\cup\Omega _{2})\) and \(u_{1},\,u_{2}\in E^{2}(\Omega _{2})\) satisfy the following boundary relation:
\begin{equation}\label{Eq etoile}
 u_{0}^{int}-(\zeta _{1}+\zeta _{2})u_{1}^{int}+\zeta _{1}\zeta _{2}u_{2}^{int}~=~u_{0}^{ext}\quad \text{a.e. on}\ \partial\Omega _{2}.
\end{equation}
Indeed, in this example, the polynomials in \Cref{eq:polynome1-condition-aux-bords,eq2:polynome1-condition-aux-bords} become 
\[
\prod_{j\in N_{int}(\lambda)}(1-d_j(\lambda)z)~=~(1-\zeta_1(\lambda)z)(1-\zeta_2(\lambda)z)
~=~1-(\zeta_1(\lambda)+\zeta_2(\lambda))z+\zeta_1(\lambda)\zeta_2(\lambda) z^2;
\]\[
\prod_{j\in N_{ext}(\lambda)}(1-d_j(\lambda)z)~=~1,
\]
and then \Cref{eq:condition-bords-H2prime-73343TT} gives exactly \Cref{Eq etoile}.
\par\smallskip
\(\bullet\) Suppose that \(g\) is orthogonal to \(H_{-}(T_{F})\). Then, by \Cref{Prop:EigenvectorDenseYakuJFA},  \(g\perp H_{\Omega _{1}}(T_{F})\), so that \(u_{0}=0\) on \(\Omega _{1}\). Hence \(u_{0}^{e}=0\) almost everywhere on \(\partial\Omega _{2}\) and  \Cref{Eq etoile} becomes
\[
u_{0}^{int}-(\zeta _{1}+\zeta _{2})u_{1}^{int}+\zeta _{1}\zeta _{2}u_{2}^{int}~=~0\quad \text{a.e. on}\ \partial\Omega _{2}.
\]
The functions \(\zeta _{1}\) and \(\zeta _{2}\) admit bounded analytic extensions to 
\(\Omega _{2}\setminus[2,3]\) given by
\begin{align*}
 \zeta _{1}(\lambda )&~=~\exp\left(\dfrac{1}{5}\Bigl [ \log|\lambda -2|+i\,\arg_{(0,2\pi )}\,(\lambda -2)\Bigr]  \right)\quad\text{for every }\lambda \in \Omega _{2}\setminus[2,3];\\
 \zeta _{2}(\lambda )&~=~\exp\left(\dfrac{1}{5}\Bigl [ \log|\lambda -2|+i\,\arg_{(2\pi,4\pi )}\,(\lambda -2)\Bigr]  \right)\quad\text{for every }\lambda \in \Omega _{2}\setminus[2,3].
\end{align*}
which are such that for every \(x\in[2,3]\),
\begin{align*}
 \lim_{\stackrel{y\to 0}{y>0}}\zeta _{1}(x+i\,y)&~=~(x-2)^{1/5}~\neq~(x-2)^{1/5}e^{i\frac{2\pi }{5}}~=~\lim_{\stackrel{y\to 0}{y<0}}\zeta _{1}(x+i\,y);\\
 \lim_{\stackrel{y\to 0}{y>0}}\zeta _{2}(x+i\,y)&~=~(x-2)^{1/5}e^{i\frac{2\pi }{5}}~\neq~(x-2)^{1/5}e^{i\frac{4\pi }{5}}~=~\lim_{\stackrel{y\to 0}{y<0}}\zeta _{2}(x+i\,y).
\end{align*}
The functions \(u_{0}\), \((\zeta _{1}+\zeta _{2})u_{1}\), and \(\zeta _{1}\zeta _{2}u_{2}\) belong to \(E^{2}(\Omega _{2}\setminus[2,3])\) and since 
\(u_{0}-(\zeta _{1}+\zeta _{2})u_{1}+\zeta _{1}\zeta _{2}u_{2}\) vanishes on a subset of positive measure of the boundary of \(\Omega _{2}\setminus [2,3]\), we have 
\(u_{0}-(\zeta _{1}+\zeta _{2})u_{1}+\zeta _{1}\zeta _{2}u_{2}=0\) on \(\Omega _{2}\setminus[2,3]\).
\par\smallskip
If either \(u_{1}\) or \(u_{2}\) is identically zero on \(\Omega _{2}\), the same argument as in \Cref{SubSection:Example2Circles} shows that \(u_{0}=0\) on \(\Omega _{2}\) as well. So henceforward, we suppose that \(u_{1}\) and \(u_{2}\) are not identically zero on \(\Omega _{2}\). For all \(x\in [2,3)\) except possibly countably many, we have \(u_{1}(x)\neq 0\) and \(u_{2}(x)\neq 0\). The function \(u_{0}=(\zeta _{1}+\zeta _{2})u_{1}-\zeta _{1}\zeta _{2}u_{2}\) is continuous at the point \(x\) (because $u_0\in E^2(\Omega_2)$), and taking limits of \(u_{0}(x+iy)\) and \(u_0(x-iy)\) as \(y\to 0\), \(y>0\), we have
\begin{multline*}
 (x-2)^{1/5}\Bigl (1+e^{i\frac{2\pi }{5}} \Bigr)u_{1}(x)-(x-2)^{2/5}e^{i\frac{2\pi }{5}}u_{2}(x)\\
 = ~(x-2)^{1/5}\Bigl (e^{i\frac{2\pi }{5}}+e^{i\frac{4\pi }{5}} \Bigr)u_{1}(x)-(x-2)^{2/5}e^{i\frac{6\pi }{5}}u_{2}(x).
\end{multline*}
Hence
\(
u_{1}(x)-(x-2)^{1/5}e^{i\frac{2\pi }{5}}u_{2}(x)=e^{i\frac{4\pi }{5}}u_{1}(x)-(x-2)^{1/5}e^{i\frac{6\pi }{5}}u_{2}(x)
\)
so that 
\begin{equation}\label{Eq double etoile}
 \Bigl (1-e^{i\frac{4\pi }{5}} \Bigr)u_{1}(x)~=~(x-2)^{1/5}\Bigl (e^{i\frac{2\pi }{5}}-e^{i\frac{6\pi }{5}} \Bigr)u_{2}(x)\quad \textrm{ for every } x\in [2,3].
\end{equation}
Set
\[
w(\lambda)~=~\exp\left( \frac{1}{5}\Bigl [\log|\lambda-2|+i\,\arg_{\,(-\pi ,\pi )}(\lambda-2) \Bigr] \right)\quad \text{ for every } \lambda\in\Omega _{2}\setminus[1,2]\]and
\( v_{1}(\lambda)~=~\Bigl (1-e^{i\frac{4\pi }{5}} \Bigr)u_{1}(\lambda),\quad  v_{2}(\lambda)~=~\Bigl (e^{i\frac{2\pi }{5}}-e^{i\frac{6\pi }{5}} \Bigr)u_{2}(\lambda)\) for every $ \lambda\in\Omega_2$.
Then \(v_{1}\), \(v_{2}\), and \(w\) are analytic functions on \(\Omega _{2}\setminus[1,2]\). Using the uniqueness principle and \Cref{Eq double etoile}, we deduce that \(v_{1}=w\,.\,v_{2}\) on \(\Omega _{2}\setminus [1,2]\). Now, since \(v_{2}\) is not identically zero on \(\Omega _{2}\), there exists 
\(x\in(1,2)\) with \(v_{2}(x)\neq 0\), and so \(w\) admits an analytic extension to a neighborhood of \(x\). But this contradicts the fact that 
\[
\lim_{\stackrel{y\to 0}{y>0}}w(x+iy)~=~(x-2)^{1/5}e^{i\frac{\pi }{5}}~\neq~ (x-2)^{1/5}e^{-i\frac{\pi }{5}}~=~\lim_{\stackrel{y\to 0}{y<0}}w(x+iy).
\]
So finally, we obtain that \(u_{1}=u_{2}=0\) on \(\Omega _{2}\), so that \(u_{0}=0\) on \(\Omega _{2}\). This means that \(g\perp H_{\Omega _{2}}(T_{F})\), and so, by \Cref{Prop:EigenvectorDenseYakuJFA}, \(g=0\). We finally deduce that \(H_{-}(T_{F})=H^{2}\).
\par\smallskip
\(\bullet\) Suppose, lastly, that \(g\) is orthogonal to \(H_{+}(T_{F})\). Since 
\(\bigl (\mathbb{C}\setminus\overline{\mathbb{D}} \bigr)\cap\Omega _{1}\neq\varnothing \), and \(\bigl (\mathbb{C}\setminus\overline{\mathbb{D}} \bigr)\cap\Omega _{2}\neq\varnothing \), \(g\) is orthogonal to \(H_{\Omega _{1}}(T_{F})\) and \(H_{\Omega _{2}}(T_{F})\). So \(g=0\) and \(H_{+}(T_{F})=H^{2}\). So \(T_{F}\) is hypercyclic on \(H^{2}\).
\end{example}   
\subsection{An equivalence}
\begin{figure}[ht]
\includegraphics[page=39,scale=.76]{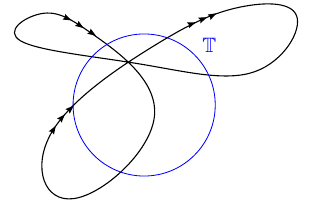}
\caption{}
\label{Fig879ljkkjh}
\end{figure}

When
$\sigma(T_F)\setminus F(T_F)=\overset{\circ}{\sigma(T_F)}$, 
 \Cref{Th:CNforHC-Intersection-JFA,Th:SimpleConseqHc-JFA} can be combined to provide a characterization of hypercyclicity.
 Observe that $\sigma(T_F)\setminus F(\mathbb T)$ coincides with the interior of $\sigma(T_F)$ if and only if the condition $F(\mathbb T)=\partial\sigma(T_F)$ is satisfied. 
 We have the following theorem:
 
\begin{theorem}\label{th-quelquechose}
Let $F$ satisfy the assumptions of \Cref{Th:YakubovichJFA}. Suppose that $F(\mathbb T)=\partial\sigma(T_F)$. Then the following assertions are equivalent:
\begin{enumerate}
    \item $T_F$ is hypercyclic on $H^2$;
    \item $O\cap\mathbb T\neq\varnothing$ for every connected component $O$ of $\overset{\circ}{\sigma(T_F)}$.
\end{enumerate} 
\end{theorem}
\Cref{Fig879ljkkjh} illustrates a situation where \Cref{th-quelquechose} can be applied.

\begin{figure}[ht]
\includegraphics[page=40,scale=1]{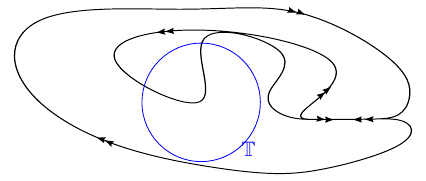}
                \caption{}
                \label{Fig39gkhjm}
                \end{figure}

\Cref{th-quelquechose}
 can be applied in particular when $F(\mathbb T)$ is a Jordan curve on which the point $F(e^{i\theta})$ travels several times.  
We can also deal with the case where $\sigma(T_F)\setminus F(\mathbb T)$ is connected:
 
\begin{theorem}\label{Th:1compCNS-JFA}
Let $F$ satisfy assumptions \ref{H1}, \ref{H2'} and \ref{H3}. If $\sigma(T_F)\setminus F(\mathbb T)$ is connected, then the following assertions are equivalent:
\begin{enumerate}
    \item $T_F$ is hypercyclic on $H^2$;
    \item $\overset{\circ}{\sigma(T_F)}\cap\mathbb T\neq\varnothing$.
\end{enumerate}    
\end{theorem}
\Cref{Fig39gkhjm} illustrates a situation where \Cref{Th:1compCNS-JFA} can be applied.

\begin{proof}
Let us begin by observing the following fact: set $\Omega=\sigma(T_F)\setminus F(\mathbb T)$. If  $\overset{\circ}{\sigma(T_F)}\cap\mathbb T\neq\varnothing$, then $\overset{\circ}{\sigma(T_F)}\cap\mathbb D\neq\varnothing$. Since, in this case, $\Omega=\overset{\circ}{\sigma(T_F)}\setminus F(\mathbb T)$, it follows that $\mathbb D$ has to intersect also $\Omega$. The same argument yields that $\mathbb C\setminus{\overline{\mathbb D}}$ intersects $\Omega$. Since $\Omega$ is connected, this implies that $\mathbb T\cap\Omega\neq\varnothing$.
We deduce that, under the hypothesis that $\Omega$ is connected, the necessary condition for hypercyclicity in \Cref{Th:CNforHC-Intersection-JFA} and the sufficient condition in \Cref{Th:SimpleConseqHc-JFA} are equivalent. We obtain thus the equivalence stated in \Cref{Th:1compCNS-JFA}.
\end{proof}

\section{Comparison with results of Abakumov, Baranov, Charpentier and Lishanskii}\label{Section:ABCL}
In \cite{AbakumovBaranovCharpentierLishanskii2021}, the authors study the hypercyclicity on $H^2$ of Toeplitz operators  with symbols which are meromorphic on $\mathbb D$ and continuous up to the boundary $\mathbb T$. Our aim in this section  will be to apply our previous results to this class of symbols, and thus to recover and extend some of the results of \cite{AbakumovBaranovCharpentierLishanskii2021}. More precisely, in this whole section, we consider symbols of the following form:
\begin{equation}\label{Form}
F(z)~=~R(1/z)+\phi(z),
\end{equation}
where $\phi\in A(\mathbb D)$ and $R$ is a rational function without poles in $\overline{\mathbb D}$. In other words,  the function $R$ can be written as
\[R(z)~=~P(z)+\sum_{l=1}^r\sum_{j=1}^{k_l}\frac{\alpha_{l,j}}{(z-\eta_l)^j},\]
where $P$ is a polynomial of degree $N_1$ and the poles $\eta_l\in\mathbb C\setminus\overline{\mathbb D}$ of the rational function $R$ are distinct with respective multiplicities $k_l$, $1\leq l\leq r$. Then setting $N_2=\sum_{l=1}^r k_l$ and $D=N_1+N_2$, we have $\deg(R)=D$. In other words, $F$ has exactly $D$ poles in $\mathbb D$, counted with multiplicity. We denote by $\mathcal{P}$ the set of these poles, which consists of the points $\eta_1^{-1},\dots,\eta_r^{-1}$, plus the point $0$ if the polynomial $P$ is not constant.

\subsection{A consequence of assumption \ref{H1}}
In this short subsection, we present a consequence of the smoothness property \ref{H1} for symbols of the form (\ref{Form}). It will be crucial to our study of exactly $D$-valent mappings later on.

\begin{lemma}\label{Fait:LocalInjectCercle}
Let $F$ be a symbol of the form given by \Cref{Form}, and suppose that $F$ is of class $C^{1+\varepsilon}$ on $\T$ with $F'\neq 0$ on $\T$. 
Then, for every $\zeta\in\mathbb T$, there exists $\delta>0$ such that $F$ is injective on $\overline{\mathbb D}\cap D(\zeta,\delta)$.
\end{lemma}

\begin{proof}
Since $F$ is of class $C^{1+\varepsilon}$ on $\mathbb T$, the function $\phi$ is also of class $\mathcal{C}^{1+\varepsilon}$ on $\mathbb T$. By the remark after \Cref{thm:privalov-Zygmund}, the functions $\phi$ and $\phi'$ belong to $ A(\mathbb D)$, and so, for every $\zeta\in\mathbb T$ we have
\[\lim_{\underset{|z|\le1}{z\to\zeta}}F'(z)=\lim_{\underset{|z|\le1}{z\to\zeta}}\Big(R(1/z)\Big)'+\lim_{\underset{|z|\le1}{z\to\zeta}}\phi'(z)=-\frac1{\zeta^2}R'(1/\zeta)+\phi'(\zeta)=F'(\zeta).\]
Given $\zeta\in\mathbb T$, we know that $F'(\zeta)\neq0$. By multiplying eventually by a unimodular constant, we can assume without loss of generality that $\Re   F'(\zeta)>0$. Then, by continuity, there exists $a>0$ and $\delta>0$ such that, for every $z\in\overline{\mathbb D}\cap D(\zeta,\delta)$, $\Re  F'(z)\ge a>0$. A well-known result (see for instance 
\cite[Lemma 1]{MR1501813}) implies now that $F$ is injective. 
\end{proof}

\subsection{Link with the setting of \cite{AbakumovBaranovCharpentierLishanskii2021}}
Since, for every $\lambda\notin F(\mathbb T)$, the function $F-\lambda$ has no zero nor pole on $\mathbb T$,  the argument principle implies that
\begin{equation}\label{Eq:WindVal}
\w_F(\lambda)~=~n_F(\lambda)-D \quad\text{ for all }\lambda\notin F(\mathbb T),
\end{equation}
where $n_F(\lambda)$ is defined for every $\lambda\in\C$ as the number of solutions of  the equation $F(z)=\lambda$ in $\mathbb D$, counted with multiplicity. 
This equality provides a relation between the orientation of the curve $F(\mathbb T)$ and the valence of $F$.  More precisely, $F(\mathbb T)$ is negatively wound whenever $F$ is $D$-valent on $\mathbb D$, i.e. $F$ satisfies $n_F\le D$. Note also that if $F$ is negatively wound, then $D$ is related to the $N=\max|\w_F|$ introduced in \Cref{Section:Yakubovich-court}  by the following relation: $N= D-\min_{\lambda\notin F(\mathbb T)}{n_F(\lambda)}$. In particular, if $F$ is not surjective, i.e. if $F(\overline{\mathbb D})\neq\widehat{\mathbb C}$, then $D=N$. 
\par\smallskip
In the rest of this section, we will implicitly view $F$ as a function from $\overline{\mathbb D}$ into $\widehat{\mathbb{C}}$. The hypothesis that $F$ is not surjective will always be implied by the hypothesis that we will consider, so that we will always have $D=N$.
\par\smallskip
If $\Omega$ is a connected component of $\C\setminus F(\T)$, we denote by $n_F(\Omega)$ the common value of the integers $n_F(\lambda)$, $\lambda\in\Omega$. We have $n_F(\Omega)=\w_F(\Omega)+D$. If $\Omega_{\infty}$ denotes the unbounded component of $\C\setminus F(\T)$, then $n_F(\Omega_{\infty})=D$.
\par\medskip
Let now $A$ be a subset of 
$\mathbb{D}$, and let $K$ be a positive integer.
We say that $F$ is $K$-valent on $A$ (resp. exactly $K$-valent on $A$)  if for any $\lambda\in F(A)$, the equation $F(z)=\lambda$ has at most (resp. exactly) $K$ solutions in $A$, counted with multiplicity. By \Cref{Fait:LocalInjectCercle}, $F$ is injective on a neighborhood of every point of $\mathbb T$. So, by counting once the unimodular solutions of the equation $F=\lambda$ (i.e. by considering that each such solution has ``multiplicity" one), we can extend the definition of (exact) $K$-valence to subsets $A$ of $\overline{\mathbb D}$.
Note that, taking 
$A=\mathbb D\setminus\mathcal{P}$ (or $A=\D$),
the K-valence of $F$ on $A$ (which means by definition that $n_F(\lambda)\le K$ for every $ \lambda\in \C$)  is equivalent to the property $n_F(\lambda)\le K$ for every $ \lambda\notin F(\mathbb T)$, i.e. 
$\w_F(\lambda)\le K-D$ for every $ \lambda\notin F(\mathbb T)$.
Indeed, if $n_F(\lambda)\le K$ for every $ \lambda\notin F(\mathbb T)$, then there is no point $\lambda_0\in F(\T)$ such that $n_F(\lambda_0)> K$: if it were the case, then by Rouch\'e's Theorem there would exist $\lambda\notin F(\T)$ such that $n_F(\lambda)> K$, a contradiction. So $n_F(\lambda)\le K$ for every $\lambda\in\C$.

\par\medskip
Recall that \Cref{Prop:OrientationForHC} asserts that if a Toeplitz operator with a continuous symbol $F$ is hypercyclic on $H^p$ for some $p>1$, then $F(\mathbb T)$ is negatively wound, i.e. $\w_F\le 0$ on $\mathbb C\setminus F(\mathbb T)$. When $F$ is defined as in \Cref{Form}, the relation between the winding number and the valence in \Cref{Eq:WindVal} together with \Cref{Prop:OrientationForHC} imply the following result:
\par\smallskip
\begin{proposition}\label{prop:HC->Nval}
Let $F$ be given by \Cref{Form}, and let $p>1$. If $T_F$ is hypercyclic on $H^p$, then $F$ is D-valent on 
$\mathbb D$.
\end{proposition}\par\smallskip
We observe that if $F$ is D-valent on $\mathbb D$, then \Cref{Lemma:spectral-properties-toeplitz-operators} and \Cref{Eq:WindVal} imply that 
\begin{equation}\label{spectre-Toeplitz-valence}
\sigma(T_F)~=~\left\{\lambda\in\mathbb{C}\setminus F(\mathbb{T})\,;\, n_F(\lambda)<D\right\}\cup F(\mathbb{T}).
\end{equation}
We finish this subsection by spelling out the link between the eigenvectors  $h_{\lambda, j}$ of $T_F$, for $|\w_F(\lambda)|=D$ and $0\le j\le D-1$,  and the eigenvectors used in  \cite{AbakumovBaranovCharpentierLishanskii2021}.

\par\smallskip
Consider a symbol \(F\) of the form (\ref{Form}), satisfying assumptions \ref{H1}, \ref{H2}, and \ref{H3}, and suppose that $F$ is not surjective. For \(\lambda \in\mathbb{C}\setminus F(\overline{\mathbb{D}})=\mathbb{C}\setminus\overline{F(\mathbb{D})}\) (i.e \(\lambda \) belongs to a connected component \(\Omega \in\mathcal{C}\) such that \(|\w_{F}(\Omega )|=D\)), the eigenvectors \(h_{\lambda ,j}\), \(0\le j\le D-1\), can be written as
\begin{equation}\label{eq:vecteur-propre-ABCL}
h_{\lambda ,j}(z)~=~z^{j}\cdot\dfrac{c_\lambda}{z^{N_{1}}Q(z)(F(z)-\lambda )}\quad\text{ for every } z\in\mathbb{D}
\end{equation}
where \[Q(z)~=~\prod_{l=1}^{r}(1-\eta_l z)^{k_{l}}\] and \(c_\lambda\) is a non-zero constant given by the value at \(0\) of the function \(z\longmapsto z^{N_{1}}Q(z)(F(z)-\lambda)\), which is analytic on \(\mathbb{D}\). 
Indeed, write 
\[
\phi_\lambda(z)~=~z^N(F(z)-\lambda)=\frac{z^{N_2}}{Q(z)}[z^{N_1}Q(z)(F(z)-\lambda)]\quad\text{ for every } z\in\mathbb D,
\]
and observe that the function \(z\longmapsto z^{N_{1}}Q(z)(F(z)-\lambda )\) belongs to \(A(\mathbb{D})\) and does not vanish on \(\overline{\mathbb{D}}\). In particular, one can define an analytic branch of its logarithm on $\mathbb{D}$, denoted by $\log(z^{N_{1}}Q(z)(F(z)-\lambda ))$, which also belongs to $A(\mathbb D)$. Moreover, since $N_2=\sum_{l=1}^r k_l$, we have 
\[
\frac{z^{N_2}}{Q(z)}~=~\frac{1}{\prod_{l=1}^r \left(\frac{1}{z}-\eta_l\right)^{k_l}}\cdot
\]
Observe now that for every $1\le l\le r$, the function $z\longmapsto\log(z-\eta_l)$ (where we choose here a suitable determination of the logarithm) belongs to $A(\mathbb D)$, and thus the function $z\longmapsto \log\left(\frac{1}{z}-\eta_l\right)$  belongs to $H^p_-$, where $H^p_-=\{f\in H^p\,;\, \hat{f}(n)=0\text{ for every }n\ge 0\}$.
Define now 
\[
U_\lambda(z)~=~-\sum_{l=1}^r k_l \log\left(\frac{1}{z}-\eta_l\right)+\log\left(z^{N_{1}}Q(z)(F(z)-\lambda )\right)\quad\text{ for every } z\in\mathbb T.
\]
We have
\[
e^{U_\lambda(z)}~=~z^{N_1}Q(z)(F(z)-\lambda)\frac{z^{N_2}}{Q(z)}~=~\phi_\lambda(z)\quad\text{ for every } z\in\T.
\]
Hence, according to \Cref{eq:vecteur-propre-Yaku33434ZD}, we see that 
\(
F_\lambda^+=e^{P_+U_\lambda}.
\)
As the function $\log(z^{-1}-\eta_l)$ belongs to  $H^p_-$ for every $1\leq l\leq r$, we have 
$$P_+\left(\log(z^{-1}-\eta_l)\right)~=~0.$$
Since $\log(z^{N_{1}}Q(z)(F(z)-\lambda ))$ belongs to $A(\mathbb D)$, 
$$P_+\Big(\log\left(z^{N_{1}}Q(z)(F(z)-\lambda )\right)\Big)~=~\log\left(z^{N_{1}}Q(z)(F(z)-\lambda )\right).$$
Finally, 
\[
F_\lambda^+(z)~=~e^{\log\left(z^{N_{1}}Q(z)(F(z)-\lambda )\right)}~=~z^{N_1}Q(z)(F(z)-\lambda).
\]
Taking \Cref{eve} into account eventually yields  \Cref{eq:vecteur-propre-ABCL} for the eigenvectors $h_{\lambda,j}$.
For any \(\Omega\in\mathcal{C}\) with \(n_{F}(\Omega )=0\), i.e. with $\w_F(\Omega)=-D$, let 
\begin{equation}\label{h_lambda}
h_{\lambda }~=~\frac{1}{F-\lambda }\quad\text{ for every } \lambda \in\Omega.    
\end{equation} We have, for every \(0\le j\le D-1\) and $\lambda\in\Omega$,
\[
h_{\lambda ,j}(z)~=~z^{j}\cdot\dfrac{c_\lambda}{z^{N_{1}}Q(z)}\cdot h_{\lambda }(z)\quad\text{ for every } z\in\mathbb{D}.
\]
Therefore we recover the eigenvectors obtained in \cite[Section 4.1]{AbakumovBaranovCharpentierLishanskii2021}.

\subsection{The case where $F$ is exactly D-valent}\label{Section-N-valent}
In this subsection, we consider a function $F$ defined by \Cref{Form}, satisfying \ref{H1}, and which is exactly $D$-valent on
$\overline{\mathbb{D}}$. This means that for every $\lambda\in F(\overline{\D})$, the equation $F(z)=\lambda$ has exactly $D$ solutions in $\overline{\D}$. Moreover, for any $\lambda\notin F(\mathbb T)$, $n_F(\lambda)\leq D$, and \Cref{Eq:WindVal} implies that $F$ satisfies \ref{H3}. 
\par\smallskip
An important consequence of the exact $D$-valence  of $F$ on $\overline{\mathbb D}$ (and not just on $\mathbb D$) is the following: 
\par\smallskip
\begin{lemma}\label{Lemme:exacte-D-valence}
Let $F$ be given by \Cref{Form}, and suppose that $F$ is of class $C^{1+\varepsilon}$ on $\T$ with $F'\neq 0$ on $\T$ (which is true if $F$ satisfies \ref{H1} for some $p>1$). If $F$ is exactly $D$-valent on $\overline{\mathbb D}$ then $F(\T)$ does not intersect the interior of $F(\overline{\D})$. In particular,  $F(\overline{\D})\neq\widehat{\C}$.
\end{lemma}

\begin{proof}
By contradiction, suppose that there exists 
a point  $\lambda_0\in F(\mathbb T)$ and an open neighborhood of $\lambda_0$ contained in $F(\overline{\mathbb D})$.
    Let $z_1,\dots z_s\in\overline{\mathbb D}$ be the distinct solutions of $F=\lambda_0$ in $\overline{\mathbb D}$, and let $1\le i\le s$. If $z_i\in\mathbb D$, we denote by $m_i$ the multiplicity of $z_i$ as a zero of $F-\lambda_0$ (which is holomorphic in a neighborhood of $z_i$). If $z_i\in \mathbb T$, \Cref{Fait:LocalInjectCercle} implies that $F$ is injective on a neighborhood of $z_i$, and thus, as explained above, we set $m_i=1$. Then the exact $D$-valence of $F$ on $\overline{\D}$ implies that
    \[m_1+\dots+m_s=D.\]
We now construct open  and disjoint neighborhoods $V_1,\dots, V_s$ of the points $z_1,\dots, z_s$  (in $\overline{\mathbb{D}}$) as follows. Let $1\le i\le s$. 
\par\medskip
-- If $z_i\in\mathbb D$, choose $\alpha_i>0$ sufficiently small so that $\overline D(z_i,\alpha_i)$ is contained in $\mathbb D$, contains no pole of $F$, and the only zero of $F-\lambda_0$ in  $\overline D(z_i,\alpha_i)$ is $z_i$. Let $0<\varepsilon_i<\inf_{|z-z_i|=\alpha_i}|F(z)-\lambda_0|$. By Rouché's theorem, the equation $F=\lambda$ has exactly $m_i$ solutions in $D(z_i,\alpha_i)$ for every $\lambda\in D(\lambda_0,\varepsilon_i)$.  We set $V_i=D(z_i,\alpha_i)\cap F^{-1}(D(\lambda_0,\varepsilon_i)).$
\par\smallskip
-- If $z_i\in \mathbb T$, then there exists by \Cref{Fait:LocalInjectCercle} a radius $\alpha_i>0$ such that $F$ is injective on $D(z_i,\alpha_i)\cap\overline{\mathbb D}$. We set $V_i=D(z_i,\alpha_i)\cap\overline{\mathbb D}$ in this case.
\par\medskip
Making the radii $\alpha_i$ smaller if necessary, we can assume that the sets $V_i$ are pairwise disjoint. Let $V= V_1\cup\dots\cup V_s$. Then $V$ is an open neighborhood of the set $\{z_1,\dots,z_s\}$ of all the preimages of $\lambda_0$. Since $F$ is of the form (\ref{Form}), 
$F:\overline{\mathbb D}\longrightarrow \widehat{\mathbb C}$ is continuous. So (see for instance \cite[Lemma 4.21]{Forster}), for every $\lambda\in\widehat{\mathbb C}$ and every open neighborhood $V_\lambda$ of $F^{-1}(\{\lambda\})$ in $\overline{\mathbb D}$, 
\begin{equation}\label{eq:ssdqsdqs32323}
\text{there exists an open neighborhood $U_\lambda$ of $\lambda$ in $\widehat{\mathbb C}$ such that $F^{-1}(U_\lambda)\subset V_\lambda$}.
\end{equation}
Apply this to $\lambda_0$ and $V_{\lambda_0}=V$. Then there exists an open neighborhood $U$ of $\lambda_0$ in $\widehat{\mathbb C}$ such that $F^{-1}(U)\subseteq V$.
Since there exists an open neighborhood of $\lambda_0$ contained in $F(\overline{\mathbb D})$
we can assume (taking  $U$ smaller if necessary) that $U$ is a simply connected open neighborhood of $\lambda_0$ contained in $F(\overline{\mathbb D})$.
\par\smallskip
   Let $\lambda\in U$. Then $\lambda\in F(\overline{\mathbb D})$. It follows from the exact $D$-valence of $F$ and from the construction of $V$ that for every $1\le i\le s$, the equation $F=\lambda$ has exactly $m_i$ solutions in $V_i$. So let $i_0$ be such that $z_{i_0}\in\mathbb T$ (such an index $i_0$ exists because $\lambda_0\in F(\mathbb T)$). For every $\lambda\in U$, the equation $F=\lambda$ has exactly one solution in the open subset $W$ of $\overline{\mathbb D}$ given by $W=F^{-1}(U)\cap V_{i_0}$, and so $F:W\longrightarrow U$ is a continuous bijective map. Moreover, using \eqref{eq:ssdqsdqs32323}, it can be checked that its inverse is also continuous, whence $F:W\longrightarrow U$ is an homeomorphism. Hence $F:W\setminus\{z_{i_0}\}\longrightarrow U\setminus \{\lambda_0\}$ is also an homeomorphism, but this is impossible since $W\setminus\{z_{i_0}\}$ is simply connected (because $z_{i_0}\in\partial W$) while $U\setminus \{\lambda_0\}$ is not.
\end{proof}

\Cref{Lemme:exacte-D-valence} implies also that $F(\T)\cap F(\D)=\varnothing$. Hence we obtain as a direct consequence the following fact:

\begin{corollary}\label{exact-N-valence}
Let $F$ satisfy the assumptions of \Cref{Lemme:exacte-D-valence}.
If $F$ is exactly $D$-valent on $\overline{\mathbb D}$, then for every $\lambda\in F(\T)$
the $D$ solutions in $\overline{\D}$ of the equation $F(z)=\lambda$ belong to $\T$. Thus $F(\mathbb T)$ is a Jordan curve on which the point $F(e^{i\theta})$ travels exactly $D$ times, and $\partial\sigma(T_F)=F(\mathbb T)$.
\end{corollary}

\par\smallskip
These arguments allow for a more precise description of the spectrum of $T_F$ in the following two cases : if $F$ is exactly D-valent on $\mathbb{D}$,
it follows from \Cref{Eq:WindVal,spectre-Toeplitz-valence} combined with the fact that for any $\lambda\in\mathbb C\setminus F(\mathbb T)$, either $n_F(\lambda)=0$ or $n_F(\lambda)=D$, that
\[\sigma(T_F)~=~(\mathbb{C}\setminus F(\mathbb{D}))\cup F(\mathbb{T}).\]
If $F$ is exactly D-valent on $\overline{\mathbb{D}}$, 
then  $F(\mathbb{T})\cap F(\mathbb{D})=\varnothing$ by \Cref{exact-N-valence}, and thus
\[\sigma(T_F)~=~\mathbb{C}\setminus F(\mathbb{D}).\]
As we have already seen, if $F$ is exactly $D$-valent on $\overline{\mathbb D}$, then $F$ satisfies \ref{H3}. Moreover, it follows from \Cref{exact-N-valence} that $F$ satisfies also \ref{H2'}. Finally, since by \Cref{exact-N-valence}, $F(\mathbb T)$ is a Jordan curve traveled $D$ times and $\partial\sigma(T_F)=F(\mathbb T)$, the Jordan curve Theorem implies that $\sigma(T_F)\setminus F(\mathbb T)$ is connected. Hence, applying  \Cref{Th:1compCNS-JFA}, we deduce the following result:

\begin{theorem}\label{th-ABCL-1}
Let $F$ be a symbol of the form  (\ref{Form}) which satisfies \ref{H1}. Suppose that $F$ is exactly $D$-valent on $\overline{\mathbb D}$. Then the following assertions are equivalent:\par\smallskip
\begin{enumerate}
    \item $T_F$ is hypercyclic on $H^2$;\par\smallskip
    \item $\overset{\circ}{\sigma(T_F)}\cap\mathbb T\neq\varnothing$.
\end{enumerate}    
\end{theorem}

The implication (2)$\implies$(1) in \Cref{th-ABCL-1} is proved in  \cite[Th. 1.3]{AbakumovBaranovCharpentierLishanskii2021} without the additional condition \ref{H1} on the regularity on $\mathbb T$. 
Here we obtain, under assumption \ref{H1}, a characterization of the hypercyclicity of $T_F$ on $H^2$.

\subsection{The Decreasing Valence Condition}
In this subsection, we discuss some consequences of our \Cref{Th:Simple_Int-Ext} in the case where the symbol $F$ is of the form \eqref{Form} and satisfies the three assumptions \ref{H1}, \ref{H2} and \ref{H3}. Recall that \Cref{Th:Simple_Int-Ext} asserts that if $\Omega,\Omega'$ are two adjacent components of $\sigma(T_F)\setminus F(\mathbb T)$ such that $|\w_F(\Omega)|>|\w_F(\Omega')|$, then $H_{\Omega'}(T_F) \subseteq H_\Omega(T_F)$.
When rewritten in terms of valence conditions, this statement becomes: let $\Omega$ and $\Omega'$ be two adjacent components of $\sigma(T_F)\setminus F(\mathbb T)$. If $n_F(\Omega)<n_F(\Omega')$, then $H_{\Omega'}(T_F) \subseteq H_\Omega(T_F)$. 
\par\smallskip
Given two components $\Omega$ and $\Omega'$ of $\sigma(T_F)\setminus F(\mathbb T)$, we write $\Omega\overset{D.V.}{\longrightarrow}\Omega'$ if there exists a sequence $(\Omega_j)_{0\le j\le r}$ of adjacent components of $\sigma(T_F)\setminus F(\mathbb T)$ such that $\Omega_0=\Omega,\,\Omega_r=\Omega'$, and  $n_F(\Omega_j)>n_F(\Omega_{j+1})$ for every $0\le j<r$. With the notation of \Cref{Subsection:CompMax}, this is equivalent to requiring that 
$\Omega\overset{I.W.}{\longrightarrow}\Omega'$.
\par\smallskip
In the present context, a connected component $\Omega$ of $\sigma(T_F)\setminus F(\mathbb T)$ satisfies $|\w_F(\Omega)|=D$ if and only if $n_F(\Omega)=0$. Moreover, by definition, $\lambda$ belongs to $ F(\mathbb D)$ if and only if $ n_F(\lambda)\geq 1$. Hence the union of all connected components $\Omega$ of $\mathbb C\setminus F(\mathbb T)$ such that $n_F(\Omega)=0$ is
\begin{equation}\label{chose}
\left\{\lambda\in\C\setminus F(\T)\,;\,   n_F(\lambda)=0\right\}~=~ \big(\C\setminus F(\T)\big)\setminus F(\D)=\C\setminus F(\overline{\D}).
\end{equation}
Since $F:\overline{\D}\longrightarrow\widehat{\C}$ is continuous,
$F(\overline{\D})=\overline{F(\D)}$, and so the right hand term in \Cref{chose} coincides with $\C\setminus\overline{F(\D)}$.
We have:
\begin{theorem}\label{th-7.3}
Let $p>1$, and let $F$ be 
of the form \eqref{Form}, satisfying the three assumptions \ref{H1}, \ref{H2} and \ref{H3}.
Suppose that
 for any connected component $\Omega$ of $\sigma(T_F)\setminus F(\mathbb T)$, there exists a connected component $\Omega_0$ of
 $\sigma(T_F)\setminus F(\mathbb T)$ with $n_F(\Omega_0)=0$ such that 
\[\Omega_0\cap \mathbb T~\neq~\varnothing\quad\text{and}\quad\Omega ~\overset{D.V.}{\longrightarrow}~\Omega_0.\]
Then $T_F$ is hypercyclic on $H^p$.
\end{theorem}
Under the assumptions of \Cref{th-7.3}, the maximal components of $\sigma(T_F)\setminus F(\T)$ are exactly the components of $\mathbb C\setminus F(\overline{\mathbb D})$, i.e. the components $\Omega$ of $\sigma(T_F)\setminus F(\T)$ with $|\w_F(\Omega)|=D$. So \Cref{th-7.3} is a direct consequence of \Cref{Th:HcInstersCompMax}.
\begin{remark}
Remark that when $p=2$, the assumption \ref{H2} in \Cref{th-7.3} can be weakened into supposing that \ref{H2'} holds, and that $n_e=0$ almost everywhere on $F(\T)\setminus\mathcal{O}$. This follows directly from \Cref{Th:HcInstersCompMax_JFA}.
\end{remark}
In~\cite[Th.~1.5]{AbakumovBaranovCharpentierLishanskii2021}, the authors obtain a result similar to our \Cref{th-7.3}, but allowing for two distinct components intersecting $\mathbb{D}$ and $\mathbb{C}\setminus\overline{\mathbb{D}}$, respectively; that is, $\Omega ~\overset{D.V.}{\longrightarrow}~\Omega_{\pm}$ with
\[
\Omega_{+} \cap \bigl(\mathbb{C}\setminus\overline{\mathbb{D}}\bigr) \neq \varnothing
\quad \text{and} \quad
\Omega_{-} \cap \mathbb{D} \neq \varnothing.
\]
However, their theorem requires more regularity than  our \Cref{th-7.3}, since they assume that $F$ is analytic on a neighborhood of $\mathbb T$, that the curve $F(\mathbb T)$ has finitely many self-intersection points which are simple and transverse, and that $F'$ does not vanish on $\mathbb T$.
Their result can be generalized using \Cref{Th:HcInstersCompMax-faible}, leading to the following result, where we recall that $N=\max\{|\w_F(\lambda)|\,;\, \lambda\notin F(\T)\}$.
\begin{theorem}\label{th 7.4}
Let $p>1$, and let $F$ be of the form \eqref{Form}, satisfying the assumptions \ref{H1}, \ref{H2}, \ref{H3} and \ref{H4}. Suppose that $N\ge 2$, and that 
for every connected component $\Omega$ of $\sigma(T_F)\setminus F(\mathbb T)$ with $n_F(\Omega)>0$, there exist  two connected components $\Omega_{+}$ and $\Omega_{-}$ of $\sigma(T_F)\setminus F(\mathbb T)$ with $n_F(\Omega_{+})=n_F(\Omega_{-})=0$ such that
\[\Omega_-\cap\mathbb D~\neq~\varnothing,\quad \Omega_+\cap(\mathbb C\setminus\overline{\mathbb D})~\neq~\varnothing\quad{and}\quad \Omega ~\overset{D.V.}{\longrightarrow}~\Omega_{\pm}.\]
Then $T_F$ is hypercyclic on $H^p$.
\end{theorem}
 Note that, the hypothesis of this theorem imply the existence of a component $\Omega$ with $n_F(\Omega)=0$ and thus $N=D$. Also, remark that our assumption in \Cref{th 7.4} on connected components $\Omega$ of $\sigma(T_F)\setminus F(\mathbb T)$ corresponds to what is called the \emph{Decreasing Valence Condition} (DVC) in \cite{AbakumovBaranovCharpentierLishanskii2021}. Note also that, under the assumptions of \Cref{th 7.4}, the maximal components of $\sigma(T_F)\setminus F(\T)$ are exactly the components of $\mathbb C\setminus F(\overline{\mathbb D})$, i.e. the components $\Omega$ of $\sigma(T_F)\setminus F(\T)$ with $|\w_F(\Omega)|=N=D$. So the assumption $N\ge2$ implies that there is no maximal component with $-1$ as winding number and thus
\Cref{th 7.4} is a direct consequence of 
 \Cref{Th:HcInstersCompMax-faible}. If $N=1$, then $\sigma(T_F)\setminus F(\mathbb T)=\mathbb C\setminus F(\overline{\mathbb D})$ so there is no component $\Omega$ of $\sigma(T_F)\setminus F(\mathbb T)$ with $n_F(\Omega)>0$. In this situation, every component of $\sigma(T_F)\setminus F(\mathbb T)$ is maximal and has $-1$ as winding number, so this case is solved by using \Cref{Th:CNS-IndMax1}.
\par\smallskip
\subsection{About the Increasing Argument Condition}\label{Section 7.3}
In \cite{AbakumovBaranovCharpentierLishanskii2021}, the authors consider symbols satisfying the so-called Increasing Argument Condition (IAC). This condition was considered by Solomyak in \cite{Solomyak1987}, in a study of cyclicity for Toeplitz operators with analytic symbols. When rewritten in terms of the function \(F\), it becomes:
\begin{enumerate}
\item[(IAC)] there exists \(\lambda _{0}\in\C\setminus F(\overline{\mathbb{D}})\) such that\par\smallskip
\begin{itemize}
 \item[\(\bullet\)] the set \(F(\mathbb{T})\) is a finite union of \(C^{2}\) Jordan arcs;\par\smallskip
 \item[\(\bullet\)] some continuous branch of the function \(\theta \longmapsto \arg\,(F(e^{i\theta })-\lambda _{0})\) is strictly decreasing on \([0,2\pi ]\).
\end{itemize}
\end{enumerate}
This condition (IAC) is written in \cite{AbakumovBaranovCharpentierLishanskii2021} in terms of the functions $h_\lambda$, $\lambda\in\C\setminus F(\overline{\mathbb{D}})$ (which are defined in (\ref{h_lambda})), but the two formulations are equivalent.
\par\smallskip

Suppose more generally that \(F\) is of the form \eqref{Form}, satisfies \ref{H1}, \ref{H2}, and \par\smallskip
\begin{enumerate}
 \item [(DAC)]\customlabel{DAC}{(DAC)}  there exists \(\lambda _{0}\in\mathbb{C}\setminus F(\overline{\mathbb{D}}) \) 
 for which some continuous branch of the function \(\theta \longmapsto \arg\,(F(e^{i\theta })-\lambda _{0})\) is strictly decreasing on \([0,2\pi ]\).
\end{enumerate}\par
\begin{figure}[ht]
    \includegraphics[page=42,scale=1.1]{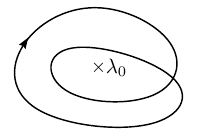}\hspace{1cm}\includegraphics[page=43,scale=.9]{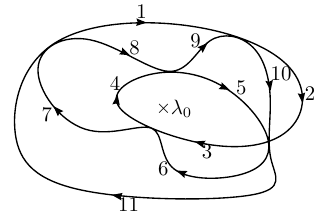}
    \caption{}\label{Fig41ldj}
\end{figure}

\smallskip
Note that, when $F$ is of the form \eqref{Form}, this condition implies that $N=D$.
\Cref{Fig41ldj} provides some examples of curves satisfying \ref{DAC}:\par\smallskip

Suppose that assumptions \ref{H1}, \ref{H2} and \ref{DAC} are satisfied.
For each \(\theta \in[0,2\pi )\), denote by \(\Delta _{\theta }\) the closed half-line having \(\lambda _{0}\) as an extremity, and making an angle \(\theta \) with the half-line \(\lambda _{0}+[0,\infty)\). For all \(\theta \) except finitely many, the function 
\(\lambda \longmapsto \w_{F}(\lambda )\) is well-defined on \(\Delta _{\theta }\setminus F(\mathbb{T})\); as \(\lambda \) goes to infinity along the line \(\Delta _{\theta }\), it takes first the value \(-D\), then the value \(-(D-1)\), \dots, until it reaches the value \(-1\), then \(0\). All values \(-k\), \(0\le k\le D\), are taken in decreasing order, and the jumps between two successive values take place when \(\lambda \) crosses the curve \(F(\mathbb{T})\). There are \(D\) such crossings.
As a consequence, we have 
\begin{fact}\label{Fact 7.3}
 Let \(F\) be of the form \eqref{Form} and satisfy the conditions \ref{H1}, \ref{H2}, and \ref{DAC}. Then 
 \begin{enumerate}[(1)]
  \item the condition \ref{H3} is automatically satisfied, i.e \(\w_{F}(\lambda )\le 0\) for every \(\lambda \in\mathbb{C}\setminus F(\mathbb{T})\);
  \item there exists exactly one component \(\Omega \in\mathcal{C}\) with 
 \(|\w_{F}(\Omega) |=D\);
 \item for any \(\Omega \in\mathcal{C}\) with \(|\w_{F}(\Omega )|<D\), \(\Omega \) is adjacent to a component \(\Omega '\in\mathcal{C}\) with \(|\w_{F}(\Omega ')|>|\w_{F}(\Omega )|\). 
 \end{enumerate}
 \end{fact}

If we denote by \(\Omega _{0}\) the unique connected component of $\sigma(T_F)\setminus F(\T)$ with \(|\w_{F}(\Omega _{0})|=D\), it follows that $\Omega_0$ is the unique maximal component of $\mathcal C$. Applying  \Cref{Th:HcInstersCompMax}, we thus obtain:

\begin{theorem}\label{Theorem 773} 
 Let \(p>1\), and let \(F\) be of the form \eqref{Form} and satisfy \ref{H1} and \ref{H2}. Suppose moreover that there exists \(\lambda _{0}\in\mathbb{C}\setminus F(\overline{\mathbb{D}})\) such that some continuous branch of the argument 
 \[\theta ~\longmapsto~\arg\,\left(F(e^{i\theta })-\lambda _{0}\right)\] is strictly decreasing on \([0,2\pi ]\). 
If \(\Omega _{0}\cap\mathbb{T}\neq\varnothing\), then \(T_{F}\) is hypercyclic on \(H^p\).
\end{theorem}
This result is obtained in \cite[Th. 1.4]{AbakumovBaranovCharpentierLishanskii2021} in the case \(p=2\) and under the additional assumption that the curve \(F(\mathbb{T})\) is piecewise \(C^{2}\)-smooth.

\begin{remark}
When $p=2$, the assumption \ref{H2} in \Cref{Theorem 773} can be weakened into supposing that \ref{H2'} holds. Indeed, the assumption  that there exists \(\lambda _{0}\in\mathbb{C}\setminus F(\overline{\mathbb{D}})\) such that some continuous branch of the argument 
 \(\theta \longmapsto\arg\,(F(e^{i\theta })-\lambda _{0})\) is strictly decreasing on \([0,2\pi ]\) implies that $n_e=0$ almost everywhere on $F(\T)\setminus\mathcal{O}$. Moreover, \Cref{Fact 7.3} remains true, since for all  \(\theta \in[0,2\pi )\) except finitely many,  the function 
\(\lambda \longmapsto \w_{F}(\lambda )\) takes first the value \(-D\), then some value strictly larger than \(-D\), \dots, until it reaches the value \(0\). The jumps between two successive values take place when \(\lambda \) crosses the curve \(F(\mathbb{T})\). There are \(D\) such crossings, counted with multiplicities. It then suffices to apply \Cref{Th:HcInstersCompMax_JFA} instead of \Cref{Th:HcInstersCompMax}.
\end{remark}

\section{Some further results and open questions}\label{Section:AutresResult+Questions}
In this final section, we consider other important properties in linear dynamics in our context of Toeplitz operators. Since the results here follow, for the most part, in a straightforward manner from theorems already proved in the previous sections, we often skip the proofs. We refer to \Cref{Subsection:LinearDyn} for all unexplained terminology.
\subsection{Supercyclicity}
\Cref{GS-S} in \Cref{Section:Rappels} allows us to immediately transpose our sufficient conditions for hypercyclicity into sufficient condition for supercyclicity of $T_F$. It suffices to replace the unit circle by a circle $r\mathbb T$ for some $r>0$ in suitable places, and to keep otherwise the same assumptions. For instance \Cref{Th:HcInstersCompMax} becomes:
\begin{theorem}
Let $p>1$, and suppose that $F$ satisfy \ref{H1}, \ref{H2} and \ref{H3}. If there exists $r>0$ such that $\Omega\cap r\mathbb T\neq\varnothing$ for every maximal component of $\sigma(T_F)\setminus F(\mathbb T)$, then $T_F$ is supercyclic on $H^p$.
\end{theorem}
In the case where $\sigma(T_F)\setminus F(\mathbb T)$ is connected, $T_F$ is always supercyclic on $H^p$. 
On the other hand, some work might be required in order to extend the necessary conditions of \Cref{Section:CNforHc} to the supercyclicity setting. We still have:
\begin{proposition}\label{Prop:OrientationForSC}
Let $F$ be continuous on $\mathbb T$. If there exists $\lambda_0\in\mathbb C\setminus F(\mathbb T)$ such that $\w_F(\lambda_0)>0$, then $T_F$ is not supercyclic.
\end{proposition}
\begin{proof}
The proof relies on the same kind of argument as that of \Cref{Prop:OrientationForHC}: if $\w_F(\lambda_0)>0$ for some $\lambda_0\in\mathbb C\setminus F(\mathbb T)$, then $\w_F(\lambda)>0$ for \textit{all} $\lambda\in\Omega_0$, where $\Omega_0$ is the unique element of $\cal C$ containing $\lambda_0$. According to \Cref{Lemma:spectral-properties-toeplitz-operators}, for all $\lambda\in\Omega_0$, we should have $\dim(\ker(T_{F}^*-\lambda))=\w_F(\lambda)$. Hence, every $\lambda\in\Omega_0$ is an eigenvalue of $T_F^*$. Since the point spectrum of the adjoint of a supercyclic operator can contain at most one element (see for instance \cite[Prop 1.26]{BayartMatheron2009}), $T_F$ cannot be supercyclic.
\end{proof}
An analogue of \Cref{Th:CNforHC-Intersection} still holds too.
\begin{theorem}\label{Th:CNforSC-Intersection}
Let $p>1$, and let $F$ satisfy \ref{H1}, \ref{H2} and \ref{H3}. If $T_F$ is supercyclic on $H^p$, then there exists $r>0$ such that every connected component of the interior ${\sigma(T_F)}$ intersects $r\mathbb T$.
\end{theorem}
\begin{proof}
Proceeding as in the proof of \cite[Th. 1.24]{BayartMatheron2009}, consider the finite set $\cal G$ of connected components of the interior of ${\sigma(T_F)}$. For each $O\in\cal G$, let $I_O=\{|z|\,;\,z\in O\}$: this is a non empty open interval in $\mathbb R^+$. So we need to prove that when $T_F$ is supercyclic, the intersection of all intervals $I_O,~O\in\cal G$, is non-empty.
\par\smallskip
If this intersection is empty, then there exist $r>0$ and two disjoint connected components $O_1,O_2\in\cal G$ such that $O_1 \subseteq r\mathbb D$ and $O_2 \subseteq\mathbb C\setminus r\overline{\mathbb D}$. Indeed, write each interval $I_O$ as $I_O=(a_O,b_O)$, $a_O<b_O$. Then, if $\bigcap_{O\in\mathcal G}I_O=\varnothing$, we have $\min_{O\in\mathcal G}b_O\le \max_{O\in\mathcal G}a_O$. It follows that there exist $O_1,O_2\in\mathcal G$ such that $b_{O_1}\le a_{O_2}$, and hence $I_{O_1}\cap I_{O_2}=\varnothing$. Taking for instance $r=b_{O_1}$, we have indeed $O_1 \subseteq r\mathbb D$ and $O_2 \subseteq\mathbb C\setminus r\overline{\mathbb D}$. Observe that the fact that there are only finitely many components $O$ to be considered here is crucial to the argument, compared to what is done in \cite[Lemma 1.25]{BayartMatheron2009}. 
\par\smallskip
Let then $O_3=\overset{\circ}{\sigma(T_F)}\setminus(O_1\cup O_2)$. So $\overset{\circ}{\sigma(T_F)}=O_1\cup O_2\cup O_3$ and $O_i\cap O_j=\varnothing$ for $i\neq j$, $1\leq i,j\leq 3$.
Using the fact that $T$ admits an $H^\infty$ functional calculus on the interior of ${\sigma(T_F)}$ (see \Cref{Cor:HinftyCalculus}), there exist three non-trivial closed $T$-invariant subspaces $M_1,M_2$ and $M_3$ of $H^p$ such that $H^p$ is the topological direct sum of $M_1$, $M_2$ and $M_3$, i.e.  $H^p=M_1\oplus M_2\oplus M_3$, and such that if we denote by $T_i$ the operator induced by $T$ on $M_i$, then $\sigma(T_i)=\overline{O_i},i=1,2$ (see the proof of \Cref{Prop:CNgeneraHC}). 

Since $T$ admits an $H^\infty$ functional calculus on $\overset{\circ}{\sigma(T_F)}$, the operator $T_i$ admits an $H^\infty(O_i)$-functional calculus, so that in particular there exists a constant $C>0$ such that $\|T_1^n\|\le Cr^n$ and $\|T_2^{-n}\|\le Cr^{-n}$ for every $n\ge0$.
\par\smallskip
Suppose now by contradiction that $T_F$ is supercyclic on $H^p$, and let $x=x_1+x_2+x_3$ be a supercyclic vector for $T_F$, with $x_i\in M_i,i=1,2,3$. Then each operator $T_i$ is supercyclic on $M_i$, with supercyclic vector $x_i$ (see \Cref{lemme-HC-direct-sum}). In particular, $x_i\neq0$ for every $i=1,2,3$.
The fact that $x$ is a supercyclic vector for $T_F$ implies that there exist a sequence $(\lambda_k)$ of complex numbers and a sequence $(n_k)$ of positive integers such that
\[\lambda_kT^{n_k}x~=~\lambda_kT_1^{n_k}x_1\oplus \lambda_kT_2^{n_k}x_2\oplus \lambda_kT_3^{n_k}x_3~\longrightarrow ~x_1\oplus0\oplus0\]
as $k\longrightarrow+\infty$. Then
\[\lambda_kr^{n_k}(r^{-n_k}T_1^{n_k}x_1)~\longrightarrow~ x_1\quad\text{ and }\quad\lambda_kr^{n_k}(r^{-n_k}T_2^{n_k}x_2)~\longrightarrow~ 0 \quad\text{ as }k\to+\infty.\]
Since $\|x_2\|\le\|T_2^{-n_k}\|\|T_2^{n_k}x_2\|\le C\,r^{-n_k}\|T_2^{n_k}x_2\|$ for every $k$, the sequence $( r^{-n_k}\|T_2^{n_k}x_2\|)_k$ is bounded away from $0$ and hence $|\lambda_k|\,r^{n_k}\to 0$ as $k\to+\infty$. But since the sequence $( r^{-n_k}\|T_1^{n_k}x_1\|)_k$ is bounded, this implies that $x_1=0$, which is not the case. Hence $T_F$ cannot be supercyclic on $H^p.$
\end{proof}
We can hence extend our necessary and sufficient conditions for hypercyclicity to the context of supercyclicity. For instance, we have:
\begin{theorem}\label{Th:CNSforSC}
Under the assumptions of \Cref{Th:CNShc2}, the following assertions are equivalent:
\begin{enumerate}
    \item $T_F$ is supercyclic on $H^p$;
    \item there exists $r>0$ such that $O\cap r\mathbb T\neq\varnothing$ for every connected component $O$ of $\overset{\circ}{\sigma(T_F)}$.
\end{enumerate}
\end{theorem}
We leave it to the reader to formulate the proper analogues of the other results of \Cref{Section:ApplicationCShc} for supercyclicity.

\subsection{Chaos, frequent hypercyclicity and ergodicity}
As explained in \Cref{Subsection:LinearDyn}, the notions of chaos, frequent hypercyclicity and ergodicity are most easily investigated via unimodular eigenvectors. And sometimes, the existence of suitably many unimodular eigenvectors is a prerequisite - this is the case for chaos, as well as for ergodicity with respect to a Gaussian measure with full support.

It follows from \Cref{Coro:A12} that many of our results can be extended to yield chaos, frequent hypercyclicity, and ergodicity. In the case considered in \Cref{Th:ShkarinHp}, where  the symbol $F$ has the form $F(e^{i\theta})=ae^{-i\theta}+b+ce^{i\theta}$, with $a,b,c\in\mathbb C$, chaos and frequent hypercyclicity of the Toeplitz operator $T_F$ acting on $H^2$ have already been investigated in \cite{chaos-Toeplitz}.

\begin{theorem}\label{TheoremT} 
Let $F$ satisfy \ref{H1}, \ref{H2} and \ref{H3} and let $p>1$. Suppose that either
\begin{enumerate}[(i)]
    \item $\sigma(T_F)\setminus F(\mathbb T)$ is connected and $\overset{\circ}{\sigma(T_F)}\cap\mathbb T\neq\varnothing$;
    
    \hspace{-1.65cm}{or}
    
    \item $\Omega\cap\mathbb T\neq\varnothing$ for every maximal component $\Omega$ of $\sigma(T_F)\setminus F(\mathbb T)$.
\end{enumerate}
Then the operator $T_F$ acting on $H^p$ is chaotic, frequently hypercyclic, and ergodic with respect to a Gaussian measure with full support.
\end{theorem}

The proof of \Cref{TheoremT} relies on the following consequence of \Cref{Coro:A12}.
\begin{proposition}\label{Prop8.4bis}
Let $F$ satisfy \ref{H1}, \ref{H2} and \ref{H3}, and let $p>1$. Suppose that $\Omega_1,\dots,\Omega_r$ are some connected components of $\sigma(T_F)\setminus F(\mathbb T)$ such that 
\[\overline{\spa}\,\big[H_{\Omega_i}(T_F)\,;\,1\le i\le r\big]~=~H^p.\]
If $\Omega_i\cap\mathbb T\neq\varnothing$ for every $1\le i\le r$ then $T_F$ is chaotic, frequently hypercyclic and ergodic with respect to a Gaussian measure with full support.
\end{proposition}
\begin{proof}
Let $U=\Omega_1\cup\dots\cup\Omega_r$. For every $0\le j<N$ and every $1\le i\le r$, define $E_j$ on $\Omega_i$ by setting $E_j(\lambda)=h_{\lambda,j}$ for every $\lambda\in \Omega_i$ if $\Omega_i  \subseteq\Omega_j^+$, and $E_j(\lambda)=0$ for every $\lambda\in \Omega_i$ if $\Omega_i\not  \subseteq\Omega_j^+$. This defines an analytic map $E_j:U\longrightarrow H^p$. Moreover the closed linear span of the vectors $E_j(\lambda),\,\lambda\in U,\,0\le j<N$ contains the subspace $H_{\Omega_i}(T_F)$ for every $1\le i\le r$. Hence
\[\overline{\spa}\,\big[E_j(\lambda)\,;\,\lambda\in U\,,\,0\le j<N\big]=H^p.\]
So the assumptions of \Cref{Coro:A12} are satisfied and \Cref{Prop8.4bis} follows.
\end{proof}
\begin{proof}[Proof of \Cref{TheoremT}] Under assumption (i), \Cref{TheoremT} is an immediate consequence of \Cref{Prop8.4bis} and of the fact that in this case, $\sigma(T_F)\setminus F(\mathbb T)$ has exactly one connected component $\Omega$ -- the interior of ${\sigma(T_F)}$ -- which satisfies $H_\Omega(T_F)=H^p$.
\par\smallskip
Under assumption (ii), we enumerate as $\Omega_1,\dots,\Omega_r$ the set of all maximal connected components of $\sigma(T_F)\setminus F(\mathbb T)$. The proof of \Cref{Th:HcInstersCompMax} shows that $\overline{\spa}\,[H_{\Omega_i}(T_F)\,;\,1\le i\le r]=H^p$, which concludes the proof in this case too.
\end{proof}
\Cref{TheoremT} can be viewed as a generalization of \Cref{Th:ConnexeHC,Th:CNS-IndMax1,Th:HcInstersCompMax} to the notions of chaos, frequent hypercyclicity and ergodicity. Here is the corresponding generalization of \Cref{Th:CNShc1,Th:CNShc2,Th:Case_wind2}.
\begin{theorem}\label{Theo:fkjeowigo}
    Under the assumptions of \Cref{Th:CNShc2}, the following assumptions are equivalent:
    \begin{enumerate}
     \item the operator $T_F$ acting on $H^p$ is chaotic, frequently hypercyclic, and ergodic with respect to a Gaussian measure with full support;
        \item the unit circle $\mathbb T$ intersects every component of $\overset{\circ}{\sigma(T_F)}$.
           \end{enumerate}
\end{theorem}
\begin{proof}
The implication (1)$\implies$(2) is given by \Cref{Th:CNShc2}, and we just have to prove that (2)$\implies$(1). 

So let us suppose that $\mathbb T$ intersects every connected component of $\overset{\circ}{\sigma(T_F)}$. Our aim is to show that for any dense subset $A$ of $\mathbb T$, we have $H_A(T_F)=H^p$. Consider a function $g\in H^q$ which vanishes on $H_A(T_F)$. Let $\Theta$ be a component of the interior of $\sigma(T_F)$. Then $\Theta\cap \mathbb T\neq\varnothing$, and the proof of \Cref{Th:CNShc2} yields that for any component $\Omega$ of $\Theta\setminus F(\mathbb T)$,  we have $H_\Omega(T_F)=H_\Theta(T_F)$. We consider now separately two cases.
\par\medskip
-- If there exists a connected component $\Omega$ of $\Theta\setminus F(\mathbb T)$ such that $\Omega\cap \mathbb T\neq\varnothing$, then $A$ has an accumulation point in $\Omega$ and thus, thanks to \Cref{Prop:EigenvectorDenseYaku}, $H_{A\cap\Omega}(T_F)=H_\Omega(T_F)$. In particular, the function $g$ vanishes also on $H_\Omega(T_F)=H_\Theta(T_F)$.
\par\smallskip
-- If such an $\Omega $ does not exist, then $\mathbb T\cap \Theta$ is contained in $ F(\mathbb T)$, so there exists a closed subarc $\gamma$ of $\T\cap(F(\mathbb T)\setminus\mathcal O)\cap \Theta$ of positive length. Consider now the exterior component  $\Omega$ of $\gamma$. Applying \Cref{Coro:Ajout29390} below to the set $\gamma\cap A$ (whose closure is $\gamma$, which has positive length) will allow us to conclude that $H_\Omega(T_F)\subseteq H_{A\cap\gamma}(T_F)$, and thus that $g$ vanishes on $H_\Theta(T_F)$.

\begin{lemma}
    \label{Coro:Ajout29390} Let $\Omega\in\mathcal C$, and let $\gamma_\Omega$  be the union of all the subarcs of $F(\mathbb T)\setminus \mathcal O$ which have $\Omega$ as exterior component.
    Let $A$ be a subset of $\gamma_\Omega$ such that the set $\overline A$ has positive length. Then
    \[H_\Omega(T_F)\subseteq H_A(T_F).\]
\end{lemma}

\begin{proof}
    Suppose that $g\in H^q$ is a function which vanishes on $H_A(T_F)$.  By \Cref{Prop:ContinuiteBord}, the function $\Phi_j^\Omega:\Omega\cup\gamma_\Omega\longrightarrow A(\mathbb D)$ defined by
\[\Phi_j^\Omega(\lambda)=\begin{cases}
h_{\lambda,j}&\text{ if }\lambda\in\Omega\\
h_{\lambda,j}^{ext}&\text{ if }\lambda\in\gamma_\Omega
\end{cases}\]
is continuous on $\Omega\cup\gamma_\Omega$. It follows that
for every point $\lambda_0\in A$ and every integer $j$ with $0\le j<|\w_F(\Omega)|$, we have that $\dual{h_{\lambda_0,j}^{ext}}{g}=0$. Indeed, $(T_F-\lambda)h_{\lambda,j}=0$ for every $\lambda\in\Omega$, and since $(T_F-\lambda)h_{\lambda,j}$ tends to $(T_F-\lambda_0)h_{\lambda_0,j}^{ext}
$ as $\lambda$ tends to $\lambda_0$, $\lambda\in\Omega$, $h_{\lambda_0,j}^{ext}$ belongs to $\ker(T_F-\lambda_0)\subseteq H_A(T_F)$. \par\smallskip
Let $(u_0,\dots,u_{N-1})=Ug$. Then  $u_j(\lambda)=\dual{\Phi_j^\Omega(\lambda)}{g}$ for all $\lambda\in\Omega$, and thus $u_j$ admits a continuous extension to $\Omega\cup\gamma_\Omega$ which is given by 
    \[u_j(\lambda)=\dual{\Phi_j^\Omega(\lambda)}{g}=\begin{cases}
        \dual{h_{\lambda,j}}{g}&~\text{if }\lambda\in\Omega\\
        \dual{h_{\lambda,j}^{ext}}{g}&~\text{if }\lambda\in\gamma_\Omega.
    \end{cases}\]
    In particular, this continuous extension satisfies $u_j(\lambda)=0$ for every $\lambda\in\overline{A}$. But  since $u_j$ belongs to $E^q(\Omega)$ and the length of $\overline{A}$ is positive, this implies that $u_j=0$ on $\Omega$. So $g$ vanishes on all the eigenspaces $\ker(T_F-\lambda)$, $\lambda\in\Omega$, and thus $g$ vanishes on $H_{\Omega}(T_F)$. 
\end{proof}

 We can now conclude the proof of \Cref{Theo:fkjeowigo}. In each one of the two cases above, we have shown that $g$ vanishes $H_\Theta(T_F)$. This being true for every component $\Theta$ of the interior of $\sigma(T_F)$, we obtain that $g=0$ and thus we have shown that $H_A(T_F)=H^p$ for every dense subset $A$ of $\mathbb T$.
As a consequence, it follows that if $D$ is any countable subset of $\T$, then $H_{\T\setminus D}(T_F)=H^p$. Hence
$T_F$ has perfectly spanning unimodular eigenvectors. So $T_F$ acting on $H^p$ is frequently hypercyclic, and ergodic with respect to a Gaussian measure with full support by \Cref{Th:TheoremA13,Theo:A19}.  Applying the same reasoning to the set $A$ of all $n$-th roots of unity, we also deduce that $T_F$ is chaotic. 
\end{proof}

 \subsection{From an exterior component to an interior component}
 
 One of the most obvious problems which arises from our work is the following: is it really necessary to add an assumption to be able to go from an exterior component to an interior component?
 
 \begin{question}\label{Question2}
 Suppose that $F$ satisfies \ref{H1}, \ref{H2} and \ref{H3} and let $\Omega,\Omega'\in \cal C$ be two adjacent components of $\sigma(T_F)\setminus F(\mathbb T)$ such that $|\w_F(\Omega)|<|\w_F(\Omega')|$. Is it true that $H_{\Omega'}(T_F) \subseteq H_\Omega(T_F)$?
 \end{question}
 If \Cref{Question2} had an affirmative answer, then \Cref{Th:CNShc1,Th:CNShc2,Th:Case_wind2}  would hold in much greater generality. In this regard, here are two examples which give food for thought.
 Remark that the approach in \Cref{ExampleC1} (resp. in \Cref{ExampleC2}) is very similar to that of the example discussed in \Cref{SubSection:Example2Circles} (resp. to \Cref{Exemple:2cercle-JFA}).

 \begin{example}\label{ExampleC1}
 Consider the map $F:\mathbb T\to\mathbb C$ defined in the following way:
 \[
 F(e^{i\theta})~=~\begin{cases}
 3+3e^{-5i\theta/3}&\text{if }0\le\theta<6\pi/5\\
 \frac92+\frac32e^{-10i\theta/3}&\text{if }6\pi/5\le\theta<9\pi/5\\
 \frac{11}2+\frac12e^{-10i\theta}&\text{if }9\pi/5\le\theta<2\pi.
 \end{cases}
 \]
 \Cref{Fig43klel} gives a picture of the curve $F(\mathbb T)$.
\begin{figure}[ht]
\includegraphics[page=45,scale=.85]{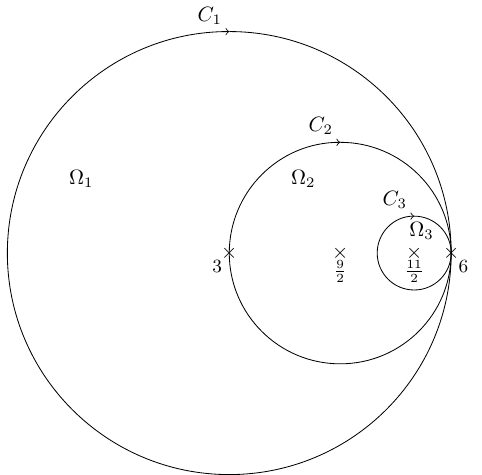}
             \caption{}
             \label{Fig43klel}
             \end{figure}

 The respective radii of the circles $C_1, C_2$ and $C_3$ are $3,\,3/2$ and $1/2$. We have 
 $$\w_F(\Omega_1)~=~-1,\quad\w_F(\Omega_2)~=~-2\quad\text{ and}\quad\w_F(\Omega_3)~=~-3.$$
 Among the connected components of $\sigma(T_F)\setminus F(\mathbb T)$, $\mathbb T$ intersects only $\Omega_1$. In order to prove that $T_F$ is hypercyclic (on $H^p,p>1)$, starting from $u_0\in E^q(\Omega_1\cup\Omega_2\cup\Omega_3),\,u_1\in E^q(\Omega_2\cup\Omega_3)$ and $u_2\in E^q(\Omega_3)$ with $u_0=0$ on $\Omega_1$ and the properties
 \begin{enumerate}[(i)]
     \item\label{i} $u_0^{int}-\zeta u_1^{int}=0$ almost everywhere on $C_2$ (since $u_0=0$ on $\Omega_1$, $u_0^{ext}=0$ on $C_2$)
     \item $u_0^{int}-\zeta u_1^{int}=u_0^{ext}$ almost everywhere on $C_3$
     \item $u_1^{int}-\zeta u_2^{int}=u_1^{ext}$ almost everywhere on $C_3$,
 \end{enumerate}
we need to be able to deduce that $u_0,u_1$ and $u_2$ are identically zero.
  However, we can use neither \Cref{Th:Complique_Ext-Int} nor \Cref{Th:Complique-Jordan_version} to deduce that $u_0=u_1=0$ on $\Omega_2$, as neither of their assumptions are satisfied. Indeed, $\partial\Omega_2=C_2\cup C_3$ is not a Jordan curve and there does not exist an open neighborhood $V$ of $\lambda_0=6$ (which is the only point of self-intersection of $F(\mathbb T)$) which satisfies $V\cap\partial \Omega_2=V\cap\partial\Omega_1\cap\partial\Omega_2$. 
 Nonetheless, working directly with the explicit expression of $\zeta$ on $C_2$, it is possible to infer what we want from (i), (ii) and (iii).

 Let $\lambda\in \partial C_2\setminus\{6\}$. In order to determine $\zeta(\lambda)=e^{i\theta(\lambda)}$, we need to solve the equation
 \[\lambda~=~\frac92+\frac32e^{10i\theta(\lambda)/3},\quad\text{with}~\theta(\lambda)\in\left(-\frac{9\pi}5,-\frac{6\pi}5\right).\]
 We obtain
 \(\theta(\lambda)~=~\frac3{10}\arg_{(-6\pi,-4\pi)}\left(\frac{2\lambda-9}3\right),\)
 and thus 
 \[
 \zeta(\lambda)~=~\exp\left[\frac{3i}{10}\arg_{(-6\pi,-4\pi)}\left(\frac{2\lambda-9}3\right)\right]\quad\text{ for every }\lambda\in\partial C_2\setminus\{6\}.
 \]
So $\zeta$ admits an analytical continuation to $\mathbb C\setminus [9/2,\infty)$, hence in particular on $\Omega_2\setminus[9/2,5)$, given by 
 \[\zeta(\lambda)~=~e^{\frac3{10}\log_{(-6\pi,-4\pi)}\left(\frac{2\lambda-9}3\right)}\quad\text{ for every } \lambda\in \Omega_2\setminus[9/2,5),\]
 where $\log_{(-6\pi,-4\pi)}$ is the analytic determination of the logarithm with imaginary part in $(-6\pi,-4\pi)$. Note that this extension is such that, for every $x\in(9/2,5)$ we have
\begin{equation}\label{Eq:ExampleC1-lim}
\lim_{\underset{y>0}{y\to0}}\zeta(x+iy)~=~\left|\frac{2x-9}3\right|^{3/10}e^{-9i\pi/5}~\neq~ \left|\frac{2x-9}3\right|^{3/10}e^{-6i\pi/5}~=~\lim_{\underset{y<0}{y\to0}}\zeta(x+iy).
\end{equation}
 Since $\zeta$ is bounded on $\Omega_2\setminus[9/2,5)$, the function $\zeta u_1$ belongs to $E^q(\Omega_2\setminus[9/2,5))$, and so $u_0-\zeta u_1\in E^q(\Omega_2\setminus[9/2,5))$. Since $u_0^{int}-\zeta u_1^{int}=0$ almost everywhere on $C_2$, which is a subset of $\partial (\Omega_2\setminus[9/2,5))$ with positive measure, we deduce that $u_0=\zeta u_1$ on $\Omega_2\setminus[9/2,5)$.
 \par\smallskip
 Suppose now that $u_1$ is not identically zero on $\Omega_2$. Then there exists $x\in(9/2,5)$ such that $u_1(x)\neq0$. Hence $\zeta$ admits an analytic extension to a neighborhood of $x$ in $\Omega_2$, given by $\zeta(z)=u_0(z)/u_1(z)$ for $|z-x|<r$, where $r>0$ is small enough. But this contradicts \Cref{Eq:ExampleC1-lim}. Hence $u_1=0$ on $\Omega_2$ and $u_0=0$ on $\Omega_2$ too. By the same argument (or by invoking either Theorem \ref{Th:Complique_Ext-Int} or Theorem \ref{Th:Complique-Jordan_version}), we deduce from (ii) and (iii) that $u_1=u_2=u_3=0$ on $\Omega_3$.
 So $T_F$ is indeed hypercyclic on $H^p$. 
 \end{example}
 Observe that our argument here depends in a crucial way of the fact that $\zeta$ does not admit any continuous extension to the whole domain $\Omega_2$, which is ensured because the interval $[9/2,6)$ is not contained in $\Omega_3$. But what happens if $[9/2,6) \subseteq\Omega_3$? We consider this situation in the next example.
 \begin{example}\label{ExampleC2}
Let $F:\mathbb T\to\mathbb C$ be defined by
 \[
 F(e^{i\theta})~=~\begin{cases}
 3+3e^{-11i\theta/6}&\text{if }0\le\theta<12\pi/11\\
 \frac92+\frac32e^{-11i\theta/3}&\text{if }12\pi/11\le\theta<18\pi/11\\
5-e^{-11i\theta/2}&\text{if }18\pi/11\le\theta<2\pi.
 \end{cases}
 \]
 In \Cref{Fig44gljwpo} the circle $C_3$ has radius $1$, so that the picture of $F(\mathbb T)$ looks like this:
 \begin{figure}[ht]
\includegraphics[page=46,scale=.85]{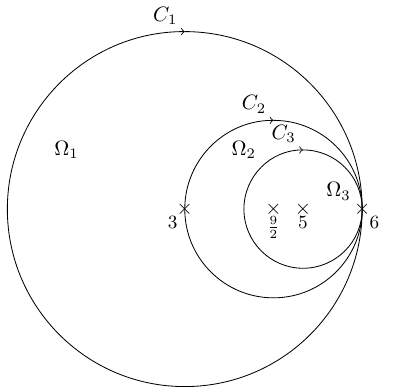}
            \caption{}
            \label{Fig44gljwpo}
            \end{figure}

 Denoting by $\zeta_2$ and $\zeta_3$ the map $1/F^{-1}$ on $\partial C_2\setminus\{6\}$ and $\partial C_3\setminus\{6\}$ respectively, we have:
 \[
 \zeta_2(\lambda)~=~\exp\left[\frac{3i}{11}\arg_{(-6\pi,-4\pi}\left(\frac{2\lambda-9}{3}\right)\right]\quad\text{ for every }\lambda\in\partial C_2\setminus\{6\},
 \]
 and 
 \[
 \zeta_3(\lambda)~=~\exp\left[\frac{2i}{11}\arg_{(-11\pi,-9\pi)}(5-\lambda)\right]\quad\text{ for every }\lambda\in\partial C_3\setminus\{6\}.
 \]
 So $\zeta_2$ and $\zeta_3$ admits analytic extensions to $\mathbb C\setminus [9/2,\infty)$ and $\mathbb C\setminus [5,\infty)$ respectively,  given by 
 \[\zeta_2(\lambda)~=~\exp\left[\frac{3i}{11}\log_{(-6\pi,-4\pi)}\left(\frac{2\lambda-9}{3}\right)\right]\quad\text{ for every } \lambda\in \mathbb C\setminus [9/2,\infty),\] 
 and
 \[\zeta_3(\lambda)~=~
 \exp\left[\frac{2i}{11}\log_{(-11\pi,-9\pi)}(5-\lambda)\right]\quad\text{ for every } \lambda\in\mathbb C\setminus [5,\infty).
 \]
 Moreover, for every $x\in (9/2,6)$, the limits
 $$\lim_{\underset{y>0}{y\to0}}\zeta_2(x+iy)\quad\text{ and }\quad \lim_{\underset{y<0}{y\to0}}\zeta_2(x+iy)$$
 exist and are different, and similarly for every $x\in (5,6)$, the limits
 $$\lim_{\underset{y>0}{y\to0}}\zeta_3(x+iy) \quad\text{ and }\quad\lim_{\underset{y<0}{y\to0}}\zeta_3(x+iy)$$
 exist and are different.
 \par\smallskip
 We proceed in the same way as in \Cref{ExampleC1}:
 starting from three analytic  functions  $u_0\in E^q(\Omega_1\cup\Omega_2\cup\Omega_3),\,u_1\in E^q(\Omega_2\cup\Omega_3)$ and $u_2\in E^q(\Omega_3)$ with $u_0=0$ on $\Omega_1$ and
 \begin{enumerate}[(i)]
     \item $u_0^{int}-\zeta_2 u_1^{int}=0$ almost everywhere on $C_2$ (as $u_0^{ext}=0$ on $C_2$)
     \item $u_0^{int}-\zeta_3 u_1^{int}=u_0^{ext}$ almost everywhere on $C_3$
     \item $u_1^{int}-\zeta_3 u_2^{int}=u_1^{ext}$ almost everywhere on $C_3$,
 \end{enumerate}
 we wish to deduce that $u_0$, $u_1$ and $u_2$ are identically zero.
 \par\medskip
 In this example, contrary to what happens in \Cref{ExampleC1}, the relation (i) does not imply that $u_0=u_1=0$ on $\Omega_2$, but only that $u_0=\zeta_2 u_1$ on $\Omega_2$ (since $u_0-\zeta_2 u_1\in E^q(\Omega_2)$, the fact that $u_0-\zeta_2 u_1$ vanishes on a subset of positive measure of the boundary of $\Omega_2$ implies that it vanishes on $\Omega_2$). So we have $u_0^{ext}=\zeta_2u_1^{ext}$ almost everywhere on $C_3$. Using the relations (ii) and (iii), it follows that
 \[ u_0^{int}-\zeta_3u_1^{int}~=~u_0^{ext}~=~\zeta_2u_1^{ext}~=~\zeta_2(u_1^{int}-\zeta_3u_2^{int})\quad\text{a.e. on }C_3,\]
 and thus 
 \[u_0^{int}-(\zeta_2+\zeta_3)u_1^{int}+\zeta_2\zeta_3u_2^{int}~=~0\quad\text{a.e. on }C_3.\]
 Since $\zeta_2$ and $\zeta_3$ both admit analytic and bounded extensions to $\Omega_3\setminus[9/2,6)$, this yields
  \begin{enumerate}[(iv)]
      \item $u_0-(\zeta_2+\zeta_3)u_1+\zeta_2\zeta_3u_2=0~\text{on }\Omega_3\setminus[9/2,6).$
  \end{enumerate}
  Suppose that $u_2$ is not identically zero on $\Omega_3$. Then for all $x\in(9/2,6)$ except countably many, $u_2(x)\neq0$, and the function $u_0\in E^q(\Omega_3)$ which satisfies  $u_0=(\zeta_2+\zeta_3)u_1-\zeta_2\zeta_3u_2$ is continuous at the point $x$. Choose such an $x$ belonging to $(9/2,5)$. Then $\zeta_3$ is continuous at the point $x$, but $\zeta_2$ is not. Taking the limits of $u_0(x\pm iy)$ as $y\to0,y>0$, we obtain the equality
  \begin{multline*}
      \left(\left(\frac{2x-9}3\right)^{3/11}e^{-18i\pi/11}+\zeta_3(x)\right)u_1(x)-\left(\frac{2x-9}3\right)^{3/11}e^{-18i\pi/11}\zeta_3(x)u_2(x)\\
      =~\left(\left(\frac{2x-9}3\right)^{3/11}e^{-12i\pi/11}+\zeta_3(x)\right)u_1(x)-\left(\frac{2x-9}3\right)^{3/11}e^{-12i\pi/11}\zeta_3(x)u_2(x)
  \end{multline*}
which yields that
  \[\left(e^{-18i\pi/11}-e^{-12i\pi/11}\right)u_1(x)~=~\left(e^{-18i\pi/11}-e^{-12i\pi/11}\right)\zeta_3(x)u_2(x).\]
  Since $e^{-12i\pi/11}\neq e^{-18i\pi/11}$ (it is here that we use the fact that the limits of $\zeta_2(x+iy)$ and $\zeta_2(x-iy)$ as $y\longrightarrow 0$, $y>0$, are distinct), we get $u_1(x)=\zeta_3(x)u_2(x)$. This being true for almost every $x\in(9/2,5)$,  by uniqueness of the analytic extension this equality is true on $\Omega_3\setminus[5,6)$. It follows that $u_1/u_2$ is an analytic extension of $\zeta_3$ to $\Omega_3$ minus the set of zeroes of $u_2$. In particular, $\zeta_3$ admits an analytic extension to a neighborhood of some points $x\in(5,6)$, and this is a contradiction. So $u_2=0$ on $\Omega_3$.
  \par\smallskip
  It then follows from equation (iv) and from a similar continuity argument that $u_0=u_1=0$ on $\Omega_3$. Then we deduce as usual that $u_0,u_1$ and $u_2$ vanish identically, and $T_F$ is hypercyclic on $H^p$.
 \end{example}
 These two examples might point towards a positive answer to \Cref{Question2}. In any case, it would be interesting to try to include the arguments used in \Cref{ExampleC1,ExampleC2} into a more general framework.
 
 \subsection{Toeplitz operators and the Godefroy-Shapiro Criterion}
It is a rather intriguing fact that all the Toeplitz operators which are known to be hypercyclic are shown to be so thanks to the Godefroy-Shapiro Criterion.  In this short section, we highlight the following question:
 
 \begin{question}\label{question-GS}
 If $T_F$ is a hypercyclic Toeplitz operator on $H^p$, $p>1$, does $T_F$ necessarily satisfy the Godefroy-Shapiro Criterion? Are there some hypercyclic Toeplitz operators which have no eigenvalue?
 \end{question}
 \subsection{Identical curves with different parametrizations}
 An intriguing fact is that some of our sufficient conditions for the hypercyclicity of $T_F$ do not depend too precisely on the properties of the symbol $F$, but rather on the geometric properties of $F(\mathbb T)$ (\Cref{Th:HcInstersCompMax,Th:CNS-IndMax1}), while some others do (\Cref{Th:CNShc1,Th:CNShc2,Th:Case_wind2}). More precisely, if $F$ and $ \widetilde F$ satisfy \ref{H1}, \ref{H2}, \ref{H3} and are such that $F(\mathbb T)= \widetilde F(\mathbb T)$, then the assumptions of \Cref{Th:HcInstersCompMax}, for instance, are satisfied for $F$ if and only if they are satisfied for $ \widetilde F$. This is certainly not true for  \Cref{Th:CNShc1,Th:CNShc2,Th:Case_wind2}. 

 \begin{question}\label{Question:2param} Let $1<p<\infty$. Can we find two functions $F$ and $\widetilde F$ on $\mathbb T$ satisfying \ref{H1}, \ref{H2} and \ref{H3} with $F(\mathbb T)=\widetilde F(\mathbb T)$, such that $T_F$ is hypercyclic on $H^p$ but  $T_{\widetilde F}$ is not?
 \end{question}
 
 The same question can be asked when $\widetilde F$ is obtained from $F$ by a ``reshuffling'' procedure. Indeed, suppose that $F$ satisfies \ref{H1}, \ref{H2} and \ref{H3}, and let $\alpha_j,\,{1\le j\le m}$, be a collection of subarcs of $\mathbb T$ given by \ref{H2}. Suppose that there exist $z_1,\dots,z_m\in\mathbb T$ such that the arcs $\widetilde \alpha_j=z_j\alpha_j,\,1\le j\le m$, form a partition of $\mathbb T$, and such that the function $\widetilde F$ defined on $\mathbb T$ by $\widetilde F(z)=F(\overline{z_j}z),~z\in\widetilde \alpha_j,\,1\le j\le m$, also satisfies \ref{H1}, \ref{H2} and \ref{H3}. In other words $\widetilde F$ is obtained from $F$ by traveling on the arcs $F(\alpha_j)$ in a different order, while preserving smoothness and ensuring that all winding numbers $\w_{\widetilde{F}}(\lambda),\,\lambda\notin \widetilde F(\mathbb T)$, remain non positive. A natural question is then:
 \begin{question}\label{Ques:2paramReshuffling}
     Under the assumptions above, is it true that $T_{\widetilde{F}}\in B(H^p)$ is hypercyclic as soon as $T_F$ is?
 \end{question}
 Let us illustrate this on the following example.
\begin{example}
 Let $F$ be such that the curve $F(\mathbb T)$ looks like in Figure \textsc{\ref{FIG47a}}. The order in which the point $F(e^{i\theta})$ travels the curve $F(\mathbb T)$, as $\theta$ grows from $0$ to $2\pi$, is given by the numbers $1$ to $9$ in the picture. And let $ \widetilde F$ be such that $ \widetilde F(\mathbb T)$ is the same curve but with a different order for the travel of $ \widetilde F(e^{i\theta})$
 on $ \widetilde F(\mathbb T)$ as represented in Figure \textsc{\ref{FIG47b}}.
\begin{figure}[ht]
  \begin{subfigure}[b]{0.45\textwidth}
  \includegraphics[page=47,scale=.75]{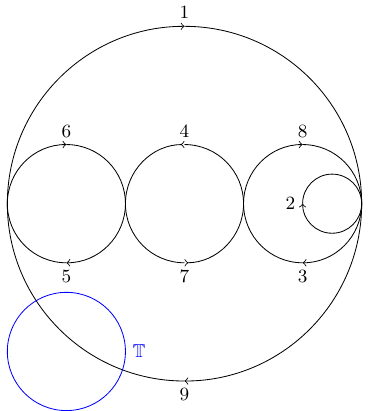}
  \caption{}\label{FIG47a}
  \end{subfigure}
  \hfill
  \begin{subfigure}[b]{0.45\textwidth}
  \includegraphics[page=48,scale=.75]{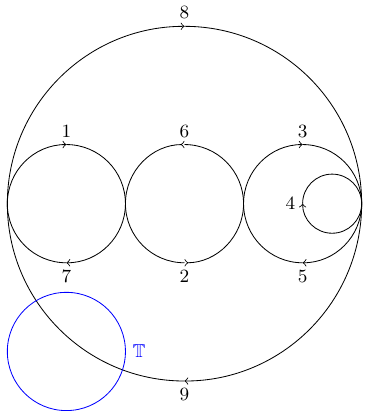}
  \caption{}\label{FIG47b}
  \end{subfigure}
             \caption{}\end{figure}

 Suppose that $\mathbb T$ intersects $\overset{\circ}{\sigma(T_F)}=\overset{\circ}{\sigma(T_{ \widetilde F})}$, and that $F$ and $ \widetilde F$ both  have an analytic extension to a neighborhood of $\mathbb T$. Then $T_F$ is hypercyclic on $H^p$ by
 \Cref{Prop:Ext-Int_AutreCas} and \Cref{Th:Simple_Int-Ext}, while our arguments do not allow us to conclude that $T_{\widetilde{F}}$ is hypercyclic. 
\end{example}
Another natural question concerns the links between the dynamical behavior of $T_F$ and $T_{\widetilde F}$ when $\widetilde F$ is obtained from $F$ by composing it with a sufficiently smooth orientation-preserving diffeomorphism of $\mathbb T$:
\begin{question}\label{Ques:2paramDiffeo}
Let $F$ be a symbol on $\T$ satisfying assumptions \ref{H1}, \ref{H2} and \ref{H3}, and let $\eta$ be an orientation-preserving diffeomorphism of $\T$ such that $\eta$ and $\eta^{-1}$ are of class $C^2$, for instance. Let $\widetilde F=F\circ\eta$. Is it true that $T_{\widetilde F} \in B(H^p)$ is  hypercyclic whenever $T_F$ is?
\end{question}

 A possible approach to this question would be via quasi-similarity. Indeed, if $T_1$ and $T_2$ are any two operators on a Banach space $X$ for which there exists $A\in \mathcal{B}(X)$ with dense range such that $AT_1=T_2A$, then $T_2$ is hypercyclic on $X$ as soon as $T_1$ is. So in particular, quasi-similar operators are simultaneously hypercyclic or non-hypercyclic (recall that $T_1,T_2\in\mathcal{B}(X)$ are quasi-similar if there exist $A,B\in \mathcal{B}(X)$ one-to-one with dense range such that $AT_1=T_2A$ and $T_1B=BT_2$). 
 \par\smallskip
 Similarity of Toeplitz operators was investigated by Clark in a series of papers \cites{Clark-1980-1,Clark-1980-2, Clark-Toeplitz-2}, in the context of rational Toeplitz operators $T_F$ such that $\sigma(T_F)\setminus F(\mathbb T)$ consists of a finite union of loops, intersecting at a finite number of points only, and by Stephenson in \cite{Stephenson1985}. It is shown in \cite{Stephenson1985} that if $F$ and $G$ are two rational symbols such that $G=F\circ \psi$ where $\psi=B_2^{-1}\circ B_1$ for some finite Blaschke products $B_1$ and $B_2$, then $T_F$ and $T_G$ are similar (as operators on $H^2$). An example is provided in \cite{Stephenson1985} of two symbols $F$ and $G$ satisfying \ref{H1}, \ref{H2} and \ref{H3}, which both map $\mathbb T$ onto a certain curve $\gamma$ consisting of three loops, but which are such that $T_F$ and $T_G$ are not similar. This example shows that the model of  \cite{Clark-1980-2} is incorrect, even for symbols which do not ``back up''. In the case of arbitrarily wound symbols, Clark provided in \cite{Clark-Toeplitz-2}  an example of two rational symbols $F$ and $G$ such that the image of $\mathbb T$ under both $F$ and $G$ is the so-called ``figure-eight loop'' but $T_F$ and $T_G$ are not similar. Note that both $F$ and $G$ satisfy \ref{H2} but only $F$ satisfies \ref{H1}: the derivative of $G$ vanishes at two points of $\mathbb T$.

 Quasi-similarity of operators in this class was studied again by Clark in \cite{Clark-quasisim}. Yakubovich also obtained in \cite{Yakubovich1991} a triangular representation for positively wound Toeplitz operators with smooth symbols, allowing him to show, in the case $p=2$, that if $F$ and $ \widetilde F$ satisfy \ref{H1}, \ref{H2} and \ref{H3} and if there exists an orientation-preserving $C^1$ diffeomorphism $\tau$ of $\mathbb T$ such that $\widetilde F=F\circ \tau$ (which is stronger than merely assuming that $F(\mathbb T)= \widetilde F(\mathbb T)$), then there exists
a bounded linear isomorphism $L$ of $H^2$ and a finite rank operator $K$ on $H^2$ such that 
 \[T_{ \widetilde F}~=~LT_FL^{-1}+K.\]
 Unfortunately, we do not know if $K=0$.
 
 \par\smallskip
 We conclude this subsection with this following observation: we have investigated dynamical properties of $T_F$ acting on the spaces $H^p,p>1$: our results depend on $p$ only in a very mild way - namely, the regularity assumption \ref{H1} requires that $F$ be of class $C^{1+\varepsilon}$, with $\varepsilon>\max(1/p,1/q)$.
 
 \begin{question}\label{question4}
 Under suitable smoothness assumptions on $F$, is it true that for any $p,r>1$, the operator $T_F$ is hypercyclic on $H^p$ if and only if it is hypercyclic on $H^r$?
 \end{question}
 
 \begin{remark}\label{rem-p-r}
  Observe that if $p>r>1$, then $T_F$ is hypercyclic on $H^r$ as soon as $T_F$ is hypercyclic on $H^p$. Indeed, we have $H^p  \subseteq H^r$ and $\|\,.\,\|_r\le \|\,.\,\|_p$. If $f\in H^p$ is a hypercyclic function for $T_F$ seen as a bounded operator on $H^p$, then $f\in H^r$. Moreover, for any analytic polynomial $q$ and any $\varepsilon>0$, there exists an integer $n$ such that $\|T_F^n f -q\|_p<\varepsilon$. It follows that $\|T_F^n f -q\|_r<\varepsilon$. Since analytic polynomials are dense in $H^r$, it follows that $f$ is a hypercyclic function for $T_F$ acting on $H^r$. Our sufficient conditions for hypercyclicity thus carry over from $H^2$ to $H^p$, $1<p<2$, and our necessary conditions from $H^2$ to $H^p$, $p>2$.
 \end{remark}

 \subsection{Boundary values of quotients of inner functions}
\begin{figure}[ht]
\includegraphics[page=49,scale=1]{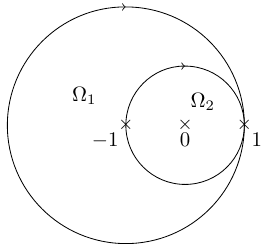}
            \caption{}\label{Fig46la}\end{figure}

One of the simplest case where \Cref{Question:2param} is open is when $F(\T)$ consists of two tangent circles. So let $F$ be a function satisfying \ref{H1} such that $F(\T)$ is represented as in \Cref{Fig46la}.

In particular, $F$ also satisfies \ref{H2} and \ref{H3}.
In \Cref{SubSection:Example2Circles}, we gave a parametrization $ \widetilde F$ of this curve for which $T_{ \widetilde F}$ became a hypercyclic Toeplitz operator. In relation with \Cref{Question:2param} above, it is natural to ask whether $T_F$ is hypercyclic for \emph{any} sufficiently smooth parametrization $F$ of the curve above. Let $\zeta=1/F^{-1}$. When trying to answer this question in this particular case, the natural thing to do is to study whether, given two functions $u,v\in E^q(\Omega_2)$, the boundary condition $u-\zeta v=0$ almost everywhere on $\mathbb T$ implies that $u=v=0$. In this situation, note that $\Omega_2=\mathbb D$ and so $E^q(\Omega_2)=H^q$. So let $u,v\in H^q$ be such that $u=\zeta v$ almost everywhere on $\mathbb T$. Since $|\zeta|=1$ almost everywhere on $\T$, one can suppose without loss of generality that the functions $u$ and $v$ are inner. Observe also that the function $ \widetilde \zeta:\theta\longmapsto\zeta(e^{i\theta})$ is $C^{1+\varepsilon}$-smooth and injective on the interval $[0,2\pi]$, with
$ \widetilde\zeta(0)\neq \widetilde\zeta(2\pi).$
\par\smallskip
This leads us to the following question concerning the boundary values of quotients of inner functions on the unit disk:

\begin{question}\label{fonctions-interieures}
Do there exist two inner functions $u, v$ on the unit disk, and a function $\xi:[0,2\pi]\to\mathbb T$ that is $C^{1+\varepsilon}$-smooth, injective (in particular, such that $\xi(0)\neq\xi(2\pi)$), such that $\xi'$ does not vanish on $[0,2\pi]$ and 
\[\xi(\theta)~=~\frac{u(e^{i\theta})}{v(e^{i\theta})}\quad\text{for almost every}~\theta\in(0,2\pi)?\]
\end{question}

Quotients of inner functions have been investigated in depth, starting from the work \cites{douglas-rudin} where it is shown that any unimodular function in $L^{\infty}(\T)$ can be uniformly approximated by quotients of inner functions. See for instance, among many others, the references \cites{bourgain, sarason, axler}.
Note that any answer to \Cref{fonctions-interieures} would be interesting:
if such functions do not exist, then, for every smooth parametrization $F$ of the curve above, the operator $T_F$ will satisfy the Godefroy-Shapiro Criterion, and thus be hypercyclic. This would also allow us, for example, to avoid the assumption regarding the existence of an analytic extension of $F$ to a neighborhood of $\Gamma$ in \Cref{Th:Complique-Jordan_version}.
In the case that the answer to \Cref{fonctions-interieures} turns out to be negative, it would be interesting to investigate whether the answer 
remains negative when we require the function $\xi$ to be 
only piecewise $C^{1+\varepsilon}$-smooth on $[0,2\pi]$. 
\par\smallskip
On the other hand, suppose that the answer to \Cref{fonctions-interieures} is affirmative, and that
such functions $u,v$ and $\xi$ do exist. Then, starting from $\xi$, we can define a parametrization $F$ of the curve above
by setting $F(\xi(\lambda))=\lambda$ for every $\lambda\in\partial\Omega_2=\mathbb T$ when the orientation of the curve $\xi(\T)$ is negative, and $F(1/\xi(\lambda))=\lambda$ for every $\lambda\in\partial\Omega_2=\mathbb T$ when the orientation of $\xi(\T)$ is positive. This defines $F$ on the subarc $\xi(\T)$ (resp. $1/\xi(\T)$) of $\T$. We then define $F$ on the whole of $\T$ in a $C^{1+\varepsilon}$-smooth way, so that $F(e^{i\theta})$ travels once over each part of the two tangent circles (except at the point $1$), and that the curve $F(\T)$ is negatively wound. 
In this way, we will obtain an operator $T_F$ that does not satisfy the Godefroy-Shapiro Criterion. Indeed, although $H_-(T_F)$ is easily seen to be equal to $H^p$, this is not the case for the subspace $H_+(T_F)$: supposing for instance that the orientation of $\xi$ is negative, define $u_0\in E^q(\Omega_1\cup\Omega_2)$ by setting $u_0=0$ on $\Omega_1$ and $u_0=u$ on $\Omega_2$, and $u_1\in E^q(\Omega_2)$ by setting $u_1=v$ on $\Omega_2$. The pair $(u_0,u_1)$ belongs to the range of the operator $U$, and there exists $g\in H^q$ such that $Ug=(u_0,u_1)$. The function $g$ is non-zero because $u_0\neq 0$ (recall that $u$ is inner).
Now, there are two possible situations: if $T_F$ if hypercyclic, then $T_F$ provides an example of a Toeplitz operator which is hypercyclic but does not satisfy the Godefroy-Shapiro Criterion, thus answering \Cref{question-GS}. If $T_F$ is not hypercyclic, then this means that the answer of \Cref{Question:2param} is negative.

\appendix
\section{Reminders}\label{Section:Rappels}
The aim of this first appendix is to recall some definitions and results which are required at various stages in the paper, either for the proofs of our main theorems or for the detailed proof of Yakubovich's \Cref{T:Yakbovich_Hp} which we provide in \Cref{Section:Yakubovich_demoHp}. We begin with some reminders in complex analysis, then present some results due to Privalov concerning the regularity of the Riesz projection on spaces of smooth functions, as well as a few facts regarding the spectral theory of Toeplitz operators. A brief presentation of Smirnov spaces on open sets which are  a finite union of simply connected components follows. Lastly, we give a brief overview of some notions in linear dynamics that appear in our work.

\subsection{Reminders in complex analysis}

\subsubsection{A consequence of Rouché's Theorem}
We first state a direct consequence of the Rouché Theorem, which is used in the proof of
\Cref{lem:description-eigenvectors}. It is certainly well-known, but we provide a short proof for completeness' sake.

\begin{lemma}\label{lem:rouche}
Let $W$ be an open subset of $\mathbb C$, let $g:W\longrightarrow \mathbb C$ be an analytic map and let $\lambda_0\in g(W)$. Assume that the equation $g(z)=\lambda_0$ has $s$ solutions $w_1,w_2,\dots,w_s$ in $W$, which are simple. Then there exist $\alpha>0$ and one-to-one analytic maps $\widetilde{d_1},\widetilde{d_2},\dots,\widetilde{d_s}$ defined on the disk $D(\lambda_0,\alpha)$, such that for every $\lambda\in D(\lambda_0,\alpha)$, the equation $g(z)=\lambda$ has exactly $s$ solutions in $W$ which are the points $\widetilde{d}_j(\lambda)$, $1\leq j\leq s$.
\end{lemma}
\begin{proof}
For every $1\leq j\leq s$, let $W_j$ be a small disk centered at $w_j$  such that $\overline{W_j} \subseteq W$ and for every $1\leq j,k\leq s$ with $j\neq k$,  $\overline{W_j}\cap \overline{W_k}=\varnothing$. Let $\alpha_j:=\min_{z\in\partial W_j}|g(z)-\lambda_0|>0$, and $\alpha:=\min_{1\leq j\leq s}\alpha_j$. Making $\alpha$ smaller if necessary, we may assume (using Rouch\'e's Theorem) that for every $\lambda\in D(\lambda_0,\alpha)$, the  set $\{z\in W\,;\,g(z)=\lambda\}$ has exactly $s$ elements.
Now, for every $\lambda\in D(\lambda_0,\alpha)$, every $1\leq j\leq s$ and every $z\in\partial W_j$, we have 
\[
|\lambda-\lambda_0|~<~\alpha_j~\leq~ |g(z)-\lambda_0|.
\]
Thus Rouch\'e's Theorem implies that the equation $g(z)=\lambda$ has exactly one solution  in $W_j$, which we denote by $\widetilde{d_j}(\lambda)$. Moreover, it follows from the Residue Theorem that
\[
\widetilde{d_j}(\lambda)~=~\frac{1}{2i\pi}\int_{\partial W_j}\frac{zg'(z)}{g(z)-\lambda}\,\mathrm{d}z,
\]
so that the function $\widetilde{d_j}$ is analytic on $D(\lambda_0,\alpha)$.
\end{proof}

\subsubsection{Lindelöf's Theorem}\label{Lindelof}
We recall in this section the following classical result due to Lindel\"of (see for instance \cite[page 89]{Garnett2007} or \cite{Lindelof}): 

\begin{theorem}\label{th-Lindelof}
Let $\eta>0$, and $p>1$. Let $f\in H^p(\mathbb D\cap \{|z-1|<\eta\})$ (that is, $|f|^p$ has a harmonic majorant on $\mathbb D\cap\{|z-1|<\eta\}$), and suppose that 
\[
\lim_{\theta\to 0^+}f(e^{i\theta})~=~\alpha\quad\text{and}\quad \lim_{\theta\to 0^-}f(e^{i\theta})~=~\beta
\]
do exists. Then $\alpha=\beta$ and
\[
\lim_{\underset{z\to 1}{z\in\D}}f(z)~=~\alpha.
\]
\end{theorem}

This result applies in particular when $f\in H^{\infty}(\D)$.
We use several times in the proofs of our main results the following direct consequence of Lindelöf's Theorem:

\begin{theorem}\label{Lindelof-Jordan}
Let $\Omega$ be a bounded Jordan domain and let $\Gamma=\partial\Omega$. Let $\lambda_0\in\Gamma$ and $f\in H^\infty(\Omega)$. Suppose that there exists a neighborhood $V$ of $\lambda_0$ such that $f$ has a continuous extension to $\overline{V\cap\Omega}\setminus\{\lambda_0\}$. 
Denote by $\Gamma_0^+$ and $\Gamma_0^-$ two disjoint sub-arcs of the boundary $\Gamma$ having $\lambda_0$ as an extremity.
\par\smallskip
If the  limit of $f(\lambda)$ as $\lambda\longrightarrow\lambda_0$, $\lambda\in \Gamma_0^+$ exists as well as the  limit of $f(\lambda)$ as $\lambda\longrightarrow\lambda_0$, $\lambda\in \Gamma_0^-$, then these two limits coincide. Denoting by $\alpha$ their common value, we have
\[\lim_{\underset{\lambda\in\overline{V\cap\Omega}\setminus\{\lambda_0\}}{\lambda\to\lambda_0}}f(\lambda)~=~\alpha.\]
\end{theorem}

\begin{proof}
Let $u:\mathbb D\to\Omega$ be a conformal map from $\D$ onto $\Omega$. Since $\Omega$ is a Jordan domain, the Carathéodory Theorem (see \cite[Th. 14.19]{Rudin1987}) implies that $u$ extends into a homeomorphism from $\overline{\mathbb D}$ onto $\overline{\Omega}$. \Cref{Lindelof-Jordan} is then a direct consequence of \Cref{th-Lindelof} applied to the function $f\circ u$.
\end{proof}

\subsubsection{Quasiconformal maps}\label{rappel:quasiconforme}
Quasiconformal functions appear naturally in the proof of \Cref{T:Yakbovich_Hp} because of the following result of Dynkin (see \cite[Th. 2]{Dynkin} and \cite[Sec. 1.3]{Dynkin2}):

\begin{theorem}\label{Dynkin}
Let $0<\varepsilon<1$, and let $F\in  C^{1+\varepsilon}(\T)$.
There exists a function $\widetilde{F}\in C^1(\C)$ such that the restriction of $\widetilde{F}$ to $\T$ is equal to $F$,
and
\begin{equation}\label{Eq:theoA4}
|\partial_{ \overline z}\widetilde{F}(z)|~\le~ C\|F\|_{ C^{1+\varepsilon}(\T)}\, \dist(z,\T)^\varepsilon\quad\text{ for every } z\in\C,
\end{equation}
where $C$ is a universal constant.
\end{theorem}
See \Cref{Section-Riesz} for the definition of 
$  C^{1+\varepsilon}(\T)$. Such a function $\widetilde{F}$ which extends $F$ and satisfies the estimate of \Cref{Dynkin} is called a pseudoanalytic
extension of $F$ to $\C$.
\par\smallskip
We will apply \Cref{Dynkin} to functions $F$ such that $F'\neq 0$ on $\mathbb T$.
According to \Cref{Eq:theoA4}, $\partial_{ \overline z}\widetilde{F}(z)=0$ for every $z\in\mathbb T$ . Moreover, since $F'$ does not vanish on $\mathbb T$,  $\partial_z\widetilde{F}$ does not vanish on $\mathbb T$. Since it is continuous on $\mathbb C$,
there exists an open neighborhood $U$ of $\mathbb T$ such that 
\[
C\|F\|_{ C^{1+\varepsilon}(\T)}\, \dist(z,\T)^\varepsilon ~\leq~ \frac{1}{2}|\partial_z \widetilde{F}(z)|
\]
for every $z\in U$. In particular, 
\begin{equation}\label{eq:conformal-map-234DSD}
|\partial_{ \overline z}\widetilde{F}(z)|~\leq~ \frac{1}{2}|\partial_z \widetilde{F}(z)|\quad \text{for every }z\in U.
\end{equation}
Recall that the Jacobian satisfies
\begin{equation}
\label{eq:jacobian-derivee-z-zbar}
J_{\widetilde{F}}(z)~=~|\partial_z\widetilde{F}(z)|^2-|\partial_{ \overline z}\widetilde{F}(z)|^2.
\end{equation}
 Then $J_{\widetilde{F}}(z)>0$ for every $z\in U$. So $\widetilde{F}$ is a local $C^1$-diffeomorphism at every point of $U$, and there exists, for every $z\in U$, an open neighborhood $U_z$ of $z$ and two positive constants $c_{1,z}$ and $c_{2,z}$ such that 
\begin{equation}\label{eq:conforme-diffeomorphism}
c_{1,z}|w-z|~\leq~ |\widetilde{F}(w)-\widetilde{F}(z)|~\leq~ c_{2,z}|w-z|   \quad\text{for every }w\in U_z.
\end{equation}
Moreover $\widetilde{F}:U_z\longmapsto \widetilde{F}(U_z)$ preserves the orientation since $J_{\widetilde{F}}(z)>0$. 
The inequality 
\begin{equation}\label{eq:sddqsdds233434f0}
|\partial_{ \overline z}g|~\leq~ k |\partial_z g|
\end{equation}
where $k$ is a constant with $0<k<1$, lies at the core of the theory of quasiconformal mappings. Indeed, if $g$ is a $C^1$-diffeomorphism between two domains $G$ and $G'$ of $\mathbb C$, the map $g$ is quasiconformal if and only if it satisfies \Cref{eq:sddqsdds233434f0} for some $0<k<1$. 

The proper setting for the definition of quasiconformal mappings is that of orientation-preserving homeomorphisms. As we will briefly need it below, we recall one of the many equivalent definitions (see for instance \cite[Ch. II, Def. B]{Ahlfors2006}). Let 
$g:G\longmapsto G'$ be an homeomorphism. Then $g$ is said to be quasiconformal if $g$ satisfies this two conditions:
\begin{enumerate}[(a)]
    \item $g$ is Absolutely Continuous on Lines (ACL) on $G$;
    \item there exists $k\in (0,1)$ such that 
    \[
    |\partial_{ \overline z}g|~\leq~ k |\partial_z g| \quad\text{a.e. on }G.
    \]
\end{enumerate}
Absolute continuity on lines means that for any rectangle $R=\{x+iy\,;\,a<x<b,\,c<y<d\}$ with $\overline{R} \subseteq G$, the function $x\longmapsto g(x+iy)$ is absolutely continuous (i.e. has bounded variation) on $(a,b)$ for almost every $y\in (c,d)$, and the function $y\longmapsto g(x+iy)$ is absolutely continuous on $(c,d)$ for almost every $x\in (a,b)$. Whenever $g$ is ACL on $G$, the partial derivatives $\partial_x g$ and $\partial_y g$ exist almost everywhere on $G$, so that (b) makes sense \cite[Ch. III, Lemma 3.1]{LehtoVirtanen1973}. Since $g$ is an homeomorphism, it follows that $g$ is differentiable almost everywhere on $G$ \cite[Ch. III, Th. 3.1]{LehtoVirtanen1973}. If (a) and (b) hold, it is even true that $g$ is regular at almost every point $z$ of $G$, i.e. $g$ is differentiable at $z$ and its Jacobian $J_g(z)$ does not vanish. See \cite[Ch. IV, Th. 1.4]{LehtoVirtanen1973} plus the explanation in \cite[Ch. IV, Sec. 5.3]{LehtoVirtanen1973} at the bottom of page 184. Since $|\partial_{ \overline z}g|\leq k |\partial_z g|$ almost everywhere on $G$ with $0<k<1$, we have 
\[
J_g(z)~=~|\partial_z g(z)|^2-|\partial_{ \overline z}g(z)|^2~\geq~ (1-k^2)|\partial_z g(z)|^2\quad\text{a.e. on }G,
\]
so that $J_g(z)\geq 0$ almost everywhere on $G$. Hence $J_g(z)>0$ almost everywhere on $G$, and $g$ is an orientation-preserving homeomorphism (see \cite[Ch. I, Sec. 1.4-1.6]{LehtoVirtanen1973} for details).  
If $g:G\longmapsto G'$ is an homeomorphism which is besides of class $ C^1$, then $g$ is clearly ACL, and thus $g$ is quasiconformal if and only if (b) holds. Let us also mention here that if $g:G\longmapsto G'$ is quasiconformal, then $g^{-1}:G'\longmapsto G$ is also quasiconformal, and that the composition of two quasiconformal mappings is quasiconformal. 

Getting back to our pseudo-analytic extension $\widetilde{F}$ of $F$ to $U$, we see that it satisfies (a) and (b). It is tempting to deduce that it is quasiconformal, but since $\widetilde{F}$ has no reason to be an homeomorphism from $U$ onto its image, it is only locally true: for every $z\in U$, the restriction of $\widetilde{F}$ to $U_z$ is a quasiconformal mapping, where $U_z$ is a neighborhood of $z$ such that $\widetilde{F}$ satisfies \Cref{eq:conforme-diffeomorphism}.

\begin{figure}[ht]
\includegraphics[page=50,scale=1.1]{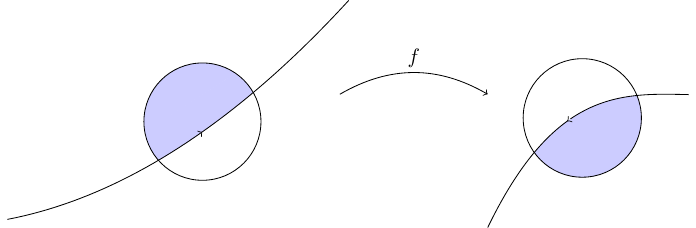}\vspace*{-1cm}
            \caption{}\label{Fig47kjl}\end{figure}

For reasons which will become clear later on in \Cref{Section:Yakubovich_demoHp}, we will need to capture the global behavior of $\widetilde{F}$ on $U$, and for this we will need the notion of quasiconformal function. A function $f:G\longmapsto \mathbb C$ is said to be a quasiconformal function if it can be written as $f=\varphi\circ g$, where $g$ is a quasiconformal mapping (in particular an homeomorphism) from $G$ onto a domain $G'$ of $\mathbb C$ and $\varphi:G'\longmapsto \mathbb C$ is a non-constant analytic function. We refer to \cite{Ahlfors2006} and \cite{LehtoVirtanen1973} for a full study on quasiconformal maps and to \cite[Ch. VI]{LehtoVirtanen1973} for an introduction of quasiconformal functions. 
It easily follows from the definition that locally, quasiconformal functions preserve the orientation, as illustrated in \Cref{Fig47kjl}. It is a particularly useful property which will be used several times in \Cref{Section:Yakubovich_demoHp}.

Now here is how one identifies quasiconformal functions \cite[Ch. VI, Th. 2.2.]{LehtoVirtanen1973}. Let $f:G\longmapsto \mathbb C$ be a non-constant function which is a generalized $L^2$-solution on $G$ of an equation  of the form
\begin{equation}\label{eq:quasiconformal-mapping-EDP}
\partial_{ \overline z}f~=~\chi\,\partial_z f,    
\end{equation}
where $\chi$ is a measurable function on $G$ with $\sup_{z\in G}|\chi(z)|<1$. Then $f$ is a quasiconformal function on $G$. That $f$ is a generalized $L^2$-solution of \Cref{eq:quasiconformal-mapping-EDP} means that $f$ is ACL on $G$ (so that $\partial_x f$ and $\partial_y f$ exist almost everywhere on $G$) and that $|\partial_x f|^2$ and $|\partial_y f|^2$ are both integrable on any compact subset of $G$ with respect to the area measure (and, of course, that \Cref{eq:quasiconformal-mapping-EDP} is satisfied). Under these assumptions, $|\partial_{ \overline z} f|\leq k |\partial_z f|$ almost everywhere on $G$, where $k=\sup_{z\in G}|\chi(z)|<1$, i.e. \Cref{eq:sddqsdds233434f0} is true. 
\par\smallskip
Suppose now that $f:G\longmapsto \mathbb C$ is a non-constant $C^1$-smooth function satisfying \Cref{eq:sddqsdds233434f0}: then $f$ is obviously ACL with partial derivatives which are measurable and square integrable on any compact subset of $G$. Suppose moreover that $J_f(z)\neq 0$ almost everywhere on $G$. Then \Cref{eq:quasiconformal-mapping-EDP} is obviously satisfied with 
\[
\chi(z)~=~\frac{|\partial_{ \overline z}f(z)|}{|\partial_z f(z)|}\quad \text{a.e. on }G.
\]
Thus one deduces from \cite[Ch. IV, Th. 2.2]{LehtoVirtanen1973} the following statement: if $f:G\longmapsto \mathbb C$ is a non-constant $ C^1$-smooth function such that $J_f(z)\neq 0$ almost everywhere on $G$ and if $|\partial_{ \overline z}f|\leq k |\partial_z f|$ almost everywhere on $G$, for some $k\in (0,1)$, then $f$ is a quasiconformal function on $G$. 
\par\smallskip
Getting back to our pseudo-analytic extension $\widetilde{F}$ of $F$ to $U$, we thus see, using \Cref{eq:conformal-map-234DSD,eq:jacobian-derivee-z-zbar}, that $\widetilde{F}$ is a quasiconformal function on $U$. So there exists a quasiconformal mapping $g$ from $U$ onto a domain $V$ of $\mathbb C$ and a non-constant analytic function $\varphi:V\longmapsto\mathbb C$ such that $\widetilde{F}=\varphi\circ g$ on $U$. Moreover, we have $J_{\widetilde{F}}(z)> 0$ for every $z\in U$, and $\widetilde{F}:U_z\longmapsto \widetilde{F}(U_z)$ is a $ C^1$-diffeomorphism, so that, in particular, $\widetilde{F}$ is injective on $U_z$. It follows that $\varphi$ is injective on the domain $g(U_z)$, and $\varphi:g(U_z)\longmapsto \varphi(g(U_z))=\widetilde{F}(U_z)$ is an analytic isomorphism. Hence, we can write $g$ as $g=\varphi^{-1}\circ \widetilde{F}$ on $U_z$, from which it follows that $g$ is a $ C^1$-diffeomorphism from $U_z$ onto $g(U_z)$. Since $g:U\longrightarrow V$ is already known to be an homeomorphism, we eventually obtain that $g$ is a $ C^1$-diffeomorphism from $U$ onto $V$. We summarize our discussion in the following theorem, which we will use as such in \Cref{Section:Yakubovich_demoHp}.
\begin{theorem}\label{thm:quasi-conformal-function}
Let $0<\varepsilon<1$, and let $F\in C^{1+\varepsilon}(\mathbb T)$ be such that $F'\neq 0$ on $\mathbb T$. There exists an open neighborhood $U$ of $\mathbb T$ and a function $\widetilde{F}\in  C^1(U)$ such that:
\begin{enumerate}[(a)]
    \item the restriction of $\widetilde{F}$ to $\mathbb T$ is equal to $F$;\par\smallskip
    \item the function $\widetilde{F}$ can be written as $\widetilde{F}=\varphi\circ g$, where $g$ is a quasiconformal $ C^1$-diffeomorphism from $U$ onto a domain $V$ of $\mathbb C$, and $\varphi$ is a non-constant analytic function on $V$. In particular, both $g$ and $g^{-1}$ are orientation preserving;\par\smallskip
    \item the Jacobian $J_{\widetilde{F}}(z)$ is positive at every point $z\in U$; hence $\widetilde{F}$ is an orientation preserving local $ C^1$-diffeomorphism at every point of $U$, and there exist, for every $z\in U$, an open neighborhood $U_z$ of $z$ and two positive constants $c_{1,z}$ and $c_{2,z}$ such that 
\begin{equation*}
c_{1,z}|w-w'|~\leq ~|\widetilde{F}(w)-\widetilde{F}(w')|~\leq~ c_{2,z}|w-w'|  \quad \text{for every } w,w'\in U_z. 
\end{equation*}
\end{enumerate}
\end{theorem}

\subsection{Properties of the Riesz projection $P_+$}\label{Section-Riesz}
Given $1<p<\infty$, let  $q$ be its conjugate exponent, i.e. $\frac{1}{p}+\frac{1}{q}=1$. Then the duality between $L^p(\mathbb T)$ and $L^q(\mathbb T)$ is given by the following duality bracket 
\begin{equation}\label{EqA7}
\dual{x}{y}~=~\frac{1}{2\pi}\int_0^{2\pi}x(e^{i\theta})y(e^{-i\theta})\,\mathrm{d}\theta,
~\end{equation}
where $x\in L^p(\mathbb T)$ and $y\in L^q(\mathbb T)$. The Hardy space $H^p$ is defined by 
\[
H^p\,=\,\left\{u\in L^p(\mathbb T)\,;\,\forall n<0,\;\hat{u}(n)=0\right\}\]
and
\[H^p_-\,=\,\left\{u\in L^p(\mathbb T)\,;\,\forall n\ge0,\;\hat{u}(n)=0\right\}.
\]
As usual, we may identify $H^p$ with the space of analytic functions $u$ on the open unit disc $\mathbb D$ such that 
\[\sup_{0\leq r<1}\|u_r\|_{L^p(\mathbb T)}~<~\infty,\]
where $u_r(z)=u(rz)$, $z\in\mathbb T$.
We also may identify the dual of $H^p$ with $H^q$ with respect to the duality bracket \eqref{EqA7}. 

For every $z\in\mathbb D$, let $k_z$ be the Cauchy kernel defined by $k_z(e^{i\theta})=(1-\overline{z}e^{i\theta})^{-1}$, $e^{i\theta}\in\T$.
We have
\[
\dual{u}{k_z}~=~\frac{1}{2\pi}\int_0^{2\pi}\frac{u(e^{i\theta})}{1-\overline{z}e^{-i\theta}}\,\mathrm{d}\theta~=~\frac{1}{2\pi}\int_0^{2\pi}\frac{e^{i\theta}u(e^{i\theta})}{e^{i\theta}-\overline{z}}\,\mathrm{d}\theta \quad\text{ for every } u\in H^p.
\]
Hence by Cauchy's formula for $H^p$ functions, we have 
\begin{equation}\label{eq:cauchy-kernel}
\dual{u}{k_z}~=~u(\overline{z})\quad \text{ for all } z\in\mathbb D\text{ and all }u\in H^p.
\end{equation}
The Riesz Theorem states that the Riesz projection $P_+:L^p(\mathbb T)\longrightarrow H^p$ defined by 
\[
P_+f(z)~=~\frac{1}{2i\pi}\int_{\mathbb T}\frac{f(\tau)}{\tau-z}\,\mathrm{d}\tau\quad\textrm{ for every } z\in\mathbb D
\]
is bounded whenever $1<p<\infty$.  Note that for every $\varphi\in L^p(\mathbb T)$ and every $v\in H^q$, we have 
\begin{equation}\label{eq:simple-P+}
\dual{\varphi}{v}~=~\dual{P_+\varphi}{v}.
\end{equation}
The Toeplitz operators that we consider are, as a general rule, associated to a $C^{1+\varepsilon}$-smooth symbol. So let us begin by recalling what this means, and by giving some useful properties related to this space of functions.
\par\smallskip
Let $0<\varepsilon<1$ and $n\ge 1$. Let $\Omega$ be a bounded subset of $\mathbb C^n$, endowed with the Euclidean norm $\|\,.\,\|_{\mathbb C^n}$, and let $Y$ be a Banach space. A function $h:\Omega\longmapsto Y$ is said to be \emph{of class $C^{\varepsilon}$}, or \emph{$C^{\varepsilon}$-smooth}, or simply  \emph{is a $C^{\varepsilon}$ function}, if there exists a constant $C>0$ such that 
\[
\|h(z_1)-h(z_2)\|_Y~\leq~ C \|z_1-z_2\|^\varepsilon_{\mathbb C^n} \quad \text{ for every } z_1,z_2\in \Omega.
\]
We denote by $C^\varepsilon(\Omega)$ the set of all such functions.
~When equipped with the norm
\begin{equation}\label{eq:NormeCeps}
\|h\|_{C^\varepsilon(\Omega)}~=~\|h\|_\infty+\sup_{z_1\neq z_2}\frac{\|h(z_1)-h(z_2)\|_Y}{\|z_1-z_2\|^\varepsilon_{\mathbb C^n}},    
\end{equation}
the space $C^\varepsilon(\Omega)$ becomes a Banach space. Given an integer $k\ge 1$, we say that a function $h: \Omega\longrightarrow\C$ is of class  $C^{k+\varepsilon}$ on $\Omega$ if $u^{(k)}$ is of class $C^{\varepsilon}$ on $\Omega$. Note that we have the following result:

\begin{lemma}\label{lemme-Cepsilon}
Let $0<\varepsilon<1$, and let $\Omega_1$ and $\Omega_2$ be bounded subsets of $\mathbb C^n$. Let the function $\phi:\Omega_1\times \Omega_2\longmapsto \mathbb C$ belong to the space $C^\varepsilon(\Omega_1\times\Omega_2)$. Write $\varepsilon$ as $\varepsilon=\gamma+\beta$, where $\gamma,\beta>0$. For every $a\in \Omega_1$ and $b\in \Omega_2$, denote by $\Phi_l(a)$ and $\Phi_r(b)$ the functions 
$\Phi_l(a)=\phi(a,\cdot)$ and $\Phi_r(b)=\phi(\cdot,b)$. 
\par\smallskip
Then 
$\Phi_l$ is a function of class $C^{\gamma}$ from $\Omega_1$ into $C^{\beta}(\Omega_2)$, and $\Phi_r$ is a function of class $C^{\beta}$ from $\Omega_2$ into $C^{\gamma}(\Omega_1)$.    
\end{lemma}
If $a$ and $b$ are two functions on a subset $\Lambda$ of $\C^d$ taking non-negative values, we write $a\lesssim b$ if there exists a positive constant $C$ such that $0\le a(w)\le Cb(w)$ for every $w\in\Lambda$.\par\smallskip
\begin{proof}
For $z_1,z_2\in \Omega_1$, we have
\begin{multline*}
\|\Phi_l(z_1)-\Phi_l(z_2)\|_{C^{\beta}(\Omega_2)}~=~
\|\Phi_l(z_1)-\Phi_l(z_2)\|_\infty \\
+\sup_{\substack{w_1,w_2\in\Omega_2\\
w_1\neq w_2}}\frac{|\phi(z_1,w_1)-\phi(z_2,w_1)-\phi(z_1,w_2)+\phi(z_2,w_2)|}{\|w_1-w_2\|_{\mathbb C^n}^{\beta}}\cdot    
\end{multline*}

Observe that, since $\phi$ is $C^\varepsilon$, we have 
\[
\|\Phi_l(z_1)-\Phi_l(z_2)\|_\infty~=~\sup_{w\in\Omega_2}|\phi(z_1,w)-\phi(z_2,w)|\leq C\|z_1-z_2\|_{\mathbb C^n}^{\varepsilon},
\]
and using the fact that $\Omega_1$ is bounded, we get 
\[
\|\Phi_l(z_1)-\Phi_l(z_2)\|_\infty~\lesssim~ \|z_1-z_2\|_{\mathbb C^n}^{\gamma}.
\]
Moreover, 
\begin{multline*}
  |\phi(z_1,w_1)-\phi(z_2,w_1)-\phi(z_1,w_2)+\phi(z_2,w_2)|\\
  \leq~|\phi(z_1,w_1)-\phi(z_2,w_1)|+|\phi(z_1,w_2)-\phi(z_2,w_2)|~\leq ~2C \|z_1-z_2\|_{\mathbb C^n}^\varepsilon.  
\end{multline*}
Similarly, we have 
\begin{multline*}
    |\phi(z_1,w_1)-\phi(z_2,w_1)-\phi(z_1,w_2)+\phi(z_2,w_2)|\\
    \leq~|\phi(z_1,w_1)-\phi(z_1,w_2)|+|\phi(z_2,w_1)-\phi(z_2,w_2)|
~\leq~ 2C \|w_1-w_2\|_{\mathbb C^n}^\varepsilon.
\end{multline*}
Therefore, 
\[
|\phi(z_1,w_1)-\phi(z_2,w_1)-\phi(z_1,w_2)+\phi(z_2,w_2)|~\leq ~2C \min(\|z_1-z_2\|_{\mathbb C^n}^\varepsilon,\|w_1-w_2\|_{\mathbb C^n}^\varepsilon).
\]
But $\min(a^\varepsilon,b^\varepsilon)\leq a^\gamma b^\beta$ for every $a,b\geq 0$, whence it follows that
\[
|\phi(z_1,w_1)-\phi(z_2,w_1)-\phi(z_1,w_2)+\phi(z_2,w_2)|~\leq~ 2C \|z_1-z_2\|_{\mathbb C^n}^{\gamma}\|w_1-w_2\|_{\mathbb C^n}^{\beta}.
\]
Thus we obtain that 
\[
\|\Phi_l(z_1)-\Phi_l(z_2)\|_{C^{\beta}(\Omega_2)}~\lesssim~ \|z_1-z_2\|_{\mathbb C^n}^{\gamma},
\]
which proves that $\Phi_l$ is a $C^{\gamma}$ function from $\Omega_1$ into $C^{\beta}(\Omega_2)$. The proof is similar for $\Phi_r$.
\end{proof}

It is a well-known result that $P_+$ is a well-defined and bounded operator on $C^\varepsilon(\mathbb T)$ and $C^{1+\varepsilon}(\mathbb T)$:\par\smallskip
\begin{theorem}[Privalov-Zygmund]\label{thm:privalov-Zygmund}
Let $0<\varepsilon<1$. Then 
\begin{enumerate}
\item $P_+$ is a bounded map from $C^\varepsilon(\mathbb T)$ into itself;
\item  $P_+$ is a bounded map from $C^{1+\varepsilon}(\mathbb T)$ into itself.
\end{enumerate}
\end{theorem}
Assertion (1) is due to Privalov. See for instance \cite[Th. 3.1.1]{CimaMathesonRoss2006}.
Assertion (2) is due to Zygmund. See \cite[Section 3.1, page 64]{CimaMathesonRoss2006} and \cite{Zygmund}. Note that (2) is a direct consequence of (1). Indeed, let $u\in C^{1+\varepsilon}(\mathbb T)$. Then $P_+u\in A(\mathbb D)$ by (1). Integrating by parts, we have that
\[(P_+u)'(z)=\frac1{2i\pi}\int_{\mathbb T}\frac{u(\tau)}{(\tau-z)^2}\mathrm d\tau=\frac1{2i\pi}\int_{\mathbb T}\frac{u'(\tau)}{\tau-z}\mathrm d\tau=P_+(u')(z)\quad\textrm{for every }z\in\D.\]
So $(P_+u)'=P_+(u')$ belongs to $C^\varepsilon(\mathbb T)$ by (1). Note also that if $f\in A(\mathbb D)$ is such that $u=f_{|\mathbb T}$ belongs to $C^{1+\varepsilon}(\mathbb T)$, then $f'=P_+(u')$ and 
its boundary limit is of class $C^\varepsilon$ on $\mathbb T$. In particular, $f'\in A(\mathbb D)$, and 
$$\lim_{\stackrel{z\to e^{i\theta}}{z\in\overline{\D}}}f'(z)=f'(e^{i\theta}).$$
See also \cite[Section 3.4]{Duren1970}.

\subsection{Carleson measures for $H^p$}\label{subsection-Carleson-measures}
Recall that a positive finite Borel measure $m$ on $\mathbb D$ is called a {\emph{Carleson measure}} if there is a constant $A>0$ such that \begin{equation}\label{eq:definition-carleson-measure}
m(S(\tau,r))~\leq~ A r\quad \text{for every}~ \tau\in\mathbb T \text{ and } 0<r<1,
\end{equation}
where
\[
S(\tau,r)~=~\{z\in\mathbb D:|z-\tau|<r\}. 
\]
The Carleson constant $C(m)$ of the measure $m$ is the infimum of the constants $A>0$ satisfying \Cref{eq:definition-carleson-measure}. A famous theorem of Carleson \cite[Th. 5.15]{MR3497010} states that if $m$ is a Carleson measure then
\begin{equation}\label{eq:Carleson-embedding2343}
\left(\int_{\mathbb D}|f(z)|^p\,\mathrm{d}m(z)\right)^{1/p}\leq~ (12e\, C(m))^{1/p}\,\|f\|_{H^p},
\end{equation}
for every $f\in H^p$ and every $1\leq p<\infty$. In other words, we can embed $H^p$ continuously into $L^p(m)$ if $m$ satisfies \Cref{eq:definition-carleson-measure}. Note that in \cite[Th. 5.15]{MR3497010}, the inequality \eqref{eq:Carleson-embedding2343} is given for $p=2$, but a standard argument using inner-outer factorization allows to extend it to other values of $1\leq p<\infty$. 
\par\smallskip
For a rectifiable curve $\gamma$ in $\mathbb C$, we also define its Carleson constant $Carl(\gamma)$ by
\[
Carl(\gamma)=\sup\left\{\frac{|\gamma\cap D(\mu,r)|}{r}:\mu\in\mathbb C;r>0\right\},
\]
where $|\gamma\cap D(\mu,r)|$ is the arc length of $\gamma\cap D(\mu,r)$. In particular, the Carleson constant of every circle (with positive radius) is $2\pi$.

In the proof of the continuity of the isomorphism involved in the construction of the model in \Cref{Section:Yakubovich_demoHp}, the following result will be needed. It is probably well-known, but we give a proof for completeness' sake. 

\begin{proposition}\label{PropA8}
Let $\Omega$ be a domain in $\mathbb C$, let $u:\Omega\to\mathbb D$ be a function of class $C^1$ on $\Omega$, and assume that there exists a constant $c>0$ such that
\begin{equation}\label{eq:PropA8-hypothesis}
\frac{1}{c}|\lambda-\mu|\leq |u(\lambda)-u(\mu)|\leq c|\lambda-\mu| \quad\text{for every }\lambda,\mu\in\Omega.
\end{equation}
Then, for every rectifiable curve $\gamma\subset\Omega$ of finite length and for every $g\in H^p$, we have $Carl(u(\gamma))\leq 2 c^2\, Carl(\gamma)$  and
\[
\int_{\gamma}|g(u(\lambda))|^p\,|\mathrm{d}\lambda|\leq 24e\,c^3\,Carl(\gamma)\,\|g\|_p^p.
\]
\end{proposition}
\par\smallskip
\begin{proof}
We may of course assume that $Carl(\gamma)<\infty$.  It follows from \Cref{eq:PropA8-hypothesis} that $u$ is a $C^1$-diffeomorphism from $\Omega$ onto its image $V\subset \mathbb D$, and its Jacobian $J_u$ satisfies
\begin{equation}\label{eq:propA8-Jacobien3E}
\frac{1}{c}\leq |J_u(\lambda)|\leq c\quad\text{for every }\lambda\in\Omega.
\end{equation}
Let $\mu\in\mathbb C$, $r>0$ such that $u(\gamma)\cap D(\mu,r)\neq\varnothing$, and let $a\in\gamma$ be such that $u(a)\in D(\mu,r)$. Observe that 
\begin{equation}\label{eq:23Zs23d3ZpropA8Carleson}
u(\gamma)\cap D(\mu,r) \subset u(\gamma\cap D(a,2cr)).
\end{equation}
Indeed, if $v\in u(\gamma)\cap D(\mu,r)$, then there exists $\lambda\in\gamma$ such that $v=u(\lambda)$ and $|v-\mu|<r$. Hence, using \Cref{eq:PropA8-hypothesis}, we have 
\[
|\lambda-a|\leq c|u(\lambda)-u(a)|\leq c(|v-\mu|+|\mu-u(a)|)\leq 2cr,
\]
which gives that $v=u(\lambda)\in u(\gamma\cap D(a,2cr))$. Now \Cref{eq:23Zs23d3ZpropA8Carleson} implies that 
\begin{equation}\label{eq:23Zs9HSD23d3ZpropA8Carleson}
|u(\gamma)\cap D(\mu,r)|\leq |u(\gamma\cap D(a,2cr))|.
\end{equation}
But, using \Cref{eq:propA8-Jacobien3E}, we have
\[
|u(\gamma\cap D(a,2cr))|~=~\int_{u(\gamma\cap D(a,2cr))}\,|\mathrm{d}z|=\int_{\gamma\cap D(a,2cr)}|J_u(\lambda)|\,|\mathrm{d}\lambda|~\leq~ c\, |\gamma\cap D(a,2cr)|.
\]
Thus, according to \Cref{eq:23Zs9HSD23d3ZpropA8Carleson} and the definition of Carleson constant of the curve $\gamma$, we finally get that
\[
|u(\gamma)\cap D(\mu,r)|~\leq~ 2c^2\,Carl(\gamma)\,r,
\]
and thus $Carl(u(\gamma))\leq 2 c^2\, Carl(\gamma)$. 

Now, let $g\in H^p$. Then, using the change of variable $z=u(\lambda)$, we have 
\begin{equation}\label{eq:23Z3ZpropA8Carleson}
    \int_\gamma |g(u(\lambda))|^p\,|\mathrm{d}\lambda|~=~\int_{u(\gamma)}|g(z)|^p |J_{u^{-1}}(z)|\,|\mathrm{d}z|~\leq~ c\,\int_{u(\gamma)}|g(z)|^p\,|\mathrm{d}z|.
\end{equation}
It remains to estimate the integral on $u(\gamma)\subset \mathbb D$. For this purpose, consider the Borel measure $\mu_\gamma$ on $\mathbb D$ defined by
\[\mu_\gamma(A)~=~|u(\gamma)\cap A| \quad\text{for every Borel subset $A$ of $\mathbb D$,}\] 
i.e. $\mu_\gamma(A)$ is the arc length measure of the Borel subset $A\cap u(\gamma)$. Then, for every Borel set $A\subset\mathbb D$, we have 
\[
\int_{\mathbb D} \mathbb{1}_A(w)\,\mathrm{d}\mu_\gamma(w)~=~\int_{u(\gamma)} \mathbb{1}_{A}(z)\,|\mathrm{d}z|.
\]
Then, by construction of the integral, for every positive measurable function $h$ on $\mathbb D$ with respect to $\mu_\gamma$, we get
\[
\int_{\mathbb D}h(w)\,\mathrm{d}\mu_\gamma(w)~=~\int_{u(\gamma)}h(z)\,|\mathrm{d}z|,
\]
which gives
\begin{equation}\label{eq:23Zsd3ZpropA8Carleson}
\int_{u(\gamma)}|g(z)|^p\,|\mathrm{d}z|~=~\int_{\mathbb D}|g(w)|^p\,d\mu_{\gamma}(w).
\end{equation}
But observe that for every $\zeta\in\mathbb T$ and every $0<r<1$, we have 
\[
\mu_\gamma(S(\zeta,r))=|u(\gamma)\cap S(\zeta,r)|\leq |u(\gamma)\cap D(\zeta,r)|\leq Carl(u(\gamma)) r\leq 2c^2 Carl(\gamma) r.
\]
Hence $\mu_\gamma$ is a Carleson measure, whose constant satisfies $C(\mu_\gamma)\leq 2c^2 Carl(\gamma)$. It follows now from \Cref{eq:Carleson-embedding2343} that 
\[
\int_{\mathbb D}|g(z)|^p\,\mathrm{d}\mu_\gamma(z)\leq 24e\,c^2\,Carl(\gamma)\,\|g\|_p^p,
\]
and \Cref{eq:23Z3ZpropA8Carleson,eq:23Zsd3ZpropA8Carleson} imply that 
\[
\int_{\gamma}|g(u(\lambda))|^p\,|\mathrm{d}\lambda|\leq 24e\,c^3\,Carl(\gamma)\,\|g\|_p^p,
\]
which ends the proof.
\end{proof}

\subsection{Toeplitz operators}\label{rappels-Toeplitz}
Given $F\in L^\infty(\mathbb T)$, the Toeplitz operator $T_F$ with symbol $F$ is defined on the Hardy space $H^p$, $p>1$, by the formula
\[
T_F(u)~=~P_+(Fu)\quad\textrm{ for every } u\in H^p.
\]
It is a well-known fact that $T_F$ is bounded from $H^p$ into itself with 
\[
\|F\|_\infty~\leq~ \|T_F\|_{\mathcal{L}(H^p)}~\leq~ c_p\|F\|_\infty,
\]
where $c_p$ is the norm of the Riesz projection on $L^p(\mathbb T)$. See 
\cite[Th. 2.7]{BottcherSilbermann1990}.
It is easy to compute the adjoint of the operator $T_F$ with respect to the duality  bracket \eqref{EqA7}:

\begin{lemma}\label{lem:adjoint-banachique}
For any $F\in L^\infty(\mathbb T)$, define $f\in L^\infty(\mathbb T)$ by setting $f(z)=F(1/z)$ for almost every $z\in\mathbb T$. Let $p>1$. Then the adjoint of the operator $T_F$ acting on $H^p$ is the operator $T_f$ acting on $H^q$. 
\end{lemma}

\begin{proof}
For every $u\in H^p$ and $v\in H^q$, using \Cref{eq:simple-P+}, we have 
\begin{eqnarray*}
\dual{u}{T_F^*v}~=~\dual{Fu}{v}&=&\frac{1}{2\pi}\int_{0}^{2\pi}F(e^{i\theta})u(e^{i\theta})v(e^{-i\theta})\,\mathrm{d}\theta\\
&=&\frac{1}{2\pi}\int_{0}^{2\pi}f(e^{-i\theta})v(e^{-i\theta})u(e^{i\theta})\,\mathrm{d}\theta\\
&=&\dual{u}{vf}~=~\dual{u}{T_f(v)},
\end{eqnarray*}
which yields that $T_F^*(v)=T_f(v)$. 
\end{proof}

The next lemma gathers some well-known results concerning spectral properties of Toeplitz operators with continuous symbols. 
\begin{lemma}\label{Lemma:spectral-properties-toeplitz-operators}
Let $F$ be a continuous function on $\T$. Consider the associated Toeplitz operator $T_F:H^p\longrightarrow H^p$, $1<p<\infty$. \par\smallskip
\begin{enumerate}[(i)]
\item  $T_F$ is a Fredholm operator if and only if $F$ does not vanish on $\mathbb T$, and in this case, its Fredholm index $j(T_F)$ satisfies 
\[
j(T_F)~=~\dim(\ker T_F)-\dim(\ker T_F^*)~=~-\w_F(0).
\]
\item We have 
\[
\sigma(T_F)~=~F(\mathbb{T})\cup \{\lambda\in\mathbb C\setminus F(\mathbb T)\,;\,\w_F(\lambda)\neq 0\}.
\]
\end{enumerate}
\end{lemma}
For a proof of \Cref{Lemma:spectral-properties-toeplitz-operators}, see for instance \cite[Th. 2.42]{BottcherSilbermann1990}.
It relies on a classical theorem due to Coburn (see \cite[Th. 2.38]{BottcherSilbermann1990}) which states the following: 

\begin{theorem}\label{Th:Coburn}
If $F\in L^\infty(\mathbb T)$ is not almost everywhere zero, then either $\text{ker }T_F=\{0\}$ or $\text{ker }T_F^*=\{0\}$.
\end{theorem}

It should be noted that  the duality pairing between $H^p$ and $H^q$ used in in \cite{BottcherSilbermann1990} is given by the duality bracket 
\begin{equation}\label{eq2:duality-pairing-sdsqdqs}
\dual{f}{g}~=~\frac{1}{2\pi}\int_0^{2\pi} f(e^{i\theta})\overline{g(e^{i\theta})}\,\mathrm{d}\theta,
\end{equation}
which is different from the one we use here.
With respect to the duality bracket (\ref{eq2:duality-pairing-sdsqdqs}), the adjoint of $T_F$ is $T_F^*=T_{\overline{F}}$. But since $\w_{\overline{F}}(0)=\w_f(0)=-\w_F(0)$, Theorems 2.42 and 2.38 from \cite{BottcherSilbermann1990} remain valid when considering the pairing \eqref{eq:duality-defn} rather than the pairing \eqref{eq2:duality-pairing-sdsqdqs}.

\subsection{Smirnov spaces $E^p(\Omega)$}\label{smirnov}
We briefly discuss in this section some standard properties of Smirnov spaces that we use in our paper. We refer to \cites{Duren1970,Khavinson1986,Priwalow1956,TumarkinHavinson1960}  and the references therein for further details.
\par\smallskip
Let $\Omega$ be a bounded simply connected domain of $\mathbb C$ with a boundary $\Gamma $ which admits a piecewise $C^1$ parametrization $\upsilon:[0,1]\longmapsto \Gamma$. We assume that $\upsilon$ is positively oriented.  Given $1<p<\infty$, we say that an analytic function $f$ on $\Omega$ belongs to the Smirnov space $E^p(\Omega)$ if there exists a sequence $(C_n)_{n\geq 1}$ of rectifiable Jordan curves  included in $\Omega$, tending to the boundary $\Gamma$ (in the sense that $C_n$ eventually surrounds each compact subdomain of $\Omega$), and such that 
\begin{equation}\label{eq-defn1-Ep}
\sup_{n\geq 1}\int_{C_n}|f(z)|^p\,|\mathrm{d}z|~<~\infty.
\end{equation}
With this definition, it is not clear that $E^p(\Omega)$ is a vector space. However, it can be proved (see \cite[page 168-169]{Duren1970}) that if $\varphi$ is a conformal map from $\mathbb D$ onto $\Omega$, and if $\Gamma_r$ is the image under $\varphi$ of the circle $\{z\in\mathbb C:|z|=r\}$, $0<r<1$, then $f\in E^p(\Omega)$ if and only if 
\begin{equation}\label{eq-defn2-Ep}
\sup_{0<r<1}\int_{\Gamma_r}|f(z)|^p\,|\mathrm{d}z|~<~\infty.
\end{equation}
Hence $E^p(\Omega)$ is a vector space and $f\in E^p(\Omega)$ if and only if $(f\circ \varphi)\cdot \varphi'^{1/p}\in H^p(\mathbb D)$ for some (all) conformal map $\varphi$ from $\mathbb D$ onto $\Omega$. 
\par
\begin{figure}[ht]
    \includegraphics[page=58]{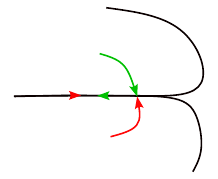}
    \caption{}\label{fig-3434ZEDSSDG}
\end{figure}

\smallskip
Moreover, since $\Gamma$ is a piecewise $C^1$ curve, each function $f\in E^p(\Omega)$ has a non-tangential limit at the point $\upsilon(t)$ for almost every $t\in [0,1]$.
Note that if the point $\upsilon(t)$ travels over a part of the curve $\partial \Omega$ in both directions (which is possible within the setting of assumption \ref{H2'}), the non-tangential limits of the function $\upsilon$ are taken with respect to the domain which remains on the left while $\upsilon(t)$ travels over the curve (and these two limits can then be different). To simplify we still denote by $f$ the non-tangential limit. The notion of non-tangential limits when we travel the curve in both directions is illustrated in \Cref{fig-3434ZEDSSDG}.
\par\smallskip

Note that the non-tangential limit cannot vanish on a set of positive measure unless $f(z)\equiv 0$ (see \cite[Th. 10.3]{Duren1970}). Furthermore, we have 
\[
\int_\Gamma |f(z)|^p\,|\mathrm{d}z|~=~\int_0^1 |f(\upsilon(t))|^p |\upsilon'(t)|\,\mathrm{d}t~<~\infty,
\]
and we can equip $E^p(\Omega)$ with the following norm
\[
\|f\|_{E^p(\Omega)}~=~\left(\int_\Gamma |f(z)|^p\,|\mathrm{d}z|\right)^{1/p}
\]
which turns $E^p(\Omega)$ into a Banach space. It is also known that each function $f\in E^1(\Omega)$ admits a Cauchy representation 
\[
f(z)~=~\frac{1}{2i\pi}\int_\Gamma \frac{f(\tau)}{\tau-z}\,\mathrm{d}\tau\quad\text{ for every } z\in\Omega,
\]
and that the integral above vanishes for all $z\in \mathbb C\setminus\overline{\Omega}$ (see \cite[Th. 10.4]{Duren1970}).
\par\smallskip
If $\Omega$ is a simply connected domain, $\varphi$ is a conformal map from $\mathbb D$ onto $\Omega$, and if $f\in E^p(\Omega)$, $f\not\equiv 0$ with $Z(f)=\{z\in\Omega:f(z)=0\}=\{z_n:n\geq 1\}$, then we have
\begin{equation}\label{eq:Blachke-Ep}
\sum_{n\geq 1}(1-|\varphi^{-1}(z_n)|)~<~\infty.
\end{equation}
Indeed, we have seen that if $f$ belongs to $E^p(\Omega)$, then  the function $F=(f\circ\varphi)\cdot \varphi'^{1/p}$ belongs to $H^p(\mathbb D)$. Moreover, $F(\varphi^{-1}(z_n))=0$ for every $n\geq 1$. Thus \Cref{eq:Blachke-Ep} is nothing but the Blaschke condition on the sequence of zeroes of a non-zero function in $H^p(\mathbb D)$. 
\par\smallskip
When $\Omega$ is a bounded open set that is a finite union of disjoint simply connected components $\Omega_1,\dots, \Omega_N$ such that $\Gamma=\partial\Omega$ 
is piecewise $C^1$, then an analytic function $f$ on $\Omega$ is said to belong to $E^p(\Omega)$ if for every $1\leq j\leq N$, the restriction $f_{|\Omega_j}$ of $f$ to $\Omega_j$ belongs to $E^p(\Omega_j)$. It can be equipped with the norm
\begin{equation}\label{eq:norm-Smirnov-fc}
\|f\|_{E^p(\Omega)}~=~\left(\sum_{j=1}^N\|f_{|\Omega_j}\|_{E^p(\Omega_j)}^p\right)^{1/p}=~\left(\sum_{j=1}^N\int_{\Gamma_j}|f_{|\Omega_j}(z)|^p\,|\mathrm{d}z| \right)^{1/p}
\end{equation}
where $\Gamma_j=\partial \Omega_j$, $1\le j\le N$. Then $E^p(\Omega)$ becomes a Banach space. 
\par\smallskip
We will also need to consider more general simply connected domains $\Omega$ within the extended complex plane $\widehat{\mathbb C}$. Let $\Omega$ be such a simply connected domain with piecewise $C^1$ boundary such that $\infty \in\Omega$.
An analytic function $f$ on $\Omega$ is said to belong to $E^p(\Omega)$ if  there exists a sequence of rectifiable Jordan curves $(\gamma_n)_{n\geq 1}$ lying in $\Omega$, tending to $\partial\Omega$, and such that 
\[
\sup_{n\geq 1}\int_{\gamma_n}|f(z)|^p\,|\mathrm{d}z|~<~\infty.
\]
Fix $a\notin\Omega$ and define a function $\Theta$ by setting $\Theta(z)=\frac{1}{z-a}$, $z\in\Omega$. Then it is not difficult to prove that 
\begin{equation}\label{Eq:smirnoviwqriou}
g\in E^p(\Omega)\text{ if and only if } h~=~\left(g\circ \Theta^{-1}\right)\cdot \left(\left(\Theta^{-1}\right)'\right)^{1/p}\in E^p\left(\widetilde{\Omega}\right),
\end{equation}
where $\widetilde{\Omega}=\Theta(\Omega)$. Then $\widetilde{\Omega}$ is a bounded simply connected domain with rectifiable boundary. And we define $E^p_0(\Omega)$ as the space of functions $f\in E^p(\Omega)$ such that $f(\infty)=0$. 

In other words, this means that, without passing by the Riemann sphere, for $\Omega\subset \mathbb C$ an unbounded domain such that $\partial \Omega$ is a   piecewise $C^1$ closed curve, we say that $f\in Hol(\Omega)$ belongs to $E^p_0(\Omega)$ if $\lim_{|z|\to\infty}f(z)=0$ and there exists a sequence of rectifiable Jordan curves $(\gamma_n)_{n\geq 1}$ lying in $\Omega$, tending to $\partial\Omega$, and such that 
\[
\sup_{n\geq 1}\int_{\gamma_n}|f(z)|^p\,|\mathrm{d}z|~<~\infty.
\]
Using \eqref{Eq:smirnoviwqriou}, it can be proved that most of the previous properties of Smirnov spaces on bounded simply connected domain extend to the space $E^p(\Omega)$ and thus to the space $E_0^p(\Omega)$. When $\Omega$ is an open set that is a finite union of simply connected components in $\widehat{\mathbb C}$, we can define $E_0^p(\Omega)$ as in the bounded case: the norm defined by \Cref{eq:norm-Smirnov-fc}  also turns it into a Banach space.
 
\subsection{Linear dynamics}\label{Subsection:LinearDyn}
We briefly recall here the definitions of some important concepts in linear dynamics which are considered in this paper. The books \cite{GrosseErdmannPeris2011} and \cite{BayartMatheron2009} are excellent sources on linear dynamics.
\subsubsection{Topological setting}
Given a separable Banach space $X$ over $\mathbb K=\mathbb R$ or $\mathbb C$, a bounded operator on $X$ is said to be \textit{hypercyclic} if there exists a vector $x\in X$ whose orbit $\{T^nx\,;\,n\ge0\}$ is dense in $X$. 
Related weaker notions are those of supercyclicity and cyclicity: $T$ is \textit{supercyclic} if there exists $x\in X$ whose projective orbit $\{\lambda T^nx\,;\,n\ge0,\lambda\in\mathbb K\}$ is dense in $X$, and $T$ is \textit{cyclic} if there exists $x\in X$ such that the linear space $\spa\,[T^nx\,;\,n\ge0]$ of the orbit of $x$ under the action of $T$ is dense in $X$. 
Clearly, hypercyclicity implies supercyclicity, which in its turn implies cyclicity.
Observe also that whenever $T$ is hypercyclic, $rT$ is supercyclic for every scalar $r\in\mathbb K\setminus\{0\}$. 
\par\smallskip
Obviously, a hypercyclic operator cannot be power bounded. Moreover, hypercyclicity entails the following spectral restriction (see for instance \cite[Th. 5.6]{GrosseErdmannPeris2011}) : if $T$ is a hypercyclic operator on a complex Banach space $X$, then every connected component of its spectrum $\sigma(T)$ meets the unit circle. Among straightforward necessary conditions for the hypercyclicity of $T\in \mathcal{B}(X)$, let us mention the fact that $T^*$ cannot have any eigenvalue. See \cite[Ch. 5]{GrosseErdmannPeris2011} for more necessary conditions for hypercyclicity. For instance, normal operators on a Hilbert space are never hypercyclic. Similar results hold for supercyclicity: if $T\in \mathcal{B}(X)$ is supercyclic, then there exists $r\ge0$ such that $r\mathbb T$ intersects each connected component of the spectrum of $T$ (see \cite[Th. 1.24]{BayartMatheron2009}). Using the equivalence between hypercyclicity and topologically transitivity, it is not difficult to see that when $T$ is invertible, then $T$ is hypercyclic if and only if $T^{-1}$ is hypercyclic (see \cite[Cor. 1.3]{BayartMatheron2009}).
\par\smallskip
We spell out the following simple fact as a Lemma.
\begin{lemma}\label{lemme-HC-direct-sum}
Let $T$ be a bounded operator on a Banach space $X$, and suppose that $X$ is decomposed as a topological direct sum $X=M_1\oplus M_2$, where $M_i$ is a non-trivial closed $T$-invariant subspace of $X$, $i=1,2$. Denote by $T_i$ the operator induced by $T$ on $M_i$. Let $x=x_1+x_2$, $x_i\in M_i$, $i=1,2$. If $x$ is a hypercyclic (resp. supercyclic) vector for $T$, then $x_i$ is a hypercyclic (resp. supercyclic) vector for $T_i$, $i=1,2$.  
\end{lemma}
\begin{proof}
Suppose that $x$ is hypercyclic. Then the set $\{T^nx_1+T^nx_2\,;\,n\ge 0\}$ is dense in $X$. Denoting by $P_i$ the bounded projection onto $M_i$, it follows that the set $\{P_i T^nx_i\,;\,\ge 0\}$ is dense in $M_i$. But $P_iT^nx_i=T_i^n x_i$ for all $n\geq 0$, so that $x_i$ is a hypercyclic vector for $T_i$. The case where $x$ is supercyclic is similar. 
\end{proof}
\par\smallskip
One of the most useful tools for proving the hypercyclicity of some concrete classes of operators is the Godefroy-Shapiro Criterion, which requires the existence of a large supply of eigenvectors associated to eigenvalues of modulus smaller than $1$ and larger than $1$, respectively. It is a straightforward consequence of the so-called Hypercyclicity Criterion (see for instance \cite[Ch. 3]{GrosseErdmannPeris2011}).

\begin{theorem}\label{GS-H}
Let $X$ be a complex separable Banach space and let $T\in \mathcal{B}(X)$. If the two subspaces
\[H_-(T)~=~\overline{\spa}\,[\ker(T-\lambda)\,;\,|\lambda|<1]
\quad\text{and}\quad H_+(T)~=~\overline{\spa}\,[\ker(T-\lambda)\,;\,|\lambda|>1]
\]
are equal to $X$, then $T$ is hypercyclic. 
\end{theorem}

The next statement provides a sufficient condition for supercyclicity in terms of eigenvectors.
\begin{theorem}\label{GS-S}
Let $X$ be a complex separable Banach space and let $T\in \mathcal{B}(X)$. If for some $r>0$, the two subspaces
\[H_{r,-}(T)~=~\overline{\spa}\,[\ker(T-\lambda)\,;\,|\lambda|<r]\quad
\text{and}\quad H_{r,+}(T)~=~\overline{\spa}\,[\ker(T-\lambda)\,;\,|\lambda|>r]
\]
are equal to $X$,
then $T$ is supercyclic. 
\end{theorem}
\Cref{GS-S} is a trivial consequence of \Cref{GS-H}: under the assumption of \Cref{GS-S}, the operator $\frac1rT$ satisfies the assumptions of \Cref{GS-H}, so $\frac1rT$ is hypercyclic, and so $T$ is supercyclic.
\par\smallskip
Another important notion in linear dynamics is that of chaos: a bounded operator $T\in \mathcal{B}(X)$ is said to be \textit{chaotic} if $T$ is hypercyclic and has a dense set of periodic points (a vector $x\in X$ is \textit{periodic} for $T$ if there exists an integer $n\ge1$ such that $T^nx=x$). Again, there is a very useful version of the Godefroy-Shapiro Criterion for chaos \cite{GodefroyShapiro1991}:

\begin{theorem}\label{GS-C}
Let $X$ be a complex separable Banach space and let $T\in \mathcal{B}(X)$. If the three subspaces of $X$ given by $H_+(T),~H_-(T)$ and
\[H_0(T)~=~\overline{\spa}\,[\ker(T-\lambda)\,;\,\lambda~\text{is an $n$-th root of unity},~n\ge1]\]
are equal to $X$, then $T$ is chaotic.
\end{theorem}

Since the set of periodic points for $T\in \mathcal{B}(X)$ is the linear vector space spanned by the eigenvectors of $T$ associated to eigenvalues which are roots of unity, a necessary condition for the chaoticity of $T$ is that $H_0(T)=X$.

\subsubsection{Measure-theoretic setting}
The few basic results on hypercyclicity that we have presented above already highlight the importance of eigenvectors in the study of linear dynamics from the topological point of view. Eigenvectors are even more important in the measure-theoretic setting. An excellent reference for the concepts recalled in this subsection is \cite[Ch. 5]{BayartMatheron2009}.
\par\smallskip
Given a (complex) separable Banach space $X$ and $T\in \mathcal{B}(X)$, the general aim of linear dynamics in the measure-theoretical setting is to endow $X$ with a $T$-invariant Borel probability measure $m$, and to study the properties of the measure-preserving dynamical system $(X,\cal B,m;T)$. Here $\cal B$ denotes the $\sigma$-algebra of Borel subsets of $X$. The measure $m$ is \textit{$T$-invariant} if $m(T^{-1}(A))=m(A)$ for every $A\in \cal B$, and $T$ is \textit{ergodic with respect to $m$} if the only sets $A\in\cal B$ such that $T^{-1}A=A$ are those with $m(A)=0$ or $m(A)=1$. The study of the dynamics of $T$ with respect to a probability measure $m$ is especially interesting when $m$ has \textit{full (topological) support}, i.e. when $m(U)>0$ for every non-empty open subset $U$ of $X$.
\par\smallskip
Gaussian measures play a specific role in the study of the dynamics of bounded operators from a measure-theoretical point of view. A \textit{Gaussian measure} $m$ on $X$ is a Borel probability measure such that every functional $x^*\in X^*$, when considered as a complex random variable on $X$, has symmetric complex Gaussian distribution.
\par\smallskip
In the Hilbertian setting, there is a very neat characterization of operators admitting an invariant ergodic measure with full support, given in terms of \textit{unimodular eigenvectors} {\cites{Fl, BG}}. Before we state this characterization as \Cref{cas-hilbertien}, we recall that $T\in \mathcal{B}(X)$ has \textit{perfectly spanning unimodular eigenvectors} if the following holds:

for \textit{any} countable subset $D$ of the unit circle $\mathbb T$, 
\[\overline{\spa}\,[\ker(T-\lambda)\,;\,\lambda\in\mathbb T\setminus D]~=~X.\]

This property admits several equivalent formulations (see {\cite{G}} for more details).

\begin{theorem}\label{cas-hilbertien}
Let $H$ be a complex separable Hilbert space, and let $T\in\cal B(H)$. The following assertions are equivalent:
\begin{enumerate}
    \item $T$ admits a $T$-invariant Gaussian probability measure with full support with respect to which $T$ is ergodic;
    \item $T$ has perfectly spanning unimodular eigenvectors.
\end{enumerate}
\end{theorem}
Outside the Hilbertian setting, the implication (2)$\implies$(1) in \Cref{cas-hilbertien} remains true in full generality \cite{BM2}.

\begin{theorem}\label{Th:TheoremA13}
Let $X$ be a complex separable space, and let $T\in \mathcal{B}(X)$. If $T$ has perfectly spanning unimodular eigenvectors, then $T$ admits an invariant Gaussian probability measure with full support with respect to which $T$ is ergodic.
\end{theorem}

One of the interests of such results is that whenever $T$ is ergodic with respect to an invariant probability-measure $m$ with full support, $T$ is \textit{frequently hypercyclic}: there exists a vector $x\in X$ (a frequently hypercyclic vector for $T$) such that for every non empty open subset $U$ of $X$, the set
\[N_T(x,U)~=~\{n\ge0\,;\,T^nx\in U\}\]
has positive lower density, i.e. 
\[\liminf_{N\to\infty}\frac1N\text{card}\{0\le n<N\,;\,T^nx\in U\}>0.\]
This follows from Birkhoff's Pointwise Ergodic Theorem. Actually, $m$-almost every vector $x\in X$ is frequently hypercyclic for $T$.
Hence we have

\begin{theorem}\label{Theo:A19}
Let $X$ be a complex separable Banach space, and let $T\in \mathcal{B}(X)$. If $T$ has perfectly spanning unimodular eigenvectors, then $T$ is frequently hypercyclic. 
\end{theorem}

We summarize some of the results presented above in a last statement, which provides us with a concrete criterion for proving chaos or ergodicity of large classes of operators.
\begin{theorem}\label{Coro:A12}
Let $X$ be a complex separable Banach space, and let $T\in \mathcal{B}(X)$. Suppose that there exist a open subset $U$ of $ \mathbb C$ and a finite or countable family $(E_i)_{i\in I}$ of analytic maps from $U$ into $X$ such that 
\begin{enumerate}[(a)]
    \item $TE_i(\lambda)=\lambda E_i(\lambda)$ for every $\lambda\in U$ and every $i\in I$;
    \item $\overline{\spa}\,[E_i(\lambda)\,;\,\lambda\in U,i\in I]=X$;
    \item every connected component of $U$ intersects the unit circle.
\end{enumerate}
Then $T$ has the following properties:
\begin{enumerate}
    \item $T$ admits an invariant Gaussian probability measure with full support with respect to which $T$ is ergodic;
    \item $T$ is frequently hypercyclic, hence hypercyclic;
    \item $T$ is chaotic.
\end{enumerate}
\end{theorem}

\begin{proof}
In order to prove assertions (1) and (2), it suffices by \Cref{Th:TheoremA13} to show that $T$ has perfectly spanning unimodular eigenvectors. So let $D$ be a countable subset of $\mathbb T$. For each connected component $O$ of $U$, the set $\mathbb T\cap O$ contains a non-empty open subarc $\Gamma_O$ of $\mathbb T$ by (c). Since $\Gamma_O\setminus D$ has accumulation points in $O$, the uniqueness principle for analytic functions implies that 
\[\overline{\spa}\,[E_i(\lambda)\,;\,\lambda\in \Gamma_O\setminus D,i\in I]~=~\overline{\spa}\,[E_i(\lambda)\,;\,\lambda\in O,i\in I].\]
Hence the closed subspace
\[\overline{\spa}\,[E_i(\lambda)\,;\,\lambda\in (\mathbb T\cap U)\setminus D,i\in I]\]
contains the linear span of all the subspaces 
\(\overline{\spa}\,[E_i(\lambda)\,;\,\lambda\in O,i\in I]\)
where $O$ varies over all connected components of $U$. It follows from (b) that
\[\overline{\spa}\,[E_i(\lambda)\,;\,\lambda\in (\mathbb T\cap U)\setminus D,i\in I]~=~X,\]
and thus by (a), 
\[\overline{\spa}\,[\ker(T-\lambda);\lambda\in\mathbb T\setminus D]~=~X.\]
Hence $T$ has perfectly spanning unimodular eigenvectors. 

The strategy to show assertion (3) is exactly the same, observing that if $E$ denotes the set of all roots of unity in $\mathbb T$, $\Gamma_O\cap E$ has accumulation points in $\mathbb T$ for every connected component $O$ of $U$, so that 
\[\overline{\spa}\,[E_i(\lambda)\,;\,\lambda\in \Gamma_O\cap E,i\in I]~=~\overline{\spa}\,[E_i(\lambda)\,;\,\lambda\in O,i\in I].\]
Hence $\overline{\spa}\,[\ker(T-\lambda)\,;\,\lambda\in E]=X$, and $T$ has a dense set of periodic points. Since we already know, by (2), that $T$ is hypercyclic, this yields that $T$ is chaotic.
\end{proof}

We finish this section by recalling an important sufficient condition for the cyclicity of a bounded operator on a Banach space, which is crucial in the works \cite{Yakubovich1991}, \cite{Yakubovich1993} and \cite{Yakubovich1996} to show the cyclicity of certain Toeplitz operators. This condition originates from \cite{Nikolski-OTAA89}. The brief proof that we present here is inspired from \cite{Grivaux-JOT2005}. For each polynomial $p\in\mathbb C[z]\setminus\{0\}$, we denote by $Z(p)$ the set of its roots.
\begin{lemma}\label{lem:cyclicite-sophie-nikolski}
Let $T$ be a bounded operator on a complex separable Banach space $X$, and let $\mathcal{D}=\{p\in\mathbb C[z]\setminus\{0\}\,;\, Z(p)\cap\sigma_p(T^*)=\varnothing\}$. Suppose that the vector space
\[
D~:=~\bigcup_{p\in\mathcal{D}}\ker p(T)
\]
is dense in $X$. Then $T$ is cyclic.
\par\smallskip
In particular, if the linear span of the eigenspaces $\ker (T-\lambda)$, $\lambda\not\in \sigma_p(T^*)$, is dense in $X$, then $T$ is cyclic.
\end{lemma}

\begin{proof}
By the Baire Category Theorem, it suffices to show that for any non-empty open subsets $U,V$ of $X$, there exists a polynomial $q\in\mathbb C[z]$ such that $q(T)(U)\cap V\neq\varnothing$. 
Let $u\in U\cap D$, and let $p\in\mathcal{D}$ be a non-zero polynomial such that $p(T)u=0$. Since no root of $p$ belongs to $\sigma_p(T^*)$, the operator $p(T)$ has dense range. Given $v\in V$, for any $\varepsilon>0$, there exists a vector $x\in X$ such that $\|p(T)x-v\|<\varepsilon$. Let now $\delta>0$ be such that the open ball $B(u,\delta)$ is contained in $U$, and let $r>\frac{2}{\delta}\|x\|$. Then $\|\frac{x}{r}\|<\delta$. Let now $ \widetilde u:=u+\frac{x}{r}$. Then $ \widetilde u\in U$, and if we set $q:=rp$, we have 
\[
q(T) \widetilde u~=~rp(T)\frac{x}{r}~=~p(T)x.
\]
So $\|q(T) \widetilde u-v\|<\varepsilon$, so that $q(T) \widetilde u\in V$ if $\epsilon$ is sufficiently small. Hence $q(T)(U)\cap V\ne\varnothing$, which is the result we were looking for.
\end{proof}

\section{Proof of the $H^p$ version of Yakubovich's model theory for Toeplitz operators}\label{Section:Yakubovich_demoHp}
In this appendix, we give a complete proof of the main result of Yakubovich \cite{Yakubovich1991} in the $H^p$ setting. We also prove some results which are true in the more general framework of \cite{Yakubovich1996}. 
\par\smallskip
\textbf{Caution :} In this \Cref{Section:Yakubovich_demoHp}, the letter $F$ will denote a \textbf{positively wound symbol} given by \ref{H3'}, while in the previous sections it was mainly denoting a negatively wound symbol, i.e. a symbol which satisfies \ref{H3}. This choice, which is of course questionable, is motivated by the fact that we preferred to stick to the notation of \cite{Yakubovich1991}, so as not confuse the reader. We will explain in Subsection \ref{mot-de-la-fin} how to connect the notations and results from \Cref{Section:Yakubovich_demoHp} to those of the rest of the paper.
\par\smallskip
Let $1<p<\infty$ and let $q$ be its conjugate exponent i.e. $1/p+1/q=1$. 
In the forthcoming  appendix, we consider the Toeplitz operator $T_F$, with symbol $F$, defined on the Hardy space $H^p$.  We recall here the assumptions on the symbol that will be used in the sequel (which are almost the same as in \Cref{Section-intro}).\par\smallskip
\begin{enumerate}[]
    \item[(H1)] $F$ belongs to the class $C^{1+\varepsilon}(\mathbb T)$ for some $\varepsilon>\max(1/p,1/q)$, and its derivative $F'$ does not vanish on $\mathbb T$;\par\smallskip
    \item[(H2)] there exist real numbers $\theta_0<\theta_1<\dots<\theta_m=\theta_0+2\pi$ such that, letting $\alpha_j$ be the open arc $\alpha_j=\{e^{i\theta};\,\theta_j<\theta<\theta_{j+1}\}$, we have
     \begin{enumerate}[(a)]
       \item $F$ is injective on  each arc $\alpha_j,~0\le j \le m-1$;
        \item for every $i\neq j,~0\le i,j\le m-1$, the sets $F(\alpha_i)$ and $F(\alpha_j)$ are disjoint;\par\smallskip
    \end{enumerate}
    \item[(H2')]  there exist real numbers $\theta_0<\theta_1<\dots<\theta_m=\theta_0+2\pi$ such that, letting $\alpha_j$ be the open arc $\alpha_j=\{e^{i\theta};\,\theta_j<\theta<\theta_{j+1}\}$, we have
     \begin{enumerate}[(a)]
       \item $F$ is injective on  each arc $\alpha_j,~0\le j \le m-1$;
        \item for every $i\neq j,~0\le i,j\le m-1$, $F(\alpha_i)=F(\alpha_j)$ or the sets $F(\alpha_i)$ and $F(\alpha_j)$ are disjoint;\par\smallskip
    \end{enumerate}
    \item[(H3')]\customlabel{H3'}{(H3')} for every $\lambda\in\mathbb C\setminus F(\mathbb T)$, $\w_F(\lambda)\ge0$, where $\w_F(\lambda)$ denotes the winding number of the curve $F(\mathbb T)$ around $\lambda$.
\end{enumerate}\par\smallskip
Denote by $\cal O  \subseteq F(\mathbb{T})$ the set of all the images by $F$ of the extremities of the arcs $\alpha_j$.
\par\smallskip
The derivative $F'$ of $F$ at a point $z_0$ of $\mathbb T$ will be understood as the limit as $z\rightarrow z_0$, $z\in\T$ of the quotients $(F(z)-F(z_0))/(z-z_0)$. Writing $z_0=e^{i\theta_0}$, $\theta_0\in [0,2\pi)$, it could also be understood as the derivative at the point $\theta_0$ of the function $\theta\longmapsto F(e^{i\theta})$ defined on $ [0,2\pi)$ (the derivative of this last function is non-zero on $ [0,2\pi)$ if and only if $F'$ is non-zero on $\T$).

\subsection{Some notations}\label{piece of notation}
\begin{figure}[ht]
    \includegraphics[page=56,scale=1]{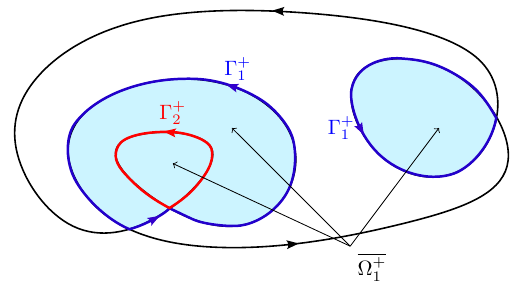}
    \caption{}\label{fig-49-notations}
\end{figure}

We now set an important piece of notation which will be involved in the statement of the boundary relations for eigenvectors of $T_F^*$. 
Let $N=\max\{\w_F(\lambda)\,;\,\lambda\in\C\setminus F(\T)\}$. We set
 \[
 \Omega_j^+~=~\{\lambda\in \mathbb C\setminus F(\mathbb T)\,;\, \w_F(\lambda)>j\}\quad\text{ for every } 0\le j\le N-1,\]
and
 \[
 \Gamma_j^+~=~\partial \overline{\Omega_j^+}
 \]
for every $0\le j\le N-1$.
We also define $\Gamma_N^+=\varnothing$ and  $\Gamma=F(\T)$. Observe that $\Gamma=\partial\Omega_0^+$. The definition of $\Gamma_j^+$ is illustrated in \Cref{fig-49-notations}.

A word of caution is in order here: the set $\Gamma_j^+$ is defined as the boundary of the closure of the set $\Omega_j^+$ (and not as the boundary of $\Omega_j^+$). Suppose that assumption \ref{H2} is satisfied. A point $\lambda\in\Gamma\setminus\mathcal{O}$ belongs to $\Gamma_j^+$ if and only if its interior and exterior components $\Omega_{int}$ and $\Omega_{ext}$, defined in \Cref{subsection:JFA-model}, are such that $\w_F(\Omega_{int})=j+1$ and $\w_F(\Omega_{ext})=j$. Hence $\Omega_{int}$ is a connected component of $\Omega_j^+$ and
$\Omega_{ext}$ is a connected component of $\C\setminus (F(\T)\cup \Omega_j^+)$.
We  have
 $$ \Gamma\cap \overline{\Omega_{j}^+}~=~\partial \Omega_{j}^+~=~\bigcup_{k=j}^N\Gamma_k^+$$ for every $0\le j\le N-1$.
Two different sets $\Gamma_j^+$ can intersect in a finite set of points at most, and
$$\Gamma~=~\bigcup_{j=0}^{N-1}\Gamma_j^+.$$ In other words the family $\Gamma_j^+$, $0\le j\le N-1$, is a partition of $\Gamma$ up to a finite set of points.
We also define for $0\leq j\leq N$ 
\[
\Omega_j~=~\{\lambda\in\mathbb C\setminus F(\mathbb T)\,;\,\w_F(\lambda)=j\}.
\]
Note that $\Omega_j^+=\displaystyle\bigcup_{k\ge j+1}\Omega_k$.
\par\smallskip
Let $r>1$ and $0\le j \le N-1$. Recall that for $u\in E^r(\Omega)$ (where $\Omega$ satisfies the assumptions of \Cref{smirnov}), we still denote by $u$ its non-tangential limit on the boundary which exists almost everywhere (in the sense given in \Cref{smirnov}). 
According to \cref{eq:norm-Smirnov-fc}, for $u\in E^r(\Omega_j^+)$, we have
\begin{equation}\label{eq:norme-E-r-3434}
\|u\|_{E^r(\Omega_j^+)}^r
~=~\sum_{k=j+1}^N\|u_{|\Omega_k}\|_{E^r(\Omega_k)}^r
~=~\sum_{k=j+1}^N\|u_{|\Omega_k}\|_{L^r(\partial\Omega_k)}^r.
\end{equation}

\begin{figure}[ht]
\includegraphics[page=52]{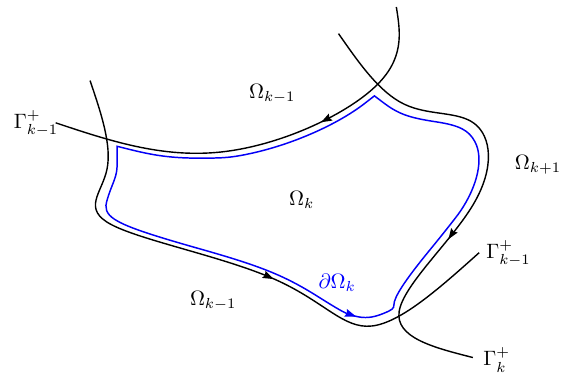}
\caption{}\label{fig-orientation50}\end{figure}

When $F$ satisfies \ref{H2}, this expression can be simplified a bit.
Note that when we travel along the curve $\Gamma_{k-1}^+$ endowed with the parametrization $F$, the set $\Omega_k$ is on the left side, and when we travel along the curve $\Gamma_{k}^+$, it is on the right side. So in order to have a positive orientation for $\partial\Omega_{k}$, which is the union of $\Gamma_{k-1}^+$ and $\Gamma_{k}^+$, we write $\partial\Omega_{k}$ as $\partial\Omega_{k}=\Gamma_{k-1}^+\bigcup (\widetilde{\Gamma}_{k}^+)$, where $\widetilde{\Gamma}_{k}^+$ is the same curve as $\Gamma_{k}^+$, but with the opposite orientation. See \Cref{fig-orientation50}.

On $\Gamma_k^+$ (still under hypothesis \ref{H2}), recall that the interior and exterior boundary values, defined in  \Cref{subsection:JFA-model}, are such that, for almost every $\lambda\in\Gamma_k^+$, 
\begin{equation}\label{Eq:limIntExt-k}
    u^{int}(\lambda)~=\lim_{\substack{\mu\to\lambda\\ \mu\in\Omega_{k+1}}}u(\mu) \quad \text{and}\quad u^{ext}(\lambda)~=\lim_{\substack{\mu\to\lambda\\ \mu\in\Omega_{k}}}u(\mu).
\end{equation}
This means that $u^{int}$ is the non-tangential limit on $\Gamma_k^+$ of $u_{|\Omega_{k+1}}$ and $u^{ext}$ is the non-tangential limit on $\Gamma_k^+$ of $u_{|\Omega_{k}}$. 
\par\smallskip
Thus, under assumption \ref{H2}, \Cref{eq:norme-E-r-3434} can be rewritten as
\begin{eqnarray}
\|u\|_{E^r(\Omega_j^+)}^r
&=&\sum_{k=j+1}^N\|u_{|\Omega_k}\|_{L^r(\partial\Omega_k)}^r\notag\\
&=&\sum_{k=j+1}^N\|u^{int}\|_{L^r(\Gamma_{k-1}^+)}^r~+~\sum_{k=j+1}^N\|u^{ext}\|_{L^r(\Gamma_{k}^+)}^r.\label{Eq:NormeErLimIntExt}
\end{eqnarray}

Finally, note that for every function $u\in E^1(\Omega_k)$, we have 
\begin{equation}\label{blablaesdsds232323}
\int_{\partial\Omega_k}u(\lambda)\,\mathrm{d}\lambda~=~\int_{\Gamma_{k-1}^+}u^{int}(\lambda)\,\mathrm{d}\lambda -\int_{\Gamma_{k}^+}u^{ext}(\lambda)\,\mathrm{d}\lambda.  
\end{equation}

\subsection{Construction of eigenvectors}
Let us now assume that $F$ satisfies the assumptions \ref{H1} and \ref{H3'}. Since $F$ is continuous, according to assertion $(ii)$ of \Cref{Lemma:spectral-properties-toeplitz-operators}, it follows that
\[
\sigma(T_F)~=~F(\mathbb T)\cup\{\lambda\in \mathbb C\setminus F(\mathbb T):\w_F(\lambda)>0\}.
\]\par\smallskip
Now, according to assertion $(i)$ of \Cref{Lemma:spectral-properties-toeplitz-operators}, for every $\lambda\in\sigma(T_F)\setminus F(\mathbb T)$, we have
\[
j(T_F)~=~\dim(\ker(T_F-\lambda))-\dim(\ker(T_F^*-\lambda))~=~-\w_F(\lambda),
\]

and then, by Coburn's theorem, 
\[\ker(T_F-\lambda)~=~\{0\}\quad\text{and}\quad\dim(\ker(T_F^*-\lambda))~=~\w_F(\lambda).\]
\par\smallskip
In this section, we will construct a basis of eigenvectors for $T_F^*-\lambda=T_f-\lambda$, with $f(z)=F(1/z),~z\in\mathbb T$.
\par\smallskip
 For $\lambda\in\mathbb C\setminus  F(\mathbb T)$ and $m=\w_f(\lambda)$, let $$\varphi_\lambda(z)~=~z^{-m}(f(z)-\lambda)\quad \textrm{ for every } z\in\T.$$ 
 Note that $0\notin \varphi_\lambda(\mathbb T)$ and that $\w_{\varphi_\lambda}(0)=0$. Consider now the function $u_\lambda$ defined on $\T$ by
\[
u_\lambda(e^{is})~=~\int_0^s ie^{it}\frac{\varphi_\lambda'(e^{it})}{\varphi_\lambda(e^{it})}~\mathrm dt\quad\textrm{ for every } s\in [0,2\pi).
\]
For every $\lambda\in\mathbb C\setminus F(\mathbb T)$, the function $z\longmapsto u_\lambda(z)$ is of class $C^{1+\varepsilon}$ on $\mathbb T$, and for every $s\in [0,2\pi)$, the function $\lambda\mapsto u_\lambda(e^{is})$ is analytic on $\mathbb C\setminus F(\mathbb T)$. Moreover, $\varphi_\lambda(e^{is})=\varphi_\lambda(1)e^{u_\lambda(e^{is})}$  (and $u_\lambda+\log(\varphi_\lambda(1))$ is a continuous logarithm of $\varphi_\lambda$ on $\T$). Define a function $v_\lambda$ on $\T$ by 
setting
\begin{equation}\label{Eq1}
    v_\lambda(z)~=~(P_+ u_\lambda)(z)~=~\frac1{2i\pi}\int_{\mathbb T}\frac{u_\lambda(\tau)}{\tau-z}\,\mathrm d\tau\quad\textrm{ for every } z\in\mathbb D.
\end{equation}
\begin{lemma}\label{Lem:analytic-en-dehors-de-la-courbe}
With the previous notation, we have 
\begin{enumerate}
\item for all $\lambda\in\mathbb C\setminus F(\mathbb T)$, $v_\lambda\in A(\mathbb D)$;
\item for all $z\in \mathbb D$, the map $\lambda\mapsto v_\lambda(z)$ is analytic on $\mathbb C\setminus F(\mathbb T)$.
\end{enumerate}
\end{lemma}
\begin{proof}
According to \Cref{thm:privalov-Zygmund}, the projection $P_+$ maps $C^{1+\varepsilon}(\mathbb T)$ into itself, whence  $v_\lambda$ belongs to $A(\mathbb D)$. 
Assertion (2) follows immediately from \Cref{Eq1}.
\end{proof}

For every $\lambda\in \C\setminus F(\T)$, define the functions $f_\lambda^+$ and $f_\lambda^-$ on $\mathbb D$ and $\mathbb C\setminus\overline{\mathbb D}$ respectively by setting
\[f_\lambda^+~=~(f(1)-\lambda)e^{v_\lambda}\quad\text{and}\quad f_\lambda^-~=~e^{P_-u_\lambda},  
\]
where $P_-=I-P_+$. Note (this will be important later on) that the definition of $f_\lambda^+$ does not depend on the choice of the continuous branch of the logarithm of $\varphi_\lambda$ that we make. 
The following equality holds on $\T$:
\[
f_\lambda^+f_\lambda^-~=~(f(1)-\lambda)e^{u_\lambda}~=~\varphi_\lambda = z^{-m}(f-\lambda).
\]
Moreover, $f_\lambda^+$ belongs to  $A(\mathbb D)$ and $\lambda\mapsto f_\lambda^+(z)$ is analytic on $\mathbb C\setminus F(\mathbb T)$. Observe also that $1/f_\lambda^+$ belongs to $A(\mathbb D)$ (because $f_\lambda^+$ belongs to $A(\mathbb D)$ and does not vanish on $\overline{\mathbb D}$). 

\begin{lemma}\label{lemma:P_+-P_-}
Let $a\in L^2(\mathbb T)$ and set $b(\tau)=a(1/\tau)$, $\tau\in\T$. Then for all $z\in\mathbb D$, we have $P_+b(z)=\hat a(0)+P_-a(1/z)$.
\end{lemma}

\begin{proof}
Set $a_n=\hat a(n)$ for every $n\in\Z$, so that  $a(\tau)=\sum_{n=-\infty}^{+\infty}a_n\tau^n$ on $\T$. Then $b(\tau)=\sum_{n=-\infty}^{+\infty}a_n\tau^{-n}$ and thus
\[P_+b(z)~=\sum_{n=-\infty}^0a_nz^{-n}~=~a_0~+\sum_{n=-\infty}^{-1}a_nz^{-n}~=~\hat a(0)+P_-a(1/z)\quad\text{ for every }z\in\D.\qedhere\]
\end{proof}

Similarly, starting from the function $\psi_\lambda$ defined on $\T$ by $\psi_\lambda(z)=z^{-n}(F(z)-\lambda)$ with $n=\w_F(\lambda)$, we may construct functions $U_\lambda$ and $V_\lambda$ associated to $F$, and then $F_\lambda^+$ and $F_\lambda^-$, satisfying the following relations:
\[
z^{-n}(F-\lambda)~=~(F(1)-\lambda) e^{U_\lambda},\quad V_\lambda~=~P_+U_{\lambda},
\]
and
\begin{equation}\label{Eq:Flambda+-}
F_\lambda^+~=~(F(1)-\lambda) e^{V_\lambda},\quad F_\lambda^-~=~e^{P_-U_\lambda}.
\end{equation}
We thus obtain that
\begin{equation}\label{Eq:B2}
F-\lambda~=~z^n F_\lambda^+ F_\lambda^-\quad \text{on }\mathbb T,
\end{equation}
Recall that (see \cite{Brown-Halmos}), $T_{F_1F_2}=T_{F_1}T_{F_2}$ if and only if $\overline{F_1}$ or $F_2$ belongs to $H^\infty$. In particular if $F_2\in H^\infty$, we have $T_{F_1F_2}(g)=T_{F_1}(F_2g)$, for every $g\in H^p$.
So \Cref{Eq:B2} gives, taking into account the fact that $F_\lambda^+\in A(\mathbb D)$ and $F_\lambda^{-}\in L^\infty(\mathbb T)$, the equality
\begin{equation}\label{eq:formule-T-F-Lambda}
T_{F-\lambda}(g)~=~T_{F_\lambda^-}(z^n F_\lambda^+ g)\quad\textrm{ for every } g\in H^p.
\end{equation}
\begin{lemma}\label{Lem:lien-Flamba-et-flambda+}
For every $\lambda\in\mathbb C\setminus F(\mathbb T)$, we have
\[
F_\lambda^-(z)~=~\frac{f_\lambda^+(1/z)}{f_\lambda^+(0)}\quad \text{ for every } z\in\mathbb C\setminus\overline{\mathbb D}.
\]
\end{lemma}

\begin{proof}
Note that $\w_f(\lambda)=-\w_F(\lambda)$ and $U_\lambda(e^{is})=u_\lambda(e^{-is})$ for every $s\in [0,2\pi)$. Then, according to \Cref{lemma:P_+-P_-}, for $z\in\mathbb D$ we have $P_+u_\lambda(z)=\widehat{U_\lambda}(0)+P_-U_\lambda(1/z)$. Therefore, for every $z\in\mathbb C\setminus\overline{\mathbb D}$, we get
\[
e^{P_-U_\lambda(z)}~=~\frac{e^{P_+u_\lambda(1/z)}}{e^{\widehat{U_\lambda}(0)}}\cdot
\]
In other words, taking the definitions of $F_\lambda^-$ and $f_\lambda^+$ into account, we have
\[
F_\lambda^-(z)~=~\frac{f_\lambda^+(1/z)}{(f(1)-\lambda)e^{\widehat{U_\lambda}(0)}}\quad \text{ for every } z\in\mathbb C\setminus\overline{\mathbb D}.
\]
Observe now that $\widehat{U_\lambda}(0)=\widehat{u_\lambda}(0)=P_+u_\lambda(0)$, which gives 
\[
F_\lambda^-(z)~=~\frac{f_\lambda^+(1/z)}{f_\lambda^+(0)}\quad \text{ for every } z\in\mathbb C\setminus\overline{\mathbb D}.\qedhere
\]
\end{proof}

We have already mentioned the fact that $f_\lambda^+$ and $F_\lambda^+$ belong to $A(\mathbb D)$. It follows immediately from \Cref{Lem:lien-Flamba-et-flambda+} that $f_\lambda^-$ and $F_\lambda^-$ belong to $A(\mathbb C\setminus\overline{\mathbb D})$.
Using this construction, we can now give an explicit description of the point spectra of both $T_F$ and $T_F^*=T_f$. For all $\lambda\in \sigma(T_F)\setminus F(\mathbb T)$ and $0\le k<\w_F(\lambda)$, we introduce the  functions $h_{\lambda,k}$ defined on $\D$ by
\begin{equation}\label{Annexe-defn-hlambdak}
h_{\lambda,k}(z)~=~z^k\dfrac{f_\lambda^+(0)}{f_\lambda^+(z)} \quad\text{ for every } z\in\D.
\end{equation}
These functions belong to $H^q$ (indeed, they are in $A(\mathbb D)$). Observe that this definition of $h_{\lambda,k}$ can be extended to any $\lambda\in \C\setminus F(\T)$ and any $k\ge0$, but only the $h_{\lambda,k}$ with $\lambda\in \sigma(T_F)\setminus F(\mathbb T)$ and $0\le k<\w_F(\lambda)$ are eigenvectors of $T_F^*$. 
\par\smallskip
\begin{lemma}\label{lem:description-eigenvectors}
Suppose that $F$ satisfies \ref{H1} and \ref{H3'}, and let $p>1$. Then, for every $\lambda\in \sigma(T_F)\setminus F(\mathbb T)$, we have
\[
\ker (T_F^*-\lambda)~=~\spa\,\big[h_{\lambda,k},\,0\le k<\w_F(\lambda)\big].
\]
\end{lemma}
\par\smallskip
\begin{proof}
For $\lambda\in \sigma(T_F)\setminus F(\mathbb T)$, let $0\leq k<n$, where $n=\w_F(\lambda)=-\w_f(\lambda)$.
Since $f-\lambda=f_\lambda^+z^{-n}f_\lambda^-$, and $z^{k-n}f_\lambda^-\in H^q_-=\ker(P_+)$ we have
\begin{multline*}
(T_F^*-\lambda)h_{\lambda,k}~=~(T_f-\lambda)h_{\lambda,k}~=~P_+\left[(f-\lambda)h_{\lambda,k}\right]\\
=~P_+\left[f_\lambda^+z^{-n}f_\lambda^-z^k\frac{f_\lambda^+(0)}{f_\lambda^+}\right]~=~P_+\left[f_\lambda^+(0)z^{k-n}f_\lambda^-\right]~=~0.
\end{multline*}
Hence we have $h_{\lambda,k}\in \text{ker }(T_F^*-\lambda)$ for every $0\leq k<n$. 
Since $\dim(\ker(T_F^*-\lambda))=\w_F(\lambda)=n$, this concludes the proof of \Cref{lem:description-eigenvectors}. 
\end{proof}

Our aim in the next section is to provide another representation of $f_\lambda^+$ (this is \Cref{Lemme-propriete-varphi} below), which will be crucial in order to obtain boundary relations for the eigenvectors of the operator $T_F^*$. 

\subsection{Another representation of the functions $f_\lambda^+$}\label{section-another-representationflambdaplus}
Assume that $F$ satisfies \ref{H1}, \ref{H2'} and \ref{H3'}.
By \Cref{thm:quasi-conformal-function}, the function $F$ can be extended to a $C^{1}$-smooth function $\widetilde{F}$ on an open neighborhood $\widetilde{U}$ of $\T$, whose Jacobian $J_{\widetilde{F}}(z)$ is positive at every point $z\in \widetilde{U}$, and which can be written as
 $\widetilde{F}=\varphi\circ g$, where $g$ is a quasiconformal $ C^1$-diffeomorphism from $\widetilde{U}$ onto a domain $V$ of $\mathbb C$, and $\varphi$ is a non-constant analytic function on $V$.
Let $U$ be an open subset of $\C\setminus\{0\}$ such that $\mathbb{T}  \subseteq U  \subseteq \overline{U} \subseteq \widetilde{U}$, and denote by $W$ the open set $W=g(U)$. Note that $\overline{W} \subseteq V$ because $g$ is an homeomorphism. Observe also that 
\begin{equation}\label{eq:jac-borne}
\sup_{z\in \overline{U}}|J_{\widetilde{F}}(z)|~<~\infty.
\end{equation}
\par\smallskip
Fix $\lambda_0\in F(\mathbb T)$. The equation $\varphi(z)=\lambda_0$ has a finite number of solutions $w_j$, $1\leq j\leq s$, in $W$. Moreover, the function $\widetilde{F}$ is locally a $ C^1$ diffeomorphism on $\widetilde{U}$, and thus is locally injective on $\widetilde{U}$. It follows that $\varphi=\widetilde{F}\circ g^{-1}$ is also locally injective on $V$, and thus the solutions $w_j$, $1\leq j\leq s$, are simple.  By \Cref{lem:rouche},
there exist $\alpha>0$ and one-to-one analytic maps $\widetilde{d_1},\widetilde{d_2},\dots,\widetilde{d_s}$  on  $D(\lambda_0,\alpha)$ such that for every $\lambda\in D(\lambda_0,\alpha)$, the solutions of the equation $\varphi(z)=\lambda$ in $W$ are exactly the $s$ points $\widetilde{d}_j(\lambda)$, $1\leq j\leq s$.
Define $d_j=\varphi^{-1}\circ\widetilde{d_j}$, $1\le j\le s$. For every $\lambda\in D(\lambda_0,\alpha)$ and every $1\leq j\leq s$, we have then
\[
\widetilde{F}(d_j(\lambda))~=~\lambda.
\]
\par\smallskip
\begin{notation}
For each $\lambda\in D(\lambda_0,\alpha)$,  denote by $N(\lambda)$ the set $N(\lambda)=\{1\leq j\leq s:|d_j(\lambda)|<1\}$, and denote by $n(\lambda)$ its cardinal. 
\end{notation}
If $d_j(\lambda_0)\notin \mathbb T$, taking $\alpha$ smaller if necessary, we can assume that $d_j(D(\lambda_0,\alpha))\cap\mathbb T=\varnothing$. 
\par\smallskip
Moreover, for every $j\neq k$, $d_j(\lambda_0)\neq d_k(\lambda_0)$, so we can also assume that
\begin{equation}\label{eq:sdfsfjsdlfsdj1323}
\overline{d_j(D(\lambda_0,\alpha))}\cap \overline{d_k(D(\lambda_0,\alpha))}~=~\varnothing.
\end{equation}
Denote by $J$ the set of indices $j\in\{1,\dots,s\}$ such that $d_j(\lambda_0)\in\mathbb T$. Note that $j$ belongs to $J$ if and only if  $d_j(D(\lambda_0,\alpha))\cap\mathbb T\neq \varnothing$. 

\par\smallskip
\begin{remark}\label{casse-pieds}
Let us point out the following consequence of these assumptions: suppose that $\lambda\in\mathcal{O}$ is an ``essential" self-intersection point of the curve $F(\T)$, in the sense that whatever the choice of the arcs $\alpha_j$ in assumption \ref{H2} (or \ref{H2'}), $\lambda$ is the image by $F$ of one of the extremities of one of the arcs $\alpha_j$. If $\lambda$ belongs to $D(\lambda_0,\alpha)$, then necessarily $\lambda=\lambda_0$. Indeed, if $\lambda\neq\lambda_0$, then the equation $F(z)=\lambda$ has strictly more solutions in $\T$ than the equation $F(z)=\lambda_0$ (see the proof of \Cref{exact-N-valence}). Hence there exists an index $j_0$ such that 
$d_{j_0}(\lambda_0)\notin \mathbb T$ while $d_{j_0}(\lambda)\in \mathbb T$, which is not possible.
\end{remark}
\par\smallskip
\begin{notation}\label{NotationB6}
We then define $D_{int}(\lambda_0)$ as the intersection of $D(\lambda_0,\alpha)$ with the interior connected component with respect to the arc $\Gamma_{\lambda_0}=D(\lambda_0,\alpha)\cap F(\mathbb T)$ and $D_{ext}(\lambda_0)$ as the intersection of $D(\lambda_0,\alpha)$ with the exterior connected component with respect to $\Gamma_{\lambda_0}$ (see \Cref{subsection:JFA-model} for the definition of the interior and exterior components). Note that if $\Gamma_{\lambda_0}$ is included in the boundary of exactly one component $\Omega$ of $\sigma(T_F)\setminus F(\mathbb T)$, this is both the interior and the exterior component. So, we denote by $D_{int}(\lambda_0)$ and $D_{ext}(\lambda_0)$ the two connected components of $\Omega\cap D(\lambda_0,\alpha)$.
Let 
\[
N_{int}(\lambda_0)~=~\{j\in J:d_j(D_{int}(\lambda_0)) \subseteq\mathbb D\},\quad N_{ext}(\lambda_0)~=~\{j\in J:d_j(D_{ext}(\lambda_0)) \subseteq\mathbb D\},
\]
and set also
$n_{i}(\lambda_0)=\text{card}(N_{int}(\lambda_0))$ and $n_{e}(\lambda_0)=\text{card}(N_{ext}(\lambda_0))$.
\end{notation}
\par\smallskip
In other words, we partition the set $\{1,\dots,s\}$ as follows:
\begin{center}
\begin{tikzpicture}
\draw[-](-.3,0.3)node[above]{$\{1,\dots,s\}$};
\draw[->](-.3,0.3)--(-.3,-1)node[below]{$J=\{j\,;\,d_j(\lambda_0)\in\mathbb T\}$};
\draw[->](-1,0.3)--(-4.5,-1);\draw(-5.4,-1)node[below]{$
\begin{aligned}
N(\lambda_0)&=\{j\,;\,|d_j(\lambda_0)|<1\}\\
&=\{j\,;\,d_j(D(\lambda_0,\alpha)) \subseteq\mathbb D\}
\end{aligned}
$};
\draw[->](.4,0.3)--(3.9,-1);\draw(4.35,-1)node[below]{$\begin{aligned}
&\{j\,;\,|d_j(\lambda_0)|>1\}\\
=\{j\,;\,&d_j(D(\lambda_0,\alpha)) \subseteq\mathbb C\setminus\overline{\mathbb D}\}
\end{aligned}$};
\draw[->](-.8,-1.8)--(-4,-3.5);\draw(-4.6,-3.5)node[below]{$\begin{aligned}
N_{int}(\lambda_0)
&=\{j\in J\,;\, d_j(D_{int}(\lambda_0)) \subseteq\mathbb D\}\\
&=\{j\in J\,;\, d_j(D_{ext}(\lambda_0)) \subseteq\mathbb C\setminus\overline{\mathbb D}\}
\end{aligned}$};
\draw[->](.2,-1.8)--(3.5,-3.5);\draw(3.2,-3.5)node[below]{$\begin{aligned}
N_{ext}(\lambda_0)
&=\{j\in J\,;\, d_j(D_{ext}(\lambda_0)) \subseteq\mathbb D\}\\
&=\{j\in J\,;\, d_j(D_{int}(\lambda_0)) \subseteq\mathbb C\setminus\overline{\mathbb D}\}
\end{aligned}$};
\end{tikzpicture}
\end{center}
\par\smallskip
In particular, we have 
\begin{equation}\label{eq:N343}
N(\lambda)~=~N(\lambda_0)\cup N_{int}(\lambda_0)\quad \text{ for every }\lambda\in D_{int}(\lambda_0)
\end{equation}
and
\begin{equation}\label{eq:N3432}
N(\lambda)~=~N(\lambda_0)\cup N_{ext}(\lambda_0)\quad \text{ for every }\lambda\in D_{ext}(\lambda_0).
\end{equation}\par\smallskip
\begin{lemma}\label{lemme:relation-wind-int-ext}
The definitions of $n_i(\lambda_0)$ and $n_e(\lambda_0)$ given in \Cref{NotationB6} coincide with the ones introduced in \Cref{section:JFA}, and moreover we have
\begin{equation}\label{eq:LemmeB7}
w_i(\lambda_0)-n_i(\lambda_0)~=~w_e(\lambda_0)-n_e(\lambda_0) \quad\text{ for every } \lambda_0\in F(\mathbb T)\setminus\mathcal O,
\end{equation}
where $w_i(\lambda_0)$ and $w_e(\lambda_0)$ denote the interior and exterior limits of  $\w_F$ at the point $\lambda_0$.
In particular, the function $\lambda\mapsto\w_F(\lambda)-n(\lambda)$ is constant on $D(\lambda_0,\alpha)\setminus F(\mathbb T)$. 
\end{lemma}

\begin{proof}
Since the functions $d_j$ are quasi-conformal, they preserve the orientation. Thus the integers $n_{i}(\lambda_0)$ and $n_{e}(\lambda_0)$ (as defined in \Cref{NotationB6}) can be linked  to the number of crossings at $\lambda_0$ in either direction along the curve $F(\mathbb T)$. Indeed, since when traveling along the unit circle in the positive orientation, the disk is always located on the left side, it follows that if $j$ belongs to  $N_{int}(\lambda_0)$, then $D_{int}(\lambda_0)$ will remain on the left side of the portion $F(\gamma_j)$ of the curve $F(\T)$, where $\gamma_j$ is a small arc in $\T$ centered at  $d_j(\lambda_0)$. In other words, $n_i(\lambda_0)$ corresponds to the number of times we cross $\lambda_0$ while keeping the interior component on the left side (and the exterior component on the right side). Similarly, $n_e(\lambda_0)$ corresponds to the number of times we cross $\lambda_0$ while keeping the exterior component on the left side (and the interior component on the right side). \Cref{figure-explication-ne-ni} below gives an example where $n_i=3$ and $n_e=1$.
\begin{figure}[ht]
\includegraphics[page=51,scale=1.1]{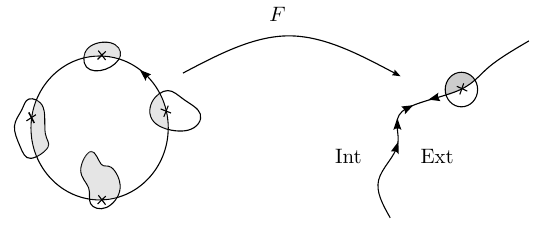}
            \caption{$n_i=3,~n_e=1$}
            \label{figure-explication-ne-ni}
\end{figure}
\par\smallskip
Recall the geometrical interpretation of the winding number (mentioned for instance in \cite{Yger2001}) that we also used in \Cref{subsection:JFA-model}: adding $1$ to the winding number of a curve at a point can be interpreted as turning one more time around the point while keeping it on the left. So, we have 
\[
n_i-n_e~=~w_i-w_e\quad\text{ on } F(\T)\setminus\mathcal{O}.
\]
Since $F$ satisfies \ref{H3'}, we have that $w_i\ge w_e\ge 0$ on $F(\T)\setminus\mathcal{O}$, and thus, $n_i\ge n_e$, proving that the two definitions of \Cref{NotationB6} and \Cref{section:JFA} coincide.
\par\smallskip
Let us now prove the second assertion of \Cref{lemme:relation-wind-int-ext}.
Since $\w_F$ is constant on connected components of $\sigma(T_F)\setminus F(\mathbb T)$ and since, according to \Cref{eq:N343,eq:N3432}, $\lambda\mapsto n(\lambda)$ is constant on $D_{int}$ as well as  on $D_{ext}$, it is sufficient to prove that given $\lambda_i\in D_{int}$ and $\lambda_e\in D_{ext}$, we have
\begin{equation}\label{EqB9}
\w_F(\lambda_i)-n(\lambda_i)~=~\w_F(\lambda_e)-n(\lambda_e).
\end{equation}
By \Cref{eq:N343,eq:N3432}, we have 
\[
n_i(\lambda_0)~=~n(\lambda_i)-n(\lambda_0)\quad\text{and}\quad n_e(\lambda_0)~=~n(\lambda_e)-n(\lambda_0).
\]
\Cref{EqB9} follows immediately by replacing $n_i(\lambda_0)$ and $n_e(\lambda_0)$ in \Cref{eq:LemmeB7}.
\end{proof}
Let us now prove the key lemma of this section:

\begin{lemma}\label{Lemme-propriete-varphi}
Suppose that $F$ satisfies \ref{H1}, \ref{H2'} and \ref{H3'}. Let $\lambda_0\in F(\mathbb T)$. With $\alpha=\alpha(\lambda_0)$ as above, there exists a function $\varphi$ defined on $D(\lambda_0,\alpha)\times \overline{\mathbb{D}}$ such that
\begin{enumerate}[(i)]
\item  for every $\gamma,\beta>0$ such that $\gamma+\beta=\varepsilon$, the functions $\varphi_l:\lambda\mapsto \varphi(\lambda,\cdot)$ and $1/\varphi_l$ are $C^\gamma$ functions from $D(\lambda_0,\alpha)$ into $C^\beta(\overline{\mathbb D})$;
\par\smallskip
\item if $\lambda\in D(\lambda_0,\alpha)\setminus F(\mathbb{T})$, then
\[
f_\lambda^+(z)~=~\varphi(\lambda,z)\prod_{j\in N(\lambda)}(z-d_j(\lambda)^{-1})\quad\text{ for every }z\in\overline{\mathbb{D}}.
\]
\end{enumerate}
\end{lemma}

\begin{proof}
$(i)$ Let $r_0:=\w_F(\lambda_1)+s-n(\lambda_1)$, where $\lambda_1$ is any element of $D(\lambda_0,\alpha)\setminus F(\mathbb{T})$. Note that, by \Cref{lemme:relation-wind-int-ext}, $r_0$ is well defined and does not depend of the point $\lambda_1\in D(\lambda_0,\alpha)\setminus F(\mathbb{T})$. Observe that for every $\lambda\in D(\lambda_0,\alpha)$ and every $1\leq j\leq s$, we have $d_j(\lambda)\in U$. In particular, $d_j(\lambda)\neq 0$.  Now, define a function $\psi$ on  $ D(\lambda_0,\alpha)\times \mathbb T$ by setting, for every $(\lambda,z)\in D(\lambda_0,\alpha)\times \mathbb T$
\[
\psi(\lambda,z)~=~
 z^{r_0}(f(z)-\lambda)\prod_{j=1}^s(z-d_j(\lambda)^{-1})^{-1}\quad\text{ if } z\neq \frac{1}{d_j(\lambda)}\text{ for every }1\le j\le s
 \]
 and
 \[ \psi(\lambda,z)~=~
 z^{r_0}f'(d_l(\lambda)^{-1})\prod_{\substack{j=1\\j\neq l}}^s(z-d_j(\lambda)^{-1})^{-1}\quad\text{ if } z=\frac{1}{d_l(\lambda)}\text{ for some }1\le l\le s.
\]
The function $\psi$ does not vanish on $D(\lambda_0,\alpha)\times \mathbb T$. Since the functions $d_j$ are of class $ C^1$ on $D(\lambda_0,\alpha)$ and $f$ is of class $C^{1+\varepsilon}$ on $\mathbb T$, it is not difficult to check that the function $\psi$ is of class $C^{\varepsilon}$ on $D(\lambda_0,\alpha)\times \mathbb T$. Observe also that $\w_{\psi(\lambda,\cdot)}(0)=0$ for every $\lambda\in D(\lambda_0,\alpha)$. Indeed, for every 
$\lambda\in D(\lambda_0,\alpha)\setminus F(\T)$, none of the points $d_j(\lambda)^{-1}$, $1\le j\le s$, belongs to $\T$ (since $\widetilde F(d_j(\lambda))=\lambda$), and thus
\[\w_{\psi(\lambda,\cdot)}(0)~=~r_0+\w_f(\lambda)+\sum_{j=1}^s\w_{(z-d_j(\lambda)^{-1})^{-1}}(0).\]
Now,
\[\w_{(z-d_j(\lambda)^{-1})^{-1}}(0)~=~\begin{cases}
0&\text{ if } |d_j(\lambda)|<1\\
-1&\text{ if } |d_j(\lambda)|>1,\\
\end{cases}\]
and hence
\[\w_{\psi(\lambda,\cdot)}(0)~=~r_0+\w_f(\lambda)-(s-n(\lambda)).\]
Now, by \Cref{lemme:relation-wind-int-ext}, $r_0=\w_F(\lambda)+s-n(\lambda)$, and it follows that $\w_{\psi(\lambda,\cdot)}(0)=0$. This equality being true for every $\lambda\in D(\lambda_0,\alpha)\setminus F(\T)$, it remains true for $\lambda\in D(\lambda_0,\alpha)\cap F(\T)$ by continuity.
Using this, one can  construct as above a determination of the logarithm of $\psi$, denoted by $\log\psi$, which is $C^\varepsilon$ on $D(\lambda_0,\alpha)\times \mathbb T$. For every $\lambda\in D(\lambda_0,\alpha)$, set $v_\lambda=\log(\psi(\lambda,\cdot))$. According to \Cref{lemme-Cepsilon}, $v_\lambda\in C^{\beta}(\mathbb T)$, and
$\|v_\lambda-v_\mu\|_{C^{\beta}(\mathbb T)}
 \lesssim |\lambda-\mu|^{\gamma}$.
Set
\[
\varphi(\lambda,\cdot)~=~e^{P_+v_\lambda}.
\]
Recall that $P_+$ is a bounded map from $C^{\beta}(\mathbb T)$ into itself (see \Cref{thm:privalov-Zygmund}). 
For every $\lambda,\mu\in D(\lambda_0,\alpha)$, we have
\[
 \|P_+v_\lambda-P_+v_\mu\|_{C^{\beta}(\mathbb T)}
 ~\lesssim~ \|v_\lambda-v_\mu\|_{C^{\beta}(\mathbb T)}
 ~\lesssim~ |\lambda-\mu|^{\gamma}.
\]
Therefore, we conclude that $\varphi_l$ is of class $C^{\gamma}$ from $D(\lambda_0,\alpha)$ into $ C^\beta(\mathbb T)$, and so into $C^\beta(\overline{\mathbb D})$ by the Maximum Principle.  Since $\varphi$ does not vanish, the function $1/\varphi_l$ is also $C^{\gamma}$ from $D(\lambda_0,\alpha)$ into $ C^\beta(\overline{\mathbb D})$.
\par\smallskip
$(ii)$ Let $\lambda\in D(\lambda_0,\alpha)\setminus F(\mathbb T)$. Then, for any $0\le j \le s$, $d_j(\lambda)\notin \mathbb T$. So for every $z\in\mathbb T$, we have 
\begin{eqnarray*}
z^{-\w_f(\lambda)}(f(z)-\lambda)&=&z^{\w_F(\lambda)-r_0}\psi(\lambda,z)\prod_{j=1}^s (z-d_j(\lambda)^{-1})\\
&=&z^{n(\lambda)-s}\psi(\lambda,z)\prod_{j=1}^s (z-d_j(\lambda)^{-1})\\
&=&z^{n(\lambda)-s}\psi(\lambda,z)\prod_{j\in N(\lambda)}(z-d_j(\lambda)^{-1}) \prod_{j\notin N(\lambda)}(z-d_j(\lambda)^{-1})\\
&=&\psi(\lambda,z)\prod_{j\in N(\lambda)}(z-d_j(\lambda)^{-1}) \prod_{j\notin N(\lambda)}(1-z^{-1}d_j(\lambda)^{-1}).
\end{eqnarray*}
For $\lambda\in D(\lambda_0,\alpha\setminus F(\mathbb T)$ and $z\in\mathbb T$, let
\[
g_\lambda(z)~=~\log\psi(\lambda,z)~+\sum_{j\in N(\lambda)}\log(z-d_j(\lambda)^{-1})~+\sum_{j\notin N(\lambda)}\log(1-z^{-1}d_j(\lambda)^{-1}).
\]
Then the function $g_\lambda$ is continuous on $\mathbb T$ and $$e^{g_\lambda(z)}~=~z^{-\w_f(\lambda)}(f(z)-\lambda)\quad\text{for every }z\in\mathbb T.$$ 
This implies that $f_\lambda^+=e^{P_+g_\lambda}$ (as mentioned above, the definition of $f_\lambda^+$ as
$f_\lambda^+=e^{P_+u_\lambda}$ does not depend of the continuous branch of the logarithm of
$z^{-m}(f-\lambda)$ on $\T$ that we consider). Observe that for every $j\in N(\lambda)$, $\log(z-d_j(\lambda)^{-1})$ belongs to $H^p$, whence it follows that $P_+\log(z-d_j(\lambda)^{-1})=\log(z-d_j(\lambda)^{-1})$. On the other hand, for every $j\notin N(\lambda)$, $\log(1-z^{-1}d_j(\lambda)^{-1})$ belongs to $H^p_-$, whence $P_+\log(1-z^{-1}d_j(\lambda)^{-1})=0$. Thus, for every $z\in\overline{\mathbb D}$,
\begin{eqnarray}
f_\lambda^+(z)&=& e^{P_+\log\psi(\lambda,z)}\exp\left(\sum_{j\in N(\lambda)}\log(z-d_j(\lambda)^{-1})\right)\\
&=& \varphi(\lambda,z) \prod_{j\in N(\lambda)}(z-d_j(\lambda)^{-1}),
\end{eqnarray}
which proves $(ii)$.
\end{proof}
\begin{remark}\label{remark:bounded-flambdaplus}
It follows from \Cref{Lemme-propriete-varphi} that the function $(\lambda,z)\longmapsto f_\lambda^+(z)$ is bounded on $D(\lambda_0,\alpha)\times \overline{\mathbb D}$. 
\end{remark}

\subsection{Boundary relations for eigenvectors}
Our aim here is to prove a boundary relation between the interior and exterior values of the eigenvectors. We start with the functions $f_\lambda^+$. See \Cref{subsection:JFA-model} for the definitions of the interior and exterior boundary values.
\begin{lemma}\label{lemme:condition-aux-bords}
Let $F$ be a symbol which satisfies \ref{H1}, \ref{H2'} and \ref{H3'}. For every $z\in\mathbb D$, the function $\lambda\longmapsto f_{\lambda}^+(z)^{-1}$ admits interior and exterior boundary values almost everywhere on $F(\mathbb T)$, satisfying the following relation:
\begin{equation}\label{eq1:relation-aux-bords}
\left(\prod_{j\in N_{int}(\lambda)}(z-d_j(\lambda)^{-1})\right){f_{\lambda,int}^+(z)}^{-1}~=~\left(\prod_{j\in N_{ext}(\lambda)}(z-d_j(\lambda)^{-1})\right){f_{\lambda,ext}^+(z)}^{-1}.
\end{equation}
In particular, if $F$ satisfies \ref{H2} rather than the more general assumption \ref{H2'}, we have
\begin{equation}\label{eq2:relation-aux-bords}
(z-\xi(\lambda)^{-1}){f^+_{\lambda,int}(z)}^{-1}~=~{f^+_{\lambda,ext}(z)}^{-1}\quad\text{for almost every } \lambda\in F(\mathbb T),
\end{equation}
where $\xi$ is defined as $\xi=F^{-1}$ on $ F(\mathbb T)\setminus\mathcal O$.
\end{lemma}

In order to simplify the notation, we write
\(f_{\lambda,int}^+\) rather than \((f_\lambda^+)^{int}\) and 
\(f_{\lambda,ext}^+\) rather than \((f_\lambda^+)^{ext}\).

\begin{proof}
Let $\lambda\in F(\mathbb T)\setminus \mathcal O$. 
Since for $\lambda_i\in D_{int}(\lambda)$ and $\lambda_e\in D_{ext}(\lambda)$, we have $N(\lambda_i)=N(\lambda)\cup N_{int}(\lambda)$ and $N(\lambda_e)=N(\lambda)\cup N_{ext}(\lambda)$ (see \Cref{eq:N343,eq:N3432}),  we can write
\[
f_{\lambda_i}^+(z)~=~\varphi(\lambda_i,z)\prod_{j\in N(\lambda)}(z-d_j(\lambda_i)^{-1})\prod_{j\in N_{int}(\lambda)}(z-d_j(\lambda_i)^{-1})\quad\text{ for every }z\in\D,
\]
and similarly for $f_{\lambda_e}^+$,
\[
f_{\lambda_e}^+(z)~=~\varphi(\lambda_e,z)\prod_{j\in N(\lambda)}(z-d_j(\lambda_e)^{-1})\prod_{j\in N_{ext}(\lambda)}((z-d_j(\lambda_e)^{-1})\quad\text{ for every }z\in\D.
\]
According to \Cref{Lemme-propriete-varphi}, the function $\varphi(\cdot,z)$ is continuous on $D(\lambda_0,\alpha)$ as well as the functions $d_j$. Hence, letting $\lambda_i$ and $\lambda_e$ tend to $\lambda$ in the equalities above gives that, for almost all $\lambda\in F(\mathbb T)$ and all $z\in\mathbb D$, we have 
\[
f_{\lambda,int}^+(z)~=~\varphi(\lambda,z)\prod_{j\in N(\lambda)}(z-d_j(\lambda)^{-1})\prod_{j\in N_{int}(\lambda)}(z-d_j(\lambda)^{-1})
\]
and 
\[
f_{\lambda,ext}^+(z)~=~\varphi(\lambda,z)\prod_{j\in N(\lambda)}(z-d_j(\lambda)^{-1})\prod_{j\in N_{ext}(\lambda)}(z-d_j(\lambda)^{-1}).
\]
Comparing these two relations gives \Cref{eq1:relation-aux-bords}.
\par\smallskip
Now assume that $F$ satisfies \ref{H2}, and let $\lambda\in F(\mathbb T)\setminus \mathcal O$. In this case the set $J$ is reduced to one point, which we call $j_0$, and $N_{int}(\lambda)=\{j_0\}$ while $N_{ext}(\lambda)=\varnothing$. Then $\xi(\lambda)=d_{j_0}(\lambda)$ and thus \Cref{eq1:relation-aux-bords} implies \Cref{eq2:relation-aux-bords}.\end{proof}

Recalling that $N=\max\{\w_F(\lambda)\,;\,\lambda\in\C\setminus F(\T)\}$, we set
 \[
 \Omega_j^+~=~\{\lambda\in \mathbb C\setminus F(\mathbb T)\,;\, \w_F(\lambda)>j\}\quad\text{ for every }
0\le j\le N-1.\]

As a direct consequence of \Cref{lemme:condition-aux-bords}, we obtain:
\begin{corollary}\label{cor:condition-aux-bords-hlambdaj}
Let $F$ satisfy \ref{H1}, \ref{H2} and \ref{H3'}. Let $j\geq 0$. Then, for every $z\in\mathbb D$, for almost every $\lambda\in \overline{\Omega_{j+1}^+}\cap F(\mathbb T)=\partial\Omega_{j+1}^+$, we have
\begin{equation}\label{eq:3434RDCS}
h_{\lambda,j}^{int}(z)~-~\xi(\lambda)h_{\lambda,j+1}^{int}(z)~=~h_{\lambda,j}^{ext}(z).
\end{equation}

\end{corollary}
Note that the definition of the $h_{\lambda,j}$ by $h_{\lambda,j}(z)=z^jf_\lambda^+(0)/f_\lambda^+(z)$ can be extended for every $j\ge0$ and \textit{every} $\lambda\in\mathbb C\setminus F(\mathbb T)$. Then \Cref{eq:3434RDCS} remains true for every $j\ge0$ and \textit{for  every} $\lambda\in F(\mathbb T)\setminus\mathcal O$.  
\begin{proof}
Using \Cref{eq2:relation-aux-bords}, for $\lambda\in F(\mathbb T)\setminus \cal O$ and $z\in\mathbb D$, we have
\begin{eqnarray*}
h_{\lambda,j}^{ext}(z)&=&z^j\frac{f_{\lambda,ext}^+(0)}{f_{\lambda,ext}^+(z)} ~=~z^j \frac{f_{\lambda,int}^+(0)}{-\xi(\lambda)^{-1}}\,\frac{z-\xi(\lambda)^{-1}}{f_{\lambda,int}^+(z)} \\
&=&z^j(1-\xi(\lambda)z)\frac{f_{\lambda,int}^+(0)}{f_{\lambda,int}^+(z)}\\
&=&z^j \frac{f_{\lambda,int}^+(0)}{f_{\lambda,int}^+(z)}-\xi(\lambda) z^{j+1}\frac{f_{\lambda,int}^+(0)}{f_{\lambda,int}^+(z)}\\
&=&h_{\lambda,j}^{int}(z)-\xi(\lambda)h_{\lambda,j+1}^{int}(z).
\end{eqnarray*}
\end{proof}

\begin{remark}\label{remark-conditions-bords-vect-H2prime3434A4}
Under the hypothesis \ref{H2'}, the functions  $h_{\lambda,j}$ also satisfy a boundary relation which can be written down rather explicitly. Indeed, we have as above $$h_{\lambda,k}^{ext}(z)~=~z^k f_{\lambda,ext}^+(0)f_{\lambda,ext}^+(z)^{-1} \quad\text{ and }\quad h_{\lambda,k}^{int}(z)~=~z^k f_{\lambda,int}^+(0)f_{\lambda,int}^+(z)^{-1}$$
for almost every $\lambda\in F(\mathbb T)$.
But according to \Cref{eq1:relation-aux-bords}, we have
\[
\frac{\left(\displaystyle\prod_{j\in N_{int}(\lambda)}(z-d_j(\lambda)^{-1})\right){f_{\lambda,int}^+(z)}^{-1}}{\left(\displaystyle\prod_{j\in N_{int}(\lambda)}(-d_j(\lambda)^{-1})\right){f_{\lambda,int}^+(0)}^{-1}}~=~\frac{\left(\displaystyle\prod_{j\in N_{ext}(\lambda)}(z-d_j(\lambda)^{-1})\right){f_{\lambda,ext}^+(z)}^{-1}}{\left(\displaystyle\prod_{j\in N_{ext}(\lambda)}(-d_j(\lambda)^{-1})\right){f_{\lambda,ext}^+(0)}^{-1}},
\]
from which it follows that for almost every $\lambda\in F(\T)$ and every $z\in\mathbb D$
\[
\left(\prod_{j\in N_{int}(\lambda)}(1-d_j(\lambda)z)\right)h_{\lambda,k}^{int}(z)~=~\left(\prod_{j\in N_{ext}(\lambda)}(1-d_j(\lambda)z)\right)h_{\lambda,k}^{ext}(z).
\]
Let us now introduce coefficients $a_p(\lambda)$ and $b_\ell(\lambda)$ such that 
\begin{equation}\label{eq:polynome1-condition-aux-bords}
\prod_{j\in N_{int}(\lambda)}(1-d_j(\lambda)z)~=~\sum_{p=0}^{n_i(\lambda)}a_p(\lambda)z^p
\end{equation}
and
\begin{equation}\label{eq2:polynome1-condition-aux-bords}
\prod_{j\in N_{ext}(\lambda)}(1-d_j(\lambda)z)~=~\sum_{\ell=0}^{n_{e}(\lambda)}b_\ell(\lambda)z^\ell.
\end{equation}
Using the fact that $z^p h_{\lambda,k}^{int}(z)=h_{\lambda,k+p}^{int}(z)$ and $z^p h_{\lambda,k}^{ext}(z)=h_{\lambda,k+p}^{ext}(z)$ for every $k,p\ge 0$, we obtain a boundary relation of the form
\begin{equation}\label{eq:relation-aux-bords-H2prime34455SD}
\sum_{p=0}^{n_{i}(\lambda)}a_p(\lambda)h_{\lambda,k+p}^{int}(z)~=\sum_{\ell=0}^{n_{e}(\lambda)}b_\ell(\lambda)h_{\lambda,k+\ell}^{ext}(z)\quad \text{ for all } z\in\D ~\text{and a.e.}~\lambda\in F(\mathbb T).
\end{equation}
\end{remark}
\par\smallskip
Proceeding in a similar way, we can show that the functions $F_\lambda^+$ defined in \Cref{Eq:Flambda+-} satisfy the boundary relation
\begin{equation}\label{eq23ZE:relation-aux-bords}
\prod_{j\in N_{ext}(\lambda)}(z-d_j(\lambda)){F_{\lambda,int}^+(z)}^{-1}~=\prod_{j\in N_{int}(\lambda)}(z-d_j(\lambda)){F_{\lambda,ext}^+(z)}^{-1}\quad\text{ a.e. on } F(\T). 
\end{equation}
In particular, when $F$ satisfies \ref{H2}, we have for almost every $\lambda\in F(\T)$ that
\begin{equation}\label{eq33434A4:relation-aux-bords}
F_{\lambda,int}^+(z)^{-1}~=~(z-\xi(\lambda))F_{\lambda,ext}^+(z)^{-1}
\quad\text{ for every } z\in\D.
\end{equation}
\par\smallskip

\subsection{Eigenvectors on the curve $F(\mathbb T)$}
In this subsection, we present an application of \Cref{Lemme-propriete-varphi} to the construction of eigenvectors of $T_F^*$ associated to some eigenvalues on the curve $F(\mathbb T)$. For any connected component $\Omega$ of $\sigma(T_F)\setminus F(\mathbb T)$, we denote by $\gamma_\Omega$ the union of the subarcs of $F(\mathbb T)\setminus \mathcal O$ which have $\Omega$ as the exterior component. The following proposition is a direct consequence of \Cref{Lemme-propriete-varphi}.

\begin{proposition}\label{Prop:ContinuiteBord}
Let $F$ satisfy \ref{H1}, \ref{H2} and \ref{H3'}. Let $\Omega$ be a connected component of $\sigma(T_F)\setminus F(\mathbb T)$, $j\ge0$ and let $\Phi_j^\Omega:\Omega\cup\gamma_\Omega\longrightarrow A(\mathbb D)$ be the function defined by
\[\Phi_j^\Omega(\lambda)=\begin{cases}
h_{\lambda,j}&\text{ if }\lambda\in\Omega\\
h_{\lambda,j}^{ext}&\text{ if }\lambda\in\gamma_\Omega.
\end{cases}\]
Then $\Phi_j^\Omega$ is well-defined and continuous on $\Omega\cup\gamma_\Omega$.
\end{proposition}

\begin{proof}
We will use in this proof the notation of the previous subsections.
It is clear that $\Phi_j^\Omega$ is well-defined and continuous on $\Omega$. Let $\lambda_0\in\gamma_\Omega$. 
Since $F$ satisfies \ref{H2}, it follows that the set $J$ is reduced to one point $j_0$. Since $d_{j_0}$ preserves local orientations, it follows from the definition of $\gamma_\Omega$ that 
    \[d_{j_0}(\Omega\cap D(\lambda_0,\alpha))\subset \mathbb C\setminus\overline{\mathbb D}.\]
    So $N(\lambda)=N(\lambda_0)$ for every $\lambda \in\Omega\cap D(\lambda_0,\alpha)$, and thus for all such $\lambda$'s we have
    \[h_{\lambda,j}(z)=z^j\frac{f_\lambda^+(0)}{f_\lambda^+(z)}=z^j\frac{\varphi(\lambda,0)}{\varphi(\lambda,z)}\prod_{j\in N(\lambda_0)}(1-d_j(\lambda)z)^{-1}~\text{for all}~z\in \mathbb D\]
by \Cref{Lemme-propriete-varphi}.
Fix now $\lambda_1\in D(\lambda_0,\alpha)\cap \gamma_\Omega$. 
Then $d_j(\lambda)\longrightarrow d_j(\lambda_1)$ when $\lambda\to\lambda_1$, $\lambda \in\Omega\cap D(\lambda_0,\alpha)$, and assertion (i) of \Cref{Lemme-propriete-varphi} implies that $\|\varphi_l(\lambda)-\varphi_l(\lambda_1)\|_{\infty}\longrightarrow0$ and $\|1/\varphi_l(\lambda)-1/\varphi_l(\lambda_1)\|_{\infty}\longrightarrow0$ when $\lambda\to\lambda_1$.  
Since \[h_{\lambda_1,j}^{ext}(z)=z^j\frac{\varphi(\lambda_1,0)}{\varphi(\lambda_1,z)}\prod_{j\in N(\lambda_0)}(1-d_j(\lambda_1)z)^{-1}~\text{for all}~z\in \mathbb D,\]
it follows that $h_{\lambda_1,j}^{ext}$ belongs to $A(\mathbb D)$ and that $\|h_{\lambda,j}-h_{\lambda_1,j}^{ext}\|_\infty\longrightarrow0$ when $\lambda\to\lambda_1$, $\lambda\in\Omega$.
We conclude that $\phi_j^\Omega$ is continuous at the point $\lambda_1$ for every $\lambda_1\in D(\lambda_0,\alpha)\cap \gamma_\Omega$, and this terminates  the proof of \Cref{Prop:ContinuiteBord}. 
\end{proof}
As a direct consequence, we recover that every $\lambda_0\in F(\mathbb T)\setminus \partial\sigma(T_F)$ which is not a self-intersection point is an eigenvalue for $T_F^*$ (which is a consequence of a result of Ahern and Clark \cite{AhernClark1985}), and that the functions $h_{\lambda_0,j}^{ext}$ with $0\le j<w_e(\lambda_0)$ are associated eigenvectors. More precisely, we have the following result:

\begin{corollary}Let $1<p<\infty$, let $F$ satisfy \ref{H1}, \ref{H2} and \ref{H3'} and let $\Omega$ be a connected component of $\sigma(T_F)\setminus F(\mathbb T)$.
    For every  $\lambda_0\in\gamma_\Omega$, we have
    \[\spa\Big[h_{\lambda_0,j}^{ext},\,0\le j\le\w_F(\Omega)\Big]\subseteq \ker(T_F^*-\lambda_0).\]
\end{corollary}
\begin{proof}
Since the functions $h_{\lambda,j}$  converge to $ h_{\lambda_0,j}^{ext}$ in $A(\mathbb D)$ as $\lambda\to\lambda_0$, $\lambda\in\Omega$, by \Cref{Prop:ContinuiteBord}, the convergence also holds in $H^q$, and thus
\[\lambda_0h_{\lambda_0,j}^{ext}=\lim_{\underset{\lambda\in\Omega}{\lambda\to\lambda_0}}\lambda h_{\lambda,j}=\lim_{\underset{\lambda\in\Omega}{\lambda\to\lambda_0}}T_F^* h_{\lambda,j}=T_F^*h_{\lambda_0,j}^{ext}.\qedhere\]
\end{proof}

\subsection{The operators $U$ and $V$}
\subsubsection{Definition and continuity of $U$}\label{Subsection:ContinuiteU}
Recall that $N=\max\{\w_F(\lambda)\,;\,\lambda\in\C\setminus F(\T)\}$, and let $0\le j<N$. Given a function $g\in H^p$, denote by $u_j$ the function defined on $\Omega_j^+$ by
\begin{equation}\label{eq:B42}
u_j(\lambda)~=~\dual{g}{h_{\lambda,j}}~=~\frac{1}{2\pi}\int_0^{2\pi}g(e^{i\theta})h_{\lambda,j}(e^{-i\theta})\,\mathrm{d}\theta\quad\textrm{ for every } \lambda\in \Omega_j^+,
\end{equation}
and set
\[
Ug~=~(u_j)_{0\leq j\leq N-1}.
\]

The main concern of this section is to show that the operator $U$ is the isomorphism which implements the model for the operator $T_F$ (\Cref{model}). 
We first state a few elementary properties of $U$.
\begin{proposition}\label{petites proprietes}
Let $F$ satisfy \ref{H1}, \ref{H2'} and \ref{H3'}. Then the operator $U$ defined above satisfies the following properties:
\begin{enumerate}
    \item for every $g\in H^p$, $0\le j\le N-1$ and $\lambda\in \Omega_j^+$, we have $(UT_Fg)_j(\lambda)=\lambda (Ug)_j(\lambda)$;
    \item let $g\in H^p$ and $\lambda\in\sigma(T_F)\setminus F(\mathbb T)$. Then $g$ vanishes on $\ker (T_f-\lambda)$ if and only if $(Ug)_j(\lambda)=0$ for every $0\le j<\w_F(\lambda)$;
    \item let $g\in H^p$ and $0\le j\le N-1$. Then the function $u_j=(Ug)_j$ is analytic on $\Omega_j^+$;
    \item let $z\in\mathbb D,~0\le j\le N-1$ and $\lambda\in\Omega_j^+$. Then $h_{\lambda,j}(z)=(Uk_{\overline z})_j(\lambda)$.
\end{enumerate}
\end{proposition}
\begin{proof}
Let $u_j(\lambda)=(Ug)_j(\lambda)=\dual{g}{h_{\lambda,j}}$. It follows from \Cref{lem:description-eigenvectors} that for every $\lambda\in \Omega_j^+$, we have
\[
(UT_Fg)_j(\lambda)~=~\dual{T_Fg}{h_{\lambda,j}}~=~\dual{g}{T_F^*h_{\lambda,j}}~=~\lambda \dual{g}{h_{\lambda,j}}~=~\lambda u_j(\lambda),
\]
which gives (1). Assertion
(2) follows immediately from the fact that $\ker (T_F^*-\lambda)=\spa\,[h_{\lambda,k}\,;\,0\le k<\w_F(\lambda)]$ (see \Cref{lem:description-eigenvectors}). Assertion
(3) is a consequence of the fact that $\lambda\mapsto h_{\lambda,j}$ is analytic on $\Omega_j^+$ and of the integral representation \eqref{eq:B42}.
Lastly, assertion
(4) follows from \eqref{eq:cauchy-kernel}. 
\end{proof}

A first consequence of  \Cref{petites proprietes} is that the range of $U$ is included in $\bigoplus \text{Hol}(\Omega_j^+)$ ($\text{Hol}(\Omega_j^+)$ being the set of analytic functions on $\Omega_j^+$). Note that, by \Cref{Lemme-propriete-varphi}, the function $\lambda\mapsto h_{\lambda,j}(z)$ belongs to  $A(\Omega)$, for every  connected component $\Omega$ of $\Omega_j^+$ and every $z\in\mathbb D$ (recall that $A(\Omega)$ is the space of analytic functions on $\Omega$ with a continuous extension to $\overline\Omega$). In particular, for every $0\le j<N$, and every $z\in\mathbb D$, the function $\lambda\mapsto h_{\lambda,j}(z) $ belongs to $E^p(\Omega_j^+)$. This means that for every $z\in\mathbb D$, the $N$-tuple $Uk_z$ lies in the space $\bigoplus E^p(\Omega_j^+)$. In the next theorem, we prove that all of the range of $U$ is included in $\bigoplus E^p(\Omega_j^+)$, and that $U$ is continuous from $H^p$ into this space.

\begin{theorem}\label{lemme-Uborne}
Let $F$ satisfy \ref{H1}, \ref{H2'} and \ref{H3'}, and let $p>1$. Then the linear map $U$ is well-defined and bounded from $H^p$ into 
$\bigoplus\limits_{0\leq j\leq N-1}E^p(\Omega_j^+)$.
\end{theorem}
The strategy of the proof is to show first that the range of $U$ is included in $\bigoplus E^p(\Omega_j^+)$, and then the continuity of $U$ will follow from the Closed Graph Theorem. 

\begin{proof}
Let $g\in H^p$. We wish to show that for every $0\le j\le N-1$, the function $u_j=(Ug)_j$ belongs to $E^p(\Omega_j^+)$. We know that it is analytic on $\Omega_j^+$. Going back to the definition of $E^p(\Omega_j^+)$, we see that we need to show that for every connected component $\Omega$ of $\mathbb C\setminus F(\mathbb T)$ contained in $\Omega_j^+$,  there exists a sequence of rectifiable Jordan curves $(\gamma_n)_n$ contained in $\Omega_j^+$ and tending to $\partial \Omega$ such that  
\[
\sup_{n\geq 1}\int_{\gamma_n}|u_j(\lambda)|^p|\,|\mathrm{d}\lambda|~<~\infty.    
\]
Observe that 
\[
\int_{\gamma_n}|u_j(\lambda)|^p\,|\mathrm{d}\lambda|~=~\int_{\gamma_n}|\dual[\big]{g}{z^j h_{\lambda,0}}|^p\,|\mathrm{d}\lambda|~=~\int_{\gamma_n}|\dual[\big]{P_+(z^{-j} g)}{h_{\lambda,0}}|^p\,|\mathrm{d}\lambda|.
\]
Since $\Omega_j^+ \subseteq \Omega_0^+$, it is sufficient to check that whatever the choice of the function $g\in H^p$, the associated function $u_0$ belongs to $E^p(\Omega_0^+)$. Namely, we need to show that for every connected component $\Omega$ of $\mathbb C\setminus F(\mathbb T)$ contained in $\Omega_0^+$, there exists a sequence of rectifiable Jordan curves $(\gamma_n)_n$ contained in $\Omega_0^+$ and tending to $\partial \Omega$ such that  
\begin{equation}\label{eq1:Smirnov-bounded-U}
    \sup_{n\geq 1}\int_{\gamma_n}|u_0(\lambda)|^p|d\lambda|~<~\infty.    
\end{equation}
\par\smallskip
For every $\lambda\in  F(\mathbb T)$, there exists a real number $\alpha(\lambda)>0$ such that \Cref{Lemme-propriete-varphi} is satisfied.
Since $F(\mathbb T)$ is compact, there exists a finite set of points $\lambda_1,\dots,\lambda_K$ of $F(\mathbb T)$ and associated radii $\alpha_1=\alpha(\lambda_1),\dots,\alpha_K=\alpha(\lambda_K)>0$ such that $F(\mathbb T) \subseteq \bigcup\limits_{k=1}^K D(\lambda_k,\alpha_k)$.
Remark that by \Cref{casse-pieds}, the ``essential" self-intersection points in $\mathcal{O}$ must all be contained in the set $\{\lambda_1,\dots,\lambda_K\}$.
Let $\Omega$ be a connected component of $\mathbb C\setminus F(\mathbb T)$ contained in $\Omega_0^+$.

By \Cref{thm:quasi-conformal-function}, $F$ has an extension $\widetilde F$ on a neighborhood of $\mathbb T$ with the following properties:  there exist a positive constant  $c$, points $\zeta_1,\dots,\zeta_Q$ in $\mathbb T$ and positive radii $r_1,\dots, r_Q$ such that $\mathbb T\subset \cup_{i=1}^R D(\zeta_i,r_i)$ and for every $i=1,\ldots, Q$, the map $\widetilde F:D(\zeta_i,r_i)\longrightarrow \widetilde F(D(\zeta_i,r_i))$ is a $C^1$ diffeomorphism that preserves the local orientation and is such that, for every $z,w\in D(\zeta_i,r_i)$,
\[\frac1c|z-w|\le |\widetilde F(z)-\widetilde F(w)|\le c|z-w|.\]
As already mentioned in \Cref{subsection-Carleson-measures}, we have $Carl(r\mathbb T)=2\pi$ for every $r>0$. Thus, by \Cref{PropA8}, we get
\begin{equation}\label{Eq:kjo}
    Carl\Big(\widetilde F\big(r\mathbb T\cap D(\zeta_i,r_i)\big)\Big)\le 2c^2 Carl\big(r\mathbb T\cap D(\zeta_i,r_i)\big)\le 2c^2 Carl(r\mathbb T)\le 4\pi c^2.
\end{equation}
Let $(\varepsilon_n)$ be a sequence of positive numbers decreasing to  $0$ sufficiently fast to ensure that $(1+\varepsilon_n)\mathbb T\cap (1-\varepsilon_n)\mathbb T$ is contained in $\cup_{i=1}^N D(\zeta_i,r_i)$, and
let, for each $n$, $\gamma_n$ be a Jordan curve such that 
\[\gamma_n\subset\Omega\,\cap\,\widetilde{F}\big((1+\varepsilon_n)\mathbb T\cup(1-\varepsilon_n)\mathbb T\big).\]
and \[\gamma_n \subseteq\bigcup\limits_{k=1}^K D(\lambda_k,\alpha_k).\]
The smoothness of $\widetilde F$ implies that $\gamma_n$ is rectifiable and tends to $\partial \Omega$ with respect to the Hausdorff distance. Moreover, by \eqref{Eq:kjo}, we have
\begin{align*}
    Carl(\gamma_n)&\le Carl\Big(\widetilde{F}\big((1+\varepsilon_n)\mathbb T\cup(1-\varepsilon_n)\mathbb T\big)\Big)\\
    &\le \sum_{i=1}^NCarl\Big(\widetilde F\big((1+\varepsilon_n)\mathbb T\,\cap\,D(\zeta_i,r_i)\big)\Big)+Carl\Big(\widetilde F\big((1-\varepsilon_n)\mathbb T\,\cap\,D(\zeta_i,r_i)\big)\Big)\\&\le 8\pi Nc^2.
\end{align*}
In particular the Carleson constants of the curves $\gamma_n$ are uniformly bounded, and their length is bounded by $8\pi NRc^2$, where $R>0$ is sufficiently large to ensure that that $\Omega\subset R\mathbb D.$
\par\smallskip
For every $n\geq 1$ and $1\leq k\leq K$, denote by $\gamma_{n,k}$ the set $\gamma_{n,k}=\gamma_n \cap D(\lambda_k,\alpha_k)$.
In order to prove \Cref{eq1:Smirnov-bounded-U}, it is sufficient to check that for every $1\leq k\leq K$, we have
\begin{equation}\label{eq:sdsds183}
\sup_{n\geq 1}\int_{\gamma_{n,k}}|u_0(\lambda)|^p|d\lambda|~<~\infty.    
\end{equation}
So fix $1\le k\le K$ and let $d_1,\dots,d_s$ be the functions defined in \Cref{section-another-representationflambdaplus} for the corresponding disk $D(\lambda_k,\alpha_k)$ (the reader should pay attention to the fact that the $d_j$'s depend on $\lambda_k$). 
Observe that the set $N(\lambda)=\{1\leq j\leq s\,;\,|d_j(\lambda)|<1\}$ is constant on $\Omega\cap D(\lambda_k,\alpha_k)$. We will denote by $N(\Omega)$ this set and by $n(\Omega)$ its cardinal. Using a partial fraction decomposition, for every $\lambda\in D(\lambda_k,\alpha_k)$, we can write \[
\prod_{j\in N(\Omega)}\frac1{z-d_j(\lambda)^{-1}}~=\sum_{j\in N(\Omega)} \frac{c_j(\lambda)}{z-d_j(\lambda)^{-1}},
\]
where 
\begin{equation}\label{eq:coeff-fraction-rationnelle}
c_j(\lambda)~=\prod_{\substack{\ell\neq j\\\ell\in N(\Omega)}}\frac1{d_j(\lambda)^{-1}-d_{\ell}(\lambda)^{-1}}\cdot
\end{equation}
Recall that $h_{\lambda,0}=f_{\lambda}^+(0)/f_{\lambda}^+$. For every $\lambda\in \gamma_{n,k}$, we have 
\begin{eqnarray*}
u_0(\lambda)&=&\dual{g}{{f_{\lambda}^+(0)}/{f_{\lambda}^+}} \\
&=& f_\lambda^+(0)\sum_{j\in N(\Omega)} c_j(\lambda) \dual[\Big]{g}{\varphi(\lambda,z)^{-1}(z-d_j(\lambda)^{-1})^{-1}},
\end{eqnarray*}
where $\varphi$ is the function given by  \Cref{Lemme-propriete-varphi}.
Define a function $\eta_j$ on $D(\lambda_k,\alpha_k)\times\T$ by setting 
\[
\eta_j(\lambda,z)~=~\frac{\varphi(\lambda,z)^{-1}-\varphi(\lambda,|d_j(\lambda)|d_j(\lambda)^{-1})^{-1}}{z-d_j(\lambda)^{-1}}, \quad \lambda\in D(\lambda_k,\alpha_k), \, z\in\T.
\]
Then, we have
\begin{align*}
\dual[\Big]{g}{\varphi(\lambda,z)^{-1}(z-d_j(\lambda)^{-1})^{-1}}=~&\dual{g}{\eta_j(\lambda,\cdot)}\\&\quad+\dual[\Big]{g}{\varphi(\lambda,|d_j(\lambda)|d_j(\lambda)^{-1})^{-1}(z-d_j(\lambda)^{-1})^{-1}}\\
=~&\dual{g}{\eta_j(\lambda,\cdot)}\\&\quad+\varphi(\lambda,|d_j(\lambda)|d_j(\lambda)^{-1})^{-1}\dual[\Big]{g}{(z-d_j(\lambda)^{-1})^{-1}}\\
=~&\dual{g}{\eta_j(\lambda,\cdot)}\\&\quad-d_j(\lambda)\varphi(\lambda,|d_j(\lambda)|d_j(\lambda)^{-1})^{-1}g(d_j(\lambda)),
\end{align*}
where the last equation follows from \Cref{eq:cauchy-kernel} and the fact that $d_j(\lambda)\in\mathbb D$ for every $j\in N(\Omega)$.
Then 
\[
\left|\dual[\big]{g}{\varphi(\lambda,z)^{-1}(z-d_j(\lambda)^{-1})^{-1}}\right|~\leq~ \|g\|_{H^p}\|\eta_j(\lambda,\cdot)\|_{L^q(\mathbb T)}+ |g(d_j(\lambda))|\sup_{D(\lambda_k,\alpha_k)\times\mathbb T}|1/\varphi|.
\]
Note that, by \Cref{Lemme-propriete-varphi}, the function $\lambda\mapsto f_\lambda^+(0)$ belongs to $A(\Omega\cap D(\lambda_k,\alpha_k))$ and $\sup_{D(\lambda_k,\alpha_k)\times\mathbb T}|1/\varphi|$ is finite. Moreover, by \Cref{eq:sdfsfjsdlfsdj1323,eq:coeff-fraction-rationnelle}, there exists a constant $C>0$ such that for every $\lambda\in D(\lambda_k,\alpha_k)$ and for every $1\leq j\leq n_k$, we have $|c_j(\lambda)|\leq C$.  Hence there exist two positive constants $C_1$ and $C_2$ such that, for every $\lambda\in \Omega\cap D(\lambda_k,\alpha_k)$,
\begin{equation}\label{eq23REZRE-continuite-U-estimation-u_0}
|u_0(\lambda)|~\leq~ C_1 \|g\|_{H^p}\!\!\! \sum_{j\in N(\Omega)}\|\eta_j(\lambda,\cdot)\|_{L^q(\mathbb T)}~+~C_2\!\!\! \sum_{j\in N(\Omega)}|g(d_j(\lambda))|.
\end{equation}
\par\smallskip
We state separately the next estimate needed to conclude the proof. Let us point out that it is here where the condition $\varepsilon>1/p$ comes into play. The condition $\varepsilon>1/q$ will in its turn be needed in the proof of \Cref{lemme-Vborne}.

\begin{fact}\label{faitC3}
There exists a constant $C_3>0$ such that $\|\eta_j(\lambda,\cdot)\|_{L^q(\mathbb T)}\leq C_3$.
\end{fact}

\begin{proof} Let $\beta$ be such that $1/p<\beta<\varepsilon$, and set $\gamma=\varepsilon-\beta>0$. According to \Cref{Lemme-propriete-varphi}, the function $1/\varphi_l$ is  of class $C^\gamma$ from $D(\lambda_k,\alpha_k)$ into $C^\beta(\overline{\mathbb D})$. By \Cref{eq:NormeCeps}, there exists $M>0$ such that 
\[
\sup_{\lambda\in D(\lambda_k,\alpha_k)}\|1/\varphi(\lambda,\cdot)\|_{C^\beta(\overline{\mathbb D})}~\leq~ M.
\]
Hence we have, for every $\lambda\in D(\lambda_k,\alpha_k)$ and for every $z_1,z_2\in\overline{\mathbb D}$, 
\[
|\varphi(\lambda,z_1)^{-1}-\varphi(\lambda,z_2)^{-1}|~\leq~ M |z_1-z_2|^\beta.
\]
Thus, for every $z\in\mathbb T$,
\[
|\eta_j(\lambda,z)|~\leq~ M \frac{|z-|d_j(\lambda)|d_j(\lambda)^{-1}|^\beta}{|z-d_j(\lambda)^{-1}|}\cdot
\]
It can be easily checked that for every $z\in\mathbb T$ and $w\in\mathbb C^*$, we have $\left|z-\frac{w}{|w|}\right|\leq 2|z-w|$, which gives 
\[
|\eta_{j}(\lambda,z)|~\leq~ M 2^{\beta}|d_j(\lambda)|^{1-\beta}\frac{1}{|1-d_j(\lambda)z|^{1-\beta}}~\leq~ 2M \frac{1}{|1-d_j(\lambda)z|^{1-\beta}}\cdot
\]
Define now, for every $\delta\in\R$, a function $J_\delta$ on $\D$ by setting
\[
J_\delta(z)~=~\int_{\mathbb T}\frac{|\mathrm d\tau|}{|1-\tau z|^{1+\delta}}\quad \textrm{ for every } z\in\mathbb D.
\]
By \cite[Th. 1.7]{HedenmalmKorenblumZhu2000}, $J_\delta$ is bounded on $\mathbb D$ if $\delta<0$. 
Take $\delta=q(1-\beta)-1$; since $\beta>p^{-1}=1-q^{-1}$, we have $\delta<0$.  Thus we get 
\[
\|\eta_j(\lambda,\cdot)\|_{L^q(\mathbb T)}^q~\leq~ (2M)^q \int_{\mathbb T}\frac{1}{|1-d_j(\lambda)\tau|^{q(1-\beta)}}\,|\mathrm{d}\tau|~=~(2M)^q J_{\delta}(d_j(\lambda))\leq C_3^q,
\]
where $C_3= (2M)\left(\sup_{z\in\D}|J_\delta(z)|\right)^{1/q}<+\infty$.  This proves \Cref{faitC3}.
\end{proof}

Getting back to our initial estimate of $|u_0(\lambda)|$, by combining \Cref{eq23REZRE-continuite-U-estimation-u_0} and \Cref{faitC3}, we obtain that, for every $\lambda\in D(\lambda_k,\alpha_k)$,
\[
|u_0(\lambda)|~\leq~ C_1C_3 n(\Omega)\|g\|_{H^p}~+~C_2\sum_{j\in N(\Omega)}|g(d_j(\lambda))|.
\]
Since 
\[
\sup_n\int_{\gamma_n}|d\lambda|~=~\sup_n|\gamma_n|~<~\infty,
\]
in order to prove \Cref{eq:sdsds183} (and hence \Cref{eq1:Smirnov-bounded-U}), it remains to show that 
\begin{equation}\label{eq2:Smirnov-bounded-U}
    \sup_{n\geq 1}\int_{\gamma_{n,k}}|g(d_j(\lambda))|^p\,|\mathrm{d}\lambda|~<~\infty.
\end{equation}
Let $c>0$ be such that 
\[
\frac{1}{c}|\lambda-\mu|\leq |d_j(\lambda)-d_j(\mu)|\leq c|\lambda-\mu|\quad\text{for every }\lambda,\mu\in D(\lambda_k,\alpha_k).
\]
According to \Cref{PropA8}, we have 
\[
\int_{\gamma_{n,k}}|g(d_j(\lambda))|^p\,|\mathrm{d}\lambda|~\leq~24e\,c^3\,Carl(\gamma_{n,k})\,\|g\|_p^p,
\]
and since 
\[
\sup_{n\geq 1}Carl(\gamma_{n,k})\leq \sup_{n\geq 1}Carl(\gamma_n)<\infty,
\]
we get \Cref{eq2:Smirnov-bounded-U}, whence \Cref{eq1:Smirnov-bounded-U}.
\par\smallskip
Finally, we can conclude that the function $u_j$ belongs to $E^p(\Omega_j^+)$ for every $0\leq j\leq N-1$, and thus the $N$-tuple $Ug=(u_j)_{0\le j\leq N-1}$ belongs to $\bigoplus E^p(\Omega_j^+)$ for every $g\in H^p$.
\par\smallskip
It remains to prove that $U$ is a bounded operator from $H^p$ into $\bigoplus E^p(\Omega_j^+)$. This is a consequence of the Closed Graph Theorem. Suppose that $(g_n)_n$ is a sequence of functions in $H^p$ such that $g_n\to 0$ in $H^p$ and $Ug_n\longrightarrow u=(u_j)_{0\leq j\le N-1}$ in $\bigoplus E^p(\Omega_j^+)$. Since convergence in $E^p(\Omega_j^+)$ implies pointwise convergence on $\Omega_j^+$,  we have
\[
\dual{g_n}{h_{\lambda,j}}~\longrightarrow~ u_j(\lambda)\quad \text{as }n\to\infty \text{ for every } 0\leq j\le N-1 \text{ and every } \lambda\in \Omega_j^+.
\]
But since $g_n\to 0$ in $H^p$, we also have that $\dual{g_n}{h_{\lambda,j}}\to 0$, and thus $u_j(\lambda)=0$ for every $0\leq j\le N-1$ and every $\lambda\in \Omega_j^+$. In other words, $u_j=0$ for every $0\leq j\le N-1$, ie. $u=0$. By the Closed Graph Theorem, this proves that $U$ is bounded.
\end{proof}

Recall that $h_{\lambda,j}(z)=(Uk_{\bar z})_j(\lambda)$ for every $z\in\D$ and every $\lambda\in\Omega_j^+$, $0\leq j\leq N-1$. Moreover, for every $z\in\D$, the map $\lambda\longmapsto h_{\lambda,j}(z)$ is in $E^p(\Omega_j^+)$. When assumption \ref{H2} holds, the functions $h_{\lambda,j}$ satisfy the boundary relation of \Cref{eq:3434RDCS}. Using the density of the linear span of the Cauchy kernels $k_z$, $z\in \mathbb D$, and the continuity of $U$, we will now prove that the range of $U$ is contained in the subspace $\mathcal{E}_F^p$ consisting of $N$-tuples $u=(u_j)_{0\leq j\leq N-1}$ in $\bigoplus E^p(\Omega_j^+)$ satisfying the following boundary relations: for every $0\leq j< N-1$, we have
\begin{equation}\label{eq:definition-espace-modele}
u_j^{int}-\xi u_{j+1}^{int}~=~u_j^{ext}\quad \text{a.e. on }\partial\Omega_{j+1}^+~=~\Gamma\cap \overline{\Omega_{j+1}^+}~=~\bigcup_{k=j+1}^N\Gamma_k^+,
\end{equation}
where $\Gamma=F(\T)$ and $\Gamma_k^+=\partial\overline{\Omega_k^+}$ for every $0\le k\le N-1$ (see \Cref{piece of notation}).
\par\smallskip 

\begin{lemma}
Let $F$ satisfy \ref{H1}, \ref{H2} and \ref{H3'}. Then
\begin{enumerate}
\item[(1)] $\mathcal{E}_F^p$ is a closed subspace of $\bigoplus\limits_{0\leq j\leq N-1}E^p(\Omega_j^+)$;
\item[(2)] $UH^p \subseteq \mathcal{E}_F^p$. 
\end{enumerate}
\end{lemma}

\begin{proof}
(1) Since the norm on $\oplus E^p(\Omega_j^+)$ is given by a sum of $L^p$ norms of boundary values, the map $u=(u_0,\dots,u_{N-1})\mapsto u_k^{int}-\xi u_{k+1}^{int}-u_k^{ext}$ is continuous from $\oplus E^p(\Omega_j^+)$ into $L^p(\Omega_{k+1}^+)$ for each $0\le k<N$, and thus $\mathcal E^p_F$ is clearly closed in $\oplus E^p(\Omega_j^+)$.
\par\smallskip
(2) For every $z\in\mathbb D$ and every $0\le j<N$, we have $(Uk_{\overline z})_j(\lambda)=h_{\lambda,j}(\lambda)$. So by \Cref{cor:condition-aux-bords-hlambdaj}, $Uk_{\overline z}\in\mathcal E^p_F$. Thus (2) is a direct consequence of the completeness of the functions $k_{\overline z}$ in $H^p$ and of the continuity of $U$ given by \Cref{lemme-Uborne}.
\end{proof}

The aim of the remaining sections is to show that under the hypothesis \ref{H2}, the operator $U$ is an isomorphism from $H^p$ onto $\cal E_F^p$. We  start by introducing another operator $V$ which will be useful to construct the inverse of $U$. We will prove that this operator is also bounded under the more general hypothesis \ref{H2'}. Let us  mention here that in the case where \ref{H2'} holds, it is possible, using \Cref{eq:relation-aux-bords-H2prime34455SD}, to show that the range of $U$ is also included in some subspace of $\bigoplus E^p(\Omega_j^+)$ defined by some boundary relations. 

\subsubsection{Definition and continuity of $V$}
For every $\lambda\in \mathbb C\setminus F(\mathbb T)$ and every integer $m\geq n=\w_F(\lambda)$, let us define
\begin{equation}\label{eq:definition-Glambda}
G_{\lambda,m}(z)~=~-z^{m-n} F_{\lambda}^+(z)^{-1}\quad \textrm{ for every } z\in\mathbb D.
\end{equation}

Recall that $N=\max\{\w_F(\lambda)\,;\, \lambda\in\C\setminus F(\T)\}$. For every $g\in H^q$ and every $0\le j\le N-1$, denote by $v_j$ the function defined on $\widehat{\mathbb C}\setminus (F(\mathbb T)\cup\Omega_j^+)$ by
\[
v_j(\lambda)~=~\dual{G_{\lambda,j}}{g}\quad \textrm{ for every } \lambda\in \widehat{\mathbb C}\setminus (F(\mathbb T)\cup\Omega_j^+),
\]
and by $Vg$ the $N$-tuple
\[
Vg~=~(v_j)_{0\le j \le N-1}.
\]

It follows from \Cref{eq33434A4:relation-aux-bords}  that when $F$ satisfies \ref{H2},  we have
\begin{equation}\label{eq:Glambda-equation-au-bord}
G_{\lambda,j+1}^{ext}(z)-\xi(\lambda)G_{\lambda,j}^{ext}(z)~=~G_{\lambda,j+1}^{int}(z)\quad \text{a.e. on }F(\T)\setminus\partial\Omega_{j+1}^+~=~\bigcup_{k=0}^j\Gamma_k^+,
\end{equation}
for every $z\in\mathbb D$ and for every $0\leq j\leq N-1$. Indeed, for every $z\in\D$, $0\leq k\leq j$, and for almost every $\lambda\in\Gamma_k^+$, we have 
\[
G_{\lambda,j}^{ext}(z)~=~\lim_{\substack{\mu\to\lambda\\ \mu\in\Omega_{k}}}G_{\mu,j}(z)~=~\lim_{\substack{\mu\to\lambda\\ \mu\in\Omega_{k}}}\left(-z^{j-k}F_\mu^+(z)^{-1}\right)~=~-z^{j-k}F_{\lambda,ext}^{+}(z)^{-1},
\]
and 
\[
G_{\lambda,j+1}^{int}(z)~=~\lim_{\substack{\mu\to\lambda\\ \mu\in\Omega_{k+1}}}G_{\mu,j+1}(z)~=~\lim_{\substack{\mu\to\lambda\\ \mu\in\Omega_{k+1}}}\left(-z^{j-k}F_\mu^+(z)^{-1}\right)~=~-z^{j-k}F_{\lambda,int}^{+}(z)^{-1}.
\]
Recall now that \cref{eq33434A4:relation-aux-bords} gives
\[
F_{\lambda,int}^+(z)^{-1}~=~(z-\xi(\lambda))F_{\lambda,ext}^+(z)^{-1}
\quad\text{ for every } z\in\D,
\]
from which it follows that, for every $0\leq k\leq j$ and for almost every $\lambda\in \Gamma_k^+$ 
\begin{eqnarray*}
G_{\lambda,j+1}^{int}(z)&=&-z^{j-k}(z-\xi(\lambda))F_{\lambda,ext}^{+}(z)^{-1}\\
&=&(z-\xi(\lambda))G_{\lambda,j}^{ext}(z)\\
&=&G_{\lambda,j+1}^{ext}(z)-\xi(\lambda)G_{\lambda,j}^{ext}(z),
\end{eqnarray*}
which proves \Cref{eq:Glambda-equation-au-bord}. 
\par\smallskip
\begin{theorem}\label{lemme-Vborne} Let $F$ satisfy \ref{H1}, \ref{H2'} and \ref{H3'}. The linear map $V$ is well-defined and bounded from $H^q$ into the space
$\bigoplus\limits_{0\le j \le N-1}E_0^q(\widehat{\mathbb C}\setminus (F(\mathbb T)\cup\Omega_j^+))$.
\end{theorem}

\begin{proof}
The proof is quite similar to that of \Cref{lemme-Uborne}. The assumption $\varepsilon>1/q$ is necessary in the present proof for exactly the same reason as the assumption $\varepsilon>1/p$ was needed in the proof of \Cref{faitC3}. We leave the details to the reader. The only point which requires an additional explanation is the following: when showing that each function $v_j$ belongs to $E_0^q(\widehat{\mathbb C}\setminus (F(\mathbb T)\cup\Omega_j^+))$, we need to check that 
$v_j(\lambda)\to 0$ as $|\lambda|\to\infty$. So let $\lambda$ be such that $|\lambda|\geq 2\|F\|_\infty$. In particular, $\w_F(\lambda)=0$. Then for every $0\leq j\leq N-1$, we have 
\[
v_j(\lambda)~=~\dual{G_{\lambda,j}}{g}~=~-\dual{z^jF_\lambda^+(z)^{-1}}{g}.
\]
So
\(
|v_j(\lambda)|\leq \|z^j F_\lambda^+(z)^{-1}\|_p \|g\|_q,
\)
and we need to bound the quantities $\|z^j F_\lambda^+(z)^{-1}\|_p$ for $|\lambda|\geq 2\|F\|_\infty$ in a suitable way. 
\par\smallskip
We have
$F_\lambda^{+}(z)^{-1}=\frac{1}{F(1)-\lambda}e^{-V_\lambda(z)}$, where
\[
V_\lambda(z)~=~P_+U_\lambda(z)~=~\frac{1}{2i\pi}\int_{\mathbb T}\frac{U_\lambda(\tau)}{\tau-z}\,\mathrm{d}\tau\quad \textrm{ for every } z\in\D
\]
and
\[
U_\lambda(e^{is})~=~\int_0^s ie^{it}\frac{\psi_\lambda'(e^{it})}{\psi_\lambda(e^{it})}\,\mathrm{d}t\quad \textrm{ for every } s\in [0,2\pi).
\]
Since we have here $\psi_\lambda(z)=F(z)-\lambda$, $z\in\T$, 
\[
U_\lambda(e^{is})~=~\int_0^s ie^{it}\frac{F'(e^{it})}{F(e^{it})-\lambda}\,\mathrm{d}t\quad\text{ for every } s\in [0,2\pi).
\]
Hence for $|\lambda|\geq 2\|F\|_\infty$ and $s,\theta\in [0,2\pi)$, we have 

\begin{equation}\label{eq:sdffdf2323434SDS}
|U_\lambda(e^{is})-U_\lambda(e^{i\theta})|~=~\left|\int_\theta^s ie^{it}\frac{F'(e^{it})}{F(e^{it})-\lambda}\,\mathrm{d}t\right|
~\leq~ \frac{\|F'\|_\infty}{\|F\|_\infty}\,|s-\theta|
~\leq~ C_1 |e^{i\theta}-e^{it}|,
\end{equation}
where the constant $C_1$ does not depend of $\lambda$. 

Now, for $t\in [0,2\pi)$ and $\lambda\in\mathbb C\setminus F(\mathbb T)$, denote by $v_{\lambda,t}(x)=V_\lambda(xe^{it})$, $x\in [0,1)$. Then, for every $0<r<1$ and $t\in [0,2\pi)$, we have 
\[
V_\lambda(re^{it})~=~V_\lambda(0)+\int_0^r v_{\lambda,t}'(x)\,\mathrm{d}x
\]
with
\[
v_{\lambda,t}'(x)~=~\frac{e^{it}}{2\pi}\int_0^{2\pi}\frac{e^{i\theta}U_\lambda(e^{i\theta})}{(e^{i\theta}-xe^{it})^2}\,\mathrm{d}\theta.
\]
Since
\[
\frac{1}{2\pi}\int_0^{2\pi}\frac{e^{i\theta}}{(e^{i\theta}-xe^{it})^2}\,\mathrm{d}\theta~=~0\quad\text{ for every } x\in (0,1),
\]
we deduce that
\[
v_{\lambda,t}'(x)~=~\frac{e^{it}}{2\pi}\int_0^{2\pi}\frac{e^{i\theta}(U_\lambda(e^{i\theta})-U_\lambda(e^{it}))}{(e^{i\theta}-xe^{it})^2}\,\mathrm{d}\theta,
\]
and thus, by \Cref{eq:sdffdf2323434SDS}, we have
\[
|v_{\lambda,t}'(x)|~\leq~ \frac{C_1}{2\pi}\int_0^{2\pi}\frac{|e^{i\theta}-e^{it}|}{|e^{i\theta}-xe^{it}|^2}\,\mathrm{d}\theta.
\]
Now, the estimate $|1-e^{i(t-\theta)}|\leq 2 |1-xe^{i(t-\theta)}|$ yields that

\begin{equation}\label{eq1:bornitude-V-Smirnov-non-borne-3435}
|v_{\lambda,t}'(x)|~\leq~ \frac{C_1}{\pi}\int_0^{2\pi}\frac{1}{|1-xe^{i(t-\theta)}|}\,\mathrm{d}\theta.
\end{equation}
In order to estimate the integral on the right-hand side, we consider separately the cases $0<x\leq 1/2$ and $1/2<x\le 1$.
For $0<x\leq 1/2$, we have $|1-xe^{i(t-\theta)}|\geq 1-x\geq 1/2$, whence 
\begin{equation}\label{eq2:bornitude-V-Smirnov-non-borne-3435}
|v_{\lambda,t}'(x)|~\leq~ 4C_1.
\end{equation}
On the other hand, according to \cite[Lemma 1.12.3]{CimaMathesonRoss2006}, there exists a constant $C_2$ independent of $x\in (1/2,1)$ such that 
\[
\int_0^{2\pi}\frac{1}{|1-xe^{i(t-\theta)}|}\,\mathrm{d}\theta ~\leq~ -C_2 \log(1-x)\quad\text{ for every } t\in[0,2\pi)\text{ and } x\in (1/2,1),
\]
whence it follows that
\begin{equation}\label{eq3:bornitude-V-Smirnov-non-borne-3435}
|v_{\lambda,t}'(x)|~\leq~ -C_3 \log(1-x) \quad \text{ for every } x\in (1/2,1).
\end{equation}
Therefore, putting together the estimates \Cref{eq2:bornitude-V-Smirnov-non-borne-3435,eq3:bornitude-V-Smirnov-non-borne-3435}, we get 
\begin{eqnarray*}
|V_\lambda(re^{it})|&\leq & |V_\lambda(0)|+\int_0^r |v_{\lambda,t}'(x)|\,\mathrm{d}x\\
&\leq & |V_\lambda(0)|+\int_0^{1/2} |v_{\lambda,t}'(x)|\,\mathrm{d}x+\int_{1/2}^1 |v_{\lambda,t}'(x)|\,\mathrm{d}x\\
&\leq & |V_\lambda(0)|+2C_1-C_2 \int_{1/2}^1 \log(1-x)\,\mathrm{d}x\\
&=& |V_\lambda(0)|+2C_1-C_2 \int_{0}^{1/2} \log(x)\,\mathrm{d}x.
\end{eqnarray*}
We are now close to our goal.
Observe that 
\[
|U_\lambda(e^{is})|~\leq~ \int_0^s \frac{|F'(e^{it})|}{|F(e^{it})-\lambda|}\,\mathrm{d}t~\leq~ \frac{\|F'\|_\infty}{\|F\|_\infty}2\pi =:C_3,
\]
and it follows from this estimate that
\[
|V_\lambda(0)|~=~|\widehat{U_\lambda}(0)|~\leq ~\|U_\lambda\|_\infty~\leq~ C_3.
\]
Finally, if we take $C_4=C_3+2C_1-C_2\int_0^{1/2}\log(x)\,\mathrm{d}x<\infty$, we see that for every $z\in\mathbb D$, and every $|\lambda|\geq 2\|F\|_\infty$, we have $$|V_\lambda(z)|~\leq~ C_4.$$ Therefore 
\[
|F_\lambda^+(z)^{-1}|~=~\frac{1}{|F(1)-\lambda|}e^{-\mathrm{Re}(V_\lambda(z))}~\leq~ \frac{1}{|F(1)-\lambda|}e^{|V_\lambda(z)|}~\leq~ \frac{e^{C_4}}{|F(1)-\lambda|},
\]
where the constant $C_4$ is independent of $\lambda$. We thus get 
\[
|v_j(\lambda)|~\leq~ \|z^j F_\lambda^+(z)^{-1}\|_p \|g\|_q~\leq~ \frac{e^{C_4}}{|F(1)-\lambda|} \|g\|_q,
\]
from which we conclude that $v_j(\lambda)\to 0$ as $|\lambda|\to \infty$. This terminates the proof of the boundedness of $V$. 
\end{proof}
\subsubsection{Left invertibility of $U$}
{\bf {From now on, we will assume that $F$ satisfies \ref{H1}, \ref{H2} and \ref{H3'}.}}
\par\smallskip
In order to show that $U$ is left-invertible as an operator from $H^p$ into $\mathcal{E}_F^p$, we
will embed the range of $U$ into the space $L^p(\Gamma)$ and the range of $V$ into the space $L^q(\Gamma)$, where $\Gamma=F(\mathbb T)=\displaystyle\bigcup_{j=0}^{N-1}\Gamma_j^+$.
\par\smallskip
Given $u=(u_j)_{0\le j \le N-1}\in \bigoplus_{0\le j \le N-1}E^p(\Omega_j^+)$, we define a function $\pi_{int}(u)$ on $\Gamma$ by setting 
\[
\pi_{int}(u)(\lambda)~=~u_j^{int}(\lambda)\quad\text{ for a.e }\lambda\in\Gamma_j^+,\; 0\le j\le N-1,
\]
and similarly for $v=(v_j)_{j\geq N}\in \bigoplus_{0\le j \le N-1}E_0^q(\widehat{\mathbb C}\setminus (F(\mathbb T)\cup\Omega_j^+))$, we define  a function $\pi_{ext}(v)$ on $\Gamma$ by setting 
\[
\pi_{ext}(v)(\lambda)~=~v_j^{ext}(\lambda)\quad\text{ for a.e. } \lambda\in \Gamma_j^+,\; 0\le j\le N-1.
 \]
The duality between $L^p(\Gamma)$ and $L^q(\Gamma)$ is given by
\[
\dualg{\widetilde{u}}{\widetilde{v}}~=~\int_{\Gamma} \widetilde{u}(\lambda) \widetilde{v}(\lambda)\,\mathrm{d}\lambda\quad \textrm{ for every } \widetilde{u}\in L^p(\Gamma),\; \widetilde{v}\in L^q(\Gamma).
\]

\begin{lemma}\label{piint-piext}
Let $\Gamma=F(\mathbb T)$. The operator $\pi_{int}$ is bounded from $\bigoplus\limits_{0\le j \le N-1}E^p(\Omega_j^+)$ into $L^p(\Gamma)$ and the operator $\pi_{ext}$ is bounded from $\bigoplus\limits_{0\le j \le N-1}E_0^q(\widehat{\mathbb C}\setminus (F(\mathbb T)\cup\Omega_j^+))$ into $L^q(\Gamma)$. Moreover, we have
\begin{equation}\label{eq:prod-scalaire-L2GAMMA23232}
    \dualg{\pi_{int}u}{\pi_{ext}v}~=~\frac{1}{2i\pi}\sum_{j=0}^{N-1}\int_{\Gamma_j^+}u_j^{int}(\lambda)v_j^{ext}(\lambda)\,\mathrm{d}\lambda.
\end{equation}
for every $u\in \bigoplus\limits_{0\le j \le N-1}E^p(\Omega_j^+)$ and every $v\in \bigoplus\limits_{0\le j \le N-1}E_0^q(\widehat{\mathbb C}\setminus (F(\mathbb T)\cup\Omega_j^+))$.
\end{lemma}

\begin{proof}
Let $u=(u_j)_{0\le j \le N-1}\in \bigoplus E^p(\Omega_j^+)$. Note that according to Equation \eqref{Eq:NormeErLimIntExt}, we have 
\[\|u_j^{int}\|_{L^p(\Gamma_j^+)}^p~\le~ \|u_j\|_{E^p(\Omega_j^+)}^p.\]
Thus we get
\[\|\pi_{int}u\|_{L^p(\Gamma)}^p~=~\sum_{j=0}^{N-1}\|u_j^{int}\|_{L^p(\Gamma_j^+)}^p~\le~ \sum_{j=0}^{N-1}\|u_j\|_{E^p(\Omega_j^+)}^p~=~\|u\|^p,\] where the norm of $u$ is taken in the space ${\bigoplus E^p(\Omega_j^+)}$.
This shows that the operator $\pi_{int}$ is bounded. A similar argument shows that $\pi_{ext}$ is bounded too.
\par\smallskip
Now, since $\bigcup_{j=0}^{N-1}\Gamma_j^+$ is a partition of $\Gamma$ (up to a finite set of points) we have
\begin{eqnarray*}
\dualg{\pi_{int}u}{\pi_{ext}v}&=&\frac{1}{2i\pi}\sum_{j=0}^{N-1}\int_{\Gamma_j^+}\pi_{int}(u)(\lambda)\pi_{ext}(v)(\lambda)\,\mathrm{d}\lambda\\
&=&\frac{1}{2i\pi}\sum_{j=0}^{N-1}\int_{\Gamma_j^+}u_j^{int}(\lambda)v_j^{ext}(\lambda)\,\mathrm{d}\lambda,
\end{eqnarray*}
which gives \Cref{eq:prod-scalaire-L2GAMMA23232} and ends the proof of the lemma.
\end{proof}

Now we can prove that $U$ is left-invertible on $H^p$. 
\begin{theorem}\label{dualite-L2-Gamma}
Let $F$ satisfy \ref{H1}, \ref{H2} and \ref{H3'}. For every $g_1\in H^p$ and every $g_2\in H^q$, we have 
\begin{equation}\label{eq:sdqsf6233SDD00000}
\dualg{\pi_{int}Ug_1}{\pi_{ext}Vg_2}~=~\dual{g_1}{g_2}.
\end{equation}
In particular, we have 
\[
(\pi_{ext}V)^*\pi_{int}U~=~I_{H^p}.
\]
\end{theorem}

\begin{proof}
Since $\pi_{int}$, $\pi_{ext}$, $U$ and $V$ are bounded operators on their respective domains by \Cref{piint-piext} and \Cref{lemme-Uborne,lemme-Vborne} respectively, it is sufficient to check \Cref{eq:sdqsf6233SDD00000} for two Cauchy kernels $g_1$ and $g_2$ (remember that the linear span of Cauchy kernels $k_z$, $z\in\D$, is dense both in $H^p$ and in  $H^q$). So let $g_1=k_{\overline{z_1}}$ and $g_2=k_{\overline{z_2}}$, where $z_1,z_2\in\mathbb D$. By \Cref{eq:cauchy-kernel}, we have $Uk_{\overline{z_1}}=(u_j)_{0\le j \le N-1}$ with 
\[
u_j(\lambda)~=~\dual{k_{\overline{z_1}}}{h_{\lambda,j}}~=~h_{\lambda,j}(z_1)\quad \text{ for every } \lambda\in \Omega_j^+,\; 0\le j\le N-1
\]
and $Vk_{\overline{z_2}}=(v_j)_{j\geq N}$ with 
\[
v_j(\lambda)~=~\dual{G_{\lambda,j}}{k_{\overline{z_2}}}~=~G_{\lambda,j}(z_2)\quad \text{ for every }  \lambda\in \widehat{\mathbb C}\setminus(F(\mathbb T)\cap \Omega_j^+),\; 0\le j\le N-1.
\]
As already noted in \Cref{piece of notation}, the interior component at a point $\lambda\in \Gamma_j^+$ is a connected component of
 $\Omega_j^+$, and the exterior component is a connected component of
$ \C\setminus (F(\T)\cup \Omega_j^+)$. Thus
we have $u_j^{int}(\lambda)=h_{\lambda,j}^{int}(z_1)$ and $v_j^{ext}(\lambda)=G_{\lambda,j}^{ext}(z_2)$. By \Cref{eq:prod-scalaire-L2GAMMA23232}, we have 
\begin{equation}\label{eq:sdsd8343ZDS}
    \dualg{\pi_{int}Uk_{\overline{z_1}}}{\pi_{ext}Vk_{\overline{z_2}}} ~=~\frac{1}{2i\pi}\sum_{j=0}^{N-1}\int_{\Gamma_j^+}h_{\lambda,j}^{int}(z_1)G_{\lambda,j}^{ext}(z_2)\,\mathrm{d}\lambda.
\end{equation}
Recalling that the functions $h_{\lambda,j}$ are defined for $\lambda\in\Omega_j^+$ as
$h_{\lambda,j}(z)=z^j f_{\lambda}^+(0)/f_{\lambda}^+(z)$, $z\in\D$, we observe that this definition can in fact be extended to any $\lambda\in \C\setminus F(\T)$ (although the resulting functions are not eigenvectors of $T_F^*$ anymore) and they also satisfies the boundary relation \Cref{{eq:3434RDCS}} almost everywhere also on $\Gamma$.

Recall now, that for $0\leq j\leq N$, we denote by $\Omega_j$ the set 
\[
\Omega_j~=~\{\lambda\in\mathbb C\setminus F(\mathbb T)\,;\,\w_F(\lambda)=j\}.
\]
Since $\w_F$ is locally constant on $\mathbb C\setminus F(\mathbb T)$, the sets $\Omega_j$ are open sets and they form a partition of $\mathbb C\setminus F(\mathbb T)$. We can now define an analytic function on $\mathbb C\setminus F(\mathbb T)$  by setting
\[
\Phi(\lambda)~=~h_{\lambda,j}(z_1)G_{\lambda,j}(z_2)\quad \textrm{ for every } \lambda\in\Omega_j,\,0\le j\le N.
\]
For $0\leq j\leq N-1$ and $\lambda\in\Gamma_j^+$, according to \Cref{Eq:limIntExt-k}, we have 
\[
\Phi^{int}(\lambda)~=~\lim_{\substack{\mu\to\lambda\\ \mu\in\Omega_{j+1}}}\Phi(\mu)~=~\lim_{\substack{\mu\to\lambda\\ \mu\in\Omega_{j+1}}} h_{\mu,j+1}(z_1)G_{\mu,j+1}(z_2)~=~h_{\lambda,j+1}^{int}(z_1)G_{\lambda,j+1}^{int}(z_2),
\]
and similarly
\[
\Phi^{ext}(\lambda)~=~\lim_{\substack{\mu\to\lambda\\ \mu\in\Omega_{j}}}\Phi(\mu)~=~\lim_{\substack{\mu\to\lambda\\ \mu\in\Omega_{j}}} h_{\mu,j}(z_1)G_{\mu,j}(z_2)~=~h_{\lambda,j}^{ext}(z_1)G_{\lambda,j}^{ext}(z_2).
\]
Observe that, for almost every $\lambda\in\Gamma_j^+$, we have 
\[
h_{\lambda,j}^{int}(z_1)(1-z_1\xi(\lambda))~=~h_{\lambda,j}^{ext}(z_1),
\]
whence 
\[
h_{\lambda,j+1}^{int}(z_1)~=~\frac{z_1}{1-z_1\xi(\lambda)}h_{\lambda,j}^{ext}(z_1).
\]
On the other hand, \Cref{eq:Glambda-equation-au-bord} gives that for almost every $\lambda\in\Gamma_j^+$, we also have 
\[
G_{\lambda,j+1}^{int}(z_2)~=~(z_2-\xi(\lambda))G_{\lambda,j}^{ext}(z_2).
\]
Thus 
\[
\Phi^{int}(\lambda)~=~\frac{z_1(z_2-\xi(\lambda))}{1-z_1\xi(\lambda)}h_{\lambda,j}^{ext}(z_1)\cdot G_{\lambda,j}^{ext}(z_2),
\]
and we obtain that
\begin{equation}\label{eq:sqdsqdsqd!!hsdsdsFTDZ}
\Phi^{int}(\lambda_0)~=~\frac{z_1(z_2-\xi(\lambda_0))}{1-z_1\xi(\lambda_0)}\Phi^{ext}(\lambda_0)\quad\text{ for a.e. }\lambda_0\in\Gamma.    
\end{equation}
For $0\leq j\leq N-1$, let now $(K_n)_n$ be a sequence of compact sets contained in $\Omega_{j+1}$ such that $K_n$ converges to $\Omega_{j+1}$. Let $\gamma_n=\partial K_n$, and consider the positive orientation on $\gamma_n$; that is, when we travel along the curve $\gamma_n$, the compact set $K_n$ remains on the left side. 
\par\smallskip
By Cauchy's Theorem, since $\Phi$ is analytic on $\Omega_{j+1}$, we have
\[
\frac{1}{2i\pi}\int_{\gamma_n}\Phi(\lambda)\,\mathrm{d}\lambda~=~0.
\]
Moreover, using \Cref{blablaesdsds232323} as $n\to \infty$, we have 
\[
\frac{1}{2i \pi}\int_{\gamma_n}\Phi(\lambda)\,\mathrm{d}\lambda~\longrightarrow~ \frac{1}{2i \pi}\int_{\Gamma_j^+}\Phi^{int}(\lambda)\,\mathrm{d}\lambda-\frac{1}{2i\pi}\int_{\Gamma_{j+1}^+}\Phi^{ext}(\lambda)\,\mathrm{d}\lambda,
\]
whence
\begin{equation}\label{egalite}
\frac{1}{2i \pi}\int_{\Gamma_{j+1}^+}\Phi^{ext}(\lambda)\,\mathrm{d}\lambda-\frac{1}{2i\pi}\int_{\Gamma_j^+}\Phi^{int}(\lambda)\,\mathrm{d}\lambda~=~0,
\end{equation}
where the integral over $\Gamma_{j+1}^+$ should be replaced by $0$ when $j=N-1$ (recall that $\Gamma_{N}^+=\varnothing$).
Plugging \Cref{eq:sqdsqdsqd!!hsdsdsFTDZ} into \Cref{egalite},  we get that for every $0\leq j\leq N-1$,
\[
\frac{1}{2i\pi}\int_{\Gamma_{j+1}^+}\Phi^{ext}(\lambda)\,\mathrm{d}\lambda-\frac{1}{2i\pi}\int_{\Gamma_j^+}\frac{z_1(z_2-\xi(\lambda))}{1-z_1\xi(\lambda)}\Phi^{ext}(\lambda)\,\mathrm{d}\lambda~=~0.
\]
Remarking that $\dfrac{z_1(z_2-\xi(\lambda))}{1-z_1\xi(\lambda)}=\dfrac{z_1z_2-1}{1-z_1\xi(\lambda)}+1$, we obtain
\[
\frac{1}{2i\pi}\int_{\Gamma_{j+1}^+}\Phi^{ext}(\lambda)\,\mathrm{d}\lambda-\frac{1}{2i\pi}\int_{\Gamma_{j}^+}\Phi^{ext}(\lambda)\,\mathrm{d}\lambda~=~\frac{z_1z_2-1}{2i\pi}\int_{\Gamma_j^+}\frac{\Phi^{ext}(\lambda)}{1-z_1\xi(\lambda)}\,\mathrm{d}\lambda.
\]
Summing over $0\leq j\leq N-1$ yields
\begin{equation}\label{egalite2}
\frac{1-z_1z_2}{2i\pi}\sum_{j=0}^{N-1}\int_{\Gamma_j^+}\frac{\Phi^{ext}(\lambda)}{1-z_1\xi(\lambda)}\,\mathrm{d}\lambda~=~\frac{1}{2i\pi}\int_{\Gamma_0^+}\Phi^{ext}(\lambda)\,\mathrm{d}\lambda.
\end{equation}
Let us now compute the right-hand side integral in \Cref{egalite2}.
Since the winding number of the curve $\Gamma_0^+$ at the point $\infty$ is equal to $-1$, the Residue Theorem implies that 
\[
\frac{1}{2i\pi}\int_{\Gamma_0^+}\Phi^{ext}(\lambda)\,\mathrm{d}\lambda~=~-\text{Res}_\infty(\Phi),
\]
where $\text{Res}_\infty(\Phi)$ is the residue of $\Phi$ at infinity.
Let us now check that  $\text{Res}_\infty(\Phi)=-1$.  

Going back to the definitions of the functions $h_{\lambda,j}$ and $G_{\lambda,j}$, we see that $\Phi(\lambda)$ can also be written as
\[
\Phi(\lambda)~=~-z_1^{\w_F(\lambda)}f_{\lambda}^+(0) f_\lambda^+(z_1)^{-1}F_\lambda^+(z_2)^{-1}.
\]
Now, for $|\lambda|$ sufficient large, using the definitions of $f_\lambda^+$ and $F_\lambda^+$, we have 
\[
\Phi(\lambda)~=~-\frac{1}{F(1)-\lambda} \exp\left(-\frac{z_1}{2i\pi}\int_{\mathbb T}\frac{u_\lambda(\tau)}{\tau(\tau-z_1)}\,\mathrm{d}\tau-\frac{1}{2i\pi}\int_{\mathbb T}\frac{U_\lambda(\tau)}{\tau-z_2}\,\mathrm{d}\tau\right).
\]
Hence we deduce that  for $|\lambda|$ sufficiently small,
\begin{equation}\label{egalite3}
-\frac{1}{\lambda^2}\Phi\left(\frac{1}{\lambda}\right)~=~\frac{1}{\lambda}\frac{1}{\lambda F(1)-1} \exp\left(-\frac{z_1}{2i\pi}\int_{\mathbb T}\frac{u_{\frac{1}{\lambda}}(\tau)}{\tau(\tau-z_1)}\,\mathrm{d}\tau-\frac{1}{2i\pi}\int_{\mathbb T}\frac{U_{\frac{1}{\lambda}}(\tau)}{\tau-z_2}\,\mathrm{d}\tau\right).
\end{equation}
Note now that
\[
u_{\frac{1}{\lambda}}(e^{is})~=~\int_0^s i e^{it}\frac{f'(e^{it})}{f(e^{it})-\frac{1}{\lambda}}\,\mathrm{d}t~=~\lambda \int_0^s \frac{ie^{it}f'(e^{it})}{\lambda f(e^{it})-1}\,\mathrm{d}t,
\]
and that we have a similar formula for $U_{\frac{1}{\lambda}}(e^{is})$. Hence the function 
\[
\lambda~\longmapsto~ -\frac{z_1}{2i\pi}\int_{\mathbb T}\frac{u_{\frac{1}{\lambda}}(\tau)}{\tau(\tau-z_1)}\,\mathrm{d}\tau-\frac{1}{2i\pi}\int_{\mathbb T}\frac{U_{\frac{1}{\lambda}}(\tau)}{\tau-z_2}\,\mathrm{d}\tau
\] 
is analytic in a neighborhood of zero and vanishes at zero. Therefore it follows from \Cref{egalite3} that
\[
\text{Res}_\infty(\Phi)~=~\text{Res}_0\left(-\frac{1}{\lambda^2}\Phi\left(\frac{1}{\lambda}\right)\right)~=~-1.
\]
Thus \Cref{egalite2} reads now as
\begin{equation}\label{egalite4}
\frac{1-z_1z_2}{2i\pi}\sum_{j=0}^{N-1}\int_{\Gamma_j^+}\frac{\Phi^{ext}(\lambda)}{1-z_1\xi(\lambda)}\,\mathrm{d}\lambda~=~1.
\end{equation}
As already noted, for almost every $\lambda\in\Gamma_j^+$, we have $\Phi^{ext}(\lambda)=h_{\lambda,j}^{ext}(z_1)G_{\lambda,j}^{ext}(z_2)$. On the other hand, according to \Cref{cor:condition-aux-bords-hlambdaj}, we have $(1-\xi(\lambda)z)h_{\lambda,j}^{int}(z)=h_{\lambda,j}^{ext}(z)$, whence it follows that
\[
\frac{\Phi^{ext}(\lambda)}{1-z_1\xi(\lambda)}~=~(1-z_1\xi(\lambda))^{-1}h_{\lambda,j}^{ext}(z_1)G_{\lambda,j}^{ext}(z_2)~=~h_{\lambda,j}^{int}(z_1)G_{\lambda,j}^{ext}(z_2).
\]
Thus
\[
\frac{1-z_1z_2}{2i\pi}\sum_{j=0}^{N-1}\int_{\Gamma_j^+}h_{\lambda,j}^{int}(\lambda)(z_1)G_{\lambda,j}^{ext}(z_2)\,\mathrm{d}\lambda~=~1
\]
by \Cref{egalite4}.
Combining this with \Cref{eq:sdsd8343ZDS,eq:cauchy-kernel}, we can conclude that 
\[
\dualg{\pi_{int}Uk_{\overline{z_1}}}{\pi_{ext}Vk_{\overline{z_2}}}~=~\frac{1}{1-z_1z_2}~=~\dual{k_{\overline{z_1}}}{k_{\overline{z_2}}},
\]
which proves the result we were looking for.
\end{proof}

\Cref{dualite-L2-Gamma} implies that when $F$ satisfies assumptions \ref{H1}, \ref{H2} and \ref{H3'}, the operator $U$ is an isomorphism onto its range.

\subsubsection{Right invertibility of $U$}
 Recall that $Ran(U)\subseteq \cal E_F^p$ and that  we proved in the previous subsection that $U$ is one to one from $H^p$ into $\cal E_F^p$. Our goal here is to prove that $Ran(U)=\cal E_F^p$, i.e. that the operator $U$ is surjective from $H^p$ to $\cal E_F^p$. To prove this, we will check that the left-inverse of $U$ obtained in the previous subsection is also the right-inverse of $U$. This statement is given by the following theorem.
 
\begin{theorem}\label{lem-inverse-a-droite}
Assume that $F$ satisfies \ref{H1}, \ref{H2} and \ref{H3'}. Then for every $u\in \mathcal{E}_F^p$, we have ${U(\pi_{ext}V)^*\pi_{int}}(u)=u$.
\end{theorem}

In order to prove this theorem, the first step is to describe explicitly $(\pi_{ext}V)^*\pi_{int}$, which is done in the next result.

\begin{lemma}\label{lem:calcul-isomorphisme-inverse-U}
For every $u=(u_j)_{0\leq j\leq N-1}\in \bigoplus_{j=0}^{N-1}E^p(\Omega_j^+)$ and every $z\in\mathbb D$, we have 
\[
((\pi_{ext}V)^*\pi_{int}u)(z)~=~\sum_{j=0}^{N-1}\frac{1}{2i\pi}\int_{\Gamma_j^+}u_j^{int}(\lambda)G_{\lambda,j}^{ext}(z)\,\mathrm{d}\lambda.
\]
\end{lemma}

\begin{proof}
Fix $z\in\mathbb D$. According to \Cref{eq:cauchy-kernel}, we have 
\[((\pi_{ext}V)^*\pi_{int}u)(z)~=~
\dual{k_{ \overline z}}{(\pi_{ext}V)^*\pi_{int}u}~=~\dualg{\pi_{ext}Vk_{ \overline z}}{\pi_{int}u}.
\]
As we already observed in the proof of \Cref{dualite-L2-Gamma},  we have 
\[
(\pi_{ext}Vk_{ \overline z})(\lambda)~=~G_{\lambda,j}^{ext}(z) \quad\text{ for every }\lambda\in\Gamma_j^+.
\]
 Hence
\[
((\pi_{ext}V)^*\pi_{int}u)(z)~=~\sum_{j=0}^{N-1}\frac{1}{2i\pi}\int_{\Gamma_j^+}u_j^{int}(\lambda)G_{\lambda,j}^{ext}(z)\,\mathrm{d}\lambda.\qedhere
\]
\end{proof}

In order to check that for every $u\in \mathcal{E}_F^p$, we have ${U(\pi_{ext}V)^*\pi_{int}}(u)=u$, we need to introduce a family of polynomials orthogonal to eigenvectors of $T_F^*$.
\par\smallskip
For every $\lambda\in\mathbb C\setminus\Gamma$, and every integer $m\geq 0$, define
\begin{equation}\label{eq:definition-rhomlambda}
\rho_{\lambda,m}~=~P_+(z^m F_\lambda^{-}),
\end{equation}
where we recall that $F_\lambda^{-1}$ is the function defined in \Cref{Eq:Flambda+-}. 
By \Cref{eq:definition-Glambda,eq:formule-T-F-Lambda}, we have for every $m\ge\w_F(\lambda)=n$
\begin{eqnarray}\label{edqs:3E3434E34D}
(\lambda-T_F)G_{\lambda,m}&=&(T_F-\lambda)(z^{m-n}(F_\lambda^+)^{-1})
~=~T_{F_\lambda^{-}}(F_\lambda^+z^nz^{m-n}(F_\lambda^+)^{-1})\\
&=&P_+(z^m F_\lambda^{-})~=~\rho_{\lambda,m}.\notag
\end{eqnarray}

Fix $m\ge0$. By \Cref{Lem:lien-Flamba-et-flambda+}, we have
\begin{equation}\label{eq122EZ:formule-rholambdam}
\rho_{\lambda,m}~=~P_+\left(z^m \frac{f_\lambda^+(1/z)}{f_\lambda^+(0)}\right)~=~P_+(z^m h_{\lambda,0}(z^{-1})^{-1}), 
\end{equation}
where we still denote by $h_{\lambda,0}$ the function $h_{\lambda,0}(z)=f_\lambda^+(0)/f_\lambda^+(z)$, $z\in\D$, even in the case where $\w_F(\lambda)=0$.
Since $h_{\lambda,0}^{-1}$ belongs to $A(\mathbb D)$, we can expand $h_{\lambda,0}^{-1}$ as a series
 \[
h_{\lambda,0}(z)^{-1}~=~\sum_{j=0}^\infty p_j(\lambda) z^j,
\]
where the sequence $(p_j(\lambda))_{j\geq 0}$ belongs to $\ell^2(\mathbb N)$. Then
\begin{equation}\label{eq:rho-polynome}
\rho_{\lambda,m}(z)~=~P_+\left(z^m\sum_{j=0}^{\infty}p_j(\lambda)z^{-j}\right)~=~\sum_{\ell=0}^m p_{m-\ell}(\lambda)z^\ell.
\end{equation}
In particular, we see that for fixed $\lambda\in\C\setminus\Gamma$ and $m\ge0$, the function $z\longmapsto \rho_{\lambda,m}(z)$ is a polynomial of degree $m$ (observe that $p_0(\lambda)=h_{\lambda,0}(0)^{-1}=1$). 

\begin{lemma}\label{rho-orthogonalite}
For every $\lambda\in\mathbb C\setminus\Gamma$, every $m\geq 0$ and every $j\geq 0$, we have 
\[
\dual{z^jh_{\lambda,0}}{\rho_{\lambda,m}}~=~\delta_{j,m},
\]
where $\delta_{j,m}=1$ if $j=m$ and $\delta_{j,m}=0$ if $j\neq m$.
\end{lemma}

\begin{proof}
According to \Cref{eq122EZ:formule-rholambdam}, we have
\begin{eqnarray*}
\dual{z^jh_{\lambda,0}}{\rho_{\lambda,m}}&=&\dual{z^j h_{\lambda,0}}{P_+(z^m h_{\lambda,0}(z^{-1})^{-1})}\\
&=&\dual{z^j h_{\lambda,0}}{z^m h_{\lambda,0}(z^{-1})^{-1}}\\
&=&\frac{1}{2\pi}\int_0^{2\pi}e^{ij\theta}h_{\lambda,0}(e^{i\theta})e^{-im\theta}h_{\lambda,0}(e^{i\theta})^{-1}\,\mathrm{d}\theta\\
&=&\frac{1}{2\pi}\int_0^{2\pi} e^{i(j-m)\theta}\,\mathrm{d}\theta=\delta_{j,m}.
\end{eqnarray*}

\end{proof}
Observe that, according to \Cref{Lem:analytic-en-dehors-de-la-courbe}, for a fixed $z\in\mathbb D$, the function $\lambda\longmapsto \rho_{\lambda,m}(z)$ is analytic on 
$\mathbb C\setminus\Gamma$. Moreover, since 
\begin{equation}\label{eq:relation-pj-h}
p_j(\lambda)~=~\dual{h_{\lambda,0}(z)^{-1}}{z^j}~=~\dual[\Big]{\frac{f_\lambda^+}{f_\lambda^+(0)}}{z^j},
\end{equation}
using \Cref{remark:bounded-flambdaplus} and Lebesgue's dominated convergence theorem, we see that 
for each component $\Omega$ of $\widehat{\mathbb C}\setminus\Gamma$, the functions $p_j$ admit continuous extensions to $\partial\Omega$. 
\par\smallskip
According to \Cref{cor:condition-aux-bords-hlambdaj} (recall that $F$ satisfies \ref{H2}), we have 
\[
h_{\lambda,0}^{ext}(z)^{-1}(1-z\xi(\lambda))~=~h_{\lambda,0}^{int}(z)^{-1}\quad {\text{a.e. on }\Gamma.}
\]
Remark that this is exactly the boundary relation of \Cref{cor:condition-aux-bords-hlambdaj} written for $j=0$, except for the fact that in \Cref{cor:condition-aux-bords-hlambdaj} it is valid on $\partial\Omega_1^+$ only. But extending $h_{\lambda,0}$ to $\C\setminus \Gamma$ as we did above, we see that it is actually valid almost everywhere on $\Gamma$.
\par\smallskip
Hence by \Cref{eq:relation-pj-h}, we have for every $0\leq j\leq N-1$ that
\[
p_{j+1}^{ext}-\xi p_j^{ext}~=~p_{j+1}^{int}\quad  \text{ a.e. on }\Gamma.
\]
Plugging this equality into \Cref{eq:rho-polynome}, we get 
\begin{equation}\label{eq:equation-aux-bords-rholambdam}
\rho_{\lambda,j+1}^{ext}(z)-\xi(\lambda)\rho_{\lambda,j}^{ext}(z)~=~\rho_{\lambda,j+1}^{int}(z),\; z\in\D\;\;\text{ for a.e. } \lambda\in \Gamma.
\end{equation}
Moreover, since $p_l$ belongs to $A(\Omega)$ for every connected component $\Omega$ of $\widehat{\mathbb C}\setminus\Gamma$, \Cref{eq:rho-polynome} implies  that for a fixed $z\in\D$, the function $\lambda\longmapsto \rho_{\lambda,j}(z)$ lies in $E^q(\Omega_0^+)$  for every $j\ge0$. 

Since $\rho_{\lambda,0}(z)= p_0(\lambda)= h_{\lambda,0}(0)^{-1}=1$ for every $z\in\D$, we also have have 
\begin{equation}\label{eq:equation-bord-rho-0}
\rho_{\lambda,0}^{int}~=~\rho_{\lambda,0}^{ext}.
\end{equation}
\begin{lemma}\label{lem:cauchy-uint-vext}
Let $F$ satisfy \ref{H1}, \ref{H2} and \ref{H3'}. Let $u\in \mathcal{E}_F^p$ and $z\in\mathbb D$. Fix also $0\leq n\leq N-1$ and a point $\mu\in\Omega_{n+1}$. Then
\[
\sum_{m=0}^{N-1}\frac{1}{2i\pi}\int_{\Gamma_m^+}\frac{u_m^{int}(\lambda)\rho_{\lambda,m}^{ext}(z)}{\lambda-\mu}\,\mathrm{d}\lambda~=~\sum_{j=0}^n u_j(\mu)\rho_{\mu,j}(z).
\]
\end{lemma}
\begin{proof}
Fix $z\in\mathbb D$ and $\mu\in\Omega_{n+1}$.
In order to simplify the notation, we write $\widetilde{v}_j(\lambda)=\rho_{\lambda,j}(z)$ and (although these formula do not necessarily correspond to boundary values of functions in a Smirnov space) $$v_j^{ext}(\lambda)~=~\dfrac{{\widetilde{v}_j}^{ext}(\lambda)}{\lambda-\mu}, \quad v_j^{int}(\lambda)~=~\dfrac{{\widetilde{v}_j}^{int}(\lambda)}{\lambda-\mu},  \quad w_j^{int}~=~u_j^{int}v_j^{int}\quad \text{ and } \quad w_j^{ext}~=~u_j^{ext}v_j^{ext}.$$

Since 
$u\in \mathcal{E}_F^p$,  we have  for every $0\leq j\leq N-1$,
\[
u_j^{int}~=~u_j^{ext}+\xi u_{j+1}^{int}\quad \text{ a.e. on }\Gamma\cap \overline{\Omega_{j+1}^+}=\bigcup_{k=j+1}^N\Gamma_k^+,
\]
and, according to \Cref{eq:equation-aux-bords-rholambdam},
\[
v_{j+1}^{ext}-v_{j+1}^{int}~=~\xi v_j^{ext}\quad \text{a.e. on }\Gamma.
\]
Hence  we have almost everywhere on $\Gamma\cap \overline{\Omega_{j+1}^+}$
\[
u_j^{int}v_j^{ext}~=~u_j^{ext}v_j^{ext}+\xi v_j^{ext}u_{j+1}^{int}~=~u_j^{ext}v_j^{ext}+(v_{j+1}^{ext}-v_{j+1}^{int})u_{j+1}^{int},
\]
which gives 
\[
u_j^{int}v_j^{ext}-u_{j+1}^{int}v_{j+1}^{ext}+w_{j+1}^{int}~=~w_{j}^{ext} \quad \text{a.e. on }\Gamma\cap \overline{\Omega_{j+1}^+}.
\]
Let $0\le m\le N-1$.
Summing these equalities for $0\leq j\leq m-1$, and remembering that $\Omega_m^+\subseteq\Omega_j^+$, we obtain 
\[
u_0^{int}v_0^{ext}-u_{m}^{int}v_m^{ext}+\sum_{j=1}^m w_{j}^{int}~=~\sum_{j=0}^{m-1}w_{j}^{ext}\quad \text{a.e. on }\Gamma\cap \overline{\Omega_m^+},
\]
that is 
\[
u_{m}^{int}v_m^{ext}~=~u_0^{int}v_0^{ext}+\sum_{j=1}^m w_{j}^{int}-\sum_{j=0}^{m-1}w_{j}^{ext}\quad \text{a.e. on }\Gamma\cap \overline{\Omega_m^+}.
\]
Denote by $I$ the integral
\[
I~=~\sum_{m=0}^{N-1}\frac{1}{2i\pi}\int_{\Gamma_m^+}\frac{u_m^{int}(\lambda)\rho_{m,\lambda}^{ext}(z)}{\lambda-\mu}\,\mathrm{d}\lambda.
\]
Then 
\begin{eqnarray*}
I&=&\sum_{m=0}^{N-1}\frac{1}{2i\pi}\int_{\Gamma_m^+}u_0^{int}(\lambda)v_0^{ext}(\lambda)\,\mathrm{d}\lambda+\sum_{m=0}^{N-1}\frac{1}{2i\pi}\sum_{j=1}^m \int_{\Gamma_m^+}w_j^{int}(\lambda)\,\mathrm{d}\lambda\\
&&~\hspace{5cm}-\sum_{m=0}^{N-1}\frac{1}{2i\pi}\sum_{j=0}^{m-1} \int_{\Gamma_m^+}w_j^{ext}(\lambda)\,\mathrm{d}\lambda\\
&=&\sum_{m=0}^{N-1}\frac{1}{2i\pi}\int_{\Gamma_m^+}u_0^{int}(\lambda)v_0^{ext}(\lambda)\,\mathrm{d}\lambda+\sum_{j=1}^{N-1}\sum_{m=j}^{N-1}\frac{1}{2i\pi}\int_{\Gamma_m^+}w_j^{int}(\lambda)\,\mathrm{d}\lambda\\ &&~\hspace{5cm}-\sum_{j=0}^{N-2}\sum_{m=j+1}^{N-1}\frac{1}{2i\pi}\int_{\Gamma_m^+}w_j^{ext}(\lambda)\,\mathrm{d}\lambda.
\end{eqnarray*}
By subsituting in the last sum $m$ by $m-1$ and taking into account that $\Gamma_N^+=\emptyset$, we obtain
\begin{multline*}
I~=~\sum_{m=0}^{N-1}\frac{1}{2i\pi}\int_{\Gamma_m^+}(u_0^{int}(\lambda)v_0^{ext}(\lambda)-w_0^{int}(\lambda))\,\mathrm{d}\lambda\\
+~\sum_{j=0}^{N-1}\sum_{m=j}^{N-1}\frac{1}{2i\pi}\int_{\Gamma_m^+}w_j^{int}(\lambda)\,\mathrm{d}\lambda-\sum_{j=0}^{N-2}\sum_{m=j}^{N-2}\frac{1}{2i\pi}\int_{\Gamma_{m+1}^+}w_j^{ext}(\lambda)\,\mathrm{d}\lambda.
\end{multline*}
Now, observe that $v_0^{ext}=v_0^{int}$ a.e. (according to \cref{eq:equation-bord-rho-0}), whence it follows that $u_0^{int}v_0^{ext}-w_0^{int}=u_0^{int}v_0^{ext}-u_0^{int}v_0^{int}=0$ a.e., so again using that $\Gamma_N^+=\varnothing$, we have
\begin{align*}
I~&=~\sum_{j=0}^{N-1}\sum_{m=j}^{N-1}\frac{1}{2i\pi}\int_{\Gamma_m^+}w_j^{int}(\lambda)\,\mathrm{d}\lambda-\sum_{j=0}^{N-2}\sum_{m=j}^{N-2}\frac{1}{2i\pi}\int_{\Gamma_{m+1}^+}w_j^{ext}(\lambda)\,\mathrm{d}\lambda\\
&=~\sum_{j=0}^{N-1}\sum_{m=j}^{N-1}\frac{1}{2i\pi}\left(\int_{\Gamma_m^+}w_j^{int}(\lambda)\,\mathrm{d}\lambda-\int_{\Gamma_{m+1}^+}w_j^{ext}(\lambda)\,\mathrm{d}\lambda\right).
\end{align*}
Observe now that $u_j\widetilde{v_j}\in E^1(\Omega_m^+)$ for $m\geq j$. Then by Cauchy's formula for functions in $E^1(\Omega_m^+)$ (see \cite[Th. 10.4]{Duren1970}), we get that
\begin{align*}
\frac{1}{2i\pi}\left(\int_{\Gamma_m^+}w_j^{int}(\lambda)\,\mathrm{d}\lambda-\int_{\Gamma_{m+1}^+}w_j^{ext}(\lambda)\,\mathrm{d}\lambda\right)&~=~\frac{1}{2i\pi}\left(\int_{\Gamma_m^+}\frac{u_j^{int}(\lambda){\widetilde{v_j}}^{int}(\lambda)}{\lambda-\mu}\,\mathrm{d}\lambda\right.\\ &\hspace{2.5cm}-\left.\int_{\Gamma_{m+1}^+}\frac{u_j^{ext}(\lambda){\widetilde{v_j}}^{ext}(\lambda)}{\lambda-\mu}\,\mathrm{d}\lambda\right)\\
&~=~\frac{1}{2i\pi}\int_{\partial\Omega_{m+1}}\frac{u_j(\lambda)\widetilde{v_j}(\lambda)}{\lambda-\mu}\,\mathrm{d}\lambda\\
&~=~\begin{cases}
0&\text{if }m\neq n\\
u_j(\mu)\widetilde{v_j}(\mu)&\text{if }m=n.
\end{cases}
\end{align*}
Finally, we get
\[
I~=~\sum_{j=0}^{N-1}\sum_{m=j}^{N-1}\frac{1}{2i\pi}\left(\int_{\Gamma_m^+}w_j^{int}(\lambda)\,\mathrm{d}\lambda-\int_{\Gamma_{m+1}^+}w_j^{ext}(\lambda)\,\mathrm{d}\lambda\right)~=~\sum_{j=0}^n u_j(\mu)\widetilde{v_j}(\mu),
\]
which is the equality we wanted to prove.
\end{proof}
We will now check the right invertibility of $U$. 

\begin{proof}[Proof of \Cref{lem-inverse-a-droite}]
Observe that, according to \Cref{lem:description-eigenvectors}, for every $\mu\in\sigma(T_F)\setminus F(\mathbb T)$, for every $\lambda\in\mathbb C\setminus F(\mathbb T)$, and every $0\leq m<\w_F(\mu)$, we have
\[
T_{F-\lambda}^*h_{\mu,m}~=~(\mu-\lambda)h_{\mu,m}.
\]
Hence, for $0\le m<\w_F(\mu)$, $\ell\geq \w_F(\lambda)$ and $\lambda\neq \mu$, we get 
\begin{eqnarray*}
\dual{G_{\lambda,\ell}}{h_{\mu,m}}&=&\frac{1}{\mu-\lambda}\dual{G_{\lambda,\ell}}{T_{F-\lambda}^*h_{\mu,m}}\\
&=&\frac{1}{\mu-\lambda}\dual{T_{F-\lambda}G_{\lambda,\ell}}{h_{\mu,m}}\\
&=&\frac{1}{\lambda-\mu}\dual{ \rho_{\lambda,\ell}}{h_{\mu,m}},
\end{eqnarray*}
the last equality following from Equation \eqref{edqs:3E3434E34D}.

Let now $u=(u_j)_{0\leq j\leq N-1}$ be an element of $\mathcal{E}_F^p$. Write $U((\pi_{ext}V)^*\pi_{int}u)=(r_j)_{}$, with $r_j\in E^p(\Omega_j^+)$ for every $0\leq j\leq N-1$. By the definition of $U$, we have  
\[
r_j(\mu)~=~\dual{(\pi_{ext}V)^*\pi_{int}u}{h_{\mu,j}}~=~\dualg{\pi_{int}u}{\pi_{ext}Vh_{\mu,j}}\quad\text{ for every }\mu\in \Omega_j^+.
\]
Now, observe that if we write $Vh_{\mu,j}=(v_\ell)_{0\leq \ell\leq N-1}$, we have, by the definition of $V$,
\[
v_\ell(\lambda)~=~\dual{G_{\lambda,\ell}}{h_{\mu,j}}~=~\frac{1}{\lambda-\mu}\dual{\rho_{\lambda,\ell}}{h_{\mu,m}} \quad \text{ for every }\lambda\in\widehat{\mathbb C}\setminus(F(\mathbb T)\cup\Omega_\ell^+).
\]
Thus, for almost all $\lambda\in \Gamma_\ell^+$, we have
\[
\pi_{ext}(Vh_{\mu,j})(\lambda)~=~v_\ell^{ext}(\lambda)~=~\frac{1}{\lambda-\mu}\dual{\rho^{ext}_{\lambda,\ell}}{h_{\mu,m}},
\]
which yields, by \Cref{lem:calcul-isomorphisme-inverse-U}, that
\[
r_j(\mu)~=~\frac{1}{2i\pi}\sum_{\ell=0}^{N-1}\int_{\Gamma_\ell^+}\frac{u_\ell^{int}(\lambda)}{\lambda-\mu}\dual{\rho_{\lambda,\ell}^{ext}}{h_{\mu,j}}\,\mathrm{d}\lambda.
\]
Observe that the function
\[
z~\longmapsto ~\int_{\Gamma_\ell^+}\frac{u_\ell^{int}(\lambda)\rho_{\lambda,\ell}^{ext}(z)}{\lambda-\mu}\,\mathrm{d}\lambda
\]
belongs to $H^p$ (because it is  actually  a polynomial). Thus, by continuity of the bilinear form $\dual{\cdot}{\cdot}$, we get 
\begin{eqnarray*}
r_j(\mu)
&=&\frac{1}{2i\pi}\sum_{\ell=0}^{N-1}\int_{\Gamma_\ell^+}\frac{u_\ell^{int}(\lambda)}{\lambda-\mu}\dual{\rho_{\lambda,\ell}^{ext}}{h_{\mu,j}}\,\mathrm{d}\lambda\\
&=&\dual[\Bigg]{\sum_{\ell=0}^{N-1}\int_{\Gamma_\ell^+}\frac{u_\ell^{int}(\lambda)\rho_{\lambda,\ell}^{ext}}{\lambda-\mu}\,\mathrm{d}\lambda}{h_{\mu,j}}.
\end{eqnarray*}

Now,  
it follows from \Cref{lem:cauchy-uint-vext} that, for every $z\in\mathbb D$ and every $\mu\in\Omega_k$, with $k\ge j+1$, we have
\[
\sum_{\ell=0}^{N-1}\int_{\Gamma_\ell^+}\frac{u_\ell^{int}(\lambda)\rho_{\lambda,\ell}^{ext}(z)}{\lambda-\mu}\,\mathrm{d}\lambda~=~\sum_{\ell=0}^{k-1}u_\ell(\mu)\rho_{\mu,\ell}(z).
\]
Hence
\[
r_j(\mu)~=~\sum_{\ell=0}^{k-1} u_\ell(\mu)\dual{\rho_{\mu,\ell}}{h_{\mu,j}}.
\]
Finally, \Cref{rho-orthogonalite} implies that 
\[
r_j(\mu)~=~u_j(\mu)\quad \text{for every }\mu\in\Omega_k,\,k\ge j+1.
\]
In other words, $r_j(\mu)=u_j(\mu)$ for every $\mu\in \Omega_j^+$, which means that $r_j=u_j$, and we have shown that $U(\pi_{ext}V)^*\pi_{int}(u)=u$.
\end{proof}
We now have everything we need to prove the main result of \cite{Yakubovich1991}.
\subsubsection{Summary}
We denote by $M_\lambda$ the multiplication by the independent variable on $\mathcal{E}_F^p$, that is, if $u=(u_j)_{0\leq j\leq N-1}\in \mathcal{E}_F^p$, then 
\[
M_\lambda u~=~(v_j)_{0\leq j\leq N-1},\quad \text{with }v_j(\lambda)~=~\lambda u_j(\lambda)\quad \textrm{ for every } \lambda\in\Omega_j^+.
\]
It is not difficult to see that $M_\lambda$ is a bounded operator on $\mathcal{E}_F^p$. More generally, if $\varphi$ is a bounded analytic function on the interior $\overset{\circ}{\sigma(T_F)}$ of the spectrum of $T_F$, then the multiplication by $\varphi$ is a bounded operator on $\mathcal{E}_F^p$ satisfying $\|M_\varphi\|\leq \|\varphi\|_\infty$. Indeed, if $u=(u_j)_{0\leq j\leq N-1}\in \mathcal{E}_F^p$, then we obviously have $\varphi u_j\in E^p(\Omega_j^+)$ and $\|\varphi u_j\|_{E^p(\Omega_j^+)}\leq \|\varphi\|_\infty \|u_j\|_{E^p(\Omega_j^+)}$.  Moreover, since $\left(\Gamma\cap\overline{\Omega_{j+1}^+}\right)\setminus\mathcal O\subset \overset{\circ}{\sigma(T_F)}$, we get immediately that $\varphi^{int}=\varphi^{ext}$ almost everywhere on $\Gamma\cap\overline{\Omega_{j+1}^+}$, $0\le j\le N-1$, and it is clear that $\varphi u_j$ satisfies the boundary conditions in \Cref{eq:definition-espace-modele}. Thus we deduce that $\varphi u\in \mathcal{E}_F^p$ and 
\[
\|\varphi u\|^p_{\mathcal{E}_F^p}~=~\sum_{j=0}^{N-1}\|\varphi u_j\|_{E^p(\Omega_j^+)}^p~\leq ~\|\varphi\|_\infty^p\sum_{j=0}^{N-1}\| u_j\|_{E^p(\Omega_j^+)}^p~=~\|\varphi\|_\infty^p\|u\|_{\mathcal{E}_F^p}^p,
\]
meaning that $M_\varphi$ is bounded on $\mathcal{E}_F^p$ with $\|M_\varphi\|\leq \|\varphi\|_\infty$. 

\begin{theorem}\label{model}
Let $F$ satisfy \ref{H1}, \ref{H2} and \ref{H3'}, and let $p>1$. 
 Then 
\begin{enumerate}
\item The operator $U$ is an isomorphism from $H^p$ onto $\mathcal{E}_F^p$.
\item We have 
\[
UT_F~=~M_\lambda U.
\]
\item $T_F$ admits an $H^\infty(\overset{\circ}{\sigma(T_F)})$ functional calculus, and there exists a constant $C>0$ such that
\[
\|\varphi(T_F)\|~\leq~ C\|\varphi\|_\infty\quad\text{ for every } \varphi\in H^\infty(\overset{\circ}{\sigma(T_F)}).
\]
\end{enumerate}
\end{theorem}
\begin{proof}
The fact that $U$ is an isomorphism from $H^p$ onto $\mathcal{E}_F^p$ follows immediately from \Cref{dualite-L2-Gamma,lem-inverse-a-droite}. The fact (2) has already been proved in \Cref{petites proprietes} (1).
Now, given a bounded analytic function $\varphi$ on $\Omega=\overset{\circ}{\sigma(T_F)}$, we can define $\varphi(T_F)=U^{-1}M_\varphi U$. Then $\varphi(T_F)$ is well-defined on $H^p$, and we have 
\[
\|\varphi(T_F)\|~\leq~ \|U\| \|U^{-1}\|\|M_\varphi\|~\leq~  \|U\| \|U^{-1}\| \|\varphi\|_\infty.
\]
It is then not difficult to check that $\varphi\longmapsto \varphi(T_F)$ defines an $H^\infty(\Omega)$ functional calculus.
\end{proof}

\subsection{From \Cref{Section:Yakubovich_demoHp} to \Cref{Section:Yakubovich-court}}\label{mot-de-la-fin}
In this very last section, we explain the change of notation which allows us to go from the setting of \Cref{Section:Yakubovich_demoHp} to the setting of \Cref{Section:Yakubovich-court}.
\par\smallskip
Let $F$ be a symbol satisfying assumptions \ref{H1}, \ref{H2} and \ref{H3} - in particular, $F$ is \textbf{negatively wound}; these are the standing assumptions in the first sections of the paper. Let $f\in L^\infty(\mathbb T)$ be defined as $f(z)=F(1/z),\, z\in \mathbb T$. Then $f$ satisfies assumptions \ref{H1}, \ref{H2} and \ref{H3'}. If $T_F$ is seen as acting on $H^p,\,p>1$, then $T_f=T_F^*$ acts on $H^q$. \Cref{model}  applies to $f$, allowing us to define the operator $U$ from $H^q$ into $\mathcal{E}_f^q$, where $\mathcal{E}_f^q$ is the set of $N$-tuples of functions $(u_j)_{0\le j\le N-1}$ with $u_j\in E^q(\Omega_j^+)$ satisfying the boundary relations
\[u_j^{int}-\zeta u_{j+1}^{int}~=~u_j^{ext}\quad\text{ a.e. on }\partial\Omega_{j+1}^+,\]
where $\zeta=f^{-1}$ on $f(\mathbb T)\setminus \mathcal O$, $\mathcal O$ denoting the set of self-intersection points of the curve $f(\mathbb T)$ (which is the same as the set of self-intersection points of the curve $F(\mathbb T)$). Since $f(\zeta(\lambda))=\lambda$ for every $\lambda\in F(\mathbb T)\setminus \mathcal O$, $F(1/\zeta(\lambda))=\lambda$, i.e. $\zeta=1/F^{-1}$ on $F(\mathbb T)\setminus \mathcal O$ (cf. the definition of $\zeta$ at the end of \Cref{Subsection:IntExt}). The space $\mathcal{E}^q_f$ is denoted by $E_F^q$ in \Cref{Section:Yakubovich-court}.
\par\smallskip
We have then $UT_f=M_\lambda U$, where $M_\lambda$ is the multiplication operator by $\lambda$ on $\mathcal{E}^q_f$, i.e. $UT_F^*=M_\lambda U$: this is exactly the relation $T_F^*=U^{-1}M_\lambda U$ given in \Cref{T:Yakbovich_Hp}. 
\par\smallskip
The eigenvectors $h_{\lambda,j}$ of $T_F$ defined by \Cref{eve} are exactly those given by \Cref{Annexe-defn-hlambdak}, since the function $F$ from \Cref{Section:Yakubovich-court} plays the role of the function $f$ from \Cref{Section:Yakubovich_demoHp}.
\par\smallskip
The same correspondence is in force when the symbol $F$ from \Cref{section:JFA} is supposed to satisfy assumption \ref{H2'} rather than \ref{H2} (besides assumptions \ref{H1} and \ref{H3}, of course).

\bibliographystyle{plain}
\bibliography{biblio}
\end{document}